\documentclass[12pt]{amsbook}
\usepackage{graphicx} 	%is this ok ??
\usepackage{color}	%for compatibility with my figures files, even if
\input{amssym.def}
\input{amssym}

\def\R{\mathbb R}
\def\N{\mathbb N}
\def\T{\mathbb T}

\def\Z{\mathbb Z}

\def\E{\mathbb E}
\def\Cplx{\mathbb C}  %ATTENTION !!!

\def\E{\mathcal E}
\def\L{\mathcal L}
\def\C{\mathcal C}
\def\<{\langle}
\def\>{\rangle}
\def\a{{\rm a}}

\newcommand{\BS}{\operatorname{BS}}

\def\vol{\mathop{\rm vol}}
\def\grad{\mathop{\rm grad}}

\def\iint{\int\!\int}

\def\Nc{{N_c}} %number of commutators
\def\sle{\preccurlyeq}

\def\H{{\mathcal H}} %in that particular paper, some Hilbert spaces
\def\V{{\mathcal V}} %
\def\K{{\mathcal K}} %

\def\I{{\rm (I)}}
\def\II{{\rm (II)}}
\def\III{{\rm (III)}}

\def\Mf0{{M^f_{\rho\, 0\, 1}}}

\def\sym{{\rm sym}}

\def\ker{\mathop{\rm Ker}}

\def\var{\varepsilon}

\def\tilde{\widetilde}
\def\pa{\partial}

\def\Om{\Omega}
\def\ov{\overline}
\def\cal{\mathcal}

\def\hat{\widehat}

\def\tilde{\widetilde}

\def\Id{{\rm Id}\,}
\def\loc{{\rm loc}}

\def\llangle{\left\langle}
\def\rrangle{\right\rangle}
\def\Blangle{\Bigl\langle}
\def\Brangle{\Bigr\rangle}

\def\dr#1#2{\left. \frac{d}{d#1} \right |_{#2}}
\def\2dr#1#2#3{\left. \frac{d^2{#1}}{d{#2}^2} \right |_{#3}}
\def\d2#1#2{\frac{d^2{#1}}{d{#2}^2}}
\def\ip#1#2{{\llangle #1, #2 \rrangle}}
\def\<{\langle}
\def\>{\rangle}
\def\derpar#1#2{\frac{\partial#1}{\partial#2}}

\def\iip#1#2{{\left (\!\left ( #1 , #2 \right ) \!\right )}}

\def\rhs{right-hand side}

\def\wrt{with respect to }

\def\iff{if and only if }

\def\dps{\displaystyle}

\def\med{\medskip}

\def\sm{\smallskip}

\def\begeq{\begin{equation}}
\def\endeq{\end{equation}}
\def\begar{\begin{eqnarray}}
\def\endar{\end{eqnarray}}
\def\begar*{\begin{eqnarray*}}
\def\endar*{\end{eqnarray*}}
\def\begal{\begin{align}}
\def\endal{\end{align}}
\def\begal*{\begin{align*}}
\def\endal*{\end{align*}}

\newtheorem{Thm}{Theorem}
\newtheorem{Lem}[Thm]{Lemma}

\newtheorem{Cor}[Thm]{Corollary}
\newtheorem{Prop}[Thm]{Proposition}

\newtheorem{Ass}{Assumption}
\newtheorem{Modprob}[Thm]{Model Problem}

\newtheorem*{Thm*}{Theorem}
\newtheorem*{Lem*}{Lemma}
\newtheorem*{Conj*}{Conjecture}
\newtheorem*{Cor*}{Corollary}
\newtheorem*{Def*}{Definition}
\newtheorem*{Prop*}{Proposition}
\newtheorem*{Exo*}{Exercise}
\newtheorem*{Exs*}{Examples}
\newtheorem*{Ex*}{Example}
\newtheorem*{Rk*}{Remark}
\newtheorem*{Rks*}{Remarks}

\theoremstyle{definition}
\newtheorem{Def}[Thm]{Definition}
\newtheorem{Rk}[Thm]{Remark}

\newtheorem{Ex}[Thm]{Example}

\begin{document}

\title{Hypocoercivity}

\author{C\'edric Villani}

\address{Unit\'e de Math\'ematiques Pures et Appliqu\'ees \\ 
UMR CNRS 5669 \\
Ecole Normale Sup\'erieure de Lyon \\
46 all\'ee d'Italie \\
F-69364 Lyon Cedex 07 \\
FRANCE} 

\email{cvillani@umpa.ens-lyon.fr}

\date{September 1, 2006} 

\begin{abstract}
This memoir attempts at a systematic study of convergence to
stationary state for certain classes of degenerate diffusive equations, 
by means of well-chosen Lyapunov functionals. 
Many open problems and possible directions
for future research are discussed.

\end{abstract}

\maketitle
\vfill

\cleardoublepage

\thispagestyle{empty} \bigskip\bigskip\bigskip\bigskip
\begin{flushright} {\em \Large }
\end{flushright}\vfill \pagebreak

\cleardoublepage
\setcounter{page}{5}		%%%changed from p.4 to p.5

\setcounter{chapter}{0}
%\dominitoc			%now useless
\tableofcontents

\chapter*{Introduction}
In many fields of applied mathematics, one is led to study 
dissipative evolution equations involving 
(i) a degenerate dissipative operator, and 
(ii) a conservative operator presenting certain symmetry properties, 
such that the combination of both operators implies convergence to a
uniquely determined equilibrium state. Typically, the dissipative part is not
coercive, in the sense that it does not admit a spectral gap; instead, it may possess
a huge kernel, which is \emph{not} stable under the action of the conservative part.
This situation is very similar to problems encountered in the theory
of hypoellipticity, in which
the object of study is not convergence to equilibrium, but regularity. 
By analogy, I shall use the word {\bf hypocoercivity},
suggested to me by Thierry Gallay, to describe this phenomenon. 
This vocable will be used somewhat loosely in general,
and in a more precise sense when occasion arises.

Once the existence and uniqueness of a steady state has been established
(for instance by direct computation, or via an abstract theorem such as
Perron-Frobenius), there are plenty of soft tools to prove \emph{convergence}
to this steady state. It is much more tricky and much more instructive
to find estimates about \emph{rates of convergence}, 
and this is the question which will be addressed here.

Both hypoellipticity and hypocoercivity often occur together in the study of 
linear diffusion generators satisfying H\"ormander's bracket condition. It is for such 
equations that theorems of exponentially fast convergence to equilibrium were first established
via probabilistic tools~\cite{MSH:ergodic:02,talay:stochHamilt:02,reybelletthomas:anharmonic:00,reybelletthomas:exp:02}, taking their
roots in the Meyn-Tweedie theory of the asymptotic behavior of Markov chains.
Some of these studies were motivated by the study of finite-dimensional approximations
of randomly forced two-dimensional Navier-Stokes equations~\cite{EMS:NS:01,MSH:ergodic:02,mattingly:NS:02};
since then, the theory has been developed to the extent that it can deal with truly 
infinite-dimensional systems~\cite{hairermattingly}.
In all these works,
exponential convergence is established, but there are no quantitative estimates of
the rate. Moreover, these methods usually try to capture information about path behavior,
which may be useful in a probabilistic perspective, but is more than what we need.

Analytical approaches can be expected to provide more precise results; they have
been considered in at least three (quite different, and complementary) settings:
\sm

- For {\bf nonlinear equations} possessing a distinguished {\bf Lyapunov functional} 
(entropy, typically), robust methods, based on functional inequalities, time-derivative
estimates and interpolation, have been developed to establish convergence
estimates in $O(t^{-\infty})$, i.e. faster than any inverse power of time.
These methods have been applied to the linear (!) Fokker-Planck equation~\cite{DV:FP:01}, 
the Boltzmann equation~\cite{DV:boltz:05}, and some variants arising in the context
of kinetic theory~\cite{CCG:linear:03,FNS:linear:04}.
So far, they rely crucially on strong regularity a priori estimates.
\sm

- For {\bf linear hypoelliptic equations} enjoying some structural properties, more specific 
methods have been developed to prove (ideally) exponential convergence to equilibrium
with explicit bounds on the rate. Up to now, this approach has been mainly
developed by H\'erau and Nier~\cite{heraunier:FP:04}, 
Eckmann and Hairer~\cite{eckmannhairer:hypo:03},
Helffer and Nier~\cite{helfnier:witten:05}, for second-order differential operators in
H\"ormander's form (a sum of squares of derivations, plus a derivation).
It uses pseudo-differential operators, and a bit of functional calculus; it can be seen
as an extension of Kohn's celebrated method for the study of hypoellipticity of H\"ormander 
operators. In fact, the above-mentioned works establish hypoellipticity at the same time
as hypocoercivity, by considering functional spaces with polynomial 
weights in both Fourier space and
physical space. After a delicate spectral analysis, they localize the spectrum inside a
cusp-like region of the complex plane, and then deduce the exponential convergence to
equilibrium. Again, in some sense these methods capture more than needed,
since they provide information on the whole spectrum.
\sm

- Finally, Yan Guo recently developed a new method~\cite{guo:landau:02}, 
which he later pushed forward with 
Strain~\cite{guostrain:almostexp:06,guostrain:exp}, 
to get rates of convergence for nonlinear kinetic equations in a close-to-equilibrium
regime. Although the method is linear in essence, it is based on robust functional
inequalities such as interpolation or Poincar\'e inequalities; so it is in some
sense intermediate between the two previously described lines of research. 
\sm

The goal of this memoir is to start a systematic study of hypocoercivity in its
own right. The {\bf basic problem considered here} consists in identifying
general structures in which the interplay between a ``conservative'' part
and a ``degenerate dissipative'' part lead to convergence to equilibrium.

With respect to the above-mentioned
works, the novelty of the approach explored here resides in its abstract nature
and its simplicity. In particular, I~wish to convey the following two messages:
\sm

1. Hypocoercivity is related to, but distinct from hypoellipticity, and in many
situations can be established quantitatively independently of regularity issues,
or after regularity issues have been settled.
\sm

2. There are some general and simple techniques, based on very elementary 
but powerful algebraic tricks, by which one can often reduce a 
mysterious hypocoercive situation to a much more standard coercive one.
\sm

There are three parts in this memoir:
\sm

Part~I focuses on the particular case of operators which 
(as in~\cite{heraunier:FP:04,eckmannhairer:hypo:03,helfnier:witten:05})
can be written in ``H\"ormander form'' $A^*\! A+B$, where $A$ and $B$
are possibly unbounded operators on a {\em given Hilbert space}. 
These results have been applied to several models, such as the kinetic
Fokker--Planck equation, the linearized 
Landau--Lifschitz--Gilbert--Maxwell model
in micromagnetism~\cite{capella}, and a model problem
for the stability of the Oseen vortices~\cite{gallgall}.
\sm

Part~II, by far the shortest, 
remains at a linear level, but considers operators
which cannot necessarily be written in the form $A^*A+B$,
at least for ``tractable'' operators $A$ and $B$.
In this part I~shall give an abstract version of a
powerful hypocoercivity theorem recently established by
Mouhot and Neumann~\cite{mouhotneumann:06}, explain why
we cannot be content with this theorem, and give
some suggestions for research in this direction.
\sm

In Part~III I~shall consider fully nonlinear equations,
in a {\em scale of Sobolev-type spaces}, in presence of
a ``good'' Lyapunov functional. In this setting I~shall
obtain results that can apply to a variety of nonlinear
models, {\em conditionally to smoothness bounds}.
In particular I~shall simplify the proof of
the main results in~\cite{DV:boltz:05}.
\sm

Though these three settings are quite different, and far from being
unified, there is a unity in the methods that will be used:
{\bf construct a Lyapunov functional by adding carefully
chosen lower-order terms to the ``natural'' Lyapunov
functional}. This simple idea will turn out to be quite
powerful.
\sm

The method will be presented in a rather systematic and abstract way. 
There are several motivations for this choice of presentation.
First, the methods are general enough and can be applied
in various contexts.
Also, this presentation may be pedagogically relevant, by emphasizing the most
important features of the problem. Last but not least, most of the time 
I~really had to set the problems in abstract terms,
to figure out a way of attacking it. 
\sm

No attempt will be made here for a {\em qualitative} study of the
approach to equilibrium, but I~believe this is a very rich topic,
that should be addressed in detail in the future. 
One of the main outcomes of my work with 
Laurent Desvillettes~\cite{DV:boltz:05} was the prediction that
solutions of the Boltzmann equation, while approaching equilibrium,
would oscillate between ``close to hydrodynamic'' and
``close to homogeneous'' states. To some extent, this guess was
in contradiction with a commonly accepted idea according to which 
the large-time behavior should be dominated by the hydrodynamic regime;
nevertheless these oscillations have been spectacularly confirmed in
numerical simulations by Francis Filbet. Further developments
can be found in~\cite{FMP:NlogN}; the results obtained by
numerical simulations are so neat that they demand
a precise explanation.
\sm

Research in the area of hypocoercivity is currently developing fast 
thanks to the efforts of several other researchers such as
Thierry Gallay, Fr\'ed\'eric H\'erau, Cl\'ement Mouhot, and others.
I~expect that further important results will soon be available
thanks to their efforts, and that this memoir really can be considered
as a starting point of a much more developed theory.

%%% the following are only there, to avoid numbering inside the
%%% Introduction as 0.x.

%\numberwithin{section}{} 
\numberwithin{equation}{section} 
%\numberwithin{figure}{chapter} 	
\def\thepart{\Roman{part}}
%%% now come parts 1 to 3

%\include{part1-1}

\part{$L=A^*A+B$} \label{partAAB}

In this part I~shall study the convergence to equilibrium for degenerate
linear diffusion equations where the diffusion operator takes the
abstract form $A^*\! A+B$, $B^*=-B$. 

The main abstract theorem makes crucial use of {\bf commutators},
in the style of H\"ormander's hypoellipticity theorem.
In its simplest form, it reduces the problem of convergence to
equilibrium for the non-symmetric, non-coercive operator $A^*A+B$, 
to that of the symmetric, possibly coercive operator $A^*A+[A,B]^*[A,B]$.
If the latter operator is not coercive, then one may consider
iterated commutators $[[A,B],B]$, $[[[A,B],B],B]$, etc.

One of the first main results (Theorem~\ref{hypocomult})
can be informally stated as follows:
{\em Let $A=(A_1,\ldots,A_m)$, $B^*=-B$, and $L=A^*A+B$ be
linear operators on a Hilbert space $\H$. Define iterated commutators
$C_j$ and remainders $R_j$ ($1\leq j\leq \Nc$) by the identities
$\dps C_0=A, \quad [C_j,B] = C_{j+1} + R_{j+1}\quad (j\leq \Nc), 
\quad C_{\Nc+1}=0.$
If $\sum_{j=0}^{\Nc} C_j^*C_j$ is coercive, and 
the operators $[A,C_k]$, $[A^*,C_k]$, $R_k$
satisfy certain bounds, then $\| e^{-tL} \|_{{\cal H}^1\to {\cal H}^1}
= O(e^{-\lambda t})$, where the ``Sobolev'' space $\H^1$ is defined
by the Hilbert norm $\|h\|_{{\cal H}^1}^2 = \|h\|^2 + \sum \|C_j h\|^2$.}
\med

The key ingredient in the proof is the construction of an auxiliary
Hilbert norm, which is equivalent to the $\H^1$ Hilbert norm, but has
additional ``mixed terms'' of the form $\ip{C_jh}{C_{j+1}h}$.

Applied to the
kinetic Fokker--Planck equation, these theorems will yield results of
convergence to equilibrium that are both more general and more precise
than previously known estimates.

After this ``abstract'' $L^2$ framework, a ``concrete'' $L\log L$
framework will be considered, leading to results of convergence
for very general data (say finite measures).

My reflexion on this subject started during the preparation of
my Cours Peccot at the Coll\`ege de France (Paris), in June 2003,
and has crucially benefited from interactions with many people. 
The first draft of the proof of Theorem~\ref{thmsimple}
occurred to me while I was struggling to understand the results 
of Fr\'ed\'eric H\'erau and Francis Nier~\cite{heraunier:FP:04} 
about kinetic Fokker--Planck equations.
The construction of the anisotropic Sobolev norm was partly inspired by
the reading of papers by Yan Guo~\cite{guo:landau:02} and Denis 
Talay~\cite{talay:stochHamilt:02}; 
although their results and techniques are quite
different from the ones in the present paper, they were the first to draw 
my attention to the interest of
introducing mixed terms such as $\nabla_vf\cdot\nabla_xf$. 
Denis also showed me a useful trick for getting long-time estimates 
on the moments of certain hypoelliptic diffusion equations, 
which is based on the construction of an adequate quadratic form.

Apart from the above-mentioned people, I was lucky enough 
to have fruitful discussions on the subject with Bernard Helffer,
Laurent Desvillettes, 
Luc Rey-Bellet, Jean-Pierre Eckmann, Martin Hairer, Cl\'ement Mouhot, 
Stefano Olla and Piere-Louis Lions, as well as with Christian Schmeiser 
and Denis Serre, who both suggested a relation between
my results and Kawashima's condition in the theory of hyperbolic 
systems of conservation laws. 
\vfill\pagebreak

\section{Notation} \label{secnot}

\subsection{Basic notation} \label{subnot}

Let $\H$ be a separable (real or complex) Hilbert space, 
to be thought of as $L^2(\mu)$, where $\mu$ is some equilibrium measure; 
$\H$ is endowed with a norm $\|\cdot\|$ coming from a scalar 
(or Hermitian) product $\ip{\cdot}{\cdot}$.

Let $\V$ be a finite-dimensional Hilbert space (say $\R^m$ or $\Cplx^m$,
depending on whether $\H$ is a real or complex Hilbert space).
Typically, $\V$ will be the space of those
variables on which a certain diffusion operator acts.
The assumption of finite dimension covers all cases that will be considered
in applications, but it is not essential.

Let $A:\H\to\H\otimes\V\simeq \H^m$ 
be an unbounded operator with domain $D(A)$, and 
let $B:\H\to\H$ be an unbounded antisymmetric operator with domain $D(B)$:
\[ \forall h,h' \in D(B),\qquad \ip{Bh}{h'} = - \ip{h}{Bh'}.\]
I shall assume that there is a dense topological vector space ${\cal S}$ in $\H$ such that
${\cal S}\subset D(A)\cap D(B)$ and $A$ (resp. $B$) continuously sends ${\cal S}$ into
${\cal S}\otimes \V$ (resp. ${\cal S}$); this assumption is here only to 
guarantee that all the computations that will be performed (involving a finite
number of applications of $A$, $A^*$ and $B$) are authorized.
As a typical example, ${\cal S}$ would be the 
Schwartz space ${\cal S}(\R^N)$ of $C^\infty$ 
functions $f:\R^N\to\R$ whose derivatives of arbitrary order decrease 
at infinity faster than all inverse polynomials; but it might be a much
larger space in case of need.

For a given linear operator $S$, I shall denote by $\|S\|$ its operator norm:
\[ \|S\| = \sup_{h\neq 0} \frac{\|Sh\|}{\|h\|} = \sup_{\|h\|,\|h'\|\leq 1} 
\ip{Sh}{h'}.\]
If there is need to emphasize that $S$ is considered as a linear operator between
two spaces $\H_1$ and $\H_2$, the symbol $\|S\|$ may be replaced by
$\|S\|_{\H_1\to\H_2}$.

The norm $\|A\|$ of an array of operators $(A_1,\ldots,A_m)$ is
defined as $\sqrt{\sum_i \|A_i\|^2}$; the norm of a matrix-valued
operator $(A_{jk})$ by $\sqrt{\sum \|A_{jk}\|^2}$; etc.

The identity operator $X\to X$, viewed as a linear mapping, will always be
denoted by $I$, whatever its domain. Often a multiplication operator
(mapping a function $f$ to $fm$, where $m$ is a fixed function)
will be identified with the multiplicator $m$ itself.

Throughout the text, the real part will be denoted by $\Re$.

\subsection{Commutators}

In the sequel, commutators involving $A$ and $B$ will play a crucial role.
Since $A$ takes its values in $\H\otimes\V$ and $B$ is only defined in $\H$,
some notational convention should first be made precise, since 
$[A,B]$, for instance, does not a priori make sense.
I~shall resolve this issue by just tensorizing with the identity:
$[A,B]=AB-(B\otimes I)A$ is an unbounded operator $\H\to\H^\V$.
In a more pedestrian writing, $[A,B]$ is the row of operators
$([A_1,B],\ldots,[A_m,B])$. Then $A^2$ stands for the matrix of operators
$(A_jA_k)_{j,k}$, $[A,[A,B]]$ for $([A_j,[A_k,B]])_{j,k}$, etc.
One should be careful about matrix operations made on the
components: For instance, $[A,A^*]$ stands for $([A_j,A_k^*]_{j,k})$,
which is an operator $\H\to\H\otimes\V\otimes\V$, while
while $[A^*,A]$ stands for $\sum_j [A_j^*,A_j]$, which is an operator
$\H\to\H$. Also note that $[A,A]$ stands for the array $([A_j,A_k])_{j,k}$,
and is therefore not necessary equal to~0. When there is a risk of confusion, 
I~shall make the notation more explicit. 

\subsection{Relative boundedness} \label{relbound}

Let $S$ and $T$ be two unbounded linear
operators on a Hilbert space $\H$, and let $\alpha\geq 0$; then the operator $S$
is said to be $\alpha$-bounded relatively to $T$ if $D(T)\subset D(S)$,
and
\[ \forall h\in D(S), \qquad \|Sh\|\leq \alpha \|Th\|;\]
or equivalently, $S^*\! S\leq \alpha T^*T$. If $S$ is $\alpha$-bounded with
respect to $T$ for some $\alpha\geq 0$, then $S$ is said to be bounded relatively
to $T$. This will be sometimes abbreviated into
\[ S \sle T. \]
Note that $S$ and $T$ need not take values in the same space.
Of course, boundedness relative to $I$ is just plain boundedness.

This notion can be generalized in an obvious way into relative boudedness
with respect to a family of operators: An operator $S$ is said to be $\alpha$-bounded 
relatively to $T_1,\ldots,T_k$ if $\cap D(T_j)\subset D(S)$, and
\[ \forall h\in D(S), \qquad \|Sh\| \leq \alpha 
\bigl(\|T_1h\| + \ldots + \|T_kh\|\bigr). \]
If such an $\alpha$ exists, then $S$ is said to be bounded relatively
to $T_1,\ldots,T_k$, and this will naturally be abbreviated into
\[ S\sle T_1,\ldots, T_k. \]

\subsection{Abstract Sobolev spaces} 

The study of partial differential equations often relies on Sobolev spaces,
especially in a linear context. If one thinks of the Hilbert space $\H$ as
a (weighted) $L^2$ space, there is a natural abstract definition of ``Sobolev norm''
adapted to a given abstract coercive symmetric operator $L=A^*\! A$:
define the $\H^k$-Sobolev norm $\|\cdot\|_{\H^k}$ by
\[ \|h\|^2_{\H^k} := \|h\|^2 + \sum_{\ell=1}^k \|A^\ell h\|^2.\]

Here is a generalization: When some operators $C_0,\ldots, C_N$ are given 
(playing the same role as derivation operators along orthogonal directions in $\R^n$),
one can define
\begeq\label{H1C} 
\|h\|^2_{\H^1}:= \|h\|^2 + \sum_{j=0}^{N}\|C_j h\|^2,\qquad
\|h\|^2_{\H^k}:=\|h\|^2 + \sum_{\ell=0}^k \sum_{j=0}^{N} \|(C_j)^\ell h\|^2. 
\endeq
Of course, there is an associated scalar product, which will be denoted by
$\ip{\cdot}{\cdot}_{\H^1}$, or $\ip{\cdot}{\cdot}_{\H^k}$.

%A ``homogeneous'' Sobolev semi-norm can also be defined, for instance
%\begeq\label{Hdot} 
%\|h\|^2_{\dot{\H}^1}:= \sum_{j=0}^{N}\|C_j h\|^2.
%\endeq

\subsection{Calculus in $\R^n$}

Most of the examples discussed below take place in $\R^n$; then I shall use standard
notation from differential calculus: $\nabla$ stands for the gradient operator, and
$\nabla\cdot$ for its adjoint in $L^2(\R^n)$, which is the divergence operator.

\begin{Ex} Let $x=(x_1,\ldots,x_n)$ and $v=(v_1,\ldots,v_n)$ stand for two
variables in $\R^n$. Let $A=\nabla_v$, then $\nabla_v^2$ is the usual Hessian
operator with respect to the $v$ variable, which can be identified with
the matrix of second-order differential operators $(\partial^2/\partial v_j\partial v_k)$
($j,k \in \{1,\ldots,n\}$). Similary, if $a$ and $b$ are smooth scalar
functions, then $[a\nabla_v, b\nabla_x]$ is the matrix of
differential operators $[a\pa_{v_j}, b\pa_{x_k}]$.
\end{Ex}

The scalar product of two vectors $a$ and $b$ in $\R^n$ or $\Cplx^n$ will be denoted either 
by $\ip{a}{b}$ or by $a\cdot b$. The norm of a vector $a$ in $\R^n$ or $\Cplx^n$
will be denoted simply by $|a|$, and the Hilbert-Schmidt norm of an $n\times n$
matrix $M$ (with real or complex entries) by $|M|$. 

The usual Brownian process in $\R^n$ will be denoted by $(B_t)_{t\geq 0}$.

The notation $H^k$ will stand for the usual Sobolev space in $\R^n$: explicitly,
$\|u\|^2_{H^k} = \sum_{j\leq k} \|\nabla^j u\|_{L^2}^2$. Sometimes I~shall use
subscripts to emphasize that the gradient is taken only with respect to certain variables;
and sometimes I~shall put a reference measure if the reference measure is not
the Lebesgue measure. For instance,
$\|u\|^2_{H^1_v(\mu)} = \|u\|^2_{L^2(\mu)} + \|\nabla_v u\|^2_{L^2(\mu)}$.

\section{Operators $L=A^*\! A+B$}

For the moment we shall be concerned with linear operators of the form
\begeq\label{L} L:= A^*\! A + B, \qquad B^*=-B,\endeq
to be thought as the negative of the generator of a certain semigroup $(S_t)_{t\geq 0}$
of interest: $S_t = e^{-tL}$.
%
%Again, it is assumed that $D(L)$ is dense. 
%
(Of course, up to regularity issues, any linear operator $L$ with nonnegative
symmetric part can be written in the form~\eqref{L}; 
but this will be interesting only if $A$ and $B$ are ``simple enough''.) 
Here below I~have gathered some properties of $L$ which can be expressed
quite simply in terms of $A$ and $B$.

\subsection{Dirichlet form and kernel of $L$}

Introduce
\[ \K:=\ker L, \qquad \Pi:= \text{orthogonal projection on $\K$}, \qquad
\Pi^\bot = I-\Pi.\]

\begin{Prop} \label{lemkerL} With the above notation,

(i) $\dps \forall h\in D(A^*\! A)\cap D(B), \quad \Re\,\ip{Lh}{h} = \|Ah\|^2$;

(ii) $\dps \K = \ker A \cap \ker B$.
\end{Prop}

\begin{proof} The proof of (i) follows at once from the identities
\[ \ip{A^*\! Ah}{h} = \ip{Ah}{Ah} = \|Ah\|^2,\qquad 
\Re\,\ip{Bh}{h}=0.\]

It is clear that $\ker A\cap \ker B \subset \K$. Conversely,
if $h$ belongs to $\K$, then $0 = \Re\ip{Lh}{h} = \|Ah\|^2$,
so $h\in \ker A$, and then $Bh=Lh-A^*\! Ah=0$. This concludes the proof of (ii).
\end{proof}

\subsection{Nonexpansivity of the semigroup}

Now it is assumed that one can define a semigroup $(e^{-tL})_{t\geq 0}$,
i.e. a mapping $(t,h)\longmapsto e^{-tL}h$, continuous as a function of both
$t$ and $h$, satisfying the usual rules $e^{0L}=\Id$, $e^{-(t+s)L}=e^{-tL}e^{-sL}$ for $t,s\geq 0$
(semigroup property), and
\[ \forall h\in D(L),\qquad
\dr{t}{t=0^+} e^{-tL}h = - Lh.\]
As an immediate consequence, for all $h\in D(A^*\! A)\cap D(B)$,
\[ \left. \frac12\, \frac{d}{dt}\right |_{t=0^+} \|e^{-tL}h\|^2 = -\Re\,\ip{Lh}{h} = 
-\|Ah\|^2\leq 0.\]
This, together with the semigroup property, the continuity of the semigroup and the
density of the domain, implies that the semigroup is
nonexpansive, i.e. its operator norm at any time is bounded by~1:
\[ \forall t\geq 0\qquad \|e^{-tL}\|_{\H\to\H} \leq 1.\]

\subsection{Derivations in $L^2(\mu)$}

In most examples considered later, the Hilbert space $\H$ takes the form
$L^2(\mu_\infty)$, for some equilibrium measure $\mu_\infty=\rho_\infty(x)\,dx$ 
on $\R^n$, with density $\rho_\infty$ with respect to Lebesgue measure;
$\V=\R^m$, $A=(A_1,\ldots, A_m)$, and the $A_j$'s and $B$ are derivations
on $\R^n$, i.e. there are vector fields $a_j(x)$ and $b(x)$ on $\R^n$ such that
\[ A_j h = a_j\cdot \nabla h, \qquad Bh = b\cdot \nabla h.\]
To write things symbolically, there is an $m\times n$ matrix $\sigma=\sigma(x)$
such that
\[ A = \sigma\nabla.\]

Below are some useful calculation rules in that context.
It will be assumed that everything is smooth enough: For instance
$\rho_\infty$ lies in $C^2(\R^n)$ and it is positive everywhere; 
and $\sigma,b$ are $C^1$. The notation $\sigma^*$ will denote 
the transpose (adjoint) of $\sigma$.

\begin{Prop} \label{propcalc} With the above notation and assumptions,
\sm

(i) $\dps B^*=-B\Longleftrightarrow \nabla\cdot (b \rho_\infty) =0$;
\sm

(ii) $\dps A^*g = -\nabla\cdot(\sigma^* g) - \bigl\<\nabla\log\rho_\infty, \sigma^* g\>$.
\end{Prop}

\begin{Rk} As a consequence of Proposition~\ref{propcalc}(ii), 
the linear second-order operator $-L = \sum A_j^2 - (B+\sum c_jA_j)$ 
has the H\"ormander form (a sum of squares of
derivations, plus a derivation). The form $A^*\! A +B$ is however
much more convenient for the purpose of the present study
--- just as in~\cite{eckmannhairer:hypo:03}. 
\end{Rk}

\begin{proof}[Proof of Proposition~\ref{propcalc}]
By polarization, the antisymmetry of $B$ is equivalent to
\[ \forall h\in \H, \qquad \ip{Bh}{h} =0. \]
But
\begin{align} \ip{Bh}{h} = \int_{\R^n} (b\cdot\nabla h)h\,\rho_\infty 
& = \frac12 \int_{\R^n} b\cdot\nabla (h^2)\rho_\infty \nonumber \\
& = - \frac12 \int_{\R^n} h^2 \nabla\cdot (b\rho_\infty).
\label{divbinf}
\end{align}
If $\nabla\cdot(b\rho_\infty)=0$, then the integral 
in~\eqref{divbinf} vanishes. If on the other 
hand $\nabla\cdot(b\rho_\infty)$ is not identically zero, one can find some $h$
such that this integral is nonzero. This proves statement (i).
\sm

To prove (ii), let $g:\R^n\to\R^m$ and $h:\R^n\to\R$, 
then $\ip{A^*g}{h}$ coincides with
\[ \ip{g}{Ah} = \int_{\R^n} g\cdot (\sigma\nabla h)\,\rho_\infty =
- \int (\sigma^* g)\cdot\nabla h \,\rho_\infty = 
- \int_{\R^n} \nabla\cdot(\sigma^* g\rho_\infty) h \]
\[ = - \int_{\R^n} \nabla\cdot (\sigma^* g)h\,\rho_\infty 
- \int_{\R^n} \sigma^* g\cdot (\nabla\log\rho_\infty) h\,\rho_\infty,\]
where the identity $\nabla\rho_\infty = (\nabla\log\rho_\infty)\rho_\infty$
was used. This proves (ii).
\end{proof}

The following proposition deals with the range of applicability for diffusion
processes.

\begin{Prop}\label{rulesdiffeq} 
Let $\sigma\in C^2(\R^n;\R^{n\times m})$ and $\xi\in C^1(\R^n;\R^n)$,
and let $(X_t)_{t\geq 0}$ be a stochastic process solving the autonomous
stochastic differential equation
\[ dX_t = \sqrt{2}\,\sigma(X_t)\,dB_t + \xi(X_t)\,dt, \]
where $(B_t)_{t\geq 0}$ is a standard Brownian motion in $\R^m$.
Then
\sm

(i) The law $(\rho_t)_{t\geq 0}$ of $X_t$ satisfies the diffusion equation
\begeq\label{FP}
\derpar{\rho}{t} = \nabla\cdot (D\nabla \rho - \xi\rho), \qquad
D:=\sigma^*\sigma;
\endeq
\sm

(ii) Assume that the equation~\eqref{FP} admits an invariant measure
$\mu_\infty(dx)=\rho_\infty(x)\,dx$ (with finite or infinite mass),
where $\rho_\infty$ lies in $C^2(\R^n)$ and is positive everywhere. 
Then the new unknown $h(t,x):=\rho(t,x)/\rho_\infty(x)$ 
solves the diffusion equation
\begeq\label{eqh} 
\derpar{h}{t} = \nabla\cdot (D\nabla h) - 
\Bigl(\xi-2D\nabla\log\rho_\infty\Bigr)\cdot\nabla h,
\endeq
which is of the form $\pa_t h + Lh =0$ with $L=A^*\! A+B$, $B^*=-B$, if one defines
\begeq\label{ABdiff} 
\H:=L^2(\mu_\infty); \qquad A := \sigma\nabla; \qquad 
B := \bigl(\xi-D\nabla\log\rho_\infty\bigr)\cdot\nabla. 
\endeq
\end{Prop}

\begin{proof}
Claim (i) is a classical consequence of It\^o's formula. To prove claim (ii), write
\begin{align*} \derpar{h}{t} & = \frac{1}{\rho_\infty} \nabla\cdot 
\Bigl(D\rho_\infty\nabla h + D h \nabla\rho_\infty - \xi\rho_\infty h\Bigr) \\
& = \nabla\cdot (D\nabla h) + 2 D\nabla h\cdot\frac{\nabla\rho_\infty}{\rho_\infty}
- \xi\cdot\nabla h + \frac{h}{\rho_\infty} 
\Bigl[ \nabla\cdot(D\nabla\rho_\infty) - \nabla\cdot (\rho_\infty\xi) \Bigr].
\end{align*}
As $\rho_\infty$ is a stationary solution of~\eqref{FP}, the last term in
square brackets vanishes, which leads to~\eqref{eqh}.
Define $A$ and $B$ by~\eqref{ABdiff}. Thanks to Proposition~\ref{propcalc} (ii), 
it is easy to check that
\[ A^*\! Ah = - \nabla\cdot (D\nabla h) - D \nabla \log\rho_\infty\cdot\nabla h,\]
so $h$ indeed satisfies $\pa_t h + Lh=0$. It only remains to check that $B^*=-B$.
By Proposition~\ref{propcalc} (i), it is sufficient to check that
\[ \nabla\cdot (\rho_\infty(\xi-D\nabla\log\rho_\infty)) =
 -\nabla\cdot (D\nabla\rho_\infty - \xi\rho_\infty)\]
vanishes; but this follows again from the stationarity of $\rho_\infty$.
So the proof of Proposition~\ref{rulesdiffeq} is complete.
\end{proof}

\subsection{Example: The kinetic Fokker--Planck equation}

The following example will serve as an important application and model.
Consider a nice (at least $C^1$) 
function $V:\R^n\to\R$, converging to $+\infty$ fast enough
at infinity (say $V(x)\geq K|x|^\alpha - C$ for some positive constants
$K$ and $C$). For $x,v\in \R^n\times\R^n$, set
\[ f_\infty(x,v) := \frac{e^{-[V(x) + \frac{|v|^2}2]}}{Z}, \qquad
\mu(dx\,dv) = f_\infty(x,v)\,dx\,dv,\]
where $Z$ is chosen in such a way that $\mu$ is a probability measure.
Define
\[ \H := L^2(\mu), \quad \V:= \R^n_v, \quad A:=\nabla_v, 
\quad B:= v\cdot\nabla_x - \nabla V(x)\cdot\nabla_v,\]
\[ L:= -\Delta_v + v\cdot\nabla_v + v\cdot\nabla_x -\nabla V(x)\cdot\nabla_v.\]
The associated equation is the kinetic Fokker--Planck equation with confinement
potential $V$, in the form
\begeq\label{kFPh} 
\partial_t h + v\cdot\nabla_x h -\nabla V(x)\cdot\nabla_v h = \Delta_v h - v\cdot\nabla_v h.
\endeq

Before considering convergence to equilibrium for this model, one should first
solve analytical issues about regularity and well-posedness.
It is shown by Helffer and Nier~\cite[Section~5.2]{helfnier:witten:05} that~\eqref{kFPh} generates a
$C^\infty$ regularizing contraction semigroup in $L^2(\mu)$ as soon as $V$ itself
lies in $C^\infty(\R^n)$. To study this equation for a less regular potential $V$, 
it is always possible to regularize $V$ into a smooth approximation $V_\var$, then
perform all a priori estimates on the regularized problem, and finally pass
to the limit as $\var\to 0$. The following well-posedness theorem justifies this
procedure by forcing the convergence of the approximate solutions to the original
solution.

\begin{Thm}\label{thmuniqueFP}
Let $V\in C^1(\R^n)$, $\inf V>-\infty$, and let 
\[ E(x,v):=V(x) + \frac{|v|^2}2, \qquad \rho_\infty= e^{-E}, \qquad
\mu(dx\,dv) = \rho_\infty(x,v)\,dv\,dx.\]
Then, for all $h_0\in L^2(\mu)$, equation~\eqref{kFPh} admits a unique
distributional solution $h=h(t,x,v) \in C(\R_+;{\cal D}'(\R^n_x\times\R^n_v))
\cap L^\infty_{\loc}(\R_+;L^2(\mu))\cap L^2_\loc (\R_+; H^1_v(\mu))$,
such that $h(0,\cdot)=h_0$.
\end{Thm}

The proof of existence is a straightforward consequence of a standard approximation
procedure, the Helffer--Nier existence results, and the a priori estimate
\begin{multline*} 
\int h^2(t,x,v)\,d\mu(x,v) 
+ \int_0^t \int h^2(s,x,v)\,d\mu(x,v)\,ds \\
\int h^2(0,x,v)\,d\mu(x,v).
\end{multline*}
There is more to say about the uniqueness statement, of which the proof is 
deferred to Appendix~\ref{appuniqueFP}.
The main subtlety lies in the absence of any growth condition on $\nabla V$;
this is overcome by a localization argument inspired 
from~\cite[Proposition~5.5]{helfnier:witten:05}. Apart from that, Theorem~\ref{thmuniqueFP} 
is just an exercise in linear partial differential equations.

Many people (including me) would rather think of~\eqref{kFPh} in the form
\begeq\label{kFP} 
\partial_t f + v\cdot\nabla_x f - \nabla V(x)\cdot\nabla_v f = \Delta_v f + \nabla_v\cdot(vf),
\endeq
in which case $f$ at time $t$ can be interpreted (if it is nonnegative) as a density of 
particles, or (if it is a probability density) as the law of a random variable 
in phase space. To switch from~\eqref{kFPh} to~\eqref{kFP} it suffices to set 
$f:=f_\infty h$. This however does not completely solve the problem
because the natural assumptions for~\eqref{kFP} are much more
general than for~\eqref{kFPh}. For instance, it is natural to assume that the initial
datum $f_0$ for~\eqref{kFP} is $L^2$ with {\em polynomial} weight; 
or just $L^1$, or even a finite measure. 
Theorem~\ref{thmuniqueFPf} below yields a uniqueness result in such a
setting, however with more stringent assumptions on the initial datum.
In the next statement, $M(\R^n\times\R^n)$ stands for the space of finite measures on 
$\R^n\times\R^n$, equipped with the topology of weak convergence (against bounded
continuous functions).

\begin{Thm}\label{thmuniqueFPf}
Let $V\in C^1(\R^n)$, $\inf V>-\infty$, and let  $E(x,v):=V(x) + \frac{|v|^2}2$.
Then, for any $f_0\in L^2((1+E)\,dx\,dv)$, equation~\eqref{kFP} admits a unique
distributional solution $f=f(t,x,v) \in C(\R_+;{\cal D}'(\R^n_x\times\R^n_v))
\cap L^\infty_{\loc}(\R_+;L^2((1+E)\,dx\,dv))\cap L^2_\loc (\R_+; H^1_v(\R^n_x\times
\R^n_v))$, such that $f(0,\cdot)=f_0$.

If moreover $\nabla^2 V$ is uniformly bounded, then for all finite measure $f_0$
the equation~\eqref{kFP} admits a unique solution 
$f=f(t,x,v)\in C(\R_+; M(\R^n_x\times\R^n_v))$.
\end{Thm}

The proof of this theorem will be deferred to Appendix~\ref{appuniqueFP}.

\section{Coercivity and hypocoercivity}

\subsection{Coercivity}

\begin{Def} \label{defco}
Let $L$ be an unbounded operator on a Hilbert space $\H$, with kernel $\K$,
and let $\tilde{\H}$ be another Hilbert space
continuously and densely embedded in ${\cal K}^\bot$, endowed with a scalar product
$\ip{\cdot}{\cdot}_{\tilde{\H}}$ and a Hilbertian norm $\|\cdot\|_{\tilde{\H}}$. 
The operator $L$ is said to be $\lambda$-coercive on $\tilde{\H}$ if 
\[ \forall h\in \K^\bot\cap D(L), \quad
\Re \ip{Lh}{h}_{\tilde{\H}} \geq \lambda \|h\|^2_{\tilde{\H}},\]
where $\Re$ stands for real part. The operator $L$ is said to be coercive 
on $\tilde{\H}$ if it is $\lambda$-coercive on $\tilde{\H}$ for some $\lambda>0$.
\end{Def}

The most standard situation is when $\tilde{\H}=\K^\bot\simeq\H/\K$. 
Then it is equivalent to say that $L$ is coercive on $\K^\bot$
(which will be abbreviated into just: $L$ is coercive),
or that the symmetric part of $L$ admits a {\em spectral gap}.
\med

Coercivity properties can classically be read at the level of the semigroup
(assuming it is well-defined), as shown by the next statement:

\begin{Prop} With the same notation as in Definition~\ref{defco},
$L$ is $\lambda$-coercive on $\tilde{\H}$ \iff 
$\dps \|e^{-tL}h_0\|_{\tilde{\H}} \leq e^{-\lambda t} \|h_0\|_{\tilde{\H}}$ 
for all $h_0\in\tilde{\H}$ and $t\geq 0$.
\end{Prop}

\begin{proof} Assume by density that $h_0\in \tilde{\H}\cap D(L)$.
On one hand the coercivity implies
\[\left. \frac{d}{dt}\right|_{t=0^+} 
\|e^{-tL} h_0\|_{\tilde{\H}}^2 = -2 \Re\, \ip{L e^{-tL}h_0}{e^{-tL}h_0}
\leq -2 \lambda\, \|e^{-tL}h_0\|^2, \]
so by Gronwall's lemma
\[ \|e^{-tL} h_0\|_{\tilde{\H}}^2 \leq e^{-2\lambda t}\|h_0\|_{\tilde{\H}}^2.\]

Conversely, if exponential decay holds, then for any 
$h_0\in \tilde{\H}\cap D(L)$,
\begin{multline*} \Re\, \ip{Lh_0}{h_0} = 
\lim_{t\to 0} \frac{\|h_0\|_{\tilde{\H}}^2 - \|e^{-tL}h_0\|_{\tilde{\H}}^2}{2t}
\\
\geq \liminf_{t\to 0} \frac{(1-e^{-2\lambda t})\|h_0\|_{\tilde{\H}}^2}{2t} 
= \lambda \|h_0\|^2,
\end{multline*}
whence the coercivity.
\end{proof}

When an operator $L$ is in the form~\eqref{L}, the coercivity of $L$ 
follows from the coercivity of $A^*\! A$, at least if $B$ has a sufficiently large kernel:

\begin{Prop} \label{largekernel}
With the notation of Subsection~\ref{subnot}, if $A^*\! A$ is $\lambda$-coercive 
on $(\ker A)^\bot$ and $\ker A\subset \ker B$, then $L$ is $\lambda$-coercive on $\K^\bot$.
\end{Prop}

\begin{proof} We know that $\K=\ker(A)\cap \ker(B)=\ker(A)$, so
for any $h\in \K^\bot$, $\ip{Lh}{h} = \|Ah\|^2 \geq \lambda \|h\|^2$. \end{proof}

\begin{Ex} Apart from trivial examples where $B=0$, one can consider the
following operator from~\cite{AMTU:FP:01}:
\[ L = -(\Delta_x -x\cdot\nabla_x) - (\Delta_v -v\cdot\nabla_v) + 
(v\cdot\nabla_x - x\cdot\nabla_v)\]
on $L^2(e^{-( |v|^2+|x|^2)/2}\,dx\,dv)$.
\end{Ex}

The main problem in the sequel is to study cases in which $A^*\! A$ is coercive, but
$L$ is {\em not}, and yet there is exponential convergence to equilibrium
(i.e. to an element of $\K$) for the semigroup $(e^{-tL})$.
In view of Proposition~\ref{largekernel}, this can only happen if {\em $\ker L$
is smaller than $\ker A$}. Here is the most typical example: With the choice
$\H=L^2(\exp(-(|v|^2+|x|^2)/2)\,dv\,dx)$ again, consider
\[ L = - (\Delta_v - v\cdot\nabla_v) + (v\cdot\nabla_x - x\cdot\nabla_v). \]
Then $\ker A$ is made of functions which depend only on $x$, but
$\ker L$ only contains constants.

\subsection{Hypocoercivity}

To fix ideas, here is a (possibly misleading, but at least precise)
definition of ``hypocoercivity'' in a Hilbertian context.

\begin{Def} \label{defhypoco}
Let $\H$ be a Hilbert space, $L$ an unbounded operator on $\H$ generating
a continuous semigroup $(e^{-tL})_{t\geq 0}$, and $\tilde{\H}$ another Hilbert space,
continuously and densely embedded in $\K^\bot$, endowed with a Hilbertian norm 
$\|\cdot\|_{\tilde{\H}}$. The operator $L$ is said to be $\lambda$-hypocoercive 
on $\tilde{\H}$ if there exists a finite constant $C$ such that
\begeq\label{hypoco} 
\forall h_0 \in \tilde{\H}, \quad \forall t\geq 0 \quad
\|e^{-tL}h_0\|_{\tilde{\H}} \leq C e^{-\lambda t} \|h_0\|_{\tilde{\H}}.
\endeq
It is said to be hypocoercive on $\tilde{\H}$ if it is $\lambda$-hypocoercive on
$\tilde{\H}$ for some $\lambda>0$.
\end{Def}

\begin{Rk} \label{rknorm1}
With respect to the definition of coercivity in terms of semigroups,
the only difference lies in the appearance of the constant $C$ in
the \rhs~of~\eqref{hypoco} (obviously $C\geq 1$, apart from trivial cases; 
$C=1$ would mean coercivity). The difference between Definition~\ref{defco}
and Definition~\ref{defhypoco} seems to be all the thinner in view
of the following fact (pointed to me by Serre): Whenever one has
a norm satisfying inequality~\eqref{hypoco} for some constant $C$, 
it is always possible to find an equivalent norm (in general, not Hilbertian)
for which the same inequality holds true with $C=1$. Indeed, just choose
\[ N(h) := \sup_{t\geq 0} \Bigl(e^{\lambda t}\, \|e^{-tL}h\|\Bigr). \]
In spite of these remarks, hypocoercivity is a strictly weaker concept
that coercivity. In particular, hypocoercivity is invariant under change of
equivalent Hilbert norm on $\tilde{\H}$, while coercivity is not.
This has an important practical consequence: If one finds an
equivalent norm for which the operator $L$ is coercive, 
then it follows that it is hypocoercive. 
I~shall systematically use this strategy in the sequel.
\end{Rk}

\begin{Rk}\label{rknorm2}
It often happens that a certain space $\tilde{\H}$ is 
convenient for proving hypocoercivity, but this particular space is
much smaller than $\H$ (stated otherwise, the Hilbert norm on $\tilde{\H}$ cannot
be bounded in terms of the Hilbert norm on $\H$): typically, $\tilde{\H}$ may
be a weighted Sobolev space, while $\H$ is a weighted $L^2$ space. 
In that situation there is in general no density argument which would allow one 
to go directly from hypocoercivity on $\tilde{\H}$, to hypocoercivity on $\H$.
However, such an extension is possible if $L$ satisfies a (hypoelliptic) 
regularization estimate of the form
\begeq\label{qdmeme1} \forall t_0>0\quad \exists C(t_0) <+\infty;\quad
\forall t\geq t_0, \ \|e^{-tL}\|_{\K^\bot \to \tilde{\H}} \leq C(t_0); 
\endeq
or, more generally, if $L$ generates a semigroup for which there is 
{\em exponential decay of singularities}:
\begeq\label{qdmeme2} \begin{cases} \forall t\geq 0,\quad e^{-tL} = S_t + R_t, \\ \\
\forall t_0>0\quad \exists C(t_0)<+\infty;\quad \forall t\geq t_0, 
\ \|S_t\|_{\K^\bot\to \tilde{\H}} \leq C(t_0); \\ \\
\exists \lambda>0; \quad \forall t\geq 0,
\ \|R_t\|_{\K^\bot\to \K^\bot} \leq C e^{-\lambda t}.
\end{cases}
\endeq
Such assumptions are often satisfied in realistic models.
For instance, integral operators (generators of jump processes) usually 
satisfy~\eqref{qdmeme2} when the kernel is integrable (finite jump measure), 
and~\eqref{qdmeme1} when the kernel is not integrable. Diffusion operators of heat or 
Fokker--Planck type usually satisfy~\eqref{qdmeme1}.
\end{Rk}

\subsection{Commutators}

If the operators $A^*\! A$ and $B$ commute, then so do their exponentials,
and $e^{-tL}=e^{-tA^*\! A} e^{-tB}$.
Then, since $B$ is antisymmetric, $e^{-tB}$ is norm-preserving, and it is 
equivalent to study the convergence for $e^{-tL}$ or for $e^{-tA^*\! A}$.
On the other hand, if these operators do {\em not} commute, one can 
hope for interesting phenomena.

\begin{Prop} \label{propkercom} With the notation of Subsection~\ref{subnot}, in particular
$L=A^*\! A+B$, define recursively the iterated commutators
\[ C_0:=A, \quad C_k:= [C_{k-1},B], \]
and then $\K':= \cap_{k\geq 0} \ker C_k$. Then $\K\subset \K'$, and $\K'$ is
invariant for $e^{-tL}$.
\end{Prop}

\begin{proof} Assume that $\K\subset \ker C_0\cap\ldots\cap \ker C_j$.
Then, for all $h\in \K\cap {\cal S}$, 
\[ C_{k+1}h = C_k B h - B C_k h = C_k B h = - C_k A^*\! A h = 0.\]
Thus $\K\subset \ker C_{k+1}$. By induction, $\K$ is included in the
intersection $\K'$ of all $\ker C_j$. 

Next, if $h\in \K'\cap {\cal S}$, 
then $Lh=Bh$, so $C_k L h = C_k B h = C_{k+1} h + B C_k h =0$;
since $k$ is arbitrary, in fact $Lh\in \K'$, so $L$ leaves $\K'$ invariant,
and therefore so does $e^{-tL}$.
\end{proof}

In most cases of interest, not only does $\K'$ coincide with $\K$, but in addition
$\K'$ can be constructed as the intersection of just \emph{finitely many} kernels
of iterated commutators. Thanks to the trivial identity
\[\bigcap_{j=0}^k \ker C_j = \ker \left ( \sum_{j=0}^k C_j^*C_j\right ),\]
the condition that $\K'$ is the intersection of finitely many iterated
commutators may be reformulated as
\begeq\label{condcomm}
\text{There exists $\Nc\in\N$ such that}
\ker \left ( \sum_{j=0}^{\Nc} C_j^*C_j \right ) = \ker L.
\endeq

\begin{Ex} For the kinetic Fokker--Planck operator~\eqref{kFPh},
$\Nc=1$ will do.\end{Ex}

If the goal is to derive estimates on the rate of convergence,
it is natural to reinforce the above condition into a more quantitative one:
\begeq \label{Horm}
\sum_{k=0}^{\Nc} C_k^* C_k \quad
\text{is coercive on $\K^\bot$}.
\endeq
Condition~\eqref{Horm} is more or less an analogue of H\"ormander's ``rank $r$''
bracket condition (as explained later, $r=2\Nc+1$ is the natural convention), 
but in the context of convergence to equilibrium and
spectral gap, rather than regularization and elliptic estimates.
There is however an important difference: Here we are taking brackets always
with $B$, while in H\"ormander's condition, brackets of the form, say,
$[A_i, A_j]$ would be allowed. This modification is intentional:
in all the cases of interest known to me, there is no need to consider such
brackets for hypocoercivity problems. A basic example which will be discussed
in Appendix~\ref{appproduct} is the following: The differential operator 
\[ L:= - (x^2{\pa_y}^*\pa_y + {\pa_x}^*\pa_x),\]
although {\em not} elliptic, is {\em coercive} (not just hypocoercive)
in $L^2(\gamma)/\R$, where $\gamma$ is the gaussian measure on $\R^2$.
For this operator, brackets of the form $[\pa_x,x\pa_y]$ play a crucial
role in the regularity study, but they are not needed to establish lower bounds
on the spectral gap.

\begin{Rk} It was pointed out to me by Serre that, when $\H$ is finite-dimensional, 
condition~\eqref{Horm} is equivalent to the statement that {\em $\ker A$ does not
contain any nontrivial subspace invariant by $B$}. In the study of convergence to
equilibrium for hyperbolic systems of conservation laws, this condition is known
as {\bf Kawashima's nondegeneracy condition}~\cite{kawashima:largetime:87,hanouzetnatalini:03,ruggeriserre:stability:04}.
It is not so surprising to note that the very same condition
appears in H\"ormander's seminal 1967 paper on hypoellipticity~\cite[p.~148]{hormander:67} 
as a necessary and sufficient condition for a diffusion equation to be hypoelliptic,
when the constant matrices $A$ and $B$ respectively stand for the diffusion matrix
and the linear drift function.\footnote{At first sight, 
it seems that both problems are completely different: Kawashima's condition
is applied to systems of unknowns, while H\"ormander's example deals with
scalar equations. The analogy becomes less surprising when one notices that
for such a diffusion equation the fundamental solution, viewed as a function
of time, takes its values in the finite-dimensional space of Gaussian
distributions, so that the equation really defines a system.}
\end{Rk}

Taking iterated commutators may rapidly lead to cumbersome expressions, 
because of ``lower-order terms''. In the present context, this might
be more annoying than in a regularity context, and so it will be 
convenient to allow for perturbations in the definition of $C_k$, say
\[ [C_k,B] = C_{k+1} + R_{k+1},\]
where $R_{k+1}$ is a ``remainder term'', chosen according to the context,
that is controlled by $C_0,\ldots,C_k$. An easy and sometimes useful
generalization is to set
\[ [C_k,B] = Z_{k+1}C_{k+1} + R_{k+1},\]
where the $Z_k$'s are auxiliary operators, typically multipliers, 
satisfying certain identities.

Once the family $(C_0,\ldots, C_{\Nc})$ is secured, one can introduce the
corresponding abstract Sobolev $\H^1$ norm as in~\eqref{H1C}.
This norm will be used on $\H$, or (more often) on $\K^\bot$. On the latter
space we may also consider ``homogeneous Sobolev norms'' such as
\begeq\label{Hdot} 
\|h\|^2_{\dot{\H}^1}:= \sum_{j=0}^{\Nc}\|C_j h\|^2.
\endeq

Note that, with the above assumptions,
the orthogonal space to the kernel $\K$ in $\H^1$ {\em does not depend}
on whether we consider the scalar product of $\H$ or that of $\H^1$.
(See the proof of Theorem~\ref{hypocomult} below.)
So a natural choice for $\tilde{\H}$ will be
\[ \tilde{\H} = \H^1/\K, \]
which is $\K^\bot$ equipped with the $\H^1$ norm.

%\begin{proof} First note that (obviously) $\K\subset\H^1$. Next, consider
%$h\in \K$ and $y\in\H$; the goal is to show that $\ip{h}{h'}=0$ if and only if
%$\ip{h}{h'}_{\H^1}=0$. By assumption $A h =0$, $Bh=0$, and so $h$ lies in the
%kernel of all iterated commutators involving $A$ and $B$. By induction, using
%the assumption $C_k = [C_{k-1},B]+R_k$, where $R_k$ is bounded relatively to
%$C_0,\ldots, C_{k-1}$, we see that $h$ lies in the kernel of all $C_k$'s and 
%all $R_k$'s. As a consequence, $\ip{h}{h'}_{\H^1}$ reduces to just 
%$\ip{h}{h'}$, and the lemma is proven.
%\end{proof}

\section{Basic theorem} \label{secbasic}

In this section linear operators satisfying a ``rank-3'' condition ($c=2$ in~\eqref{Horm})
are considered. Although this is a rather simple situation, it is already of interest,
and its understanding will be the key to more complicated extensions;
so I~shall spend some time on this case. Here it will be assumed for pedagogical
reasons that the operators $A$ and $C$ commute; this assumption will be relaxed in 
the next section.

\begin{Thm}\label{thmsimple}
With the notation of Subsection~\ref{subnot}, consider a linear operator
$L=A^*\! A+B$ ($B$ antisymmetric), and define $C:=[A,B]$. Assume the existence
of constants $\alpha,\beta$ such that
\sm

(i) $A$ and $A^*$ commute with $C$; $A$ commutes with $A$
(i.e. each $A_i$ commutes with each $A_j$);
\sm

(ii) $[A,A^*]$ is $\alpha$-bounded relatively to $I$ and $A$;
\sm

(iii) $[B,C]$ is $\beta$-bounded relatively to $A$, $A^2$, $C$ and $AC$;
\sm

Then there is a scalar product $\iip{\cdot}{\cdot}$ on $\H^1/\K$, which defines
a norm equivalent to the $\H^1$ norm, such that
\begeq\label{xLx}
\forall h\in \H^1/\K,\qquad \iip{h}{Lh} \geq K \bigl(\|Ah\|^2 + \|Ch\|^2\bigr) 
\endeq
for some constant $K>0$, only depending on $\alpha$ and $\beta$.
\sm

If, in addition,
\[ \text{$A^*\! A+C^*C$ is $\kappa$-coercive} \]
for some $\kappa>0$, then there is a constant $\lambda>0$, only depending on
$\alpha,\beta$ and $\kappa$, such that
\[ \forall h\in \H^1/\K, \qquad \iip{h}{Lh} \geq \lambda \iip{h}{h}. \]
In particular, $L$ is hypocoercive in $\H^1/\K$:
\[ \|e^{-tL}\|_{\H^1/\K \to \H^1/\K} \leq c\, e^{-\lambda t} \qquad
(c<+\infty),\]
where both $\lambda$ and $c$ can be estimated explicitly in terms of upper bounds
on $\alpha$ and $\beta$, and a lower bound on $\kappa$.
\end{Thm}

Before stating the proof of Theorem~\ref{thmsimple}, I~shall provide some remarks and
further explanations.

\begin{Rk} Up to changing $\alpha$ and $\beta$, it is equivalent to 
impose (ii) and (iii) above or to impose the seemingly more general conditions:

\quad (ii') $[A,A^*]$ is $\alpha$-bounded relatively to $I$, $A$ and $A^*$,

\quad (iii') $[B,C]$ is $\beta$-bounded relatively to $A$, $A^2$, $A^*\! A$, $C$ and $AC$;

Indeed,
\[ \ip{A^*h}{A^*h} = \ip{AA^*h}{h} = \ip{A^*\! Ah}{h} + 
\sum_i\ip{[A_i,A_i^*]h}{h}, \]
so 
\[ \|A^*h\|^2 \leq \|Ah\|^2 + \|[A,A^*]h\|\,\|h\|. \]
Then assumption (ii') implies
\begin{align*}\|A^*h\|^2 & \leq \|Ah\|^2 + 
\alpha\Bigl(\|h\|^2+ \|Ah\|\,\|h\| + \|A^*h\|\,\|h\|\Bigr) \\
& \leq \|Ah\|^2 + \alpha\Bigl(\|h\|^2+ \|Ah\|\,\|h\|\Bigr) +
\frac12 \|A^*h\|^2 + \frac{\alpha^2}{2} \|h\|^2, 
\end{align*}
and then it follows that $A^*$ is bounded relatively to $I$ and $A$,
so that (ii) holds true. This also implies that $A^*\! A$ is bounded 
relatively to $A^2$ and $A$, so (iii') implies (iii). 
\end{Rk}

\begin{Rk}\label{rkii''} Assumption (ii) in Theorem~\ref{thmsimple} can be
relaxed into

\quad (ii'') $[A,A^*]A$ is relatively bounded with respect to 
$A$ and $A^2$; where by convention
\[ \bigl\| [A,A^*]Ah\bigr\|^2 = \sum_i \Bigl\|
\sum_j [A_i, A_j^*]A_jh\Bigr\|^2.\]
\end{Rk}

\begin{Rk} \label{rkweight}
Here is a crude heuristic rule explaining a bit the assumptions (i) to (iii)
above. As is classical in H\"ormander's theory, define the weights $w(O)$ of
the operators involved, by
\[ w(A)=w(A^*)=1, \qquad w(B)=2, \qquad w([O_1,O_2]) = w(O_1) + w(O_2). \]
Then rules (i) to (iii) guarantee that certain key commutators can be estimated in
terms of operators whose order is strictly less: for instance, the weight of $[B,C]$ is
$2+3=5$, and assumption (iii) states that it should be controlled by some operators,
for which the maximal weight is~4. (This rule does not however explain why $I$ is
allowed in the right-hand side of (ii), but not in (iii);
so it might be better to think in terms of Assumption (ii'')
from Remark~\ref{rkii''} rather than in terms of Assumption (ii).)
\end{Rk}

\begin{Rk} In particular cases of interest, it may be a good idea to
rewrite the proof of Theorem~\ref{thmsimple}, 
taking into account specific features of
the problem considered, so as to obtain better constants $\lambda$ and $C$. 
\end{Rk}

\subsection{Heuristics and strategy}

The proof of Theorem~\ref{thmsimple} is quite elementary; in some sense, the most 
sophisticated analytical tool on which it rests is the Cauchy--Schwarz inequality.
The argument consists in devising an appropriate Hilbertian norm on $\H^1/\K$, 
which will be equivalent to the usual norm, but will turn $L$ into 
a coercive operator. One can see an analogy with a classical,
elementary proof of a standard theorem in linear 
algebra~\cite[pp.~147-148]{arnold:EDO:mir}: If the real parts of the
eigenvalues of a matrix $M$ are all positive, then
$e^{-tM}\to 0$ (exponentially fast) as $t\to\infty$.

Define
\begeq\label{eqnorm}
\iip{h}{h} = \|h\|^2 + a\, \|Ah\|^2 + 2b\,\Re\, \ip{Ah}{Ch} + c\, \|Ch\|^2, 
\endeq
where the positive constants $a,b,c$ will be chosen later on,
in such a way that $1\gg a \gg b \gg c$. (The constant $c$ here is
not the same as the one in the conclusion of Theorem~\ref{thmsimple}.)

By polarization, this formula defines a bilinear symmetric form on $\H^1$.
By using Young's inequality, in the form
\[ \Bigl| 2b\, \ip{Ah}{Ch} \Bigr| \leq 2b\, \|Ah\|\,\|Ch\|
\leq b \sqrt{\frac{a}{c}} \|Ah\|^2 + b \sqrt{\frac{c}{a}}\,\|Ch\|^2,\]
one sees that the scalar products $\iip{\cdot}{\cdot}$ and $\ip{\cdot}{\cdot}_{\H^1}$ 
define equivalent norms as soon as $b<\sqrt{ac}$, and more precisely
\begeq\label{equivalence} 
\min(1,a,c) \left ( 1 - \frac{b}{\sqrt{ac}} \right ) \|h\|^2_{\H^1} \leq \iip{h}{h}\leq
\max(1,a,c) \left ( 1 + \frac{b}{\sqrt{ac}} \right ) \|h\|^2_{\H^1}. 
\endeq
In particular, the scalar products $\iip{\cdot}{\cdot}$ and $\ip{\cdot}{\cdot}_{\H^1}$
define equivalent norms. 

In spite of their equivalence, the scalar products $\iip{\cdot}{\cdot}$ 
and $\ip{\cdot}{\cdot}_{\H^1}$ are quite different: 
it is possible to arrange that $L$ is coercive \wrt the former, although
it is not \wrt the latter. Heuristically, one may say that the
``pure'' terms $\|h\|^2$, $\|Ah\|^2$ and $\|Ch\|^2$ will mainly feel
the influence of the symmetric part in $L$, but that the ``mixed''
term $\ip{Ah}{Ch}$ will mainly feel the influence of 
the antisymmetric part in $L$.
The following simple calculations should help understanding this. 
Whenever $Q$ is a linear operator commuting with $A$ 
(be it $I$, $A$ or $C$ in this example),
\[ \dr{t}{t=0} \|Qe^{-tA^*\! A}h\|^2 = - 2\|QAh\|^2, \]
but on the other hand
\[ \dr{t}{t=0} \ip{Ae^{-tB}h}{Ce^{-tB}h} = 
- \ip{ABh}{Ch} - \ip{Ah}{CBh}.\]
Pretend that $B$ and $C$ commute, and this can be rewritten
\begin{align*} -\ip{ABh}{Ch} - \ip{Ah}{BCh} & = -\ip{ABh}{Ch} -\ip{B^*Ah}{Ch} \\
& = -\ip{ABh}{Ch} +\ip{BAh}{Ch} \\
& = -\ip{[A,B]h}{Ch} = -\|Ch\|^2,
\end{align*}
where the antisymmetry of $B$ has been used to go from the first to the second line.
This will yield the dissipation in the $C$ direction, which the symmetric part
of $A$ was unable to provide!

\subsection{Proof of Theorem~\ref{thmsimple}} \label{proofsimple}

Introduce the norm~\eqref{eqnorm}. By Proposition~\ref{propkercom}, any $h\in \K=\ker L$ 
satisfies $Ah=0$, $Ch=0$, in which case $\iip{h}{h'}=\ip{h}{h'}_{\H^1}=\ip{h}{h'}$.
In particular, the orthogonal space $\K^\bot$ is the same for these three scalar products.
So it makes sense to choose $\tilde{\H}=\H^1/\K$.

Let us compute
\[ - \frac12\, \frac{d}{dt} \iip{e^{-tL}h}{e^{-tL}h} = \Re\,\iip{e^{-tL}h}{Le^{-tL}h};\]
if we can bound below this time-derivative by a constant multiple
of $\iip{e^{-tL}h}{e^{-tL}h}$, then the conclusion of Theorem~\ref{thmsimple}
will follow by Gronwall's lemma. By semigroup property, it is sufficient to consider
$t=0$, so the problem is to bound below $\Re\iip{h}{Lh}$ by a multiple of 
$\iip{h}{h}$. Obviously,
\begeq\label{fmlcomp}
\Re\,\iip{h}{Lh} = \Re\,\ip{h}{Lh} + a\,\I + b\,\II + c \,\III, 
\endeq
where
\begin{multline*} \I := \Re\,\ip{Ah}{ALh}, 
\qquad \II := \Re\,\ip{ALh}{Ch} + \Re\,\ip{Ah}{CLh}, \\
\qquad \III := \Re\,\ip{Ch}{CLh}. 
\end{multline*}
By Proposition~\ref{lemkerL}(i), $\Re\,\ip{h}{Lh} = \|Ah\|^2$. For each of the terms
$\I, \II, \III$, the contributions of $A^*\! A$ and $B$ will be estimated separately, 
and the resulting expressions will be denoted $\I_A$, $\I_B$, $\II_A$,
$\II_B$, etc. For consistency with the sequel, I~shall introduce the notation
\begeq\label{R1} R_2 := [C,B]. \endeq
Moreover, to alleviate notation, I~shall temporarily assume that $\H$ is a real Hilbert space;
otherwise, just put real parts everywhere.

First of all,
\begin{multline*} 
\I_B = \ip{Ah}{ABh} = \ip{Ah}{BAh} + \ip{Ah}{[A,B]h} \\
= 0 + \ip{Ah}{Ch} \geq - \|Ah\|\|Ch\|,
\end{multline*}
where the antisymmetry of $B$ was used. Then,
\[ \I_A  = \ip{Ah}{AA^*\! Ah} = 
\ip{A^2h}{A^2h} + \ip{Ah}{[A,A^*]Ah},\]
to be understood as
\[ \sum_{ij}\ \ip{A_jA_ih}{A_iA_jh} + 
\ip{A_ih}{[A_i,A_j^*]A_jh}.\]
This can be rewritten
\begin{multline*} \sum_{ij}\ \|A_i A_j h\|^2 +
\ip{[A_j,A_i]h}{A_iA_jh} + 
\ip{A_ih}{[A_i,A_j^*]A_jh}\\
\equiv \|A^2h\|^2 + \ip{[A,A]h}{A^2h}
+ \ip{Ah}{[A,A^*]Ah}.
\end{multline*}
In the present case it is assumed that $[A,A]=0$, so
the second term vanishes.
Then from the Cauchy--Schwarz inequality we have
\[ \I_A \geq \|A^2h\|^2 - \|Ah\|\,\|[A,A^*]Ah\|.\]

Next,
\begin{align*} 
\II_B & = \ip{ABh}{Ch}+ \ip{Ah}{CBh}  \\
& = \ip{ABh}{Ch} + \ip{Ah}{BCh} + \ip{Ah}{[C,B]h} \\
& = \ip{ABh}{Ch} - \ip{BAh}{Ch} + \ip{Ah}{R_2h} \\
& = \ip{[A,B]h}{Ch} + \ip{Ah}{R_2h} \\
& \geq \|Ch\|^2 - \|Ah\|\,\|R_2h\|;
\end{align*}
\begin{align*}
\II_A & = \ip{Ah}{CA^*\! Ah}+\ip{AA^*\! Ah}{Ch}\\
& = \ip{Ah}{A^*\! CAh} + \ip{A^*\! A^2h}{Ch} + \ip{[A,A^*]Ah}{Ch} \\
& = \ip{A^2h}{CAh} + \ip{A^2h}{ACh} + \ip{[A,A^*]Ah}{Ch}\\
& = 2 \ip{A^2h}{CAh} + \ip{Ch}{[A,A^*]Ah} \\
& \geq - 2 \|A^2h\| \|CAh\| - \|Ch\| \|[A,A^*]Ah\|.
\end{align*}
(Here the commutation of $C$ with both $A$ and $A^*$ was used.)

Finally,
\begin{align*}
\III_B = \ip{Ch}{CBh} & = \ip{Ch}{BCh} + \ip{Ch}{[C,B]h} \\
& = 0 + \ip{Ch}{R_2h}\\
& \geq -\|Ch\|\,\|R_2h\|;
\end{align*}

\[
\III_A  = \ip{Ch}{CA^*\! Ah} = \ip{Ch}{A^*\! CAh} = \ip{ACh}{CAh} = \|CAh\|^2
\]
(here again the commutation of $C$ with $A$ and $A^*$ was used).
\med

On the whole,
\begin{align}\label{estL}
\Re\, \iip{h}{Lh} \geq &\ \|Ah\|^2 \\
& + a\Bigl( \|A^2h\|^2 -
\|Ah\|\,\|[A,A^*]Ah\| -\|Ah\|\,\|Ch\|\Bigr)\notag \\
& + b \Bigl( \|Ch\|^2 - \|Ah\|\,\|R_2h\| - 
2 \|A^2h\| \|CAh\| - \|Ch\| \|[A,A^*]Ah\|\Bigr) \notag\\
& + c \Bigl ( \|CAh\|^2 - \|Ch\| \|R_2h\| \Bigr ). \notag
\end{align}

The assumptions of Theorem~\ref{thmsimple} imply
\[ \|[A,A^*]y\| \leq \alpha \bigl (\|y\| + \|Ay\|  \bigr), \]
\[ \|R_2 h\| \leq \beta \bigl( \|Ah\| + \|A^2h\| + \|Ch\| + \|CAh\| \bigr).\]
Plugging this into~\eqref{estL}, follows an estimate which can be conveniently recast as
\[ \Re\,\iip{h}{Lh} \ge \langle X,\, mX \rangle_{\R^4}, \]
where $X$ is a vector in $\R^4$ and $m$ is a $4\times 4$ matrix, say upper-diagonal:
\[ X := \Bigl( \|Ah\|, \|A^2h\|, \|Ch\|, \|CAh\|\Bigr), \]
\[ m := \left [ \begin{matrix} 
1-\left(a\alpha + b\beta \right) & 
- \left(a\alpha + b\beta\right) & - (a+b\alpha+b\beta+c\beta) & -b\beta \\
0 & a & -(b\alpha+c\beta) & -2b \\
0 & 0 & b-c\beta & -c\beta \\ 
0 & 0 & 0 & c
\end{matrix} \right ] \equiv
[m_{ij}]_{1\leq i,j\leq 4}\]

If the symmetric part of $m$ is definite positive, this will imply inequality~\eqref{xLx}.
Then the rest of Theorem~\ref{thmsimple} follows easily, since
the $\kappa$-coercivity of $A^*\! A+C^*C$ implies
\begin{align*} \|Ah\|^2 + \|Ch\|^2 & 
\geq \frac12 (\|Ah\|^2 + \|Ch\|^2) + \frac{\kappa}2 \|h\|^2 \\
& \geq \frac{\min(1,\kappa)}2 \|h\|_{\H^1}^2.
\end{align*}
So it all boils down now to choosing the parameters $a$, $b$ and $c$ in such
a way that the symmetric part of $m$ is positive definite, and for this it is 
sufficient to ensure that
\[\begin{cases} \forall i,\quad m_{ii}>0; \\
\forall (i,j), \quad i\neq j \Longrightarrow m_{ij} \ll \sqrt{m_{ii} m_{jj}}.
\end{cases}\]
In the sequel, the statement ``the symmetric part of $m_1$ is
greater than the symmetric part of $m_2$'' will be abbreviated into just
``$m_1$ is greater than $m_2$''.

Let $M:=\max(1,\alpha,\beta)$. Assume, to fix ideas, that
\begeq\label{1abc} 1\geq a\geq b\geq 2c.
\endeq
Then $m$ can be bounded below by
\[ \left [ \begin{matrix} 
1-2Ma & -2Ma & -4Ma & -Mb \\
0 & a & -2Mb & -2Mb \\
0 & 0 & b-Mc & -Mc \\ 
0 & 0 & 0 & c
\end{matrix} \right ].\]
If now it is further assumed that
\begeq\label{condabM}
a\leq \frac1{4M}, \qquad c\leq \frac{b}{2M}, 
\endeq
then the latter matrix can in turn be bounded below by
\[ \left [ \begin{matrix} 
1/2 & -2Ma & -4Ma & -Mb \\
0 & a & -2Mb & -2Mb \\
0 & 0 & b/2 & -Mc \\ 
0 & 0 & 0 & c
\end{matrix} \right ]\equiv [\tilde{m}_{ij}].\]

By imposing 
\begeq\label{a}
\tilde{m}_{ij} \leq \sqrt{\tilde{m}_{ii}\tilde{m}_{jj}}/2
\leq (\tilde{m}_{ii} + \tilde{m}_{jj})/4,
\endeq
it will follow
\[ \sum_{ij} \tilde{m}_{ij} X_i X_j \geq \sum_{i} \tilde{m}_{ii} X_i^2
- \frac34 \sum_{i}\tilde{m}_{ii} X_i^2 = \frac14 \sum \tilde{m}_{ii} X_i^2.\]
(The 3 in $3/4$ is because each diagonal term should participate
in the control of three off-diagonal terms.) To ensure~\eqref{a},
it suffices that
\begin{multline*} 
2Ma \leq \sqrt{\frac{a}8}, \quad 4 Ma \leq \sqrt{\frac{b}{16}}, 
\quad Mb \leq \sqrt{\frac{c}{8}}, \quad 2Mb \leq \sqrt{\frac{ab}8},\\
2Mb \leq \sqrt{\frac{ac}4}, \quad Mc\leq \sqrt{\frac{bc}8}. 
\end{multline*}
All these conditions, including~\eqref{condabM}, are fulfilled if
\begeq
a, \ \frac{b}{a}, \ \frac{c}{b} \ \leq \frac1{32\,M^2} \qquad
\frac{a^2}{b}, \ \frac{b^2}{ac}\ \leq \frac1{256\, M^2}. 
\endeq
Lemma~\ref{ll} in Appendix~\ref{toolbox} shows that it is always possible to 
choose $a,b,c$ in such a way. This concludes the proof of Theorem~\ref{thmsimple}.
$\square$

\begin{Rk} There are other possible ways to conduct these calculations.
In an early version of this work, the last term $\III_A$ was rewritten
in three different forms to create helpful terms in $\|AC^*h\|^2$ and $\|A^*\! Ch\|^2$,
at the cost of requiring additional assumptions on $[C,C^*]$.
\end{Rk}

\section{Generalization} \label{secgener}

Now I shall present a variant of Theorem~\ref{thmsimple} which covers more
general situations.

\begin{Thm}\label{hypocomult}
Let $\H$ be a Hilbert space, let $A:\H\to\H^n$ and $B:\H\to \H$ be unbounded
operators, $B^*=-B$, let $L:=A^*\! A+B$ and $\K:=\ker L$. Assume the existence of $\Nc\in\N$ and
(possibly unbounded) operators $C_0, C_1, \ldots, C_{\Nc+1}$, $R_1,\ldots,R_{\Nc+1}$
and $Z_1,\ldots,Z_{\Nc+1}$ such that
\[ C_0=A, \qquad [C_j,B] = Z_{j+1}C_{j+1} + R_{j+1}\quad (0\leq j\leq \Nc), \qquad C_{\Nc+1}=0,\]
and, for all $k\in \{0,\ldots, \Nc\}$,
\sm

(i) $[A,C_k]$ is bounded relatively to $\{C_j\}_{0\leq j\leq k}$ and $\{C_jA\}_{0\leq j\leq k-1}$;
\sm

(ii) $[C_k,A^*]$ is bounded relatively to $I$ and $\{C_j\}_{0\leq j\leq k}$;
\sm

(iii) $R_k$ is bounded relatively to $\{C_j\}_{0\leq j\leq k-1}$ and $\{C_jA\}_{0\leq j\leq k-1}$.
\sm

(iv) There are positive constants $\lambda_j$, $\Lambda_j$ such that
$\lambda_j I \leq Z_j \leq \Lambda_j I$.
\sm

Then there is a scalar product $\iip{\cdot}{\cdot}$ on $\H^1$, which defines
a norm equivalent to the $\H^1$ norm,
\[ \|h\|_{\H^1}:= \sqrt{\|h\|^2 + \sum_{k=0}^{\Nc} \|C_kh\|^2},\]
such that
\begeq\label{xLx'}
\forall h\in \H^1/\K,\qquad \Re\,\iip{h}{Lh} \geq K \sum_{j=0}^{\Nc} \|C_jh\|^2
\endeq
for some constant $K>0$, only depending on the bounds appearing implicitly
in assumptions (i)--(iii).
\sm

If, in addition,
\[ \text{$\sum_{j=0}^{\Nc} C_j^*C_j$ is $\kappa$-coercive} \]
for some $\kappa>0$, then there is a constant $\lambda>0$, only depending on
$K$ and $\kappa$, such that
\[ \forall h\in \H^1/\K, \qquad \Re\,\iip{h}{Lh} \geq \lambda\: \iip{h}{h}. \]
In particular, $L$ is hypocoercive in $\H^1/\K$: There are constants
$C\geq 0$ and $\lambda>0$, explicitly computable in terms of the bounds
appearing implicitly in assumptions (i)--(iii), and $\kappa$, such that
\[ \|e^{-tL}\|_{\H^1/\K \to \H^1/\K} \leq C e^{-\lambda t}.\]
\end{Thm}

This result generalizes Theorem~\ref{thmsimple} in several respects:
successive commutators are allowed, remainders $R_{j+1}$ and multiplicators
$Z_{j+1}$ are allowed in the identity defining $C_{j+1}$ in terms of $C_j$,
and the operators $C_0,\ldots, C_{\Nc}$ are not assumed to commute.

\begin{Rk} The same rule as in Remark~\ref{rkweight} applies to
Assumptions (i)--(iii).
\end{Rk}

\begin{Rk} Theorem~\ref{hypoellmult} in Appendix~\ref{appreg}
will show that the same structure assumptions (i)--(iv) imply an 
immediate regularization effect $\H\to\H^1$. This effect extends the
range of application of the method, allowing data which do
not necessarily lie in $\H^1$ but only in $\H$.
\end{Rk}

\begin{proof}[Proof of Theorem~\ref{hypocomult}]
The proof is an amplification of the proof of Theorem~\ref{thmsimple}. Let
\begeq\label{iipm}
\iip{h}{h}:= \|h\|^2 + \sum_{k=0}^{\Nc} \Bigl( a_k\|C_kh\|^2 + 2 \Re\, b_k 
\ip{C_kh}{C_{k+1}h}\Bigr),
\endeq
where $\{a_k\}_{0\leq k\leq \Nc+1}$ and $\{b_k\}_{0\leq k\leq \Nc}$ are families of positive
coefficients, satisfying
\begeq\label{condakbk}
0\leq k\leq N_c \Longrightarrow 
\begin{cases} a_0\leq \delta,\qquad b_k\leq \delta\, a_k, \qquad a_{k+1}\leq \delta\, b_k, \\ \\
a_k^2 \leq \delta\, b_{k-1}b_k \quad (1\leq k\leq N_c), 
\qquad b_k^2 \leq \delta\, a_k a_{k+1} \quad (0\leq k\leq N_c).
\end{cases}\endeq
The small number $\delta>0$ will be chosen later on, 
and the existence of the coefficients
$a_k$, $b_k$ is guaranteed by Lemma~\ref{ll} again.

Since $C_{\Nc+1}=0$, the last term in~\eqref{iipm}, with coefficient $b_{\Nc}$,
does not play any role. For $k\leq \Nc-1$, the inequality $b_k\leq\delta \sqrt{a_k a_{k+1}}$ 
implies
\[ b_k \Bigl|\ip{C_kh}{C_{k+1}h}\Bigr| \leq \frac{\delta a_k}{2} \|C_kh\|^2 + 
\frac{\delta a_{k+1}}2 \|C_{k+1}h\|^2. \]
Hence, for $\delta$ small enough,
\[ \iip{h}{h} \geq \|h\|^2 + \frac12 \sum_{k=0}^{\Nc} a_k \|C_kh\|^2. \]
So the norm defined by~\eqref{iipm} is indeed equivalent to
the $\H^1$ norm.
\med

The next observation is that the space $\K^\bot$ is the same, whether the orthogonality
is defined with respect to the scalar product in $\H$, the one in $\H^1$ or the one
defined by~\eqref{iipm}. To show this it is sufficient to prove
\begeq\label{toprovekerL} 
h\in \ker L \Longrightarrow \qquad \forall k\in \{0,\ldots, \Nc+1\}, \quad C_k h=0.
\endeq
This will be achieved by finite induction on $k$. Let $h\in\ker L$; 
by Lemma~\ref{lemkerL}, $Ah=0$, $Bh=0$; so~\eqref{toprovekerL} is true for $k=0$.
Assume now that $C_jh=0$ for $j\leq k$; it is obvious that also $C_jAh=0$ for
$j\leq k$; then our assumption on $R_{k+1}$ implies that $R_{k+1}h=0$. So
$C_{k+1}h= C_kBh - BC_kh - R_{k+1}h=0$. This concludes the proof of~\eqref{toprovekerL}.
\med

To prove~\eqref{xLx'} it is obviously sufficient to establish
\begeq \label{goalnow}
\Re\,\iip{h}{Lh} \geq \frac12 \|Ah\|^2 + \sum_{k=0}^{\Nc} \Bigl( \frac{a_k}2 \|C_kAh\|^2 
+ \frac{b_k}2 \|C_{k+1}h\|^2 \Bigr). 
\endeq 

As in the proof of Theorem~\ref{thmsimple} one can compute, with obvious notation,
\begeq\label{iipLx} 
\Re\,\iip{h}{Lh} = \|Ah\|^2 + \sum_{k=0}^{\Nc} \Bigl \{ a_k [\I_A^k + \I_B^k]
+ 2 b_k [ \II_A^k + \II_B^k ] \Bigr\}.
\endeq

To alleviate the notation, assume for a moment that we are working in a real Hilbert space,
so there is no need to take real parts (otherwise, just put real parts everywhere).
Explicit computations yield, for $k\leq \Nc$,
\begin{align*}
\I_B^k = \ip{C_kh}{C_kBh} & = \ip{C_kh}{[C_k,B]h} + \ip{C_kh}{BC_kh} \\
& = \ip{C_kh}{[C_k,B]h} + 0 \\
& = \ip{C_kh}{Z_{k+1}C_{k+1}h} + \ip{C_kh}{R_{k+1}h} +0 \\
& \geq - \Lambda_{k+1}\|C_kh\|\,\|C_{k+1}h\| - \|C_kh\|\,\|R_{k+1}h\|
\end{align*}
\begin{align*}
\I_A^k = \ip{C_kh}{C_kA^*\! Ah} & = \ip{C_kh}{A^*\! C_kAh} + \ip{C_kh}{[C_k,A^*]Ah} \\
& = \ip{AC_kh}{C_kAh} + \ip{C_kh}{[C_k,A^*]Ah} \\
& = \ip{C_kAh}{C_kAh} + \ip{[A,C_k]h}{C_kAh} + \ip{C_kh}{[C_k,A^*]Ah} \\
& \geq \|C_kAh\|^2 - \|C_kAh\|\,\|[A,C_k]h\| - \|C_kh\|\,\|[C_k,A^*]Ah\|
\end{align*}
and, for $k\leq \Nc-1$ (when $k=\Nc$, $\II_B^k=0$),
\begin{align*}
\II_B^k & = \ip{C_kBh}{C_{k+1}h} + \ip{C_kh}{C_{k+1}Bh} \\
& = \ip{C_kBh}{C_{k+1}h} + \ip{C_kh}{BC_{k+1}h} + \ip{C_kh}{[C_{k+1},B]h} \\
& = \ip{C_kBh}{C_{k+1}h} - \ip{C_kh}{BC_{k+1}h} + \ip{C_kh}{[C_{k+1},B]h} \\
& = \ip{[C_k,B]h}{C_{k+1}h} + \ip{C_kh}{[C_{k+1},B]h} \\
& = \ip{C_kBh}{C_{k+1}h} - \ip{BC_kh}{C_{k+1}h} + \ip{C_kh}{Z_{k+2}C_{k+2}h} 
+ \ip{C_kh}{R_{k+2}h} \\
& = \ip{[C_k,B]h}{C_{k+1}h} + \ip{C_kh}{Z_{k+2}C_{k+2}h} + \ip{C_kh}{R_{k+2}h} \\
& = \ip{Z_{k+1}C_{k+1}h}{C_{k+1}h} + \ip{R_{k+1}h}{C_{k+1}h} 
+ \ip{C_kh}{Z_{k+2}C_{k+2}h} + \ip{C_kh}{R_{k+2}h} \\
& \geq \lambda_{k+1}\|C_{k+1}h\|^2 - \|C_{k+1}h\|\,\|R_{k+1}h\| - 
\Lambda_{k+2}\|C_kh\|\,\|C_{k+2}h\| - \|C_kh\|\,\|R_{k+2}h\|
\end{align*}
\begin{align*}
\II_A^k & = \ip{C_kh}{C_{k+1}A^*\! Ah} + \ip{C_kA^*\! Ah}{C_{k+1}h} \\
& = \ip{C_kh}{[C_{k+1},A^*]Ah} + \ip{C_kh}{A^*\! C_{k+1}Ah} + \ip{A^*\! C_kAh}{C_{k+1}h} \\
   & \hspace*{80mm}    + \ip{[C_k,A^*]Ah}{C_{k+1}h} \\
& = \ip{C_kh}{[C_{k+1},A^*]Ah} + \ip{AC_kh}{C_{k+1}Ah} + \ip{C_kAh}{AC_{k+1}h} \\
   & \hspace*{80mm} + \ip{[C_k,A^*]Ah}{C_{k+1}h} \\
& = \ip{C_kh}{[C_{k+1},A^*]Ah} + \ip{C_kAh}{C_{k+1}Ah} + \ip{[A,C_k]h}{C_{k+1}Ah} \\
   & \hspace*{30mm} + \ip{C_kAh}{C_{k+1}Ah} + \ip{C_kAh}{[A,C_{k+1}]h} 
     + \ip{[C_k,A^*]Ah}{C_{k+1}h} \\
& \geq - \|C_kh\|\,\|[C_{k+1},A^*]Ah\| - \|C_kAh\|\,\|C_{k+1}Ah\| - \|C_{k+1}Ah\|\,\|[A,C_k]h\| \\
   & \hspace*{20mm}    - \|C_kAh\|\,\|C_{k+1}Ah\| 
     - \|C_kAh\|\,\|[A,C_{k+1}]h\| - \|C_{k+1}h\|\,\|[C_k,A^*]Ah\|.
\end{align*}
%horrible typography

The next step is to use the quantities 
$a_j \|C_jAh\|^2$ and $b_j\|C_{j+1}h\|^2$
to control all the remaining terms. For this I~shall apply Young's inequality,
in the form $XY\leq \var X^2 + C(\var)Y^2$ ($C(\var)=\var^{-1}/4$). In the computations below,
the dependence of $C$ on the constants $\Lambda_j$ will not be recalled.
\begin{align}
a_k [\I_A^k +\I_B^k] + b_k [\II_A^k+\II_B^k] &
\geq a_k \|C_kAh\|^2 + b_k \lambda_{k+1}\|C_{k+1}h\|^2 \\
& - \var b_{k-1} \|C_kh\|^2 - C \frac{a_k^2}{b_{k-1}} \|C_{k+1}h\|^2 \label{1}\\
& - \var b_{k-1} \|C_kh\|^2 - C \frac{a_k^2}{b_{k-1}} \|R_{k+1}h\|^2 \label{2}\\
& - \var a_k \|C_kAh\|^2 - C a_k \|[A,C_k]h\|^2 \label{3}\\
& - \var b_{k-1} \|C_kh\|^2 - C \frac{a_k^2}{b_{k-1}} \|[C_k,A^*]Ah\|^2 \label{4}\\
& - \var b_k \|C_{k+1}h\|^2 - C b_k \|R_{k+1}h\|^2 \label{5}\\
& - \var b_{k-1} \|C_kh\|^2 - C \frac{b_k^2}{b_{k-1}} \|C_{k+2}h\|^2 \label{6}\\
& - \var b_{k-1} \|C_kh\|^2 - C \frac{b_k^2}{b_{k-1}} \|R_{k+2}h\|^2 \label{7}\\
& - \var b_{k-1} \|C_kh\|^2 - C \frac{b_k^2}{b_{k-1}} \|[C_{k+1},A^*]Ah\|^2 \label{8}\\
& - \var a_k \|C_kAh\|^2 - C \frac{b_k^2}{a_k} \|C_{k+1}Ah\|^2 \label{9}\\
& - \var a_{k+1} \|C_{k+1}Ah\|^2 - C \frac{b_k^2}{a_{k+1}} \|[A,C_k]h\|^2 \label{10}\\
& - \var a_k \|C_kAh\|^2 - C \frac{b_k^2}{a_k} \|C_{k+1}Ah\|^2 \label{11}\\
& - \var a_k \|C_kAh\|^2 - C \frac{b_k^2}{a_k} \|[A,C_{k+1}]h\|^2 \label{12}\\
& - \var  b_k \|C_{k+1}h\|^2 - C \frac{b_k^2}{a_k} \|[C_k,A^*]Ah\|^2 \label{13},
\end{align}
with the understanding that~\eqref{5} to~\eqref{13} in the above are not present
when $k=\Nc$. The problem is to show that each of the terms appearing in 
lines~\eqref{1} to~\eqref{13} can be bounded below by
\begeq\label{etaexpr} 
-\var \Bigl(\|Ah\|^2 +\sum_j (a_j \|C_jAh\|^2 + b_j\|C_{j+1}h\|^2)\Bigr), 
\endeq
as soon as $\delta$ is small enough. This is true, by construction, of all the
terms appearing on the left in these lines; so let us see how to control all the terms 
on the right. In the sequel, the notation $u\ll v$ means $u\leq \eta\, v$, where $\eta>0$
becomes arbitrarily small as $\delta\to 0$.
\med

- For lines~\eqref{1}, \eqref{6}, \eqref{9} and~\eqref{11}, it is sufficient to
impose
\[ \frac{a_k^2}{b_{k-1}} \ll b_k, \qquad \frac{b_k^2}{b_{k-1}} \ll b_{k+1}, \qquad
\frac{b_k^2}{a_k} \ll a_{k+1}. \]
The first and the third of these inequalities are true by construction;
as for the second one, it follows from
\[ b_k^2 \ll a_k a_{k+1} \ll \sqrt{b_{k-1}b_k}\,\sqrt{b_kb_{k+1}} \Longrightarrow \qquad
b_k^2 \ll b_{k-1}b_{k+1}. \]
\med

- For lines~\eqref{2}, \eqref{5} and~\eqref{7}, we know that
$\|R_{k+1}h\|$ is controlled by a combination of $\|Ah\|$, $\|A^2h\|$, $\|C_1h\|$,
$\|C_1Ah\|$, \ldots, $\|C_kh\|$, $\|C_kAh\|$; hence it is sufficient to bound the coefficients
appearing in front of $R_{k+1}h$ (resp. $R_{k+2}h$) by a small multiple of $a_k$ (resp. $a_{k+1}$). 
So these terms are fine as soon as
\[ \frac{a_k^2}{b_{k-1}} \ll a_k, \qquad b_k \ll a_k, \qquad
\frac{b_k^2}{b_{k-1}} \ll a_{k+1}.\]
The second of these inequalities is true by construction, while the first and third one 
follow from
\[ a_k^2\ll b_{k-1} b_k \ll b_{k-1} a_k, \qquad b_k^2 \ll b_{k-1}b_{k+1} \ll b_{k-1} a_{k+1}. \]
\med

- For lines~\eqref{3}, \eqref{10} and~\eqref{12}, we know that
$\|[A,C_k]h\|$ is controlled by $\|Ah\|$, $\|A^2h\|$, $\|C_1h\|$, $\|C_1Ah\|$, \ldots, $\|C_kh\|$.
By a reasoning similar to the one above, it is sufficient to ensure
\[ a_k \ll b_{k-1}, \qquad \frac{b_k^2}{a_{k+1}} \ll b_{k-1}, \qquad
\frac{b_k^2}{a_k} \ll b_k. \]
The first and third of these inequalities are true by construction, 
while the second one follows from
\[ b_k^2 \ll a_k a_{k+1} \ll b_{k-1} a_{k+1}. \]
\med

- For lines~\eqref{4}, \eqref{8} and~\eqref{13}, we know that
$\|[C_k,A^*]y\|$ is controlled by $\|y\|$, $\|Ay\|$, $\|C_1y\|$, \ldots, $\|C_ky\|$,
so $\|[C_k,A^*]Ah\|$ is controlled by $\|Ah\|$, $\|A^2h\|$, $\|C_1Ah\|$, \ldots, $\|C_kAh\|$.
By a reasoning similar to the one above, it is sufficient to ensure
\[ \frac{a_k^2}{b_{k-1}} \ll a_k,\qquad \frac{b_k^2}{b_{k-1}} \ll a_{k+1}, \qquad
\frac{b_k^2}{a_k} \ll a_{k}. \]
The first and third of these inequalities are true by construction, while the second one
follows from
\[ b_k^2 \ll a_k a_{k+1} \ll b_{k-1} a_{k+1}. \]
\med

Putting all together, for all $\eta$ there is a $\delta$ such that each of the ``error terms''
which appeared in the estimates above can be bounded below by~\eqref{etaexpr}.
Then, if $N_e$ stands for the number of error terms,
\[ \Re\,\iip{h}{Lh} \geq \bigl(1-N_e\var\bigr) \left( \|Ah\|^2 + \sum_{k=0}^{\Nc} 
\Bigl( \frac{a_k}2 \|C_kAh\|^2 + \frac{b_k}2 \|C_{k+1}h\|^2 \Bigr) \right), \]
which implies~\eqref{goalnow} for $\delta$ small enough.
The proof of~\eqref{xLx'} is now complete, and the end of Theorem~\ref{hypocomult}
follows easily, as in the proof of Theorem~\ref{thmsimple}.
\end{proof}

I~shall conclude this section with a simple generalization of Theorem~\ref{hypocomult}:

\begin{Thm} \label{thmCjA}
With the same notation as in Theorem~\ref{hypocomult}, define
$C_{j+1/2}=C_jA$. Then the conclusion of Theorem~\ref{hypocomult} still holds true if 
Assumptions (i) to (iii) are relaxed as follows: There exists a constant $M$ such that

\quad (i') $\|[A,C_k]h\| \leq M \sum \|C_\alpha h\|^\theta\,\|C_\beta h\|^{1-\theta}$,
where the sum is over all couples of indices $(\alpha,\beta)$ such that
$\theta\alpha+(1-\theta)\beta \leq k+1/2$, $(\alpha,\beta) \neq (k+1/2, k+1/2)$;

\quad (ii') $\|[C_k,A^*]Ah\| \leq M \sum \|C_\alpha h\|^\theta\,\|C_\beta h\|^{1-\theta}$,
where the sum is over all couples of indices $(\alpha,\beta)$ such that
$\theta\alpha+(1-\theta)\beta \leq k+1$, $(\alpha,\beta) \neq (k+1, k+1)$;

\quad (iii') $\|R_kh\| \leq M \sum \|C_\alpha h\|^\theta\,\|C_\beta h\|^{1-\theta}$,
where the sum is over all couples of indices $(\alpha,\beta)$ such that
$\theta\alpha+(1-\theta)\beta \leq k$, $(\alpha,\beta) \neq (k, k)$.
\end{Thm}

In statements (i') to (iii'), $\theta$ may vary from one couple $(\alpha,\beta)$ 
to the other. The conditions on admissible couples $(\alpha,\beta)$ 
can be understood more easily if one remembers
Remark~\ref{rkweight}; then the weight $w(C_\alpha)$ is $2\alpha+1$. 
If one formally attributes to $[A,C_k]$, $[A^*,C_k]A$ and $R_{k}$ the 
weights $2k+2$, $2k+3$ and $2k+1$, and decides that the weight of a formal product
$\|O_1h\|^\theta\,\|O_2h\|^{1-\theta}$ is $\theta w(O_1) + (1-\theta) w(O_2)$,
then Conditions (i') to (iii') mean that each of the operators $[A,C_k]$,
$[A^*,C_k]A$ and $R_k$ can be bounded in terms of lower weights.
For instance, an estimate like
\[ \|R_5h\| \leq M \sqrt{\|C_3h\|\,\|C_5h\|}\]
is admissible, since $[w(C_3)+w(C_5)]/2=(7+11)/2<11=w(R_5)$. 

\begin{proof}[Proof of Theorem~\ref{thmCjA}]
The strategy is the same as in Theorem~\ref{hypocomult}; but now one should 
use Young's inequality in the form
\[ a^\theta b^{1-\theta} \leq \var a + C\, b, \]
and note that, if $(u_m)$ is given by Lemma~\ref{ll}, then
\[ \left[\ell \leq \frac{m+n}2, \quad (m,n)\neq (\ell,\ell)\right] \Longrightarrow\qquad
u_\ell \ll \sqrt{u_m\, u_n}. \]
Then all the estimates entering the proof of Theorem~\ref{hypocomult} can be adapted
without difficulty.
\end{proof}

\section{Hypocoercivity in entropic sense} \label{secentropic}

In this section I~shall consider the problem of convergence to equilibrium
for solutions of diffusion equations in an $L\log L$ setting. This represents
a significant extension of the results already discussed, 
because in many cases of interest, after a finite time
the solution automatically belongs to $L\log L(\mu)$, where
$\mu$ is the stationary solution, but not to 
$L^2(\mu)$.\footnote{Of course this does
not contradict the fact that it will be {\em locally} $C^\infty$.
The most basic illustration is the case of the linear Fokker--Planck equation
$\pa_t h = \Delta_v h - v\cdot\nabla_v h$: The theory of hypercontractivity
tells us that the semigroup at time $t$ is regularizing from $L^p$ to $L^q$
only after a time $\log((q-1)/(p-1))/2$, which is finite only if
$p>1$.}

For that purpose, I~shall use the same information-theoretical functionals 
as in the theory of logarithmic Sobolev inequalities: First, the
{\bf Kullback information} (or Boltzmann $H$ functional, or Shannon information),
\[ H_\mu(\nu) = \int h\log h\,d\mu,\qquad \nu = h\mu;\]
and secondly, the {\bf Fisher information}
\[ I_\mu(\nu) = \int \frac{|\nabla h|^2}{h}\,d\mu,\qquad \nu = h\mu.\]
Recall that a probability measure $\mu$ on $\R^N$ satisfies a 
{\em logarithmic Sobolev inequality} if there is a constant $\lambda>0$ such that
\[ H_\mu(\nu) \leq \frac{1}{2\lambda} I_\mu(\nu),\]
for all probability measures $\nu$ on $\R^N$ (with the convention that
$H_\mu(\nu)=I_\nu(\mu)=+\infty$ if $\nu$ is not absolutely continuous with
respect to $\mu$).

The main difference with the classical theory is that I~shall {\em distort}
the Fisher information by using a suitable field of quadratic forms; that is,
replace $\int |\nabla h|^2/h\,d\mu$ by $\int \<S\nabla h, \nabla h\>/h\,d\mu$,
where $x\to S(x)$ is a function valued in the space of quadratic forms, in
such a way that $S(x)\geq \kappa I_N$ for some $\kappa>0$, independently
of $x$. It turns out that the same algebraic tricks which worked
in a Hilbertian context will also work here, at the price of more
stringent assumptions on the vector fields: The proofs will
be based on some slightly miraculous-looking computations, which
may be an indication that there is more structure to understand.

Here below is the main result of this section. Note carefully that this is
not expressed in terms of linear operators in abstract Hilbert spaces, but
in terms of derivation operators on $\R^N$. (The theorem might possibly
be generalized by replacing $\R^N$ by a smooth manifold.)
So as not to be bothered with regularity issues,
I~shall assume here that the reference density is rapidly decaying
and that all coefficients are $C^\infty$ and have at most polynomial
growth; but of course these assumptions can be relaxed. I~shall
also assume that the solution is smooth if the initial datum
is smooth.

\begin{Thm} \label{thmLlogL}
Let $E\in C^2(\R^N)$, such that $e^{-E}$ is rapidly decreasing,
and $\mu(dX) = e^{-E(X)}\,dX$ is a probability measure on $\R^N$. 
Let $(A_j)_{1\leq j\leq m}$ and $B$ be
first-order derivation operators with smooth coefficients.
Denote by $A_j^*$ and $B^*$ their respective adjoints in $L^2(\mu)$, 
and assume that $B^*=-B$. Denote by $A$ the collection $(A_1,\ldots,A_m)$, 
viewed as an unbounded operators whose range is made of functions valued in $\R^m$.
Define 
\[ L = A^*A + B = \sum_{j=1}^m A_j^*A_j + B,\]
and assume that $e^{-tL}$ defines a well-behaved semigroup on a
suitable space of positive functions
(for instance, $e^{-tL}h$ and $\log (e^{-tL}h)$ are $C^\infty$ and all
their derivatives grow at most polynomially if $h$ is itself $C^\infty$ 
with all derivatives bounded, and $h$ is bounded below by a positive
constant).

Next assume the existence of $\Nc\geq 1$,
derivation operators $C_0,\ldots,C_{\Nc+1}$, $R_1,\ldots,R_{\Nc+1}$;
and vector-valued functions $Z_1,\ldots,Z_{\Nc+1}$ (all of them with
$C^\infty$ coefficients, growing at most polynomially, as their
partial derivatives) such that
\[ C_0 =A, \qquad [C_j,B] = Z_{j+1}\,C_{j+1} + R_{j+1}\quad
(0\leq j\leq \Nc),\qquad C_{\Nc+1}=0,\]
and

(i) $[A,C_k]$ is pointwise bounded relatively to $A$;

(ii) $[C_k,A^*]$ is pointwise bounded relatively to 
$I, \{C_j\}_{0\leq j\leq k}$;

(iii) $R_k$ is pointwise bounded with respect to $\{C_j\}_{0\leq j\leq k-1}$;

(iv) there are positive constants $\lambda_j,\Lambda_j$ such that
$\lambda_j \leq Z_j\leq \Lambda_j$;

(v) $[A,C_k]^*$ is pointwise bounded relatively to $I,A$.
\sm

Then there is a function $x\to S(x)$, valued in the space of quadratic
forms on $\R^N$, uniformly bounded,
such that if one defines
\[ {\cal E}(h) := \int h\,\log h\,d\mu + 
\int \frac{\< S \nabla h, \nabla h\>}{h}\,d\mu,\]
then one has the estimate
\[ \frac{d}{dt} {\cal E} (e^{-tL}h) \leq -\alpha\,
\frac{\< S \nabla h, \nabla h\>}{h}\,d\mu,\]
for some positive constant $\alpha>0$, which is explicitly computable
in terms of the bounds appearing implicitly in conditions (i)--(v).
\sm

If furthermore 
\sm

(a) there is a positive constant $\lambda$ such that
$\sum_k C_k^* C_k \geq \lambda I_N$, pointwise on $\R^N$;
\sm

(b) $\mu$ satisfies a logarithmic Sobolev inequality
with constant $K$;
\sm
\sm

\noindent then $S(x)$ is uniformly positive definite,
and there is a constant $\kappa>0$ such that
\[ \frac{d}{dt} {\cal E} (e^{-tL}h) \leq -\kappa\, {\cal E}(e^{-tL}h).\]
In particular,
\[ I_\mu((e^{-tL}h\,\mu) = O(e^{-\kappa t}),\qquad
H_\mu((e^{-tL}h\,\mu) = O(e^{-\kappa t}),\]
and all the constants in this estimate can be estimated explicitly
in terms of the bounds appearing implicitly in conditions (i)--(v),
and the constants $\lambda,K$.
\end{Thm}

\begin{Rk} The matrices $S(x)$ will be constructed from 
the vector fields entering the equation, by linear combinations
{\em with constant coefficients}. I~expect that for more degenerate
situations it will be useful to use varying coefficients.
\end{Rk}

\begin{Rk} \label{rkexpl}
A major difference between the assumptions of Theorem~\ref{hypocomult}
and the assumptions of Theorem~\ref{thmLlogL} is that the latter impose
{\em pointwise} bounds on $\R^N$, in the following sense.
First, $A$ is an $m$-tuple of derivation operators $(A_i)_{1\leq i\leq m}$, 
each of which can be identified with a vector field $\sigma_i$, 
in such a way that $A_i f = \sigma_i\cdot\nabla f$; so 
$\sigma=(\sigma_i)_{1\leq i\leq m}$ can be seen as a map valued in
$(m\times N)$ matrices. Then each commutator $C_k$ is also an $m$-tuple of 
derivation operators $(C_{k,j})_{1\leq j\leq m}$, so that $C_{k,j}$
has been obtained from the commutation of $C_{k-1,j}$ with $B$.
Then $[A_i,C_{k,j}]$ is represented by a vector field $\xi_{i,j,k}$;
and Assumption (i) says that $|\xi_{i,j,k}(x)|$ is bounded,
{\em for all $x$}, by $c\|\sigma(x)\|$, where $c$ is a constant.
The other pointwise conditions are to be interpreted similarly.
Let us consider for instance Assumption (ii).
Since $A_i$ is a derivation, the adjoint of $A_i$ takes the form
$-A_i + a_i$, where $a_i$ is a function; so the adjoint of $A$ is
of the form $(g_1,\ldots,g_m)\to -\sum A_i g_i + \sum a_i g_i$.
Then the commutator of $C$ with $A^*$ is the same as the commutator
of $C$ with $A$, up to an array of operators that are the multiplication
by the {\em functions} $C_j a_i$.
So the second part of Assumption (ii) really says that 
the functions $C_j a_i$ are all bounded.
Finally, note that in Assumption (a), each
$C_k$ is an $m$-tuple of derivations, so it can be identified to
a function valued in $m\times N$ matrices; and $\sum C_k^*C_k$ 
to a function valued in $N\times N$ matrices, which should be
uniformly positive definite.
\end{Rk}

\begin{Rk} Theorem~\ref{hypoellLog} below will show that the Assumptions
of Theorem~\ref{thmLlogL}
entail an immediate ``entropic'' regularization effect:
If the initial datum is only assumed to have finite entropy, then
the functional ${\cal E}$ becomes immediately finite.
This allows to extend the range of application of the method to data
which are only assumed to have finite entropy. In the case of the
Fokker--Planck equation I~shall show how to relax even this assumption
of finite entropy and treat initial data which are only
assumed to have finite moments of large enough order.
\end{Rk}

The key to the proof of Theorem~\ref{thmLlogL} is the following lemma,
which says that the computations arising in the time-differentiation
of the functional ${\cal E}$ are quite the same as the computations
arising in Theorem~\ref{hypocomult}, {\em provided that
$A$ and $C_k$ commute}.

In the next statement, I~shall use the notation
\[ \left(\frac{d}{dt} \right)_{S} {\cal F}(h)\]
for the time-derivative of the functional ${\cal F}$ along the
semigroup generated by the linear operator $-S$. More explicitly,
\[ \left(\frac{d}{dt} \right)_{S} {\cal F}(h) = 
\left.\frac{d}{dt}\right|_{t=0} {\cal F}(e^{-tS}h).\]
Moreover, when no measure is indicated this means that the
Lebesgue measure should be used.

\begin{Lem} \label{lemcalcHI}
Let $\mu(dX)=e^{-E(X)}\,dX$, $A=(A_1,\ldots,A_m)$, $B$ 
and $L=A^*A+B$ be as in Theorem~\ref{thmLlogL}.
Let $C=(C_1,\ldots,C_m)$ and $C'=(C'_1,\ldots,C'_m)$ be $m$-tuples 
of derivation operators on $\R^N$ (all of them with smooth
coefficients whose derivatives grow at most polynomially).
Then, with the notation $f=he^{-E}$, $u=\log h$, one has

\begeq\label{ddtHB} \left(\frac{d}{dt} \right)_{B}
\int h\,\log h\,d\mu = 0;
\endeq

\begeq\label{ddtHA}
- \left(\frac{d}{dt} \right)_{A^*\!A}
\int h\,\log h\,d\mu = \int \frac{|Ah|^2}{h}\,d\mu = 
\int f |Au|^2,
\endeq
where by convention $|Au|^2 = \sum_i (A_i u)^2$;

\begin{align}\label{ddtIB}
- \left( \frac{d}{dt} \right)_{B} \int \frac{\<Ch, C'h\>}{h}\,d\mu
& = \int \frac{\<Ch, \, [C',B]h\>}{h}\,d\mu
+ \int \frac{\<[C,B]h,\, C'h\>}{h}\,d\,\mu \\
& = \int f\,\bigl\<Cu,\, [C',B]u\bigr\> + 
\int f\, \bigl\<[C,B]u,\, C'u\bigr\>,
\end{align}
where by convention $\<[C,B]u,\, C'u\> = \sum_j ([C_j,B]u)(C'_j u)$;

\begin{align}\label{ddtIA}-\left(\frac{d}{dt} \right)_{A^*\!A}
\int \frac{\<Ch,C'h\>}{h}\,d\mu = & \
2 \int f \<C A u,\, C'Au\>  \\
& + \left( \int f\,\<[C,A^*]Au,\, C'u\> + 
\int f\,\<C u, \, [C',A^*]Au\>\right) \nonumber \\
& + \left( \int f\, \bigl\<CAu,\, [A,C']u\bigr\> +
\int f\, \bigl\<[A,C]u,\, C'Au\bigr\> \right) \nonumber \\
& + \int f\, Q_{A,C,C'}(u), \nonumber
\end{align}
where by convention $\<C u, \, [C',A^*]Au\> = 
\sum_{ij} (C_ju)([C'_j,A_i^*]A_iu)$, etc. and
\begin{multline}\label{QACC}
Q_{A,C,C'}(u) := [A,C]^*(Au\otimes C'u) +
[A,C']^*(Au\otimes Cu) \\
:= \sum_{ij} [A_i,C_j]^*(A_iu\, C'_ju) +
[A_i,C'_j]^*(A_iu\, C_ju).
\end{multline}
\end{Lem}

\begin{Rk}
If $[A_i,C_j]=[A_i,C'_j]=0$ for all $i,j$, then obviously
$Q_{A,C,C'}$ vanishes identically. The same conclusion holds
true if $A=C=C'$, even if $[A_i,A_j]$ is not necessarily~0.
Indeed, $[A_j,A_i]^*(A_ju\, A_iu) = - [A_i,A_j]^*(A_iu\,A_ju)$,
so $\sum_{ij} [A_i,A_j]^* (A_iu\, A_ju)=0$ by the symmetry
$i\leftrightarrow j$. I~don't know whether there are simple
general conditions for the vanishing of $Q_{A,C,C'}$, that
would encompass both $[A,C']=[A,C]=0$ and $A=C=C'$ as
particular cases.
\end{Rk}

\begin{Rk} One of the conclusions of this lemma is that the derivatives 
of the quantities $\int h\,\log h\,d\mu$ and $\int \<Ch,\,C'h\>/h\,d\mu$
can be computed just as the derivatives of the quantities
$\int h^2\,d\mu$ and $\int \<Ch,\,C'h\>\,d\mu$, if one replaces
in the final result the measure $\mu$ by $f = h\,\mu$,
and in the integrand the function $h$ by its logarithm,
{\em as long as the quantity $Q_{A,C,C'}$ vanishes}.
In the special case when $C=C'=A$, this principle is well-known in the
theory of logarithmic Sobolev inequalities, where it is stated
in terms of Bakry and \'Emery's ``$\Gamma_2$ calculus''. 
As in the theory of $\Gamma_2$ calculus, Ricci curvature should
play a crucial role here, since it is related to the commutator 
$[A,A^*]$. (In a context of
Riemannian geometry, this is what Bochner's formula is about.)
\end{Rk}

\begin{proof}[Proof of Lemma~\ref{lemcalcHI}]
The proofs of~\eqref{ddtHB} and \eqref{ddtHA} are easy and well-known,
however I~shall recall them for completeness. The proof
of~\eqref{ddtIB} will not cause any difficulty. But the proof
of~\eqref{ddtIA} will be surprisingly complicated and indirect, 
which might be a indication that a more appropriate formalism 
is still to be found.

As before, I~shall assume that the function $h$ is very smooth,
and that all the integration by parts or other manipulations needed
in the proof are well justified. To alleviate notation, I~shall
abbreviate $e^{-tL}h$ into just $h$, with the understanding that
the time dependence is implicit. 
Also $f=he^{-E}$ and $u=\log h = \log f + E$ 
will depend implicitly on the time $t$. Recall that the Lebesgue measure 
is used if no integration measure is specified.

First, by the chain-rule,
\begin{multline*} - \left(\frac{d}{dt}\right)_B \int h\log h\,d\mu = 
\int (\log h+1)\,(Bh)\,d\mu = \int B(h\log h)\,d\mu \\
= \int (B^*1) (h\log h)\,d\mu,
\end{multline*}
and this quantity vanishes since $B^*=-B$ is a derivation.
This proves~\eqref{ddtHB}.

Next,
\begin{multline*}- \left(\frac{d}{dt}\right)_{A^*A} \int h\log h\,d\mu = 
\int (\log h+1)\, (A^*\!Ah)\,d\mu = \int \<A(\log h+1),\, Ah\>\,d\mu\\
= \int \bigl\<\frac{Ah}{h},\, Ah\bigr\>\,d\mu.
\end{multline*}
This proves~\eqref{ddtHA}.

To prove~\eqref{ddtIB}, it suffices to remark that (a) the integrand can
be written as a quadratic expression of $\sqrt{h}$: Indeed, by chain rule,
\[ \int \frac{\<Ch, \,C'h\>}{h}\,d\mu =
4 \int \< C\sqrt{h},\, C'\sqrt{h}\>\,d\mu;\]
and that (b) the evolution equation for $\sqrt{h}$ along $B$ is the same
as for $h$: Indeed, $\pa_t h + Bh=0$ implies $\pa_t\sqrt{h} + B\sqrt{h}=0$.
So to compute the time-derivative in~\eqref{ddtIB} the problem reduces
to a {\em quadratic} computation:
\[ \int \<CB\sqrt{h},\, C'\sqrt{h}\>\,d\mu = 
\int \<BC\sqrt{h},\, C'\sqrt{h}\>\,d\mu
+ \int \<[C,B]\sqrt{h},\, C'\sqrt{h}\>\,d\mu.\]
Then the first term in the right-hand side vanishes since $B$ is
antisymmetric. Formula~\eqref{ddtIB} follows upon use of the
chain-rule again.

Now it only remains to establish~\eqref{ddtIA}.
Before starting the computations, let us recast the equation 
$\pa_t h + A^*\!Ah=0$ in terms of $f=he^{-E}$. 
It follows by Proposition~\ref{rulesdiffeq},
with $\rho_\infty=e^{-E}$ and $B=0$, that
\[ \pa_t f = \nabla\cdot (D (\nabla f + f\, \nabla E))
= \nabla\cdot (D f \nabla u),\]
with the diffusion matrix $D=A^*A$, or more rigorously
$\sigma^*\sigma$, where $\sigma$ is such that $Ah = \sigma (\nabla h)$.
In particular, if $v$ and $w$ are two smooth functions, then
\begeq \label{DnablaV}\<D \nabla v,\, \nabla w\> = \<Av,\, Aw\>.
\endeq
Another relation will be useful later: by explicit computation,
if $g$ is a vector-valued smooth function, then
\[ A^*g = \nabla\cdot (\sigma^*g) - \<\sigma\nabla E,\, g\>;\]
it follows that, for any real-valued smooth function $u$,
\begeq\label{Dnablau}
\nabla\cdot (D\nabla u) - \<D\nabla E,\, \nabla u\> = A^*Au.
\endeq

Next, by chain-rule,
\[ \int \frac{\<Ch,\,C'h\>}{h}\,d\mu = \int f\, \< Cu,\, C'u\>.\]
So the left-hand side of~\eqref{ddtIA} is equal to
\begin{multline}\label{123}
- \int \nabla\cdot (f\, D\nabla u)\, \<Cu,\,C'u\> - 
\int f\, \Bigl\< \frac{C \nabla\cdot (Df\nabla u)}{f},\, C'u\Bigr\>\\
- \int f\, \Bigl\< Cu,\, \frac{C' \nabla\cdot (Df\nabla u)}{f}\Bigr\>.
\end{multline}

The three terms appearing in the right-hand side
of~\eqref{123} will be considered separately. First,
by integration by parts and~\eqref{DnablaV},
\begin{align}
- \int \nabla\cdot (f D\nabla u)\, \<Cu,\,C'u\> & =
 \int f\, \Bigl\<D\nabla u,\, \nabla\<Cu,\, C'u\>_{\R^m}\Bigr\>_{\R^N}\\
& =  \int f\, \Bigl\<Au,\, A\<Cu,\, C'u\> \Bigr\>.
\end{align}

For the second term in~\eqref{123}, we use the identity
\[ \nabla\cdot (Df\nabla u) = 
f\nabla\cdot (D\nabla u) + \<D \nabla u,\, \nabla f\>
= f \nabla \cdot (D\nabla u) + f \<D \nabla u,\, \nabla \log f\>.\]
So
\begin{multline*} 
- \int f\, \Bigl\< \frac{C \nabla\cdot (Df\nabla u)}{f},\, C'u\Bigr\>
= - \int f\, \Bigl\< C\nabla\cdot (D\nabla u),\, C'u\Bigr\>
\\- \int f\, \Bigl\<C \<D \nabla u,\, \nabla \log f\>_{\R^N},\,
C'u\Bigr\>_{\R^m} 
\end{multline*}
\begin{multline*} = -\int f \Bigl\<  C\nabla\cdot (D\nabla u),\, C'u\Bigr\>
- \int f\, \Bigl\<C \<D \nabla u,\, \nabla u\>_{\R^N},\, C'u\Bigr\>_{\R^m}
\\ + \int f\, \Bigl\<C \<D \nabla u,\, \nabla E\>_{\R^N},\, C'u\Bigr\>_{\R^m}.
\end{multline*}
By combining the first and third integrals in the expression above,
then using~\eqref{Dnablau} and~\eqref{DnablaV} again, we find that
\begin{align}
& -\int f\, \Bigl\< \frac{C \nabla\cdot (Df\nabla u)}{f},\, C'u\Bigr\> 
\nonumber\\
& = -\int f\, \Bigl\< C \bigl( \nabla\cdot(D\nabla u) -
\<D\nabla E,\,\nabla u\>_{\R^N}\bigr),\, C'u\Bigr\>_{\R^m}
- \int f\, \Bigl\< C\<D\nabla u,\nabla u\>_{\R^N},\, C'u\Bigr\>_{\R^m}
\nonumber\\
& = \int f\,\< C A^*\!Au,\, C'u\> -
\int f\,\< C |Au|^2, \,C'u\> \nonumber\\
& = \int f\, \< [C,A^*]Au,\, C'u\> +
\int f\, \<A^*\!CAu,\,C'u\> - 2 \int f\, \<(C Au)\cdot (Au),\, C'u\>.
\label{lastbeforerewr}
\end{align}
(In the last term, the dot is just here to indicate the
evaluation of the matrix $CAu$ on the vector $Au$.
Also $\<A^*CAu, C'u\>$ should be understood as
$\sum_{ij} \<A_i^*C_jA_iu, C'_ju\>$.)

Now the second integral in~\eqref{lastbeforerewr} needs some rewriting.
By using the chain rule as before, and the definition of the adjoint,
\begin{align}
\int f\, \<A^*\!CAu,\,C'u\> & =
\int \< A^*\!CAu,\, h C'\log h\>\,d\mu\\
& = \int \<A^*\!CAu,\, C'h\>\,d\mu \nonumber \\ 
& = \int \<CAu,\,AC'h\>\,d\mu \nonumber\\
& = \int \<CAu,\,[A,C']h\>\,d\mu +
\int \<CAu,\, C'Ah\>\,d\mu. \label{CAuC'Ah}
\end{align}

The first term in~\eqref{CAuC'Ah} can be rewritten as
\begeq
\int \<CAu,\, [A,C']u\>\,h\,d\mu = \int f\, \< CAu,\, [A,C']u\>.
\endeq
As for the second term in~\eqref{CAuC'Ah}, since $C'$ is a derivation,
it can be recast as
\begin{align}
\int \<CAu,\,C'(hAu)\>\,d\mu  \nonumber
& = \int \<CAu,\, hC'Au\>\,d\mu + \int \<CAu,\, (C'h)\otimes Au\>\,d\mu \\
& = \int f\,\< CAu,\, C'Au\> + \int f\,\<(CAu)\cdot (Au),\, C'u\>.
\end{align}
Note that there is a partial simplification with the last term
of~\eqref{lastbeforerewr} (only partial since the coefficients
are not the same).

Of course, the expressions which we obtained for the second term 
in~\eqref{123} also hold for the third term, up to the exchange
of $C$ and $C'$. After gathering all these results, we obtain
\begin{align}
-\left(\frac{d}{dt}\right)_{A^*A} \int f\<Cu,\,C'u\> =
& \int f\, \bigl\<Au,\,A\<Cu,\,C'u\>\bigr\> \label{l1} \\
& + \left( \int f\, \< [C,A^*]Au,\, C'u\> +
\int f\,\<Cu,\, [C',A^*]Au\> \right) \label{l2} \\
& + \left( \int f\,\< CAu,\, [A,C']u\>  + 
\int f\,\<[A,C]u,\,C'Au\>\right) 
\label{l3}\\
& + 2 \int f \<CAu,\, C'Au\> \label{l4}\\
& - \left(\int f\, \<(CAu)\cdot Au,\, C'u\> +
\int f\,\<(C'Au)\cdot Cu,\, Au\> \right).\label{l5}
\end{align}

The terms appearing in~\eqref{l2}, \eqref{l3} and~\eqref{l4}
coincide with some of the ones which appear in~\eqref{ddtIA}, so
it only remains to check that the ones in~\eqref{l1} and~\eqref{l5}
add up to~\eqref{QACC}. By using the identity
\[ A\<Cu,\,C'u\> = (ACu)\cdot (C'u) + (AC'u)\cdot (Cu),\]
we see that the sum of~\eqref{l1} and~\eqref{l5} can be recast as
\begin{multline*} 
\int f\, \Bigl\< Au,\, \<(AC-CA)u,\, C'u\>\Bigr\> +
\int f\, \Bigl\< Au,\, \<(AC'-C'A)u,\, Cu\>\Bigr\>\\
= \int f\,  \Bigl\< Au,\, \<[A,C]u,\, C'u\>\Bigr\> +
\int f\, \Bigl\< Au,\, \<[A,C']u,\, Cu\>\Bigr\>;
\end{multline*}
or, more explicitly:
\begeq\label{moreexplic} 
\sum_{ij} \int f\,  (A_iu)([A_i,C_j]u)(C'_ju)+
\sum_{ij}\int f\, (A_iu)([A_i,C'_j]u)(C_ju).
\endeq
It remains to check that~\eqref{moreexplic} can be transformed
into~\eqref{QACC}. Consider for instance the first term
in~\eqref{moreexplic}, for some index $(i,j)$. Since
$[A_i,C_j]$ is a derivation,
\begin{align*} \int f\,  (A_iu)([A_i,C_j]u)(C'_ju) 
& = \int (A_iu)([A_i,C_j]h)(C'_ju)\,d\mu\\
& = \int h [A_i,C_j]^*(A_iu\,C'_ju)\,d\mu \\
& = \int f [A_i,C_j]^*(A_iu\,C'_ju).
\end{align*}
This concludes the proof of Lemma~\ref{lemcalcHI}.
\end{proof}

\begin{proof}[Proof of Theorem~\ref{thmLlogL}]
Here I~shall use the same conventions as in the proof of Lemma~\ref{lemcalcHI}.
The functional ${\cal E}$ will be searched for in the form
\[ {\cal E}(h) = \int f u + \sum_{k=0}^{\Nc}
\left( a_k\, \int f |C_k u|^2 \ + \ 2 b_k\, \int f\<C_k u,\, C_{k+1} u\>
\right).\]
In other words, the quadratic form $S$ in the statement will be
looked for in the form
\[ \<S(x) \xi,\xi\>_m  = \sum a_k \bigl|C_k(x)\xi\bigr|_{\R^m}^2 + 
2 \sum b_k \Bigl\<C_k(x)\xi, C_{k+1}(x)\xi\Bigr\>_{\R^m},\]
where $C_k$ is identified with a function valued in
$m\times N$ matrices.

If the inequalities~\eqref{condakbk} are enforced, then for $\delta$ 
small enough
\[ \< S(x)\xi,\xi\>_m \geq K \sum |C_k(x)\xi|_m^2;\]
then $S$ will be a nonnegative symmetric matrix.

Next, we consider the evolution of ${\cal E}$ along the
semigroup. As recalled in Lemma~\ref{lemcalcHI},
\[ - \frac{d}{dt} \int h\log h\,d\mu = 
+ \int \frac{|Ah|^2}{h}\,d\mu. \]

Next,
\[ \begin{cases} \dps
- \frac{d}{dt} \int \frac{|C_kh|^2}{h}\,d\mu =
2\,\bigl(\I_A^k + \I_B^k\bigr),\\ \\
\dps - \frac{d}{dt} \int \frac{\<C_k h,\, C_{k+1}h\>}{h}\,d\mu
= \II_A^k + \II_B^k, 
\end{cases}\]
where the subscript $A$ indicates the contribution of the
$A^*A$ operator, and the subscript $B$ indicates the contribution
of the $B$ operator. The goal is to show that these terms
can be handled in exactly the same way as in Theorem~\ref{hypocomult}:
Everything can be controlled in terms of the quantities
\begeq\label{qties}
\int f |C_k Au|^2 \qquad \text{and} \qquad \int f|C_k u|^2
= \int \frac{|C_k h|^2}{h}\,d\mu.
\endeq
(These integrals play the role that the quantities
$\|C_k Ah\|^2$ and $\|C_kh\|^2$ were playing
in the proof of Theorem~\ref{hypocomult}.)

The terms $\I_B^k$ and $\II_B^k$ are most easily dealt with.
By Lemma~\ref{lemcalcHI}, we just have to reproduce the
result of the computations in the proof of Theorem~\ref{hypocomult}
and divide the integrand by $h$.
So in place of
\[ \int \<C_kh,\, Z_{k+1} C_{k+1}h\>\,d\mu + 
\int \< C_kh, \, R_{k+1}h\>\,d\mu,\]
we have
\[
\I_B^k = \int \frac{\<C_kh,\, Z_{k+1} C_{k+1}h\>}{h}\,d\mu 
+ \int \frac{\< C_kh, \, R_{k+1}h\>}{h}\,d\mu.\]
Then we proceed just as in the proof of Theorem~\ref{hypocomult}:
By Cauchy--Schwarz inequality (applied here for vector-valued
functions),
\begin{multline*} \I_B^k \geq - \Lambda_k 
\sqrt{\int \frac{|C_kh|^2}{h}\,d\mu}\,
\sqrt{\int \frac{|C_{k+1}h|^2}{h}\,d\mu} \\
- \sqrt{\int \frac{|C_kh|^2}{h}\,d\mu}\,
\sqrt{\int \frac{|R_{k+1}h|^2}{h}\,d\mu}.
\end{multline*}
Then $|R_{k+1}h|$ can be bounded {\em pointwise} in terms
of $|Ah|,\ldots,|C_kh|$, so $\int |R_{k+1}h|^2/h\,d\mu$
can be controlled in terms of $\int |C_jh|^2/h\,d\mu$
for $j\leq k$.

The treatment of $\II_B^k$ is similar:
\begin{multline*} \II_B^k = 
\int \frac{\<Z_{k+1} C_{k+1}h, C_{k+1}h\>}{h}\,d\mu +
\int \frac{\<R_{k+1}h,\, C_{k+1}h\>}{h}\,d\mu \\ +
\int \frac{\<C_k h,\, Z_{k+2} C_{k+2} h\>}{h}\,d\mu +
\int \frac{\<C_kh, \, R_{k+2}h\>}{h}\,d\mu
\end{multline*}
\begin{multline*}
\qquad \geq \lambda_{k+1}\,\int \frac{|C_{k+1}h|^2}{h}\,d\mu
- \sqrt{\int \frac{|R_{k+1}h|^2}{h}\,d\mu}\,
\sqrt{\int \frac{|C_{k+1}h|^2}{h}\,d\mu} \\
- \Lambda_{k+2} \sqrt{\int \frac{|C_kh|^2}{h}\,d\mu}\,
\sqrt{\int \frac{|C_{k+2}h|^2}{h}\,d\mu} 
- \sqrt{\int \frac{|C_kh|^2}{h}\,d\mu}\,
\sqrt{\int \frac{|R_{k+2}h|^2}{h}\,d\mu}.
\end{multline*}
Then once again, one can control the functions $|R_{k+2}h|$
by $|C_jh|$ for $j\leq k+1$.

Now consider the terms coming from the action of $A^*A$.
Let us first pretend that the extra terms $Q_{A,C,C'}$ in~\eqref{ddtIA}
do not exist.
Then by Lemma~\ref{lemcalcHI} again,
\[ \I_A^k = \int f |C_kAu|^2 + \int f \<[C_k,A^*] Au, \, C_k u\>
+ \int f \<C_kAu, [A,C_k]u\>.\]
By Cauchy--Schwarz inequality (for vector-valued functions),
\begin{multline*} \I_A^k \geq \int f |C_kAu|^2 - 
\sqrt{ \int f \bigl| [C_k,A^*]Au\bigr|^2}{\sqrt{\int f |C_ku|^2}}\\
- \sqrt{\int f |C_kAu|^2}\sqrt{\int f \bigl| [A,C_k]u\bigr|^2}.
\end{multline*}
Then Assumption (iii) implies
\begin{multline} \label{assiiiimplies}
\sqrt{\int f \bigl| [C_k,A^*]Au\bigr|^2} \leq
c \Bigl( \sqrt{\int f |A^2u|^2} + \sqrt{\int f |A C_1u|^2} + \ldots \\
+ \sqrt{\int f|AC_ku|^2} \Bigr).
\end{multline}

Finally,
\begin{multline*} \II_A^k = 2 \int f\<C_kA u, C_{k+1}Au\>+
\int f \bigl\< [C_k,A^*]Au,\, C_{k+1}u\bigr\> + 
\int f\bigl\<C_k u,\, [C_{k+1},A^*]Au\bigr\> \\
+ \int f \bigl\< C_kAu, [A,C_{k+1}]u\bigr\>
+ \int f \bigl\< [A,C_k]u,\, C_{k+1}Au\bigr\>.
\end{multline*}
and this can be bounded below by a negative multiple of
\begin{multline*} - \sqrt{\int f |C_kAu|^2}\,\sqrt{\int f|C_{k+1}Au|^2}
- \sqrt{\int f |C_ku|^2}\sqrt{\int f \bigl| [C_{k+1},A^*]Au\bigr|^2}\\
- \sqrt{\int f |C_{k}Au|^2}\sqrt{\int f \bigl| [A,C_{k+1}]u\bigr|^2}
- \sqrt{\int \bigl| [A,C_k]u\bigr|^2}
\sqrt{\int f|C_{k+1}Au|^2};
\end{multline*}
then one can apply~\eqref{assiiiimplies} (as it is, and also
with $k$ replaced by $k+1$) to control the various terms above.

All in all, everything can be bounded in terms of
the integrals appearing in~\eqref{qties}, and the computations
are {\em exactly} the same as in the proof of Theorem~\ref{hypocomult};
then the same bounds as in Theorem~\ref{hypocomult} will work,
provided that the coefficients $a_k$ and $b_k$ are well chosen.
The result is
\begeq\label{prelimddtE} \frac{d}{dt} {\cal E}(h) \leq
- K \int \frac{\<S(x) \nabla h(x), \nabla h(x)\>}{h(x)}\,d\mu(x).
\endeq

Now let us see what happens if Assumptions (a) and (b) are enforced.
By assumption (a), we have $\sum |C_k(x)\xi|^2 \geq \lambda |\xi|^2$,
where $\lambda>0$; so there exists $\kappa>0$ such that
\[ {\cal E}(h) \geq \int f u + \kappa \int f |\nabla u|^2 =
\int h\log h\,d\mu + \kappa \int \frac{|\nabla h|^2}{h}\,d\mu.\]
Thus ${\cal E}$ will dominate both the Kullback information
$H_\mu(h\mu)$, and the Fisher information $I_\mu(h\mu)$.

Then, since $S$ is uniformly positive definite,
\[ \int \frac{\<S(x) \nabla h(x), \nabla h(x)\>}{h(x)}\,d\mu(x)
\geq \lambda \int \frac{|\nabla h|^2}{h}\,d\mu = \lambda I_\mu(h\mu).\]
As a consequence, by Assumption (b),
\[ \int \frac{\<S(x) \nabla h(x), \nabla h(x)\>}{h(x)}\,d\mu(x)\geq
\kappa\, H_\mu(h\mu) = \kappa \int fu\]
for some $\kappa>0$. So the right-hand side of~\eqref{prelimddtE} 
controls also $\int fu$, and in fact there is a positive constant $\kappa$
such that
\[ \frac{d}{dt} {\cal E}(h) \leq - \kappa\, {\cal E}(h).\]
Then we can apply Gronwall's inequality to conclude
the proof of Theorem~\ref{thmLlogL}.
\med

It remains to take into account the additional terms generated by
$Q_{A,C,C'}$ in~\eqref{ddtIA}. More precisely, in $\I_A^k$ we should
consider $\int f\, Q_{A,C_k,C_k} (u)$; and in $\II_B^k$ we should
handle $\int f\, Q_{A,C_k,C_{k+1}}(u)$.
So the problem is to bound also these expressions in terms
of the quantities~\eqref{qties}.

We start with the additional term in $\I_A^k$, that is,
\begeq\label{addtermIAk} 
\int f\, Q_{A,C_k,C_k}(u) =\ \int f [A,C_k]^* (Au\otimes C_ku).
\endeq
By assumption $[A,C_k]^*$ is controlled by $I$ and $A$, so there
is a constant $c$ such that
\begeq\label{cIA}
\left|\int [A,C_k]^* (Au\otimes C_ku)\right|\leq
c \left( \int f\,\bigl|A(Au\otimes C_ku)\bigr|\ 
+ \ \int f\,\bigl|Au\otimes C_ku\bigr|\right).
\endeq
Next, by the rules of derivation of products,
\begin{align*} 
A(AuC_ku) & = (A^2u)(C_ku) + (Au)(AC_ku)\\
& = (A^2u)(C_ku) + (Au)(C_kAu) + (Au)([A,C_k]u).
\end{align*}
Here as in the sequel, I~have omitted indices for simplicity; the above
equation should be understood as
$A_\ell (A_iu\,C_{k,j}u) = (A_\ell A_i u)(C_{k,j}u) +
(A_\ell u)(C_{k,j}A_iu) + (A_\ell u)([A_i,C_{k,j}]u)$.
Since by assumption $[A,C_k]$ is controlled by $\{C_j\}_{0\leq j\leq k}$,
there exists some constant $c$ such that the following pointwise bounds holds:
\[ |A(Au\otimes C_ku)| \leq 
c \left(|A^2u|\, |C_ku| + |Au|\,|C_kAu| + |Au|^2
+ \sum_{0\leq j\leq k} |Au|\,|C_ju|\right).\]
Plugging this in~\eqref{cIA} and then in~\eqref{addtermIAk}, 
then using the Cauchy--Schwarz inequality,
we end up with
\begin{multline*}
\left| \int  f\, Q_{A,C_k,C_k}(u) \right| \leq
c\left( \sqrt{\int f\, |A^2u|^2} \sqrt{\int f|C_ku|^2} +
\sqrt{\int f\,|Au|^2}\sqrt{\int f|C_kAu|^2} \right. \\ \left.
+ \sum_{j\leq k} \sqrt{\int f\,|Au|^2}\sqrt{\int f |C_ju|^2} 
+ \sqrt{\int f\,|Au|^2}\sqrt{\int f|C_ku|^2}\right).
\end{multline*}
All these terms appear in $\I_A^k$ with a coefficient. So they
can be controlled in terms of~\eqref{qties}, with the right coefficients,
as in the proof of Theorem~\ref{hypocomult}, if
\[ a_k \ll \max (\sqrt{a_0\, b_{k-1}}, \ \sqrt{a_k}, \ 1, \
\max_{j\leq k} \sqrt{b_{j-1}}, \ \sqrt{b_{k-1}}).\]
These conditions are enforced by the construction of the coefficients
$(a_j)$ and $(b_j)$.
\med

Now we proceed similarly for the additional terms in $\II_B^k$.
By repeating the same calculations as above, we find
\begin{multline*}
\left|\int f\, Q_{A,C_k,C_{k+1}}(u)\right|\leq
c \left( \sqrt{\int f\, |Au|^2}\sqrt{\int f\,|C_{k+1}u|^2} +
\sqrt{\int f\, |A^2u|^2} \sqrt{\int f\,|C_{k+1}u|^2} \right. \\
\left.
+ \sqrt{\int f\, |Au|^2}\sqrt{\int f\,|C_{k+1}Au|^2}
+ \sum_{j\leq k+1} \sqrt{\int f |Au|^2}
\sqrt{\int f|C_ju|^2} + \right. \\ \left.
+ \sqrt{\int f\,|Au|^2}\sqrt{\int f\, |C_ku|^2}
+ \sqrt{\int f\,|A^2u|^2}\sqrt{\int f\,|C_ku|^2}
+ \sqrt{\int f\, |Au|^2}\sqrt{\int f\,|C_kAu|^2}\right)
\end{multline*}
All these terms come with a coefficient $b_k$, and they are properly
controlled by~\eqref{qties} if
\[ b_k \ll \max (\sqrt{b_k},\ \sqrt{a_0\,b_k}, \max_{j\leq k+1}
\sqrt{b_{j-1}}, \ \sqrt{a_1}).\]
Again, these estimates are enforced by construction.
This concludes the proof of Theorem~\ref{thmLlogL}.
\end{proof}

\section{Application: the kinetic Fokker--Planck equation} \label{secFP}

In this section I~shall apply the preceding results to the kinetic linear 
Fokker--Planck equation, which motivated
and inspired the proof of Theorem~\ref{thmsimple} as well as previous 
works~\cite{DV:FP:01,heraunier:FP:04,helfnier:witten:05}.

The equation to be studied is~\eqref{kFPh}, which I~recast here:
\begeq\label{kFPh'} 
\derpar{h}{t} + v\cdot\nabla_x h -\nabla V(x)\cdot\nabla_v h = 
\Delta_v h - v\cdot\nabla_v h;
\endeq
and the equilibrium measure takes the form
\[ \mu(dx\,dv) = \gamma(v) e^{-V(x)}\,dv\,dx,\qquad 
\gamma(v) = \frac{e^{-\frac{|v|^2}{2}}}{(2\pi)^{n/2}},\qquad \int e^{-V}=1. \]

Let $\H:=L^2(\mu)$, $A:=\nabla_v$, 
$B:= v\cdot\nabla_x - \nabla V(x)\cdot\nabla_v$.
Then~\eqref{kFPh'} takes the form $\partial h /\partial t +L h=0$, with
$L=A^*\! A+B$, $B^*=-B$. The kernel $\K$ of $L$ is 
made of constant functions, and the space $\H^1=H^1(\mu)$ is the 
usual $L^2$-Sobolev space of order~1, with derivatives
in both $x$ and $v$ variables, and reference weight $\mu$:
\[ \|h\|_{\H^1}^2 = \int_{\R^n\times\R^n} \Bigl(|\nabla_v h(x,v)|^2
+|\nabla_x h(x,v)|^2\Bigr) \,\mu(dx\,dv).\]
By direct computation,
\[ [A,A^*] = I, \qquad C := [A,B] = \nabla_x, \qquad [A,C]=[A^*,C]=0, \]
\[ \qquad [B,C] = \nabla^2V(x)\cdot\nabla_v.\]

\subsection{Convergence to equilibrium in $H^1$}

In the present case, assumptions (i)--(iii) of Theorem~\ref{thmsimple} 
are satisfied if
\begeq\label{condthmsimpleFP} 
\text{$\nabla^2 V$ is relatively bounded by $\{I,\nabla_x\}$ in $L^2(e^{-V})$.}
\endeq
By Lemma~\ref{lemD2Vbdd} in Appendix~\ref{toolbox}, 
this is true as soon as there exists a constant $c\geq 0$ such that
\begeq \label{condD2V}
|\nabla^2 V| \leq c (1+|\nabla V|).
\endeq

The other thing that we should check is the coercivity of $A^*\! A+C^*C$,
which amounts to the validity of a Poincar\'e inequality of the form
\begeq\label{Poinctens} 
\int (|\nabla_v h|^2 + |\nabla_x h|^2)\,d\mu \geq \kappa
\left[ \int h^2\,d\mu - \left(\int h\,d\mu\right)\right]^2.
\endeq
Since $\mu$ is the tensor product of a Gaussian distribution in $\R^n_v$
(for which the Poincar\'e inequality holds true with constant~1) and of the distribution
$e^{-V}$ in $\R^n_x$, the validity of~\eqref{Poinctens} is equivalent to the
validity of a Poincar\'e inequality (in $\R^n_x$)
\begeq \label{Poincx}
\int |\nabla_x h(x)|^2\,e^{-V(x)}\,dx \geq \lambda 
\left[\int h^2\,e^{-V} - \left(\int h e^{-V}\right)\right]^2.
\endeq
This functional inequality has been studied by many many authors, and it is
natural to take it as an assumption in itself.
Roughly speaking, inequality~\eqref{Poincx} needs $V$ to 
grow ``at least linearly'' at infinity. In Theorem~\ref{thmDS} in 
Appendix~\ref{appproduct} I~recall a rather general sufficient condition 
for~\eqref{Poincx} to be satisfied; it holds true for instance
if~\eqref{condD2V} is true and $|\nabla V(x)|\to \infty$ at infinity.
Then Theorem~\ref{thmsimple} leads to the following statement:

\begin{Thm}
\label{thmFPH1}
Let $V$ be a $C^2$ potential in $\R^n$, satisfying conditions~\eqref{condD2V}
and~\eqref{Poincx}. Then, with the above notation, there are constants $C\geq 0$ 
and $\lambda>0$, explicitly computable, such that for all $h_0\in H^1(\mu)$,
\[ \left\| e^{-tL} h_0 - \int h_0\,d\mu \right\|_{H^1(\mu)} 
\leq C e^{-\lambda t} \|h_0\|_{H^1(\mu)}. \]
\end{Thm}

\begin{Rk} Conditions~\eqref{condD2V} and~\eqref{Poincx}
morally mean that the potential $V$ should grow at least linearly, and at most
exponentially fast at infinity. These conditions are more general than those imposed
by Helffer and Nier~\cite{helfnier:witten:05}
\footnote{This comparison should not hide the fact that the estimates by
Helffer and Nier were already remarkably general, and constituted
a motivation for the genesis of this paper.}
in that no regularity at order higher than~2
is needed, and there is no restriction of polynomial growth on $V$.
Here is a more precise comparison: Helffer and Nier prove exponential convergence
under two sets of assumptions: on one hand, \cite[Assumption 5.6]{helfnier:witten:05};
on the other hand, \cite[Assumption 5.7]{helfnier:witten:05} plus a spectral gap condition
which is equivalent to~\eqref{Poincx}. Both these assumptions 5.6 and 5.7 
contain~\cite[eq.(5.17)]{helfnier:witten:05}, which is stronger than~\eqref{condD2V}.
Finally, the spectral gap condition is not made explicitly in 
\cite[Assumption 5.6]{helfnier:witten:05},
but is actually a consequence of that assumption, since it implies~\eqref{DS}.
\end{Rk}

\begin{proof}[Proof of Theorem~\ref{thmFPH1}]
We already checked all the assumptions of Theorem~\ref{thmsimple},
except for the existence of a convenient dense subspace ${\cal S}$.
If $V$ is $C^\infty$, it is possible to choose
the space of all $C^\infty$ functions on $\R^n_x\times\R^n_v$ whose derivatives
of all orders vanish at infinity faster than any inverse power of 
$(1+|\nabla V|)(1+|v|)$. (Note that the operators appearing in the theorem
preserve this space because $|\nabla^2V|$ is bounded by a multiple
of $1+|\nabla V|$.)
Then there only remains the problem of approximating
$V$ by a $C^\infty$ potential, without damaging Condition~\eqref{condD2V}.
This can be done by a standard convolution argument: 
let $V_\var:=V\ast\rho_\var$, where $\rho_\var(x)=
\var^{-n}\rho(x/\var)$, and $\rho$ is $C^\infty$, supported in the unit ball,
nonnegative and of unit integral. Then, for all $\var>0$,
\[ |\nabla^2 V_\var(x)| \leq \sup_{|x-y|\leq \var} |\nabla^2 V(y)| 
\leq C \sup_{|x-y|\leq \var} (1 + |\nabla V(y)|).\]
But~\eqref{condD2V} implies that $\log (1+|\nabla V|^2)$ is $L$-Lipschitz
($L=2C$), so
\[ |x-y|\leq \var \Longrightarrow \qquad
\frac{1+|\nabla V(x)|^2}{1+|\nabla V(y)|^2} \leq e^{L\var}. \]
In particular, by~\eqref{condD2V}, $|\nabla^2V(y)|$ can be controlled in terms 
of $|\nabla V(x)|$, for $y$ close to $x$. It follows
\[ |x-y|\leq \var\Longrightarrow |\nabla V(x) - \nabla V(y)| \leq 
C (1+ |\nabla V(x)|) \var e^{L\var}.\]
As a consequence,
\[ |\nabla V_\var(x) - \nabla V(x)| \leq C (1+ |\nabla V(x)|)\var e^{L\var}.\]
From this it is easy to deduce that $1+|\nabla V_\var(x)| \geq (1-C'\var) (1+|\nabla V(x)|)$,
for some explicit constant $C'$. Then, $V_\var$ satisfies the same condition~\eqref{condD2V}
as $V$, up to replacing the constant $C$ by some constant $\tilde{C}(\var)$ which
converges to $C$ as $\var\to 0$. 

All in all, after replacing $V$ by $V_\var$, we can apply the first part of
Theorem~\ref{thmsimple} to get
\begin{multline}\label{iipfvar}
\iip{h_\var(t)}{h_\var(t)} + K \int_0^t \left(\int (|\nabla_x h_\var(s)|^2
+ |\nabla_v h_\var(s)|^2)e^{-V_\var(x)} \gamma(v)\right)\,ds\\
\leq \iip{h_0}{h_0}, 
\end{multline}
where $h_\var(t)=e^{-tL_\var}(h_0-\int h_0)$, and $K$ is a 
constant independent of $\var$.

By the uniqueness theorem of Appendix~\ref{appuniqueFP}, $h_\var(t)$ converges
to $h(t)=e^{-tL}h$, in distributional sense as $\var\to 0$.
Also, $\int h_0\, e^{-V_\var}\gamma \longrightarrow \int h_0\, e^{-V}\gamma$.
Since the left-hand side is a convex functional of $h$ and $V_\var$ converges
locally uniformly to $V$, inequality~\eqref{iipfvar} passes to the limit
as $\var\to 0$.
The Poincar\'e inequality for $e^{-V}$ and the definition of the
auxiliary scalar product guarantee the existence of $K'>0$ such that
\[ \iip{h(t)}{h(t)} + K' \int_0^t \iip{h(s)}{h(s)}\,ds \leq 
\iip{h_0}{h_0}. \]
The exponential convergence of $\iip{h(t)}{h(t)}$ to~0 follows, and the theorem
is proved.
\end{proof}

\subsection{Explicit estimates}

As a crude test of the effectiveness of the method, one can repeat the
proof of Theorem~\ref{thmsimple} on the particular example of the Fokker--Planck
equation, taking advantage of the extra structure to get more precise results. 
Using $[A,A^*]=I$, one obtains
\begeq\label{xLxsimpl} 
\iip{h}{Lh} \geq \|Ah\|^2 + a\|A^2h\|^2 +b \|Ch\|^2 + c\|CAh\|^2 - (E), 
\endeq
\begin{multline*} (E):= a (\|Ah\|^2 + \|Ah\|\,\|Ch\| ) + b (\|Ah\|\,\|R_2h\| + 2 \|A^2h\|\,\|CAh\| 
+ \|Ah\|\,\|Ch\|) \\ + c\|Ch\|\,\|R_2h\|.
\end{multline*}
Moreover, $R_2h=-[B,C]=\nabla^2V\cdot A$; to simplify computations even more, assume
that $|\nabla^2 V|\leq M$ (in Hilbert-Schmidt norm, pointwise on $\R^n$). Then
\begin{align*}
(E) & \leq a (\|Ah\|^2 + \|Ah\|\,\|Ch\|) + b\Bigl( M\|Ah\|^2 
+ 2\|A^2h\|\,\|CAh\| + \|Ah\|\,\|Ch\|\Bigr) \\
& \hspace*{80mm} + cM \|Ah\|\,\|Ch\| \\
& = (a+bM) \|Ah\|^2 + (a+b+cM) \|Ah\|\,\|Ch\| + 2b \|A^2h\|\,\|CAh\| \\
& \leq (a+bM+1/4) \|Ah\|^2 + (a+b+cM)^2 \|Ch\|^2 + a\|A^2h\|^2 + \frac{b^2}{a} \|CAh\|^2.
\end{align*}
Since $b^2/a\leq c$, the last two terms above are bounded by the terms in $\|A^2h\|^2$
and $\|CAh\|^2$ in~\eqref{xLxsimpl}; so
\[ \iip{h}{Lh}  \geq \left[ 1- (a+bM+1/4)\right] \|Ah\|^2 + \left[ b - (a+b+cM)^2\right] \|Ch\|^2 \]

On the other hand, taking into account the spectral gap assumption on $A^*\! A+C^*C$,
\[ \iip{h}{h} \le (2a+\kappa^{-1}) \|Ah\|^2 + (2c+\kappa^{-1}) \|Ch\|^2. \]
So the proof yields a convergence to equilibrium in $H^1$ 
like $O(e^{-\ov{\lambda} t})$, where
\[ \ov{\lambda}:= \sup_{(a,b,c)}\ 
\min \left( \frac{1-(a+bM+\frac14)}{2a+\kappa^{-1}}, \
\frac{b-(a+b+cM)^2}{2c+\kappa^{-1}}\right), \]
and the supremum is taken over all triples $(a,b,c)$ with $b^2\leq ac$.

In the particular (quadratic) case where $\nabla^2V$ is the identity, one has $M=1$,
$\kappa=1$; then the choice $a=b=c=0.05$ yields $\ov{\lambda}=0.025$, which is off the
true (computable) rate of convergence to equilibrium $\lambda=1/2$ 
(see~\cite[p.~238--239]{risken:book}) by a factor~20. 
Thus, even if the method is not extremely sharp,
it does yield quite decent estimates.\footnote{As a comparison, the bounds
by H\'erau and Nier~\cite[formula (4)]{heraunier:FP:04} yield a lower bound
on $\lambda$ which is around $10^{-4}$.} 
Note that the coefficients $a,b,c$ chosen 
in the end do not satisfy $c\ll b\ll a$!

\subsection{Convergence in $L^2$}

Theorem~\ref{thmFPH1} is stated for $H^1$ initial data. However, it can be combined
with an independent regularity study: Under condition~\eqref{condD2V},
one can show that solutions of~\eqref{kFPh'} satisfy the estimate
\begeq\label{estregFP} 0\leq t\leq 1\Longrightarrow\qquad
\|f(t,\cdot)\|_{H^1(\mu)} \leq \frac{C}{t^{3/2}}\|f(0,\cdot)\|_{L^2(\mu)}.
\endeq
A proof is provided in Appendix~\ref{appreg}.
Combined with Theorem~\ref{thmFPH1}, 
this estimate trivially leads to the following statement:

\begin{Thm} \label{thmFPL2}
Let $V$ be a $C^2$ potential in $\R^n$, satisfying conditions~\eqref{condD2V}
and~\eqref{Poincx}. Then, with the above notation, there are constants $C\geq 0$ 
and $\lambda>0$, explicitly computable, such that for all $h_0\in L^2(\mu)$,
\[ t\geq 1 \Longrightarrow\qquad
\left\| e^{-tL} h_0 - \int h_0\,d\mu \right\|_{H^1(\mu)} \leq C e^{-\lambda t} \|h_0\|_{L^2(\mu)}. \]
\end{Thm}

\subsection{Convergence for probability densities}

Write $e^{-V}(x)\gamma(v)  = e^{-E(x,v)}$, and set $f=e^{-E}h$, 
then the Fokker--Planck equation~\eqref{kFPh'} becomes the kinetic equation 
for the density of particles:
\begeq\label{kFPf'}
\derpar{f}{t} + v\cdot\nabla_x f - \nabla V(x)\cdot\nabla_v h =
\nabla_v\cdot (\nabla_v f + fv).
\endeq

The previous results show that there is exponential convergence to
equilibrium as soon as
\[ \int f^2 e^{E}\,dx\,dv < +\infty.\]
As an integrability estimate, this assumption is not very natural
for a probability density; as a decay estimate at infinity, it is extremely 
strong. The goal now is to establish convergence to equilibrium under
much less stringent assumptions on the initial data, maybe at the price
of stronger assumptions on the potential $V$.

An obvious approach to this problem consists in using stronger hypoelliptic
regularization theorems. For instance, it was shown by H\'erau and 
Nier~\cite{heraunier:FP:04} that if the initial datum in~\eqref{kFPh'} takes the form
is only assumed to be a tempered distribution, then the solution at later times
lies in $L^2(\mu)$, and in fact takes the form $\sqrt{e^{-E}} g$, where
$g$ is $C^\infty$ with rapid decay. Similar results can also be shown by variants of
the method exposed in Appendix~\ref{appreg}; for instance one may show that
if the initial datum belongs to a negative $L^2$-Sobolev space of order $k$
then for positive times the solution belongs to a positive $L^2$-Sobolev space
of order $k'$, whatever $k$ and $k'$. In particular, this approach works fine if 
the initial datum for~\eqref{kFPf'} is a probability measure $f_0$ satisfying
\[ \int e^{E(x,v)/2} f_0(dx\,dv) <+\infty.\]
However this still does not tell anything if we assume only polynomial moment
bounds on $f_0$. 

In the next result (apparently the first of its kind),
this problem will be solved with the help of Theorem~\ref{thmLlogL},
that is, by using an entropy approach.

\begin{Thm} \label{thmmesdata}
Assume that $V$ is $C^\infty$ with $|\nabla^j V(x)|\leq C_j$ for all $j\geq 2$;
and that the reference measure $e^{-V}$ satisfies a logarithmic Sobolev inequality.
Let $f_0$ be a probability measure with polynomial moments of all orders:
\[ \forall k\in\N,\qquad \int (1+|x|+|v|)^k\, f_0(dx\,dv) < +\infty.\]
Then the solution to~\eqref{kFPf'} is $C^\infty$ in $x$ and $v$ for
all positive times, and converges to $e^{-E}$ exponentially fast
as $t\to\infty$, in the sense that
\[ \int f(t,x,v)\log \left(\frac{f(t,x,v)}{e^{-E(x,v)}}\right)\,dx\,dv \ 
= O(e^{-\alpha t})\qquad (t\geq 1),\]
with explicit estimates.
\end{Thm}

\begin{Rk} A well-known sufficient condition for $e^{-V}$ to satisfy 
a logarithmic Sobolev inequality is $V=W+w$, where $\nabla^2 W\geq \kappa I_n$,
$\kappa>0$, and $w$ is bounded (this is the so-called
``uniformly convex + bounded'' setting).
\end{Rk}

\begin{proof}[Proof of Theorem~\ref{thmmesdata}]
Since by assumption the function $\nabla V(x)$ is Lipschitz by assumption,
it can be shown by standard techniques that the Fokker--Planck equation 
admits a unique measure-valued solution. So it is sufficient to
establish the convergence for very smooth initial data, with rates
that do not depend on the smoothness of the initial datum, and then
use a density argument.

Since $\nabla^2 V$ is bounded, the transport coefficients appearing
in~\eqref{kFPf'} are Lipschitz (uniformly for $(x,v)\in\R^n_x\times\R^n_v$),
and it is easy to show by classical estimates that all moments
increase at most linearly in time:
\[ \int (1+|v|^2+|x|^2)^{k/2} f(t,dx\,dv) = O (1+t).\]

It is shown in Appendix~\ref{appreg} that $f(t,\cdot)$ also belongs to 
all Sobolev spaces (in $x$ and $v$) for $t>0$; in fact,
estimates of the form
\[ \|f(t,\cdot)\|_{H^k_x H^\ell_v(\R^n_x\times\R^n_v)} 
= O(t^{-\beta(k,\ell)})\qquad 0< t\leq 1\]
will be established there. Then by elementary interpolation, 
$f(t,\cdot)$ lies in all weighted Sobolev spaces
for all $t\in (0,1)$: That is,
\[ \|f(t,\cdot)\|_{H^k_s} := 
\| f(t,x,v) (1+|v|^2+|x|^2)^{s/2}\|_{H^k} < +\infty.\]

It is shown in~\cite[Lemma~1]{TV:slow:00} that
$I(f)\leq C \|f\|_{H^k_s}$ for $k$ and $s$ large enough (depending on $n$),
where $I(f)$ stands for the Fisher information, $\int f |\nabla (\log f)|^2/f$.
So $f(t,\cdot)$ has a finite Fisher information (in both $x$ and $v$ variables)
for all $t>0$. Since also $f(t,\cdot)$ has all its moments bounded and
$|\nabla E|=O(1+|x|+|v|)$, we have in fact
\[ \int f \bigl|\nabla (\log f + E)\bigr|^2 = O(t^{-\gamma})\qquad 0<t\leq 1\]
for some $\gamma>0$, where the time variable is omitted in the
left-hand side. So from time $t=t_0>0$ on, the solution $f$ has 
a finite relative Fisher information with reference measure 
$\mu(dx\,dv)=e^{-E(x,v)}\,dx\,dv$.

Then we can apply Theorem~\ref{thmLlogL} with $A=\nabla_v$,
$B=v\cdot\nabla_x - \nabla V(x)\cdot\nabla_v$,
$C_1=[A,B]=\nabla_x$, $R_1=0$, $C_2=0$, $R_2=\nabla^2V(x)\cdot\nabla_v$,
$Z_j=I$. Assumptions (i), (ii), (iii) and (v) in Theorem~\ref{thmLlogL} are
automatically satisfied, and Assumption (iv) is also satisfied
since $\nabla^2 V$ is bounded. (This is the place where the boundedness
of the Hessian of $V$ is crucially used.) Since the reference measure $\mu$
is the product of $e^{-V(x)}\,dx$ (which satisfies a logarithmic Sobolev
inequality by assumption) with $\gamma(v)\,dv$ (which also satisfies
a logarithmic Sobolev inequality), $\mu$ itself satisfies a logarithmic
Sobolev inequality. So Theorem~\ref{thmLlogL} yields the estimate
\[ \int f (\log f +E) = O(e^{-\lambda (t-t_0)})\]
for $t\geq t_0$. In words: The relative entropy of the solution
with respect to the equilibrium measure converges to~0 exponentially
fast as $t\to\infty$. This concludes the proof of Theorem~\ref{thmmesdata}.
\end{proof}

\section{The method of multipliers}

A crucial ingredient in the $L^2$ treatment of the Fokker--Planck equation was
the use of the mixed second derivative $CAh=\nabla_v\nabla_x h$ to control
the error term $[B,C]h = (\nabla^2 V)\cdot\nabla_v h$. There is an
alternative strategy, which does not need to use $CA$: It
consists in modifying the quadratic form~\eqref{eqnorm} thanks to
well-chosen auxiliary operators, typically multipliers. 
In the case of the Fokker--Planck equation, this method leads to less 
general results; it is however of independent interest, and can certainly
be applied to many equations. In this section I shall present a
variant of Theorem~\ref{thmsimple} allowing for multipliers, 
and test its applicability to the Fokker--Planck equation.
Some extensions are feasible, but I~shall not consider them.

Let again $A$ and $B$ be as in Subsection~\ref{subnot}, and
$C=[A,B]$, $R_2=[C,B]$; assume that $[A,C]=0$ for simplicity.
Let $M$, $N$ be two self-adjoint, invertible nonnegative operators such that
\[ [B,M]=0, \qquad [B,N]=0, \qquad [M,N] =0\]
(these conditions can be somewhat relaxed by imposing only an adequate
control on the commutators, but this leads to cumbersome calculations).
Instead of~\eqref{eqnorm}, consider the quadratic form
\begeq\label{eqnormMN}
\iip{h}{h} = \|h\|^2 + a\|MAh\|^2 +2b \ip{MAh}{NCh} + c\|NCh\|^2. 
\endeq
By straightforward variants of the calculations performed in the
proof of Theorem~\ref{thmsimple}, one obtains
\[
\iip{h}{Lh} = \|Ah\|^2 + a \|MA^2h\|^2 + b \|\sqrt{MN}Ch\|^2 
+ c\|NCAh\|^2 - (E),\] 
where
\begin{multline*} 
- (E) :=  a\Bigl( \ip{MAh}{MCh} + \ip{MAh}{[M,A^*]A^2h} + \ip{[A,M]Ah}{MA^2h} \\
+ \ip{MAh}{M[A,A^*]Ah} \Bigr) 
\end{multline*}
\begin{multline*}
+ b\Bigl( \ip{MAh}{N[B,C]h} + 2 \ip{MA^2h}{NCAh} 
+ \ip{[A,M]Ah}{NCAh} + \ip{MAh}{[N,A^*]CAh} \\ + \ip{[A,N]Ch}{MA^2h} + 
\ip{NCh}{[M,A^*]A^2h} + \ip{NCh}{M[A,A^*]Ah} \Bigr)
\end{multline*}
\[ + c\Bigl(\ip{NCh}{NR_2h} + \ip{[A,N]Ch}{NCAh} + \ip{NCh}{[N,A^*]CAh}\Bigr).\]

It would be a mistake to use Cauchy--Schwarz inequality right now. Instead, one should
first ``re-distribute'' the multipliers $M$ and $N$ on the two factors in the scalar
products above. For instance, $\ip{MAh}{MCh}$ is first rewritten $\ip{Ah}{M^2Ch}$
since it should be controlled by (inter alia) $\|Ah\|^2$, not $\|MAh\|^2$.
To obtain the correct weights, one is sometimes led to introduce the inverses
$M^{-1}$ and $N^{-1}$. In the end,
\begin{multline*}
(E) \leq a \Bigl ( \|Ah\|\,\|M^2Ch\| + \|Ah\|\,\|M[M,A^*]A^2h\| + \|[A,M]Ah\|\,\|MA^2h\| \\ 
+ \|Ah\|\,\|M^2[A,A^*]Ah\| \Bigr) 
\end{multline*}
\begin{multline*}
+ b \Bigl ( \|Ah\|\,\|MN[B,C]h\| + 2 \|MA^2h\|\,\|NCAh\| + \|[A,M]Ah\|\,\|NCAh\|
\\ + \|Ah\|\,\|M[N,A^*]CAh\| + \|(\sqrt{MN})Ch\|\,\|(\sqrt{N/M})[M,A^*]A^2h\| 
+ \|[A,N]Ch\|\,\|MA^2h\|\\
+ \|\sqrt{MN}Ch\|\,\|\sqrt{MN}[A,A^*]Ah\| \Bigr) \end{multline*}
\begin{multline*}
+ c \Bigl( 
+ \|(\sqrt{MN}) Ch\|\,\|(N^{3/2}/M^{1/2}) R_2h\| \\ + \|[A,N]Ch\|\,\|NCAh\| 
+ \|(\sqrt{MN})Ch\|\,\|(\sqrt{N/M})[N,A^*]CAh\| \Bigr).
\end{multline*}
Of course, in the above $\sqrt{N/M}$ stands for $N^{1/2}M^{-1/2}$, etc.

Repeating the scheme of the proof of Theorem~\ref{thmsimple}, it is easy to see that
$(E)$ can be controlled in a satisfactory way as soon as, say (conditions are listed
in order of appearance and the notation of Subsection~\ref{relbound} is used),
\[ M^2 \sle \sqrt{MN}, \qquad [M,A^*]\sle I, \qquad [A,M]\sle I, \qquad
M^2[A,A^*]\sle I, \]
\[ MN[B,C]\sle A, \qquad [A,M] \sle I, \qquad M [N,A^*]\sle N, \qquad
[A,N]\sle \sqrt{MN}, \]
\[(\sqrt{N/M}) [M,A^*] \sle M, \qquad \sqrt{MN} [A,A^*] \sle I, 
\qquad (N^{3/2}/M^{1/2}) [B,C]\sle A, \]
\[ [A,N] \sle \sqrt{MN}, \qquad (\sqrt{N/M})[N,A^*] \sle N. \]

For homogeneity reasons it is natural to assume $N=M^3$. Then the above conditions
are satisfied if
\begeq\label{condMN1}
M^2[A,A^*]\sle I, \qquad M^4[B,C] \sle A, 
\endeq
\begeq\label{condMN2}
[M,A]\sle I, \qquad [M,A^*] \sle I, \qquad [M^3,A]\sle M^2, \qquad [M^3, A^*]\sle M^2.
\endeq

If these conditions are satisfied, then one can repeat the scheme of the proof
of Theorem~\ref{thmsimple}, with an important difference: instead of
$\|Ah\|^2 + \|Ch\|^2$, it is only $\|Ah\|^2 + \|(\sqrt{MN})Ch\|^2
= \|Ah\|^2 + \|M^2Ch\|^2$ which is controlled in the end.
This leads to the following theorem.

\begin{Thm} \label{thmmultsimple} 
With the notation of Subsection~\ref{subnot}, assume that
\[ [C,A]=0, \qquad [C,A^*]=0, \]
and that there exists an invertible nonnegative self-adjoint bounded operator $M$ 
on $\H$, commuting with $B$, such that conditions~\eqref{condMN1} 
and~\eqref{condMN2} are fulfilled. Define
\begeq\label{eqnormM}
\iip{h}{h} = \|h\|^2 + a\|MAh\|^2 -2b \ip{M^2Ah}{Ch} + c\|M^3Ch\|^2. 
\endeq
Then, there exists $K>0$, only depending on the bounds appearing implicitly
in~\eqref{condMN1} and~\eqref{condMN2}, such that
\[ \Re\,\iip{h}{Lh} \geq K (\|Ah\|^2 + \|M^2Ch\|^2). \]

If in addition
\begeq\label{spectralAMC} 
A^*\! A + C^* M^4 C \quad \text{admits a spectral gap $\kappa>0$},
\endeq
then $L$ is hypocoercive on $\H^1/\K$: there exists constants $C\geq 0$ and
$\lambda>0$, explicitly computable, such that
\[ \|e^{-tL}\|_{\H^1/\K\to \H^1/\K} \leq C e^{-\lambda t}. \]
\end{Thm}

As usual, it might be better in practice to guess the right multipliers and 
re-do the proof, than to apply Theorem~\ref{thmmultsimple} directly. It is also clear
that many generalizations can be obtained by combining the method of
multipliers with the methods used in Theorem~\ref{hypocomult}.
Rather than going into such developments, I~shall just
show how to apply Theorem~\ref{thmmultsimple} 
on the Fokker--Planck equation with a potential $V \in C^2(\R^n)$. 
In that case, $[A^*,A]=I$ and $[B,C]=(\nabla^2 V)A$. When $\nabla^2V$ is bounded and
$A^*\! A+C^*C$ is coercive, there is no need to introduce an auxiliary operator $M$:
the choice $M=I$ is sufficient to provide exponential convergence to equilibrium.
But a multiplier might be useful when $\nabla^2V$ is unbounded.
Assume, to fix ideas, that $V$ behaves at infinity like $O(|x|^{2+\alpha})$ for
some $\alpha>0$, and $|\nabla^2V|$ like $O(|x|^{\alpha})$; 
then it is natural to use an operator $M$ which behaves polynomially,
in such a way as to compensate the divergence of $V$. In the rest of the section,
I shall use this strategy to recover the exponential convergence for the kinetic
Fokker--Planck equation under assumptions~\eqref{newcondFP} and~\eqref{D2Vgeq} below.

Let $M$ be the operator of multiplication by $m(x,v)$, where
\[ m(x,v) := \frac{1}{\left(V_0+V(x)+ \frac{|v|^2}{2} \right)^{\frac{\alpha}{4(2+\alpha)}}}, \]
and $V_0$ is a constant, large enough that $V_0+V$ is bounded below by~1.
Since $Bm=0$ and $B$ is a derivation, it is true that $B$ commutes with $M$.
Assume that $V_0+V$ is bounded below by a multiple of $1+|x|^{\alpha+2}$; then
\[ m^4 \leq \frac{1}{(V_0+V(x))^{\frac{\alpha}{(2+\alpha)}}} \leq 
\frac{K}{(1+|x|)^\alpha} \]
for some constant $K>0$, and then $m^4(\nabla^2 V)$ is bounded, so that
$M^4[B,C]$ is relatively bounded by $A$. Finally, condition~\eqref{condMN2} reduces
to 
\[ |\nabla_v m| \leq K m, \]
which is easy to check. To summarize, conditions~\eqref{condMN1} and~\eqref{condMN2}
are fulfilled as soon as there exist constants $C\geq 0$ and $K>0$ such that
\begeq\label{newcondFP} \forall x\in\R^n,\qquad
V(x) \geq K |x|^{2+\alpha}- C, \qquad |\nabla^2 V(x)| \leq C (1+|x|^\alpha)\qquad
(\alpha>0). 
\endeq

To recover exponential convergence under these assumptions, it remains to check
the spectral gap assumption~\eqref{spectralAMC}! This will be achieved under the
following assumption: there exists a potential $W$, and constants $C\geq 0$,
$K>0$ such that
\begeq\label{D2Vgeq} \forall x\in\R^n, \qquad
|V(x)-W(x)|\leq C, \qquad \nabla^2 W(x) \geq K(1+|x|)^{-\alpha}.
\endeq

From~\eqref{newcondFP} there exists $K>0$ such that
\[ m^4(x,v) \geq \left( \frac{K}{(1+|v|)^{\frac{1}{(2+\alpha)}}}\right )
\frac{1}{(1+|x|)^\alpha}
= : \Phi(v)\Psi(x). \]
Then
\[ \nabla_x^* m^4 \nabla_x \geq \Phi (\nabla_x^* \Psi \nabla_x). \]
Now I~claim that $\nabla_x^*\Psi\nabla_x$ (where $\Psi$ is a shorthand
for the multiplication by $\Psi$) is coercive in $L^2(e^{-V}\,dx)$, or in other
words that there exists $K>0$ such that
\[ \int f e^{-V} =0 \Longrightarrow \qquad
\int \Psi |\nabla_xf|^2 e^{-V} \geq K \int f^2 e^{-V};\]
or equivalently, that there is a constant $C$ such that 
for all $f\in L^2(e^{-V})$,
\[ \int [f(x)-f(y)]^2\, e^{-V(x)} e^{-V(y)}\,dx\,dy \leq
C\int \Psi(x) |\nabla_xf(x)|^2 e^{-V(x)}\,dx. \]
Indeed, with $C$ standing for various positive constants, one can write
\begin{align*}
\int [f(x)-f(y)]^2\, e^{-V(x)} e^{-V(y)}\,dx\,dy 
& \leq C \int [f(x)-f(y)]^2\, e^{-W(x)} e^{-W(y)}\,dx\,dy \\
& \leq C \int \bigl\langle (\nabla^2W(x))^{-1} \nabla_xf(x), 
\nabla_xf(x) \bigr\rangle\,e^{-W(x)}\,dx \\
& \leq C \int (1+|x|)^\alpha |\nabla_x f(x)|^2\,e^{-V(x)}\,dx,
\end{align*}
where the passage from the first to the second line is justified by the
{\bf Brascamp-Lieb inequality}~\cite[Theorem~4.1]{braslieb:BMPL:76}.

Now it is possible to conclude: the operator $A^*\! A=\nabla_v^*\nabla_v$ is coercive on 
$L^2(\gamma)$, $\gamma$ standing for the Gaussian distribution in the $v$ variable,
and the operator $\nabla_x^*\Psi\nabla_x$ is coercive on $L^2(e^{-V})$.
Theorem~\ref{thmcoertens} in Appendix~\ref{appproduct}
shows that $A^*\! A+\Phi\nabla_x^*\Psi\nabla_x$
is coercive on $L^2(\mu)$, where $\mu$ is the equilibrium distribution for the
Fokker--Planck equation. By monotonicity, 
$A^*\! A+C^*M^4C$ also admits a spectral gap;
this was the last ingredient needed for Theorem~\ref{thmmultsimple}
to apply.

\section{Further applications and open problems}

A very nice application of Theorem~\ref{hypocomult} was recently done by
Capella, Loeschcke and Wachsmuth on the so-called
Landau--Lifschitz--Gilbert--Maxwell model arising in micromagnetism.
Under certain simplifying assumptions, the linearized version of this
model can be written 
\begeq\label{LLGM} \begin{cases}
\pa_t m = J(h-m)\\ \\
\pa_t h = - \nabla\wedge\nabla\wedge h - J(h-m)\\ \\
\nabla\cdot h = - \nabla\cdot m,
\end{cases}
\endeq
where $m:\R^3\to\R^2$ stands for the (perturbation of the) magnetization,
and $h:\R^3\to\R^3$ for the (perturbation of the) magnetic field;
moreover, $J$ is the usual symplectic operator $J[x_1,x_2,x_3]=
[-x_2,x_1]$. Obviously, the system~\eqref{LLGM} is dissipative
but strongly degenerate, since the dissipation term
$-\nabla\wedge\nabla\wedge h$ only acts on $h$,
and not even on all components of $h$. This case turns out to be particularly
degenerate since one needs {\em three} commutators to apply
Theorem~\ref{hypocomult}. For further details I~refer the reader
to the preprint by Capella, Loeschcke and Wachsmuth~\cite{capella}.
\med

Still, many issues remain open in relation to the
hypocoercivity of operators of the form $A^*A+B$. I~shall
describe four of these open problems below.

\subsection{Convergence in entropy sense for rapidly increasing potentials}

In the present paper I~have derived some first results of exponential convergence
to equilibrium for the kinetic Fokker--Planck equation based on an entropy method
(Theorem~\ref{thmmesdata}). While these results seem to be the first of their kind, 
they suffer from the restriction of boundedness imposed on the Hessian of the
potential. It is not clear how to relax this assumption in order to treat,
say, potentials that behave at infinity like a power of $|x|$ that is higher
than~2. A first possibility would be to try to adapt the method of multipliers,
but then we run into two difficulties: (a) Entropic variants of the
Brascamp-Lieb inequality do not seem to be true in general,
and are known only under certain particular restrictions on the reference
measure (see the discussion by Bobkov and 
Ledoux~\cite[Proposition~3.4]{bobkled:BM:00}); 
(b) It is not clear that there is an entropic analogue 
of Theorem~\ref{thmSGproduct}. Both problems (a) and (b) have 
their own interest.
\sm

Another option would be to try to relax the {\em local} conditions
(i)--(iv) into {\em global} (integrated) boundedness conditions,
so as to have an analogue of Lemma~\ref{lemD2Vbdd} where the reference
measure would be the solution $f$ of the Fokker--Planck equation.
This is conceivable only if $f$ satisfies some good a priori estimates
for positive times.
\med

\subsection{Application to oscillator chains} \label{applosc}

One of the motivations for the present study was the hope to revisit
the works by Eckmann, Hairer, Rey-Bellet and others on hypoelliptic
equations for oscillator chains, modelling heat 
diffusion~\cite{EPRB:anharmonic:99,reybelletthomas:anharmonic:00,reybelletthomas:exp:02,eckmannhairer:uniqueness:01,eckmannhairer:hypo:03}.
So far I~have obtained only very partial success in that direction.
If we try to apply Theorem~\ref{hypocomult} to the model, as it is described
e.g. in the last section of~\cite{eckmannhairer:hypo:03}, we find that
the assumptions of Theorem~\ref{hypocomult} apply as soon as

(a) the ``pinning potential'' $V_1$ and the ``interaction potential'' $V_2$ have
bounded Hessians;

(b) the Hessian of the interaction potential is bounded below by a positive constant;

(c) the second derivatives of the logarithm of the stationary density are bounded;

(d) the stationary measure satisfies a Poincar\'e inequality.
\med

Let us discuss these assumptions.
Assumption (a) is a bit too restrictive, since it excludes for instance the
quartic double-well potentials which are classically used in that field;
but it would still be admissible for a start; and hopefully this restriction
can be relaxed later by a clever use of the method of multipliers.
By the way, it is interesting to note that such assumptions are not covered 
by the results in~\cite{eckmannhairer:hypo:03} which need a superquadratic growth
at infinity. Next, Assumption (b) is not so surprising since (as far as I~know) it 
has been imposed by all authors who worked previously on the 
subject.\footnote{More generally, as pointed out to me by Hairer,
all existing results seem to require that 
the interaction potential does dominate the pinning potential.}
But it is a completely open problem to derive
sufficient conditions for Assumptions (c) and (d), except in the simple case
where the two temperatures of the model are equal. This example illustrates an
important remark: The range of application of Theorem~\ref{hypocomult} 
(and other theorems of the same kind) will be considerably augmented when one
has qualitative theorems about the stationary measure for nonsymmetric
diffusion processes. For instance,

- When does the stationary measure satisfy a Poincar\'e inequality?

- Can one derive bounds about the Hessian of the logarithm of its density?
\med

The first question was addressed recently in papers by R\"ockner and
Wang (see for instance~\cite{rocknerwang:poincare:04}) in the context 
of {\em elliptic} equations, and it looks like a challenging open problem to 
extend their results to hypoelliptic equations. The second question seems to 
be completely open; of course it has its intrinsic interest, since very little has been
known so far about the stationary measures constructed e.g. in~\cite{EPRB:anharmonic:99}.

\subsection{The linearized compressible Navier--Stokes system}

An extremely interesting instance of hypocoercive linear {\em system}
is provided by the linearized compressible Navier--Stokes equations
for perfect gases. In this example, the noncommutativity does not
arise because of derivation along noncommuting vector fields,
but because of the noncommutativity of the space where the unknown
takes its values.

Obtained by linearizing the nonlinear system of Section~\ref{secCNS}
around the equilibrium state $(1,0,1)$, the linearized compressible
Navier--Stokes system reads as follows:
\begeq\begin{cases} \label{lCNS}
\pa_t\rho + \nabla\cdot u = 0; \\ \\
\pa_t u + \nabla(\rho+\theta) = \mu\, \Delta u + \mu \left(1-\frac2N\right)
\nabla(\nabla\cdot u);\\ \\
\pa_t\theta + \frac2N \nabla\cdot u = \kappa\, \Delta\theta.
\end{cases}\endeq
Here $N$ is the dimension,
$(\rho,u,\theta)$ are {\em fluctuations} of the density,
velocity and temperature respectively, $\mu>0$ is the viscosity of the fluid
and $\kappa>0$ the heat conductivity.
So it is natural to define $\H = L^2(\Om; \R\times\R^N\times\R)$,
where $\Om\subset\R^N$ is the position domain, and
the target space $\R\times\R^N\times\R$ is equipped with the
Euclidean norm
\[ \Bigl\| (\rho,u,\theta)\Bigr\|^2 = \rho^2 + |u|^2 + \frac{N}2\, \theta^2,\]
which is (up to a factor $-1/2$) the quadratic approximation of the
usual entropy of compressible fluids.

Let $h=(\rho,u,\theta)$; it turns out that~\eqref{lCNS} can 
be written in the form $\pa_t h + L h=0$, where
$L=A^*A+B$, $B^*=-B$, and $A$, $B$ are quite simple: 
\begeq
\begin{cases}
Ah = \Bigl( 0, \ \sqrt{2\mu}\, \{\nabla u\}, \ \nabla\theta \Bigr)\\ \\
Bh = \Bigl( \nabla\cdot u, \ \nabla\rho+\nabla\theta, \
\frac{2}N\, \nabla\cdot u\Bigr).
\end{cases} \endeq
Here I have used the notation
\[ \{\nabla u\}_{ij} = \left( \frac12\Bigl(\derpar{u_i}{x_j} + 
\derpar{u_j}{x_i} \Bigr) - \Bigl(\frac{\nabla\cdot u}{N}\Bigr)\,
\delta_{ij}\right)\]
for the traceless symmetrized (infinitesimal) strain tensor of the fluid.

The system~\eqref{lCNS} is degenerate in two ways.
First, the diffusion on the velocity variable $u$ does not control
all directions: In general it is false that $\int |\{\nabla u\}|^2$
controls the whole of $\int |\nabla u|^2$ (see the 
discussion in~\cite{DV:boltz:05} for instance: one needs at least an 
additional control on the divergence). Secondly, there is no
diffusion on the density variable $\rho$. This suggests to consider
commutators between $\tilde{A}:h \to (0,0,\nabla\theta)$ and
$B$. After some computations one gets (in slightly sketchy notation)
\[ [\tilde{A},B] = C_1+R_1;\qquad
C_1h = \frac2N \Bigl(0, 0, \nabla\nabla\cdot u\Bigr);\qquad
R_1 = - (0, \nabla^2\theta, 0);\]
\[ [C,B] = C_2 + R_2; \qquad
C_2h=\frac2N ( 0,0, \nabla\Delta\rho);\qquad
R_2h = \frac2N (0,-\nabla^2\nabla\cdot u, \nabla\Delta\theta).\]
So the commutator $C_1$ controls the variations of the divergence of $u$,
while the iterated commutator $C_2$ controls the variations of
the density $\rho$. However, if we try to apply Theorem~\ref{hypocomult}
in this situation, we immediately run into problems to control
the remainder $R_2$, and need to modify the strategy.
This problem is tricky enough to deserve a
separate treatment, so I~shall not consider it in this memoir.

\subsection{A model problem arising in the study of Oseen vortices}

All the material in this subsection was taught to me by Gallay.
Oseen vortices are certain self-similar solutions to the
two-dimensional incompressible Navier--Stokes equation, in vorticity
formulation~\cite{gallaywayne:invariant:02,gallaywayne:stability:05}.
The linear stability analysis of these vortices reduces to the spectral
analysis of the operator $S+\alpha B$ in $L^2(\R^2)$, where
\begeq\begin{cases}
\dps S\omega = -\Delta\omega + \frac{|x|^2}{16}\, \omega - \frac{\omega}2,
\\ \\
\dps B\omega = \BS[G]\cdot\nabla\omega + 2 \BS [G^{1/2}\omega]\cdot
\nabla G^{1/2};
\end{cases} \endeq
here $\BS[\omega]$ is the velocity field reconstructed from the
vorticity $\omega$:
\[ \BS[\omega](x) = \frac1{2\pi} \int_{\R^2}
\frac{(x-y)^\bot}{|x-y|^2}\:\omega(y)\,dy,\]
and $v^\bot$ is obtained from $v$ by rotation of angle $\pi/2$;
moreover $G$ is a Gaussian distribution: $G(x) = e^{-|x|^2/4}/(4\pi)$;
and $\alpha$ is a real parameter.

The spectral study of $S+\alpha B$ turns out to be quite tricky.
In the hope of getting a better understanding, one can decompose
$\omega$ in Fourier series: $\omega = \sum_{n\in\Z}
\omega_n(r) e^{in\theta}$, where $(r,\theta)$ are standard polar
coordinates in $\R^2$. For each $n$, the operators $S$ and~$B$
can be restricted to the vector space generated by $e^{in\theta}$,
and can be seen as just operators on a function $\omega(r)$:
\[\begin{cases}\dps (S_n \omega)(r) = - \partial_r^2 \omega -
\left(\frac{r}2 + \frac1{r}  \right)\partial_r \omega -
\left( 1- \frac{n^2}{r^2}\right)\omega,\\
\dps (B_n\omega)(r) = i\,n\, (\varphi\omega - g \Omega_n);
\end{cases}\]
here $g(r) = e^{-r^2/4}/{4\pi}$,
$\varphi(r) = (1-e^{-r^2/4})/{2\pi r^2}$,
and $\Omega_n(r)$ solves the differential equation
\[ - (r\Omega')' + \frac{n^2}r \Omega = \frac{r}{2} \omega.\]

The regime $|\alpha|\to\infty$ is of physical interest and has
already been the object of numerical investigations by physicists.
There are two families of eigenvalues which are imposed by symmetry reasons;
but apart from that, it seems that all eigenvalues
converge to infinity as $|\alpha|\to\infty$, and for some
of them the precise asymptotic rate of divergence $O(|\alpha|^{1/2})$
has been established by numerical evidence.
If that is correct, this means that the ``perturbation''
of the symmetric part $S$ by the antisymmetric, lower-order
operator $\alpha B$ is strong enough to send most eigenvalues to
infinity as $|\alpha|\to\infty$. Obviously, this is again a manifestation 
of a hypocoercive phenomenon.

To better understand this stability issue,
Gallagher and Gallay suggested the following

\begin{Modprob}
Identify sufficient conditions on $f:\R\to\R$, so that
the real parts of the eigenvalues of
\[ L_\alpha : \omega\longmapsto (-\pa_x^2\omega + x^2\omega - \omega)
+ i \alpha f \omega\]
in $L^2(\R)$ go to infinity as $|\alpha|\to\infty$, and estimate this rate.
\end{Modprob}

Here is how Gallagher and Gallay partially solved this problem.
Set $\H=L^2(\R;\C)$,
$A=\pa_x\omega + x\omega$, $B\omega=(i\alpha f)\omega$.
Then $C\omega = i\alpha f'\omega$, so the operator
$A^*A+C^*C$ is of Schr\"odinger type:
\[ (A^*A + C^*C)\,\omega = (-\partial^2_x\omega + x^2\omega -\omega)
+ \alpha^2 f'^2\omega,\]
and the spectrum of $A^*A+C^*C$ can be studied via standard
semi-classical techniques. For instance, if
$f'(x)^2 = x^2/(1+x^2)^k$, $k\in\N$, then the real part of the spectrum
of $A^*A + C^*C$ is bounded below like $O(|\alpha|^{2\nu})$,
with $\nu=\min(1,2/k)$. Then a careful examination of the
proof of Theorem~\ref{thmsimple} yields a lower bound
like $O(|\alpha|^\nu)$ on the real part of the spectrum
of $A^*A+B$.

This estimate is optimal for $A^*A+C^*C$, but it is in general
{\em not so} for $A^*A+B$. For instance, if $f(x) = 1/(1+x^2)$,
then $\nu = 1/4$, but numerical simulations suggest that the growth
is like $|\alpha|^{1/2}$. This might indicate a fundamental
limitation of the techniques developed in this part, 
and motivate the development of a refined analysis.

\part{The auxiliary operator method}

In this part I~shall present an abstract hypocoercivity theorem
applying to a linear operator $L$ whose symmetric part is nonnegative,
but which does not necessarily take the form $A^*A+B$. 
Still it will be useful to decompose $L$ into its
symmetric part $S$ and its antisymmetric part $B$. Of course,
we could always define $A$ to be the square root of $S$, but this
might be an extremely complicated operator, and the assumptions
of the $A^*A+B$ Theorems might in practice be impossible to check.
Important applications arise when the operator $S$ contains
an integral part, as in the linearized Boltzmann equation.

A classical general trick in spectral theory, when one studies
the properties of a given linear operator $L$, consists in introducing
an auxiliary operator which has good commutation properties with $L$.
Here the idea will be similar, with just an important twist: 
We shall look for an auxiliary operator $A$
which ``almost commutes'' with $S$ and ``does not at all'' commute
with $B$, in the sense that the effect of the commutator $[A,B]$ 
will be strong enough to enforce the coercivity of $S + [A,B]^*[A,B]$.

With this idea in mind, I~had been looking for a hypocoercivity
theorem generalizing, say, Theorem~\ref{thmsimple}, but stumbled
on the problem of practical verification of my assumptions.
In the meantime, Cl\'ement Mouhot and Lukas Neumann found
a theorem which, while in the same spirit of Theorem~\ref{thmsimple},
has some important structural differences. The Mouhot--Neumann theorem
is quite simple and turns out to be applicable to many important cases, as
investigated in~\cite{mouhotneumann:06}; so in the sequel I~shall
only present their approach, with just slight variations and a
more abstract treatment. Then I~shall discuss the weak
points of this method, and explain why another theory still needs to be
developed, probably with slightly more sophisticated tools. 
At the time of writing, Fr\'ed\'eric H\'erau has made partial
progress in this direction, but still did not manage to get things
to work properly.
\vfill\pagebreak

\section{Assumptions} \label{secassL}

In the sequel, $\H$ is a separable Hilbert space on $\R$ or $\Cplx$,
$S$ is a nonnegative symmetric, possibly unbounded operator $\H\to\H$
and $B$ is an antisymmetric, possibly unbounded operator $\H\to\H$.
Then $A=(A_1,\ldots,A_m)$ is an array of unbounded operators $\H\to\H$.
All of these operators are defined on a common dense domain.
I~shall actually ignore all regularity issues and be content with 
formal calculations, to be considered as a priori estimates.

The same conventions as in Section~\ref{secnot} will apply.
Some of the assumptions below will involve $\sqrt{S}^{-1}U$ for various
operators $U$; of course, this is not rigorous since $\sqrt{S}$ is in
general not invertible. To make sense of these assumptions, one can
either consider them as
a priori estimates for a regularized problem in which $S$ is
replaced by an invertible approximation (something like $S+\var I$,
and one tries to get estimates which are independent of $\var$);
or supply them with the condition that $\sqrt{S}$ is invertible on
the range of $U$ (a trivial case of application is when $U=0$).

The object of interest is the semigroup generated by the operator
\[ L = S + B.\]
The next hypocoercivity theorem for $L$ will make crucial use
of the commutator of $A$ and $B$. I~shall write
\[ [A,B] = Z\,C + R,
\]
where $Z$ is bounded from above and below,
and $R$ is some ``remainder''.

Now come a bunch of commutator conditions which will be used in
Section~\ref{secmainL}. Later in Section~\ref{secsimplifiedL}
I~shall make some simplifying assumptions which will drastically
reduce the number of these conditions; but for the moment I~shall
keep the discussion at a general level.
\med

{\bf (A1)}
$\dps
\begin{cases} \text{{\bf either}\quad} 
[C,S]\sle \sqrt{S} \\
\text{{\bf or}\quad\quad}
\sqrt{S}^{-1} [C,S] \sle \sqrt{S}
\end{cases}$
\med

%{\bf (A2)}
%$\dps R\sle \sqrt{S}, \ \sqrt{S} C$
%\med
%
{\bf (A2)}
$\dps
\begin{cases} \text{{\bf either}\quad} 
(A\sle \sqrt{S} A, \ C, \ \sqrt{S}C, \ \sqrt{S}) \quad\text{and}\quad
([C,L]\sle \sqrt{S}, \ \sqrt{S}C)\\
\text{{\bf or}\quad\quad}
\sqrt{S}^{-1} [C,L] \sle \sqrt{S}, \ \sqrt{S}C
\end{cases}$
\med

{\bf (A3)}
$\dps
\sqrt{S}[A^*,C]\sle \sqrt{S},\ \sqrt{S}C,\ C, \ \sqrt{S}A$
\med

{\bf (A4)}
$\dps (\sqrt{S}A^*\sle \sqrt{S}A, \ C, \ \sqrt{S}C,\ \sqrt{S})\quad
\text{and}\quad 
(\sqrt{S}C^*\sle \sqrt{S}C,\ \sqrt{S})$
\med

{\bf (A5)}
$\dps
\begin{cases} \text{{\bf either}\quad} 
(A^*\sle \sqrt{S}A, \ C, \ \sqrt{S} C, \ \sqrt{S}) \quad\text{and}\quad
([C^*,S]\sle \sqrt{S}, \ \sqrt{S}C)\\
\text{{\bf or}\quad\quad}
(\sqrt{S}A^* \sle \sqrt{S}A, \ C, \ \sqrt{S}C, \ S)\quad\text{and}\quad
(\sqrt{S}^{-1} [C^*,S]\sle \sqrt{S},\ \sqrt{S}C)\\
\text{{\bf or}\quad}
[C^*,S]=0
\end{cases}$
\med

{\bf (A6)} $\dps R\sle \sqrt{S}, \sqrt{S}C$.
\med

{\bf (A7)} There exist constants $\kappa, \ov{c}>0$ such that for
all $h\in\H$,
\[ \<Ah,ASh\> \geq \kappa \<SAh, Ah\> - \ov{c}
\Bigl( \<Sh,h\> + \<SCh,Ch\> + \|Ch\|^2\Bigr).\]
\med

Here is a simple, but sometimes too restrictive,
sufficient condition for {\bf (A7)} to hold (the proof is
left to the reader):
\med

{\bf (A7')} 
$\dps
\begin{cases} \text{{\bf either}\quad} 
(A\sle \sqrt{S}A, \ C, \ \sqrt{S} C, \ \sqrt{S})
\quad\text{and}\quad
([S,A]\sle C, \sqrt{S}C, \sqrt{S})\\
\text{{\bf or}}\quad\quad
\sqrt{S}^{-1}[S,A] \sle \sqrt{S}, \ \sqrt{S}C, \ C
\end{cases}$

\begin{Rk} Some of the assumptions {\bf (A1)}--{\bf (A7)} can be
replaced by other assumptions involving the commutator $[A,S]$.
I~did not mention these alternative assumptions since they are
in general more tricky to check that the ones which I~chose.
In case of need, the reader can easily find them by adapting
the proof of the main theorem below.
\end{Rk}

\section{Main theorem} \label{secmainL}

\begin{Thm}[hypocoercivity for $L=S+B$] \label{thmSB}
With the same notation as in Section~\ref{secmainL},
assume that {\bf (A1)}--{\bf (A7)} are satisfied.
Further assume that
\begeq\label{HL}
\begin{cases}
(i)\qquad\exists\kappa,c>0;\quad \forall h\in (\ker L)^\bot,\\
\qquad\qquad\qquad\qquad 
\<ASh, Ah\> \geq \kappa \|Ah\|^2\ - \ c\bigl(\<Sh,h\> + \<SCh,Ch\> +
\|Ch\|^2\bigr);\\ \\
(ii)\qquad 
S + A^*SA + C^*SC + A^*A + C^*C \quad \text{is coercive on $(\ker L)^\bot$.}
\end{cases}
\endeq
Then there are constants $c,\lambda>0$, depending only on the 
constants appearing implicitly in {\bf (A1)}--{\bf (A7)} 
and~\eqref{HL}, such that
\[ \Bigl\| e^{-tL}\Bigr\|_{\tilde{\H}\to\tilde{\H}} \leq c\, e^{-\lambda t},\]
where $\tilde{\H}\subset(\ker L)^\bot$ is defined by the
Hilbert norm
\[ \|h\|_{\tilde{\H}}^2 = \|h\|^2 + \|Ah\|^2 + \|Ch\|^2.\]
\end{Thm}

\begin{Rk} Although Condition~\eqref{HL}(i) formally resembles
Assumption {\bf (A7)}, I~have preferred to state it together with
\eqref{HL}(ii) because its practical verification
often depends on a control of $\|h\|^2$ by
$\<Sh,h\> + \|Ch\|^2$.
\end{Rk}

\begin{proof}[Proof of Theorem~\ref{thmSB}]
The proof is quite similar in spirit to the proof of 
Theorem~\ref{thmsimple}, so I~shall be sketchy and
only point out the main steps in the calculations.

First note that $(\ker L)^\bot$ is stable under the evolution
by $e^{-tL}$. Indeed, if $k\in\ker L$, then
$(d/dt) \< e^{-tL}h,k\> = \<L e^{-tL}h,k\>
= \<e^{-tL}h,L^*k\>$, so it is sufficient to show that $L^*k=0$.
But $Lk=0$ implies $\<Sk,k\>=\<Lk,k\>=0$, so $k\in\ker S$
(here the nonnegativity of $S$ is essential), so $k\in\ker B$ also,
and $L^*k = (S-B)k=0$.

Next let
\[ {\cal F}(h) = \|h\|^2 + a \|Ch\|^2 + 2 b\, \Re\,
\<Ch, Ah\> + c \|Ah\|^2,\]
where $\Re$ stands for real part, and
$a, b, c$ will be chosen later in such a way that
$1\gg a\gg b\gg c>0$, $a\ll \sqrt{b}$, $b\ll \sqrt{ac}$. 
In particular, ${\cal F}(h)$ will be bounded from above and below by
constant multiples of $\|h\|_{\tilde{\H}}^2$; so to prove the theorem
it is sufficient to establish the estimate
$(-d/dt){\cal F}(e^{-tL}h) \geq {\rm const.} {\cal F}(e^{-tL}h)$.
Without loss of generality, we can do it for $t=0$ only.
In the sequel, I~shall also pretend that $\H$ is a real Hilbert space,
so I~shall not write real parts.

By direct computation,
\begin{align}\label{dFdtL}
\left.-\frac{d}{2\,dt}\right|_{t=0} {\cal F}(e^{-tL}h)\ = 
& \ \<Sh, h\> \\
& + a \< CSh,Ch\> + a \< CBh, Ch\> \nonumber\\
& + b \<CLh, Ah\> + b \<Ch, ALh\> \nonumber\\
& + c \<ASh, Ah\> + c \<ABh, Ah\>\nonumber.
\end{align}
\noindent Now we shall estimate~\eqref{dFdtL} line after line.
\sm

(1) The first line of~\eqref{dFdtL} is kept unchanged.
\med

(2) The second line of~\eqref{dFdtL} is rewritten as follows:
\begeq\label{CShCh} 
a \<CSh, Ch\> = a \<SCh, Ch\> + a \<[C,S]h, Ch\>.
\endeq
Then the second term in the right-hand side of~\eqref{CShCh}
is estimated from below, either by
$-a \|[C,S]h\|\,\|Ch\|$, or by $-a \|\sqrt{S}^{-1}[C,S]h\|\,\|\sqrt{S}Ch\|$;
By Assumption ({\bf A1}), these expressions can in turn be estimated 
from below by a constant multiple of
\[ -a \Bigl( \|\sqrt{S}h\| \,\|Ch\| + \|\sqrt{S}h\|\, \|\sqrt{S}Ch\|\Bigr).\]
(Here I~used the identity $\<Su,u\> = \|\sqrt{S}u\|^2$.)
\med

(3) The treatment of the third line of~\eqref{dFdtL} is crucial; this is where 
the added coercivity from the commutator $[A,B]$ will show up.
To handle the first term in this line, we write
\begin{align*}
\<CLh, Ah\> & =
\<LCh, Ah\> + \<[C,L]h, Ah\> \\
& = \<SCh, Ah\> + \<BCh, Ah\> + \<[C,L]h, Ah\> \\
& = \<SCh, Ah\> - \<Ch, BAh\> + \<[C,L]h, Ah\>.
\end{align*}
When we add this to the second term of the third line,
$\<Ch, ALh\> = \<Ch, ABh\> + \<Ch, ASh\>$, we obtain
\begin{align*} 
& \<Ch, (AB-BA)h\> + \<SCh, Ah\> + \<[C,L]h, Ah\> + \<Ch, ASh\>\\
= \ & \<Ch, \,(ZC + R)
h\> + \<SCh, Ah\> + \<[C,L]h, Ah\> + \<Ch, ASh\>\\
\geq \ & \kappa \|Ch\|^2 + \<Ch,Rh\> + \<SCh, Ah\> + \<[C,L]h, Ah\>
+ \<Ch, ASh\>.
\end{align*}
So there are four ``error'' terms to estimate from below:
\begeq\label{3terms}
\<Ch, Rh\>, \qquad 
\<SCh, Ah\>, \qquad \<[C,L]h, Ah\>, \qquad \<Ch, ASh\>.
\endeq
\sm

- To estimate the first term in~\eqref{3terms}, just write
\[ \<Ch, Rh\> \geq -\|Ch\|\, \|Rh\|\]
and apply Assumption {\bf (A6)}; it follows that there is a lower bound
by a constant multiple of
\[ - b \|Ch\| \, \bigl( \|\sqrt{S}h\| + \|\sqrt{S}Ch\|\bigr).\]
\sm

- To estimate the second term in~\eqref{3terms}, use the Cauchy--Schwarz
inequality:
\[ \< SCh, Ah\> \geq -\,\|\sqrt{S}Ch\|\, \|\sqrt{S}Ah\|.\]
\sm

- To estimate the third term in~\eqref{3terms}, write either
\[ \<[C,L]h, Ah\> \geq -\, \|[C,L]h\|\,\|Ah\|\]
or
\[ \<[C,L]h, Ah\> \geq -\, \|\sqrt{S}^{-1}[C,L]h\|\,\|\sqrt{S}Ah\|\]
and apply Assumption {\bf (A2)}. It results a lower bound by a constant
multiple of
\begin{multline*} 
-b \bigl( \|\sqrt{S}Ah\| + \|Ch\| + \|\sqrt{S}Ch\| + \sqrt{S}h\|\bigr)
\bigl( \|\sqrt{S}h\| + \|\sqrt{S}Ch\|\bigr) \\
-b\ \bigl( \|\sqrt{S}h\|+\|\sqrt{S}Ch\|\bigr) \|\sqrt{S}Ah\|\bigr).
\end{multline*}

\sm

- The fourth term in~\eqref{3terms} is a bit more tricky:
\begin{align*}
\<Ch, ASh\> & = \<A^*Ch, Sh\> \\
& = \<[A^*,C]h, Sh\> + \<CA^*h, Sh\> \\
& = \<[A^*,C]h, Sh\> + \<A^*h, C^*Sh\> \\
& = \<[A^*,C]h, Sh\> + \<A^*h, SC^*h\> + \<A^*h, [C^*,S]h\>.
\end{align*}
This gives rise to three more terms to estimate:
\begeq\label{threeterms}
\< [A^*,C]h, Sh\>, \quad \<A^*h, SC^*h\>, \quad \<A^*h, [C^*,S]h\>.
\endeq

- To handle the first term in~\eqref{threeterms}, write
\[ \< [A^*,C]h, Sh\> =
\< \sqrt{S} [A^*,C]h, \sqrt{S}h\> \geq
- \|\sqrt{S} [A^*,C]h\|\, \|sqrt{S}h\|;\]
then apply Assumption {\bf (A3)} to bound $\|\sqrt{S}[A^*,C]h\|$.
The result is a lower bound by a constant multiple of
\[ -b\|\sqrt{S}h\| \bigl( \|\sqrt{S}h\| + \|\sqrt{S}Ch\|
+ \|Ch\| + \|\sqrt{S}Ah\|\bigr).\]
\sm

- To bound the second term in~\eqref{threeterms}, write
\[ \<A^*h, SC^*h\> = 
\<\sqrt{S} A^*h, \sqrt{S}C^*h\> \geq
- \|\sqrt{S} A^*h\|\,\|\sqrt{S}C^*h\|;\]
then apply Assumption {\bf (A4)} to bound these two norms separately.
The result is a lower bound by a constant multiple of
\[ -b \bigl( \|\sqrt{S}Ah\| + \|Ch\| + \|\sqrt{S}Ch\| + \|\sqrt{S}h\|\bigr)
\bigl(\|\sqrt{S}Ch\| + \|\sqrt{S}h\|\bigr).\]
\sm

- To bound the last term in~\eqref{threeterms}, one possibility is to write
\[ \<A^*h, [C^*,S]h\> \geq - \|A^*h\|\,\|[C^*,S]h\|;\]
another possibility is
\[ \<A^*h, [C^*,S]h\> = \<\sqrt{S}A^*h,\, \sqrt{S}^{-1} [C^*,S]h\>
\geq - \|\sqrt{S}A^*h\|\,\|\sqrt{S}^{-1} [C^*,S]h\|.\]
Then one can apply Assumption {\bf (A5)} to control these terms.
In the end, this gives a lower bound by a constant multiple of
\[ -b \bigl( \|\sqrt{S}Ah\| + \|Ch\| + \|\sqrt{S}Ch\| + \|\sqrt{S}h\|\bigr)
(\|\sqrt{S}h\| + \|\sqrt{S}Ch\|).\]
\sm

(4) Finally, the fourth line of~\eqref{dFdtL} is handled as follows:
\begin{multline}\label{4thline} \<Ah, ASh\> + \<Ah, ABh\> = 
\alpha \<Ah, ASh\> + \beta\<Ah, ASh\> \\ + \<Ah, BAh\> + \<Ah, [A,B]h\>,
\end{multline}
where $\alpha,\beta\geq 0$ and $\alpha+\beta=1$. The first term
$\alpha\<Ah, ASh\>$ is estimated by means of Assumption {\bf (A7)};
the second term $\beta\<Ah, ASh\>$ by means of Assumption~\eqref{HL}(i);
altogether, these first two terms can be bounded below by a constant
multiple of
\[ c (\|\sqrt{S}Ah\| + \|Ah\|^2) \ - c
\bigl( \|\sqrt{S}h\| + \|\sqrt{S}Ch\| + \|Ch\|^2\bigr).\]
Then the third term $\<Ah, BAh\>$ in~\eqref{4thline} vanishes; 
and the last term $\<Ah, [A,B]h\>$ is bounded below by
$-\|Ah\|\,\|ZCh\|-\|Ah\|\,\|Rh\|$, which in view of Assumption
{\bf (A6)} can be bounded below by a constant multiple of
\[ -c\|Ah\|\,\|Ch\| - c \|Ah\|\, \bigl( \|\sqrt{S}h\| + \|\sqrt{S}Ch\|\bigr).\]
\med

Gathering up all these lower bounds, we see that
\[ -\frac{d}{2\,dt} {\cal F} \geq {\rm const.}\, \<X, mX\>,\]
where
\[ X = \Bigl( \|\sqrt{S}h\|, \, \|\sqrt{S}Ch\|,\,
\|Ch\|,\, \|\sqrt{S}Ah\|,\, \|Ah\|\Bigr),\]
$m$ is the $5\times 5$ matrix
\[ m = \left [ \begin{matrix} 
1-Mb-Mc & -Ma-Mb & -Ma-Mb & -Mb & -Mc \\
0 & a-Mb-Mc & -Mb & -Mb & -Mc \\
0 & 0 & b-Mc & 0 & -Mc \\ 
0 & 0 & 0 & c & 0 \\
0 & 0 & 0 & 0 & c
\end{matrix} \right ],\]
and $M$ is a large number depending on the bounds appearing in
the assumptions of the theorem.

Then by reasoning as in Section~\ref{secbasic} and using
Lemma~\ref{ll}, we can find coefficients $a,b,c>0$ and a constant
$\kappa>0$ such that
\begeq\label{dFdtL2}
\left.-\frac{d}{2\,dt}\right|_{t=0} {\cal F}(e^{-tL}h)\ \geq
\kappa \Bigl( \<Sh, h\> + \<SCh, Ch\> + \|Ch\|^2 + \<SAh, Ah\> +
\|Ah\|^2\Bigr).
\endeq

By Assumption~\eqref{HL}(ii), this implies the existence of $\kappa',
\kappa''>0$ such that
\begin{align*}
\left.-\frac{d}{2\,dt}\right|_{t=0} {\cal F}(e^{-tL}h)\ & \geq
\kappa' \Bigl( \|h\|^2 + \<Sh, h\> + \<SCh, Ch\> + \|Ch\|^2 + \<SAh, Ah\> +
\|Ah\|^2\Bigr) \\
& \geq \kappa'' {\cal F}(h).
\end{align*}

This concludes the proof.
\end{proof}

\section{Simplified theorem and applications} \label{secsimplifiedL}

In this section I~shall consider a simplified version of
Theorem~\ref{thmSB}.

\begin{Cor} \label{thmSBsimpl}
Let $A=(A_1,\ldots,A_m),B,S$ be linear operators on a Hilbert space $\H$,
and let $C=[A,B]$. Assume that
\[ A^*=-A, \quad B^*=-B, \quad C^*=-C, \quad S^*=S\geq 0;\]
\[ [C,A]=0,\quad [C,B]=0,\quad [C,S]=0.\]
Further assume that there exists $\kappa,c>0$ such that for all
$h\in (\ker A\cap \ker B)^\bot$,
\begin{multline}\label{AhASh} 
\<Ah, ASh\> \geq \kappa \bigl( \<SAh, Ah\> + \|Ah\|^2\bigr) \\
- \ c\bigl( \<Sh, h\> + \<SCh, Ch\> + \|Ch\|^2\bigr);
\end{multline}
and that
\begeq\label{HC'}
S+C^*C \quad \text{is coercive on $(\ker A\cap\ker B)^\bot$}.
\endeq

Then there exists $\lambda>0$ such that
\[ \|e^{-t(S+B)}\|_{\tilde{H}\to\tilde{H}} = O(e^{-\lambda t}),\]
where $\tilde{\H}\subset(\ker L)^\bot$ is defined by the
Hilbert norm
\[ \|h\|_{\tilde{\H}}^2 = \|h\|^2 + \|Ah\|^2 + \|Ch\|^2.\]
\end{Cor}

\begin{proof}[Proof of Corollary~\ref{thmSBsimpl}]
The assumptions of the theorem trivially imply assumptions
{\bf (A1)}--{\bf (A6)} from Section~\ref{secassL}.
Assumption~\ref{AhASh} is equivalent to the conjunction of
{\bf (A7)} and~\eqref{HL}(i). Finally, \eqref{HC'} is obviously
stronger than~\eqref{HL}(ii).
\end{proof}

Now let us make the link with the Mouhot--Neumann hypocoercivity
theorem~\cite[Theorem~1.1]{mouhotneumann:06}. Although the set
of assumptions in that reference is not exactly the same as in
the current section, we shall see that under a small
additional hypothesis, the assumptions in~\cite{mouhotneumann:06}
imply the present ones.

In~\cite{mouhotneumann:06},
the Hilbert space $\H$ is $L^2(\T^n_x\times\R^n_v)$,
and $A=\nabla_v$, $B=v\cdot\nabla_x$, $C=\nabla_x$;
and the operator $S$ only acts on the velocity variable $v$,
so we have indeed $A^*=-A$, $B^*=-B$, $C^*=-C$, and $C$ commutes
with $A$, $B$ and $S$.
The kernel of $L$ is similar to the kernel of $S$ (up to identifying
$v\to h(v)$ with $(x,v)\to h(v)$), and contains constant functions.
Since $C^*C=-\Delta_x$ has a spectral gap, Condition \eqref{HC'}
is equivalent to the fact that $S$ has a spectral gap
in $L^2(\R^n_v)$, which is Assumption H.3 in~\cite{mouhotneumann:06}.
So it only remains to check~\eqref{AhASh}, which will be true as
soon as
\begeq\label{AhASh'} 
\<\nabla_v h, \nabla_v Sh\> \geq \kappa 
\bigl( \<S\nabla_vh, \nabla_vh\> + \|\nabla_vh\|^2\bigr)
\ - \ c\|h\|^2.
\endeq

It is assumed in~\cite{mouhotneumann:06} that $S$, viewed as an operator
on $L^2(\R^n_v)$, can be decomposed into the difference of two 
self-adjoint operators: $S=\Lambda-K$, where $\Lambda$ is positive
definite and 
\begeq\label{LK1}
\<\nabla_v h, \nabla_v \Lambda h\> \geq 
\kappa \<\nabla_vh, \Lambda \nabla_v h\> - c \|h\|^2;
\endeq
\begeq\label{LK2}
\forall\delta>0,\quad \exists c(\delta)>0;\qquad
\<\nabla_vh, \nabla_v K h\> \leq \delta \|\nabla_vh\|^2 + c(\delta)
\|h\|^2.
\endeq
Let us further assume that {\em $K$ is compact relatively to $\Lambda$},
in the sense that
\[
\forall\var>0, \quad \exists c(\var)>0;\qquad
K \leq \var \Lambda + c(\var) I,
\]
or equivalently (since $\Lambda = S+K$)
\begeq\label{LK+}
\forall\var>0, \quad \exists c(\var)>0;\qquad
K \leq \var S + c(\var) I.
\endeq
By using~\eqref{LK1}, \eqref{LK2} and~\eqref{LK+},
and denoting by $c$ and $\kappa$ various positive constants,
one easily obtains
\begin{align*}
\<\nabla_v h, \nabla_v Sh\> & = 
\<\nabla_v h , \nabla_v \Lambda h\> - \<\nabla_vh, \nabla_v Kh\> \\
& \geq \kappa \<\nabla_v h, \Lambda\nabla_v h\>  - c\|h\|^2
- \<\nabla_v h, \nabla_v Kh\> \\
& \geq \kappa \bigl(\<\nabla_v h ,\Lambda \nabla_v h\> + \|\nabla_vh\|^2\bigr)
- c \|h\|^2 - (\kappa/2) \|\nabla_vh\|^2 - c \|h\|^2 \\
& \geq \kappa \bigl(\<\nabla_v h ,\Lambda \nabla_v h\> + \|\nabla_vh\|^2\bigr)
- c \|h\|^2 \\
& \geq \kappa \bigl( \<\nabla_v h, S \nabla_v h\> + \|\nabla_vh\|^2\bigr)
- c \bigl( \|h\|^2 + \<\nabla_v h, K \nabla_v h\> \bigr)\\
& \geq \kappa \bigl(\<\nabla_v h, S \nabla_v h\> + \|\nabla_vh\|^2\bigr)
- c \|h\|^2.
\end{align*}
This establishes~\eqref{AhASh'}.

Assumption~\eqref{LK+} is not made in~\cite{mouhotneumann:06},
but it is satisfied in all the examples discussed therein:
linear relaxation, semi-classical relaxation, linear Fokker--Planck
equation, Boltzmann and Landau equations for hard potentials.
So all these examples can be treated by means of 
Theorem~\ref{thmSBsimpl}.
I~refer to~\cite{mouhotneumann:06} for more explanations and results
about all these models. Mouhot and Neumann also use these hypocoercivity 
results to construct smooth solutions for the corresponding 
{\em nonlinear} models close to equilibrium,
thereby simplifying parts of the theory developed by
Guo, see e.g.~\cite{guo:landau:02}.

\section{Discussion and open problems}

Although it already applies to a number of interesting models,
Theorem~\ref{thmSB} suffers from several shortcomings.
Consider for instance the case when $S$ is a bounded operator
(as in, say, the linearized Boltzmann equation for Maxwellian
cross-section), and there is a force term $-\nabla V(x)\cdot\nabla_v$
in the left-hand side of the equation. Then the higher derivative
term in $[C,L]$ is $-\nabla^2 V(x) \cdot\nabla_v$, which certainly
{\em cannot} be bounded in terms of $S$ and $C$; so Assumption
{\bf (A2)} does not hold. It is likely that Theorem~\ref{thmSB} rarely
applies in practice when $[C,B]=[[A,B],B]\neq 0$.

Other problems are due to Assumption {\bf (A4)}. 
This assumption will not hold for, say,
$C=\nabla_x$ in a bounded domain $\Omega\subset\R^n$;
indeed, in a slightly informal writing, $C^*=-C + \sigma\cdot dS$,
where $\sigma$ is the outer unit normal vector on $\pa\Omega$
and $dS$ is the surface measure on $\pa\Om$. So the computation
used in the proof of Theorem~\ref{thmSB} does not seem to give
any result in such a situation.\footnote{By the way, at present there 
seems to be no really satisfactory treatment of bounded domains in
linearized kinetic theory, apart of course from the case of a periodic box.}

A last indication that Theorem~\ref{thmSB} is not fully satisfactory
is that it does not seem to contain Theorem~\ref{thmsimple} as a particular
case, although we would like to have a unified treatment of the
general case $L=S+B$ and the particular case $L=A^*A+B$.
In fact, as the reader may have noticed, the choices of coefficients
in the auxiliary functionals appearing respectively in the proof
of Theorem~\ref{thmsimple} and in the proof of Theorem~\ref{thmSB}
go in the opposite way!! Indeed, in the first case it was
$\|h\|^2 + a\|Ah\|^2 + 2b \< Ah, Ch\> + c\|Ch\|^2$ with 
$a\gg b \gg c$, while in the second case it was
$\|h\|^2 + a\|Ch\|^2 + 2b \<Ch, Ah\> + c \|Ah\|^2$.

Some playing around with the functionals suggests that these problems
can be solved only if the auxiliary operator $A$ is ``{\em comparable}''
to $\sqrt{S}$, say in terms of order of differential operators.
So if $S$ is bounded, then also $A$ should be bounded.
This suggests to modify the Mouhot--Neumann strategy in the case
when $S$ is bounded, by choosing, instead of $A=\nabla_v$, something
like $A= (I-\Delta_v+ v\cdot\nabla_v)^{-1/2}\nabla_v$.
(I wrote $\Delta_v - v\cdot\nabla_v$ rather than $\Delta_v$,
because in many cases known to me, the natural reference measure
is the Gaussian measure in $\R^n_v$.)
Then computations involve nonlocal operators and become more intricate.
I~shall leave the problem open for future research.

\part{Fully nonlinear equations}

In this part I~shall consider possibly nonlinear equations, and 
I~shall not depend on ``exact'' commutator identities.
To get significant results under such weak structure assumptions, 
I~shall assume that I~deal with
solutions that are {\em very smooth}, uniformly in time.
Moreover, I~shall only prove results of convergence like
$O(t^{-\infty})$, that is, faster than any inverse power of $t$.

As in Remark~\ref{rknorm2}, the assumption of uniform smoothness
can be relaxed as long as one has good estimates of exponential
decay of singularities, together with a stability result
(solutions depart from each other no faster than exponentially fast).
However, I~shall not address this issue here.

At the level of generality considered here, the rate $O(t^{-\infty})$
cannot be so much improved, since some cases are included for which
exponential convergence simply does not hold, even for the linearized
equation. In many situations one can still hope for rates of convergence
like $O(e^{-\lambda t^\gamma})$, as in the close-to-equilibrium
theory of the Boltzmann equation with soft potentials~\cite{guostrain:exp}. 
If a linearized study
suggests convergence like $O(e^{-\lambda t})$ or $O(e^{-\lambda t^\gamma})$
for a particular nonlinear model, then
one can try to obtain this rate of convergence by putting together
the present nonlinear analysis (which applies far from equilibrium) 
with a linearization procedure (close to equilibrium) and a 
subsequent linear study.

This part is strongly influenced by my collaborations with Laurent 
Desvillettes on the convergence to 
equilibrium for the linear Fokker-Planck equation~\cite{DV:FP:01} 
and the nonlinear Boltzmann equation~\cite{DV:boltz:05}. 
The method introduced in these papers was based on 
the study of second-order time differentiation of certain functionals;
since then it has been successfully applied to other 
models~\cite{CCG:linear:03,FNS:linear:04}. Our scheme of proof
had several advantages: It was very general, physically meaningful,
and gave us the intuition
for the strong time-oscillations between hydrodynamic and homogeneous
behavior, that were later observed numerically with a high 
accuracy~\cite{FMP:NlogN}. On the other hand, our method had two major 
drawbacks: First, the heavy amount of calculations entailed by the
second-order differentiations (especially in the presence of several 
conservation laws); and secondly, the particularly tricky analysis
of the resulting coupled systems of second-order differential inequalities.

The approach will I~shall adopt in the sequel remedies these drawbacks:
First, it only uses first-order differentiation; secondly,
it confines many heavy computations into a black box that can be
used blindly. The price to pay will be the loss of intuition in the proof.

The main result is a rather abstract theorem
stated in Section~\ref{secmain} and proven in Section~\ref{secproof}.
Then I~shall show how to use
this abstract result on various examples: the compressible Navier--Stokes
system (Section~\ref{secCNS}); the 
Vlasov--Fokker--Planck equation with smooth and small coupling
(Section~\ref{secVFP}); and the Boltzmann equation
(Section~\ref{secboltz}). 

In the case of the
Vlasov--Fokker--Planck equation to be considered, the coupling
is simple enough that all the smoothness bounds appearing in the
assumptions of the main theorem can be proven in terms of just
assumptions on the initial data. In the other cases, the
results will be conditional (depend on the validity of uniform
regularity estimates).

The hard core of the proof of the main result
was conceived during the conference 
``Advances in Mathematical Physics'' in the honor of Carlo Cercignani 
(Montecatini, September 2004). It is a pleasure to thank the organizers 
of that meeting (Luigi Galgani, Maria Lampis, Rossanna Marra, Giuseppe Toscani)
for helping to create a fruitful and pleasant atmosphere of work.
The main results were first announced two weeks later, 
in an incomplete and preliminary form, at the 
Conference ``Mathematical Aspects of Fluid and Plasma Dynamics''
(Kyoto, September 2004), beautifully organized 
by Kazuo Aoki. During the Summer of 2006, for the purpose of various 
lectures in Porto Ercole, Trieste and Xining, I~rewrote and generalized
the main theorem, and added new applications.
Additional thanks are due to Kazuo for an important remark
about the treatment of the Boltzmann equation with Maxwellian diffusive 
boundary condition.

\vfill\pagebreak

\section{Main abstract theorem} \label{secmain}

The assumptions in this section are expressed in a rather abstract 
formalism. ``Concrete'' examples will be provided later in
Sections~\ref{secCNS} to~\ref{secboltz}.

\subsection{Assumptions and main result} \label{subass}

The theorem below involves five kinds of objects:
\sm

- a family of normed spaces $(X^s,\|\cdot\|_s)_{s\geq 0}$;
the index $s$ can be thought of as a way to quantify the regularity
(smoothness, decay, etc.);
\sm

- two ``differential'' operators $B$ and $\C$,
such that $B$ is ``conservative'' and $\C$ is ``dissipative'';
\sm

- a ``very smooth'' solution $t\to f(t)$ of the equation
\[ \pa_t f + B f = \C f,\]
with values in a subset $X$ of the intersection of all the spaces $X^s$;
\sm

- a Lyapunov functional $\E$, which is dissipated by the equation
above, and admits a unique absolute minimizer $f_\infty$;
\sm

- a finite sequence of ``nested nonlinear projections''
$(\Pi_j)_{1\leq j\leq J}$; one can think that $\Pi_j$ is the
projection onto the space of minimizers of $\E$ under $J-j$
constraints, and in particular $\Pi_J$ is the map which takes
everybody to $f_\infty$.
\sm

The goal is to prove the convergence of $f(t)$ to the stationary
state $f_\infty$, and to get estimates on the rate of convergence.
\sm

I~shall make several assumptions about these various objects.
Even though these assumptions may look a bit lengthy and complicated,
I~tend to believe that they are satisfied in many natural cases.
The following notation will be used: 

- If $A$ is an operator, then the image of a function $f$ by $A$ 
will be denoted either by $A(f)$ or simply by $Af$.

- The Fr\'echet derivative of $A$, evaluated at a function $f$,
will be denoted by $A'(f)$ or $A'_f$; so $A'_f\cdot g$ stands for
the Fr\'echet derivative of $A$ evaluated at $f$ and applied to the
``tangent vector'' $g$.

- The notation $\|A'(f)\|_{X\to Y}$ stands for the norm of the linear
operator $A'(f):X\to Y$, i.e. the smallest constant $C$ such that
$\|A''(f)\cdot g\|_Y \leq C \|g\|_X$ for all $g\in X$.

- Similarly, the second (functional) derivative of $A$, 
evaluated at a function $f$,
will be denoted by $A''(f)$ or $A''_f$; so $A''_f\cdot (g,h)$
stands for the Hessian of $A$ evaluated at $f$ and applied to
the two ``tangent vectors'' $g$ and $h$. The notation $\|A''(f)\|_{X\to Y}$
stands for the smallest constant $C$ such that
$\|A'(f)\cdot (g,h)\|_{Y} \leq C \|g\|_X \|h\|_X$ for all $g,h\in X$.

\sm

\begin{Ass}[scale of functional spaces] \label{assfunct}
$(X^s,\|\cdot\|_s)_{s\geq 0}$ is a nonincreasing family
of Banach spaces such that
\sm

(i) $X^0$ is Hilbert; its norm $\|\cdot\|_0$ will be denoted by
just $\|\cdot\|$;

(ii) The injection $X^{s'}\subset X^s$ is continuous for $s'\geq s$;
that is, there exists $C=C(s,s')$ such that
\begeq
\|f\|_{s} \leq C \|f\|_{s'}.
\endeq
\sm

(iii) The family $(X^s)_{s\geq 0}$ is an interpolation family:
For any $s_0,s_1\geq 0$ and $\theta\in [0,1]$ there is a constant
$C=C(s_0,s_1,\theta)$ such that
\begeq\label{interp}
s=(1-\theta)\,s_0 + \theta\,s_1\Longrightarrow\quad
\forall f \in X^{s_0}\cap X^{s_1},\quad
\|f\|_{s}\leq C\, \|f\|_{s_0}^{1-\theta}\,\|f\|_{s_1}^\theta;
\endeq
\end{Ass}

One may think of $s$ as an index quantifying the regularity of $f$,
say the number of derivatives which are bounded in a certain norm.
In the sequel, I~shall sometimes refer informally to $s$ as an
index standing for a number of derivatives, even if it is not
necessarily so in general.

\begin{Ass}[workspaces]\label{assworkspace}
$X$ and $Y$ are two sets such that
$X\subset Y\subset \cap_{s\geq 0} X^s$; moreover,
$Y$ is convex and bounded in all spaces $X^s$.
\end{Ass}

\begin{Ass}[solution] \label{asssol}
$f\in C(\R_+;X^s)\cap C^1((0,+\infty);X^s)$ for all $s$;
moreover $f(t)\in X$ for all $t$. (In particular
$f$ is bounded in all spaces $X^s$.)
\end{Ass}

In the sequel, the notation $f_0$ will be a shorthand for $f(0)$.

\begin{Ass}[equation] \label{asseq}
$f$ solves the equation
\begeq\label{eqf}
\derpar{f}{t} + Bf = \C f,
\endeq
where
\sm

(i) $B$, $\C$ are well-defined on $Y$ and valued in a bounded
subset of $X^s$ for all $s$;

(ii) For any $s$ there is $s'$ large enough such that
$B'$ is bounded $X^{s'}\to X^s$, uniformly on $Y$;

(iii) $\C$ is Lipschitz $X^s\to X^0$, uniformly on $Y$, for $s$ large enough.
\end{Ass}

In short, $B$ and $\C$ satisfy a ``Lipschitz condition with
possible loss of derivatives''. If $X^s$ is a Sobolev space of order $s$
on a bounded domain, then any reasonable differential operator of
finite order, with smooth coefficients, will satisfy these assumptions.

\begin{Ass}[stationary state] \label{assstat}
$f_\infty$ is an element of $X$, satisfying $Bf_\infty = \C f_\infty =0$.
\end{Ass}

\begin{Ass}[projections] \label{assproj}
$(\Pi_j)_{1\leq j\leq J}$ are nonlinear operators defined on $Y$,
with $\Pi_J(Y)=\{f_\infty\}$. ($\Pi_J$ sends everybody to the
stationary state.) Moreover, for all $j\in \{1,\ldots, J\}$,

(i) $\Pi_j(X)\subset Y$, \ $\C\circ \Pi_j=0$;

(ii) $\Pi_jf_\infty= f_\infty$;

(iii) For any $s$ there is $s'$ large enough such that
$(\Pi_j)'$ and $(\Pi_j)''$ are bounded $X^{s'}\to X^s$,
uniformly on $Y$.
\end{Ass}

The last of these assumptions morally says that $\Pi_j$ is $C^2$ with
possible loss of derivatives.

\begin{Ass}[Lyapunov functional] \label{asslyap}
$\E:Y\to\R$ is $C^1$ on $Y$ viewed as a subset of $X^s$ for
$s$ large enough. For all $f$ one has
$\E(f)\geq \E(\Pi_1f)\geq \E(f_\infty)$, and more precisely

(i) For any $\var\in (0,1)$ there is $K_\var>0$ such that
for all $f\in Y$,
\begeq\label{EEEe}
\E(f) - \E(\Pi_1f) \geq K_\var \, \|f-\Pi_1f\|^{2+\var};
\endeq

(ii) For any $\var\in (0,1)$ there are $K_\var, C_\var>0$ such that
for all $f\in Y$,
\begeq\label{EEEEE}
K_\var\, \|\Pi_1f-f_\infty\|^{2+\var} \leq \E(\Pi_1f)-\E(f_\infty)\leq
C_\var \, \|\Pi_1f-f_\infty\|^{2-\var}.
\endeq
\end{Ass}

Note that $\Pi_1f$ and $f$ are bounded uniformly, so these bounds
become more and more stringent when $\var$ decreases.

\begin{Ass}[Key hypocoercivity assumptions] \label{asskey}
\

(i) $\C$ alone is dissipative, strictly out of the range of $\Pi_1$:
For any $\var>0$ there is a constant $K_\var>0$ such that for all
$f\in X$,
\begeq\label{Cdissip}
-\E'(f)\cdot (\C f) \geq K_\var\, \bigl[ \E(f) - \E(\Pi_1f)\bigr]^{1+\var};
\endeq

(ii) $\C-B$ is dissipative just as well: For any $\var>0$ there is 
$K_\var>0$ such that for all $f\in X$,
\begeq\label{CBdissip}
{\cal D}(f) := -\E'(f)\cdot (\C f- Bf) \geq 
K_\var\, \bigl[ \E(f) - \E(\Pi_1f)\bigr]^{1+\var};
\endeq

(iii) For any $k\leq J$ and for any $\var>0$ there is a constant $K_\var>0$
such that for all $f\in X$,
\begeq\label{DIdPi}
{\cal D}(f) + \sum_{j\leq k}
\Bigl\| (\Id - \Pi_j)'_{\Pi_jf}\cdot (B\Pi_jf)\Bigr\|^2 \geq
K_\var\, \bigl\| (\Pi_k - \Pi_{k+1})f\bigr\|^{2+\var}.
\endeq
\end{Ass}

\begin{Rk}[Simplified assumptions]
In many cases of application, $B$ is conservative,
in the very weak sense that $\E'(f)\cdot (Bf)=0$;
then Assumption~\ref{asskey}(ii) trivially follows from
Assumption~\ref{asskey}(i). Also most of the time, 
Assumption~\ref{asskey}(iii) will be replaced by the stronger property
\begeq\label{key'} 
\Bigl\| (\Id - \Pi_k)'_{\Pi_kf}\cdot (B\Pi_jf)\Bigr\|^2 \geq
K_\var\, \bigl\| (\Pi_k - \Pi_{k+1})f\bigr\|^{2+\var}.
\endeq
In the sequel, I~shall however discuss an important case where
none of these simplifications holds true (Boltzmann equation with 
Maxwellian diffuse boundary condition).
\end{Rk}

\begin{Rk}[Practical verification of the key conditions]
\label{rkpractic}
Often the $\Pi_j$'s are {\em nested projectors}, in the sense
that $\Pi_{j+1}\Pi_j = \Pi_{j+1}$. Then~\eqref{key'} becomes
\[ \Bigl\| (\Id - \Pi_k)'_g\cdot (Bg)\Bigr\|^2 \geq
K_\var\, \bigl\| (\Id - \Pi_{k+1}) g\bigr\|^{2+\var},\qquad
g\in \Pi_k(X).\]
So the recipe is as follows: (a) Take $g\in \Pi_k(X)$,
let it evolve according to $\pa_t g + Bg=0$;
(b) compute $\pa_t (\Pi_kg)$ at $t=0$; (c) check that
$\|Bg + \pa_t (\Pi_kg)\|$ controls $\|g-\Pi_{k+1}g\|^{1+\var}$
for any $\var>0$.
\end{Rk}

\begin{Rk}[Connection with earlier works]
To make the connection with the method used in~\cite{DV:boltz:05},
note that if $g=\Pi_kg$ at $t=0$, then, since $\|\cdot\|$ is Hilbertian,
\[ \Bigl\| (\pa_t)_{t=0} (g_t-\Pi_{k}g_t) \Bigr\|^2
= \left.\frac{d^2}{dt^2}\right|_{t=0} \|g_t - \Pi_k g_t\|^2.\]
In view of this remark, Assumption~\ref{DIdPi} can be understood
as a very abstract reformulation of the property of
``instability of hydrodynamic description'' introduced in~\cite{DV:boltz:05}.
\end{Rk}

Now comes the main nonlinear result in this memoir:

\begin{Thm} \label{thmmain}
Let Assumptions~\ref{assfunct} to~\ref{asskey} be satisfied.
Then, for any $\beta>0$ there is a constant $C_\beta$, only depending
on the constants appearing in these assumptions, on $\beta$ and on 
an upper bound on $\E(f_0)-\E(f_\infty)$, such that
\[ \forall t\geq 0,\qquad \E(f(t)) - \E(f_\infty) \leq C_\beta\, t^{-\beta}.\]
As a consequence, for all $s\geq 0$,
\[ \|f(t)-f_\infty\|_s = O(t^{-\infty}).\]
\end{Thm}

\subsection{Method of proof}

To estimate the speed of approach to equilibrium, the first natural
thing to do is to consider the rate of decay of the Lyapunov functional $\E$.
From the assumptions of Theorem~\ref{thmmain}, if
$\var>0$ is small enough then
\begeq\label{firstattempt}
\frac{d}{dt} [\E(f)-\E(f_\infty)] = 
-{\cal D}(f)\leq -K_\var \bigl[ \E(f)-\E(\Pi_1f)\bigr]^{1+\var}.
\endeq
(I have omitted the explicit dependence of $f$ on $t$.)
But the differential inequality~\eqref{firstattempt} cannot in general
be closed, since $\E(f)-\E(\Pi_1 f)$ might be much smaller
than $\E(f)-\E(f_\infty)$. It may even be the case that
$f=\Pi_1f$, yet $f\neq f_\infty$ (the dissipation vanishes).
So this strategy seems to be doomed.

In~\cite{DV:FP:01,DV:boltz:05} we solved this difficulty by coupling the
differential inequality~\eqref{firstattempt} with some second-order
differential inequalities involving other functionals.
Here on the contrary, I~shall modify the functional
$\E$ by adding some ``lower-order'' terms. So the proofs
in the present paper are based on the following auxiliary functional:
\begeq\label{calL}
\L(f) = \bigl[\E(f)-\E(f_\infty)\bigr] + 
\sum_{j=1}^{J-1} a_j\, \Bigl \< (\Id-\Pi_j)f, 
   \: (\Id-\Pi_j)'_f\cdot (Bf) \Bigr\>,
\endeq
where $\langle\cdot,\cdot\rangle$ denotes the scalar product in $X^0$,
and $a_j>0$ ($1\leq j\leq J-1$) are carefully chosen small numbers,
depending on smoothness bounds on $f$, and also on
upper and lower bounds on $\E(f)-\E(f_\infty)$.

The coefficients $a_j$ will be chosen in such a way that $\L(f)$ is
always comparable to $\E(f)-\E(f_\infty)$; still the time-derivatives
of these two quantities will be very different, and it will be
possible to close the differential inequalities defined in terms
of $\L$.

When the value of $\E(f)-\E(f_\infty)$ has substantially decreased,
then the expression of $\L$ should be re-evaluated (the coefficients
$a_j$ should be updated), so $\L$ in itself does not really define 
a Lyapunov functional. But it will act just the same: On any time-interval
where $\E(f)-\E(f_\infty)$ is controlled from above and below,
one can choose the coefficients $a_j$ in such a way that
$(d/dt)\L(f) \leq -K \L(f)^{1+\delta}$, for any fixed $\delta>0$.
This will be sufficient to control the rate of decay of
$\L$ to~0, and as a consequence the rate of decay of $\E$ to
its minimum value.

Complete proofs will be given in the next section. It is clear that
they enjoy some flexibility and can be slightly modified or 
adapted in case of need.

\section{Proof of the Main Theorem}  \label{secproof}

Theorem~\ref{thmmain} will be obtained as a consequence
of the following more precise result:

\begin{Thm}\label{thmprecise}
Let Assumptions~\ref{assfunct} to~\ref{asskey} be satisfied,
and let $E>0$ be such that
\begeq\label{EEE} 
\frac{E}2 \leq \E(f) -\E(f_\infty) \leq E.
\endeq
Let further
\begeq\label{defL} 
\L(f) = \bigl[\E(f)-\E(f_\infty)\bigr] + 
\sum_{j=1}^{J-1} a_j\, \Bigl \< (\Id-\Pi_j)f, 
   \: (\Id-\Pi_j)'_f\cdot (Bf) \Bigr\>,
\endeq
where $(a_j)_{1\leq j\leq J-1}$ are positive numbers; let $a_0=1$.
Then, 
\sm

(i) For any $\var\in (0,1)$, there 
is a constant $K>0$, depending only on $\var$ and on the
constants appearing in Assumptions~\ref{assfunct} to~\ref{asskey}
(but not on $E$) such that if $a_j \leq K E^\var$ for all $j$, then
\[ \forall f\in X,\qquad \frac{E}4 \leq \L(f) \leq \frac{5E}4.\]

(ii) There are absolute constants $\var_0, k>0$, and there
are constants $K, K'> 0$, depending only
on $\var$, on an upper bound on $\E(f)-\E(f_\infty)$ and on the constants
appearing in Assumptions~\ref{assfunct} to~\ref{asskey} such that,
if $0<\var\leq \var_0$ and 
\[ a_{j+1} \leq a_j;\qquad \frac{a_{j+1}^2}{a_j} \leq K\, 
a_{J-1}^{1+\var}\, E^{k\var},\]
for all $j\in\{0,\ldots, J-2\}$, then
\[\forall f\in X,\qquad 
\L'(f)\cdot (\C f -Bf) \leq -\: a_{J-1}\, K'\,E^{1+\var}.\]
\end{Thm}

\begin{Rk} Lemma~\ref{lemchoice} in Appendix~\ref{toolbox} 
shows that Conditions~(i) and~(ii) can be fulfilled
with $a_{J-1} \geq K_1 E^{\ell \var}$, where $\ell$ only depends on
$J$ and $k$.
\end{Rk}

\begin{Rk} In concrete situations, the explicit form of
$\L$ might be extremely complicated. In the case of the Boltzmann
equation, to be considered later on, the formula for $\L$ requires
eight lines of display.
\end{Rk}

Before starting the proof of Theorem~\ref{thmprecise},
let me make some remarks to facilitate its reading. First of all,
when uniform bounds in the $X^s$ spaces are taken for granted,
a bound from above by, say, $\|f-f_\infty\|_s^\alpha$ is better if the
exponent $\alpha$ is {\em higher}; this is somewhat contrary to what one
is used to when working on smoothness a priori estimates.

In all the sequel the exponents $s$, $s'$ and the constants
$C$, $C'$, $K$, $K'$, etc. may change from one formula to the other.
These quantities can all be computed in terms of an upper bound on
$\E(f_0)-\E(f_\infty)$, the exponents and constants
appearing in Assumptions~\ref{assfunct} to~\ref{asskey} (and for given $\var$,
they only involve a finite number of these constants and exponents).
As a general rule, the symbols $C$, $C'$, etc. will stand for constants which
should be taken {\em large} enough, while the symbols $K$, $K'$, etc. will 
stand for positive constants which should be taken {\em small} enough.

Finally, I~shall frequently use the following fact:
{\em If $\|g\|_{s'}\leq C_{s'}$ for all $s'\geq 0$, then for any $s$ and
any $\delta$ there exists a constant $C$, only depending on $C_{s'}$
for some $s'$ large enough, such that
\begeq\label{interpg}
\|g\|_s \leq C \|g\|_{s'}^{1-\delta}.
\endeq
}
To see this, it suffices to use~\eqref{interp} with
$s_0=0$, $s_1=s/\delta$, $\theta=\delta$.
In other words, it is always possible to replace the norm in some
$X^s$ by the norm in any other $X^s$, up to a arbitrarily small deterioration 
of the exponents.
\med

\begin{proof}[Proof of Theorem~\ref{thmprecise}]
To prove (i), it is sufficient to show that there exists $C$ such that
\begeq\label{boundIPij}
\Bigl| \Bigl\< (\Id-\Pi_j) f,\: 
(\Id-\Pi_j)'_f\cdot (Bf) \Bigr\> \Bigr| \leq C\, [\E(f)-\E(f_\infty)]^{1-\var}
\endeq
for all $f\in X$. Indeed, it will follow from~\eqref{EEE} that
\[ \Bigl| \Bigl\< (\Id-\Pi_j) f,\: 
(\Id-\Pi_j)'_f\cdot (Bf) \Bigr\> \Bigr| \leq \frac{2^\var C}{E^\var}
\,[\E(f)-\E(f_\infty)];\]
then if $a_j\leq KE^\var$, the definition of ${\cal L}$ (formula~\eqref{defL})
will imply
\[ (1- 2^\var JKC)\, [\E(f) - \E(f_\infty)] \leq \L(f) \leq
(1+ 2^\var JKC)\, [\E(f) - \E(f_\infty)].\]
Then the conclusion will be obtained by choosing, say, $K=1/(2^{2+\var} JC)$.
(Here $C$ is the same constant as in~\eqref{boundIPij}.)
\sm

To prove~\eqref{boundIPij}, I~shall first apply the Cauchy-Schwarz inequality,
and bound separately $\|(\Id-\Pi_j)f\|$ and
$\|(\Id-\Pi_j)'_f\cdot (Bf)\|$.
\med

\noindent{\em Bound on $\|(\Id-\Pi_j)f\|$:}
\med

By Assumption~\ref{assproj}(ii), $f-\Pi_j f = 
(f-f_\infty) - (\Pi_j f - \Pi_j f_\infty)$, so
\[ \|f-\Pi_j f\| \leq \|f - f_\infty\| + \|\Pi_j f - \Pi_j f_\infty\|.\]
By Assumption~\ref{assproj}(iii) and the convexity of $Y$,
$\Pi_j$ is Lipschitz $X^s\to X^0$ for some $s$ large enough; so
\[ \|\Pi_j f - \Pi_j f_\infty\| \leq C \|f-f_\infty\|_s.\]
Both $f$ and $f_\infty$ belong to $Y$, so by 
Assumption~\ref{asssol} they are bounded in $X^{s'}$ for all $s'$, 
and we can apply the interpolation inequality~\eqref{interpg}:
\[ \|f-f_\infty\|_s \leq C \|f-f_\infty\|^{1-\frac{\var}2}.\]
Then by Assumption~\ref{asslyap}(i)-(ii),
\[ \|f-f_\infty\|^{1-\frac{\var}2} \leq C [\E(f)-\E(f_\infty)]
^{\frac12 -\frac{\var}2}.\]
All in all,
\begeq\label{boundPij} 
\|f-\Pi_j f\| \leq C [\E(f)-\E(f_\infty)]^{\frac12 -\frac{\var}2}.
\endeq
\med

\noindent{\em Bound on $\|(\Id-\Pi_j)'_f\cdot (Bf)\|$:}
\med

By Assumption~\ref{assproj}(iii), there are constants $C$ and $s$ such that
\[ \|(\Id-\Pi_j)'_f\cdot (Bf) \| \leq C \|Bf\|_s.\]
By Assumption~\ref{asseq}(i), $Bf$ is bounded in all spaces $X^{s'}$, so
by interpolation,
\[ \|Bf\|_s \leq C \|Bf\|^{1-\frac{\var}4}.\]
It follows from Assumption~\ref{asseq}(ii) and the convexity of $Y$
that $B$ is Lipschitz $X^s\to X^0$ on $Y$; in view of 
Assumption~\ref{assstat} ($Bf_\infty=0$), this leads to
\[ \|Bf\|^{1-\frac{\var}4} 
= \|Bf - Bf_\infty\|^{1-\frac{\var}4}
\leq C \|f-f_\infty\|_s^{1-\frac{\var}4}.\]
The end of the estimate is just as before:
\[ \|f-f_\infty\|_s^{1-\frac{\var}4} 
\leq C \|f-f_\infty\|^{1-\frac{\var}2} \leq
C' [\E(f)-\E(f_\infty)]^{\frac12 - \frac{\var}2}.\]
All in all,
\[ \|(\Id-\Pi_j)'_f\cdot (Bf) \|  \leq 
C\,[\E(f)-\E(f_\infty)]^{\frac12 - \frac{\var}2}.\]
This combined with~\eqref{boundPij} establishes~\eqref{boundIPij}.
\med

Now we turn to the proof of (ii), which is considerably more tricky.
Let
\begeq\label{defDtilde}
\tilde{\cal D}(f) := - {\cal L}(f)\cdot (\C f - Bf).
\endeq
The argument will be divided in three steps.
\med

\noindent {\bf Step~1:} The estimates in this step are mainly based on
regularity assumptions.
\sm

By direct computation,
\begin{align*}
- \tilde{\cal D}(f) = - {\cal D}(f)
& + \sum_{j=1}^{J-1} a_j \Blangle (\Id-\Pi_j)'_{f}\cdot(\C f-Bf),
\: (\Id-\Pi_j)'_{f}\cdot(Bf)) \Brangle \\
& + \sum_{j=1}^{J-1} a_j  \Blangle (\Id-\Pi_j)f, \; (\Id-\Pi_j)''_f\cdot
(\C f-Bf, Bf) \Brangle \\
&+ \sum_{j=1}^{J-1} a_j  \Blangle (\Id-\Pi_j)f, \; 
(\Id-\Pi_j)'_f\cdot\bigl(B'_f\cdot(\C f-Bf)\bigr) \Brangle 
\end{align*}
\begin{align*}
= -{\cal D}(f)
&- \sum_{j=1}^{J-1} a_j \bigl\| (\Id-\Pi_j)'_f\cdot(Bf) \bigr\|^2 \\
&+\sum_{j=1}^{J-1}a_j \Blangle (\Id-\Pi_j)'_{f}\cdot(\C f),
\: (\Id-\Pi_j)'_{f}\cdot(Bf) \Brangle \\
&+ \sum_{j=1}^{J-1} a_j \Blangle (\Id-\Pi_j)f,\: (\Id-\Pi_j)''_{f}\cdot(\C f-Bf, Bf) \Brangle \\
&+ \sum_{j=1}^{J-1} a_j  \Blangle (\Id-\Pi_j)f, \; 
(\Id-\Pi_j)'_f\cdot\bigl(B'_f\cdot(\C f-Bf)\bigr) \Brangle
\end{align*}

Then by Cauchy-Schwarz inequality,
\begin{align} \label{25}
-\tilde{\cal D}(f) \leq -{\cal D}(f)
&- \sum_{j=1}^{J-1} a_j \bigl\| (\Id-\Pi_j)'_f\cdot(Bf) \bigr\|^2 \\
&+ \sum_{j=1}^{J-1} a_j \bigl\|(\Id-\Pi_j)'_{f}\cdot(\C f)\bigr\|
\bigl\| (\Id-\Pi_j)'_{f}\cdot(Bf) \bigr\| \notag \\
&+ \sum_{j=1}^{J-1} a_j 
\bigl\|(\Id-\Pi_j)f\bigr\| \bigl\|(\Id-\Pi_j)''_{f}\cdot(\C f-Bf, Bf)\bigr\|
\notag\\
&+ \sum_{j=1}^{J-1} a_j  \bigl\| (\Id-\Pi_j)f\bigr\|\; 
\bigl\|(\Id-\Pi_j)'_f\cdot\bigl(B'_f\cdot(\C f-Bf)\bigr) \bigr\|.\notag
\end{align}

By applying the inequality $ab\leq (a^2+b^2)/2$, with
$a=\|(\Id-\Pi_j)'_f\cdot(Bf)\|$ and $b=\|(\Id-\Pi_j)'_f\cdot(\C f)\|$,
we see that the second and third terms in the right-hand side of~\eqref{25}
can be bounded by
\begeq\label{12sum} 
- \frac12\sum_{j=1}^{J-1} a_j \bigl\| (\Id-\Pi_j)'_f\cdot(Bf) \bigr\|^2
+ \frac12 \sum_{j=1}^{J-1} a_j \bigl\|(\Id-\Pi_j)'_f\cdot(\C f)\bigr\|^2.
\endeq
Then we apply the Hilbertian inequality
\[ -\|a\|^2 \leq -\,\frac{\|b\|^2}2 + \|b-a\|^2\]
with $a=(\Id-\Pi_j)'_f\cdot(Bf)$ and $b=(\Id-\Pi_j)'_{\Pi_jf}\cdot(B\Pi_jf)$,
to bound~\eqref{12sum} by
\begin{multline*} 
- \frac14\sum_{j=1}^{J-1} a_j 
\bigl\|(\Id-\Pi_j)'_{\Pi_jf}\cdot (B\Pi_jf) \bigr\|^2 
+ \frac12 \sum_{j=1}^{J-1} a_j
\bigl\|(\Id-\Pi_j)'_{\Pi_jf}\cdot(B\Pi_jf) 
- (\Id-\Pi_j)'_f\cdot(Bf)\bigr\|^2 \\
+ \frac12 \sum_{j=1}^{J-1} a_j \bigl\|(\Id-\Pi_j)'_f\cdot(\C f)\bigr\|^2.
\end{multline*}
It follows, after plugging these bounds back in~\eqref{25}, that
\begeq\label{expr1}
-\tilde{\cal D}(f) \leq -{\cal D}(f)
- \frac14\sum_{j=1}^{J-1} 
a_j \bigl\|(\Id-\Pi_j)'_{\Pi_jf}\cdot (B\Pi_jf) \bigr\|^2 
+ \ \sum_{j=1}^{J-1} a_j\, (R)_j,
\endeq
where
\begin{align}
(R)_j:= \ & \frac12\, 
\Bigl\|(\Id-\Pi_j)'_{\Pi_jf}\cdot (B\Pi_jf) - (\Id-\Pi_j)'_f\cdot(Bf)\Bigr\|^2 
\label{26}\\
& + 
    \frac12\,  \bigl\|(\Id-\Pi_j)'_f\cdot (\C f)\bigr\|^2 \notag\\
& + 
   \|(\Id-\Pi_j)f\| 
   \Bigl( \bigl\|(\Id-\Pi_j)''_f\cdot(\C f - Bf, Bf) \bigr\|\notag \\
&\qquad\qquad\qquad\qquad\qquad\qquad 
   + \bigl\|(\Id-\Pi_j)'_f\cdot(B'_f\cdot(\C f-Bf))\bigr\| \Bigr). 
\notag
\end{align}

Now I~shall estimate the various terms in~\eqref{26} one after the other.
\med

\noindent{\em  First line of~\eqref{26}:}
\med

First,
\begin{multline}\label{1stline} 
\bigl\|(\Id-\Pi_j)'_{\Pi_jf}\cdot (B\Pi_j f) -
(\Id-\Pi_j)'_f\cdot (Bf)\bigr\|\leq
\bigl\|(\Id-\Pi_j)'_{\Pi_jf}\cdot(Bf - B\Pi_jf)\bigr\| \\
+ \bigl\|[(\Pi_j)'_{\Pi_jf} - (\Pi_j)'_f]\cdot (Bf)\bigr\|.
\end{multline}

By Assumption~\ref{assproj}(iii),
\[ \bigl\|(\Id-\Pi_j)'_{\Pi_jf}\cdot(Bf - B\Pi_jf)\bigr\|
\leq C \|Bf-B\Pi_jf\|_s.\]
(Here I~use the fact that $\Pi_j(X)\subset Y$.)
Also, from Assumptions~\ref{asssol} and~\ref{asseq}
(and again $\Pi_j(X)\subset Y$), $Bf$ and $B\Pi_j f$ are bounded 
in all spaces $X^{s'}$, so by interpolation
\[ \|Bf-B\Pi_jf\|_s \leq C \|Bf-B\Pi_jf\|^{1-\frac{\var}2}.\]
As a consequence of Assumption~\ref{asseq}(ii) and the convexity of $Y$,
$B$ is Lipschitz continuous $X^s\to X^0$ on $Y$, so
\[ \|Bf-B\Pi_jf\|^{1-\frac{\var}2} \leq C \|f-\Pi_jf\|_s^{1-\frac{\var}2}
\leq C \|f-\Pi_jf\|^{1-\var},\]
where the last inequality is obtained again from interpolation.
This provides a bound for the first term on the right-hand side
of~\eqref{1stline}
\sm

Next, as a consequence of Assumption~\ref{assproj}(iii), $(\Pi_j)'$ is 
Lipschitz continuous on $Y$, in the sense that for all $f,g\in Y$,
\[ \bigl\|[(\Pi_j)'_f - (\Pi_j)'_g]\cdot h\bigr\| \leq C \|f-g\|_s \|h\|_s.\]
Combining this with Assumption~\ref{asseq}, we find
\begin{multline*} \bigl\|[(\Pi_j)'_{\Pi_jf} - (\Pi_j)'_f]\cdot (Bf)\bigr\|\leq
C \|\Pi_jf - f\|_s\, \|Bf\|_s \\ \leq C' \|f-\Pi_j f\|_s\leq
C'' \|f-\Pi_jf\|^{1-\var}.
\end{multline*}
This takes care of the second term on the right-hand side
of~\eqref{1stline}. The conclusion is that the first line of~\eqref{26} 
is bounded by $O(\|f-\Pi_jf\|^{1-\var})$, for any $\var\in (0,1)$.
\med

\noindent{\em Second line of~\eqref{26}:}
\med

First, by Assumption~\ref{assproj}(iii),
\[ \bigl\|(\Id-\Pi_j)'_f\cdot (\C f) \bigr\| \leq C \|\C f\|_s.\]
By Assumption~\ref{asseq}(i), $\C f$ is bounded in all spaces $X^{s'}$,
so by interpolation:
\[ \|\C f\|_s\leq C \|\C f\|^{1-\frac{\var}2}.\]
By Assumption~\ref{assproj}(i), $\C \Pi_j f=0$; and
by Assumption~\ref{asseq}(iii), $\C$ is Lipschitz $X^s\to X^0$ on $Y$; so
\[ \|\C f\|^{1-\frac{\var}2} = \|\C f - \C \Pi_j f\|^{1-\frac{\var}2}
\leq C \|f-\Pi_j f\|^{1-\frac{\var}2}_s \leq 
C' \|f-\Pi_jf \|^{1-\var}.\]
The conclusion is that the second line of~\eqref{26} can be
bounded just as the first line, by $O(\|f-\Pi_jf\|^{1-\var})$,
for any $\var\in (0,1)$.
\med

\noindent{\em Third and fourth lines of~\eqref{26}:}
\med

By Assumption~\ref{assproj}(iii),
\[ \bigl\| (\Id-\Pi_j)''_f\cdot (\C f-Bf, Bf) \| \leq 
C \|\C f - B f\|_s\, \|Bf\|_s \leq
C' (\|Bf\|_s^2 + \|\C f\|_s^2).\]
The second term $\|\C f\|_s^2$ can be bounded by 
$O(\|f-\Pi_kf\|^{2-2\var})$, as we already saw; 
by taking $k=J$ we get a bound like $O(\|f-f_\infty\|^{2-\var})$.
As for the first term $\|Bf\|_s^2$, we saw before that 
it is also bounded like $O(\|f-f_\infty\|^{2-\var})$.
In the sequel I~shall only keep the worse bound $O(\|f-f_\infty\|^{1-\var})$.

Next, by Assumption~\ref{assproj}(iii) again,
\[ \bigl\| (\Id-\Pi_j)'_f\cdot (B'_f\cdot (\C f - Bf))\bigr\|
\leq C \|B'_f\cdot (\C f - B f)\|_s.\]
From Assumption~\ref{asseq}(ii),
\[ \|B'_f\cdot (\C f - B f)\|_s \leq C \|\C f - B f\|_{s'}
\leq C \bigl(\|\C f\|_{s'} + \|Bf\|_{s'}\bigr),\]
and as before this can be controlled by
$O(\|f-f_\infty\|^{1-\var})$.

The conclusion is that the third and fourth lines of~\eqref{26} can
be bounded by $C \|f-\Pi_jf\|\, \|f-f_\infty\|^{1-\var}$, for
any $\var\in (0,1)$.
\med

Gathering all these estimates and replacing $\var$ by $\var/2$, we
deduce that the expression in~\eqref{26} can be bounded as follows:
\[ (R)_j \leq C \
\Bigl( \|f-\Pi_j f\|^{2-\frac{\var}2} + 
\|f-\Pi_jf\|\, \|f-f_\infty\|^{1-\frac{\var}2}\Bigr).\]
As we already saw before,
\[ \|f-\Pi_jf\|^{1-\frac{\var}2} \leq C \|f-f_\infty\|^{1-\var},\]
so actually
\begeq\label{soR} (R)_j \leq C\, \|f-\Pi_jf\|\, \|f-f_\infty\|^{1-\var}.
\endeq

The temporary conclusion is that
\begin{multline}\label{tempconcl}
-\tilde{\cal D}(f) \leq - {\cal D}(f)
-\frac14 \sum_{j=1}^{J-1} a_j\,
\bigl\| (\Id-\Pi_j)'_{\Pi_jf}\cdot (B\Pi_jf) \bigr\|^2\: \\
+ \Bigl( \sum_{j=1}^{J-1} a_j\, \|f-\Pi_jf\|\Bigr)\, \|f-f_\infty\|^{1-\var}.
\end{multline}
\med

\noindent{\bf Step~2:} This step uses Assumption~\ref{asskey} crucially.
\med

By triangle inequality,
\[ \|f-\Pi_jf\| =\|\Pi_0 f - \Pi_j f\|\leq
\sum_{0\leq k\leq j-1} \|\Pi_k f - \Pi_{k+1} f\|;\]
so
\begin{align*} \sum_{1\leq j\leq J-1} a_j \|f-\Pi_j f\| 
& \leq \sum_{1\leq j\leq J-1} a_j \Bigl(
\sum_{0\leq k\leq j-1} \|\Pi_k f -\Pi_{k+1} f\|\Bigr) \\
& = \sum_{0\leq k\leq J-2} \Bigl(
\sum_{k+1\leq j\leq J-1} a_j\Bigr) \|\Pi_kf-\Pi_{k+1}f\|\\
& \leq \sum_{0\leq k\leq J-2} (Ja_{k+1})\,
\|\Pi_k f - \Pi_{k+1}f\|,
\end{align*}
where the last inequality follows from the fact that the 
sequence $(a_j)_{0\leq j\leq J-1}$ is nonincreasing.

Renaming $k$ as $j$, plugging this inequality back in~\eqref{tempconcl},
we arrive at
\begin{align}\label{new26}
-\tilde{\cal D}(f) \leq 
- {\cal D}(f) & - \frac14 \, \sum_{0\leq j\leq J-2}
a_j\, \Bigl\| (\Id-\Pi_j)'_{\Pi_jf}\cdot (B\Pi_jf)\Bigr\|^2 \\
& + C \, \sum_{0\leq j\leq J-2} a_{j+1}\,
\bigl\|\Pi_jf - \Pi_{j+1}f\bigr\|\, \|f-f_\infty\|^{1-\var}.\nonumber
\end{align}

Since $a_k\leq 1$, we can write
\begin{align}\label{akleq1}
{\cal D}(f) & = \frac{{\cal D}(f)}{2} + \frac1{2(J-1)}\,
\sum_{k=1}^{J-1} {\cal D}(f) \\
& \geq \frac{{\cal D}(f)}2 + \frac1{2(J-1)}\,
\sum_{k=1}^{J-1} a_k {\cal D}(f). \nonumber
\end{align}
On the other hand, $k\geq j \Longrightarrow a_j\geq a_k$, so
\begin{align}\label{kgeqj}
\sum_{j=1}^{J-1}a_j \Bigl\| (\Id-\Pi_j)'_{\Pi_jf}\cdot (B\Pi_jf)\Bigr\|^2
& = \frac1{J-1} \sum_{j=1}^{J-1} \sum_{k=1}^{J-1} a_j\,
\Bigl\|(\Id-\Pi_j)'_{\Pi_jf}\cdot (B\Pi_jf)\Bigr\|^2 \\
& \geq \frac1{J-1} \sum_{j=1}^{J-1}\sum_{k=1}^{J-1} a_k\,
\Bigl\|(\Id-\Pi_j)'_{\Pi_jf}\cdot (B\Pi_jf)\Bigr\|^2. \nonumber
\end{align}
From~\eqref{new26}, \eqref{akleq1} and~\eqref{kgeqj},
\begin{align*}
-\tilde{\cal D}(f) &\leq
- \frac{{\cal D}(f)}2 \\
& - \frac1{4(J-1)}\,
\sum_{k=1}^{J-1} a_k\, 
\Bigl( {\cal D}(f) + \sum_{j=1}^k \Bigl\|
(\Id-\Pi_j)'_{\Pi_jf}\cdot (B\Pi_jf)\Bigr\|^2 \Bigr) \\
& + C\, \sum_{j=0}^{J-2} a_{j+1}\,
\bigl\| \Pi_jf -\Pi_{j+1}f\bigr\|\, \|f-f_\infty\|^{1-\var}.
\end{align*}
At this point we can apply Assumption~\ref{asskey}(ii)-(iii)
and Assumption~\ref{asslyap} and we get constants $K,C$ such that
\begin{align}\label{expr3}
-\tilde{\cal D}(f) \leq & - K \Bigl( [\E(f) - \E(\Pi_1f)]^{1+\frac{\var}4} 
+ \|f-\Pi_1f\|^{2+\var} \Bigr)\\
& - K \sum_{k=1}^{J-1}a_k\, \|\Pi_kf - \Pi_{k+1}f\|^{2+\var} \nonumber\\
& + C \sum_{j=0}^{J-2} a_{j+1}\, \|\Pi_{j}f-\Pi_{j+1}f\|\, 
\|f-f_\infty\|^{1-\var}. \nonumber
\end{align}

By applying Young's inequality, in the form
\[ a X Y^{1-\var} \leq b \, \frac{X^{2+\var}}{2+\var}
\: + \: \left(\frac{a^{\frac{2+\var}{1+\var}}}{b^{\frac1{1+\var}}}\right)\,
\frac{(Y^{1-\var})^{\left(\frac{2+\var}{1+\var}\right)}}
{\left(\frac{2+\var}{1+\var}\right)},\]
with $a=a_{j+1}$, $X=\|\Pi_jf-\Pi_{j+1}f\|$, $Y=\|f-f_\infty\|$,
$b=K a_j$ in the last line of~\eqref{expr3}, and get
\begin{align}\label{expr4}
-\tilde{\cal D}(f) \leq & - K [ \E(f) - \E(\Pi_1f)]^{1+\frac{\var}4} 
- K \|f-\Pi_1f\|^{2+\var} \\
& - K \sum_{j=1}^{J-1}a_j \|\Pi_jf - \Pi_{j+1}f\|^{2+\var}\nonumber \\
& + C \sum_{0\leq j\leq J-2} \left(\frac{a_{j+1}^{2+\var}}{a_j}\right)
^{\frac1{1+\var}} \|f-f_\infty\|^{\frac{(1-\var)(2+\var)}{1+\var}}.
\nonumber
\end{align}
Since $a_{j+1}\leq 1$, we can bound trivially $a_{j+1}^{2+\var}$ by
$a_{j+1}^2$. Moreover, for $\var\leq \var_0$ small enough,
we have $(1-\var)(2+\var)/(1+\var)\geq 2-4\var$ and, so
\[ \|f-f_\infty\|^{\frac{(1-\var)(2+\var)}{1+\var}} \leq 
  C \|f-f_\infty\|^{2-4\var}.\]
Taking into account once again the fact that $a_j\geq a_{J-1}$ for
all $j$, \eqref{expr4} implies, with the convention $\Pi_0f=f$, $a_0=1$,
\begin{align}\label{expr5}
-\tilde{\cal D}(f) \leq & - K [ \E(f) - \E(\Pi_1f)]^{1+\frac{\var}4} \\
& - K a_{J-1} \left(\sum_{0\leq j\leq J-1} 
\|\Pi_jf-\Pi_{j+1}f\|^{2+\var} \right) \nonumber\\
& + C \sup_{0\leq j\leq J-2} \left(\frac{a_{j+1}^2}{a_j}\right)
^{\frac1{1+\var}} \|f-f_\infty\|^{2-4\var}.\nonumber
\end{align}
Next,
\[ \|f-f_\infty\|^{2+\var} = \|\Pi_0 f - \Pi_J f\|^{2+\var} \leq 
C \sum_{0\leq j\leq J-1} \|\Pi_jf - \Pi_{j+1}f\|^{2+\var}, \]
so from~\eqref{expr5} we deduce
\begin{multline}\label{expr6}
-\tilde{\cal D}(f) \leq - K [ \E(f) - \E(\Pi_1f)]^{1+\frac{\var}4} 
- K\, a_{J-1} \|f-f_\infty\|^{2+\var}
\\
+ C \sup_{0\leq j\leq J-2} \left(\frac{a_{j+1}^2}{a_j}\right)
^{\frac1{1+\var}} \|f-f_\infty\|^{2-4\var}.
\end{multline}
\med

\noindent{\bf Step~3:} Now a few complications will arise because we only 
have a control {\em from below} of $\E(f)-\E(f_\infty)$ in terms 
of $\|f-f_\infty\|$; so the fact that $\E(f)-\E(f_\infty)$ is of order 
$E$ does not imply any lower bound on $\|f-f_\infty\|$, and then
$\|f-f_\infty\|^{2-4\var}$ might be much, much higher 
than $\|f-f_\infty\|^{2+\var}$. To solve this difficulty,
a little additional detour will be useful.

From Assumption~\ref{assproj}(ii)-(iii) and interpolation,
\begeq\label{firstdetour} 
\|\Pi_1 f - f_\infty\|^{2+2\var} =
\|\Pi_1 f - \Pi_1 f_\infty\|^{2+2\var} 
\leq C \|f-f_\infty\|^{2+\var};
\endeq
on the other hand,
\begeq\label{seconddetour} 
\|f-f_\infty\|^{2-4\var} \leq C 
\bigl( \|f-\Pi_1f\|^{2-4\var} + \|\Pi_1 f - f_\infty\|^{2-4\var}\bigr).
\endeq
By using~\eqref{firstdetour} and~\eqref{seconddetour}
in~\eqref{expr5}, and replacing the exponent $1+\var/4$ by
the worse exponent $1+2\var$ (which is allowed since
$\E(f) - \E(\Pi_1f) \leq \E(f)-\E(f_\infty)$ is uniformly bounded), 
we obtain
\begin{multline}\label{expr7}
-\tilde{\cal D}(f) \leq - K [ \E(f) - \E(\Pi_1f)]^{1+2\var} 
- K\, a_{J-1} \|\Pi_1 f-f_\infty\|^{2+2\var}
\\
+ C \delta^{\frac1{1+\var}} \bigl(\|f-\Pi_1f\|^{2-4\var}+
\|\Pi_1f-f_\infty\|^{2-4\var}\bigr),
\end{multline}
where 
\[ \delta:= \max_{0\leq j\leq J-1}\ \frac{a_{j+1}^2}{a_j}.\]

Then from Assumption~\eqref{asslyap}(i)-(ii) (both the
upper and the lower bounds are used in (ii)),
\begin{multline}\label{expr8}
-\tilde{\cal D}(f) \leq - K [ \E(f) - \E(\Pi_1f)]^{1+2\var} 
- K\, a_{J-1} [\E(\Pi_1f) - \E(f_\infty)]^{1+2\var}
\\
+ C \delta^{\frac1{1+\var}} \,[\E(f) - \E(\Pi_1f)]^{1-3\var}
+ C\delta^{\frac1{1+\var}} \, [\E(\Pi_1f) - \E(f_\infty)]^{1-3\var}.
\end{multline}

Let us distinguish two cases:
\med

\noindent{\em First case: $\E(\Pi_1f)-\E(f_\infty) \leq
\E(f) - \E(\Pi_1f)$.}
\med

Then
\[ \E(f) - \E(f_\infty) \leq 2 [\E(f) - \E(\Pi_1f)],\]
and in particular
\begeq\label{EEPif}
\E(f) - \E(\Pi_1f) \geq \frac{E}4.
\endeq
In that case we throw away the second negative term in~\eqref{expr8},
and bound the last term by the but-to-last one:
\begin{align}
-\tilde{\cal D}(f) & \leq - K [ \E(f) - \E(\Pi_1f)]^{1+2\var} 
+ C \delta^{\frac1{1+\var}} \,[\E(f) - \E(\Pi_1f)]^{1-3\var}\notag\\
& = -K \left( 1 - \frac{C\delta^{\frac1{1+\var}}}{K}\, 
[\E(f) - \E(\Pi_1f)]^{-5\var} \right)\, 
[\E(f) - \E(\Pi_1f)]^{1+2\var}. \label{lastbef8}
\end{align}
If
\begeq\label{aaa}
\delta^\frac1{1+\var} \leq \frac{K}{2C}\, \bigl[\E(f)-\E(\Pi_1f)\bigr]^{5\var}
\endeq
(where $K$ and $C$ are the same constants as in~\eqref{lastbef8}),
then~\eqref{lastbef8} can be bounded above by
$-K' [\E(f) -\E(\Pi_1f)]^{1+2\var}$, and by~\eqref{EEPif}
this can also be bounded above by $-K'' E^{1+2\var}$.

Finally, in view of~\eqref{EEPif} again, \eqref{aaa} is satisfied if
\begeq\label{suffdelta} 
\delta^{\frac1{1+\var}} \leq K' E^{5\var},
\endeq
where $K'=4^{-5\var}K/(2C)$. Since $\var\leq 1$ and $E$ is uniformly
bounded, a sufficient condition for~\eqref{suffdelta} to hold
is $\delta\leq K'' E^{10\var}$.

\med

\noindent{\em Second case: $\E(\Pi_1f)-\E(f_\infty) \geq
\E(f) - \E(\Pi_1f)$.} In that case
\begeq\label{EEPif'}
\E(\Pi_1f) - \E(f_\infty) \geq \frac{E}4,
\endeq
and we retain from~\eqref{expr8} that
\[
-\tilde{\cal D}(f) \leq
- K\, a_{J-1} [\E(\Pi_1f) - \E(f_\infty)]^{1+2\var}
+ C\delta^{\frac1{1+\var}} \, \bigl[\E(\Pi_1f) - \E(f_\infty)\bigr]^{1-3\var};
\]
and by a reasoning similar as the one above, this is bounded above by
\[ -\, \frac{K}2\, a_{J-1}\, \bigl[ \E(\Pi_1 f) - 
\E(f_\infty)\bigr]^{1+2\var}\]
as soon as 
\[ \delta^{\frac1{1+\var}} \leq K' a_{J-1} E^{5\var}.\]
This condition is fulfilled as soon as
\[ \delta\leq K'' a_{J-1}^{1+\var} E^{10\var},\]
and, a fortiori (since $a_{J-1}\leq 1$) if
\[ \delta\leq K''' a_{J-1}^{1+2\var} E^{10\var}.\]
\med

In both cases, we have concluded that if $\var\leq \var_0$ and
\[ \frac{a_{j+1}^2}{a_j} \leq K a_{J-1}^{1+2\var} E^{10\var},\]
then
\[-\tilde{\cal D}(f) = \L'(f)\cdot (\C f-Bf) \leq - K' a_{J-1}\, E^{1+2\var}.\]
Up to the replacement of $\var$ by $\var/2$ and $\var_0$ by
$\var_0/2$, this is exactly the desired conclusion.
\end{proof}

\begin{proof}[Proof of Theorem~\ref{thmmain}]
Let $\ov{E}$ be such that
$\E(f_0) - \E(f_\infty)\leq \ov{E}$. Since $\E(f(t))$ is
a nonincreasing function of $t$,
\[\forall t\geq 0,\qquad \E(f(t)) - \E(f_\infty) \leq \ov{E}.\]

Let now $\var>0$, and $E\in (0,\ov{E}]$. Let $[t_0, t_0+T]$
be the time-interval where
\[\frac{E}2\leq \E(f(t)) - \E(f_\infty) \leq E;\]
this time-interval is well-defined (at least if $T$ is a priori allowed
to be infinite) since $\E(f(t))-\E(f_\infty)$ is
a continuous nonincreasing function.
The goal is to show that if $\var$ is small enough, then 
\begeq\label{Tleq}
T\leq C E^{-\lambda\var},
\endeq
where $\lambda$ only depends on $J$,
and $C$ may depend on $\var$ but not on $E$.
When~\eqref{Tleq} is proven, it will follow from a classical
argument that 
\begeq\label{Eftfinfty}
\E(f(t))-\E(f_\infty)=O(t^{-1/((\lambda+1)\var)}).
\endeq
Indeed, let $E_0:=\E(f_0)-\E(f_\infty)$; then $\E(f(t))-\E(f_\infty)$
will be bounded by $E_0/2^{m+1}$ after a time
\begin{align*} T_m:= & C \Bigl( E_0^{-\lambda\var} + 
\left(\frac{E_0}2\right)^{-\lambda\var}
+\left(\frac{E_0}4\right)^{-\lambda\var} + \ldots + 
\left(\frac{E_0}{2^m}\right)^{-\lambda\var} \Bigr)\\
& \leq C \bigl( \sum_{j=0}^{m} 2^{\lambda j\var}\bigr) E_0^{-\lambda\var}
\leq C' 2^{\lambda m \var} E_0^{-\lambda\var}.
\end{align*}
So $\E(f(t))-\E(f_\infty) = O(2^{-(m+1)})$ after a time
proportional to $2^{\lambda m\var}$, and~\eqref{Eftfinfty} follows 
immediately.

Then $\E(f(t)) -\E(f_\infty) = O(t^{-\infty})$ since $\var$ 
can be chosen arbitrarily small and $\lambda$ does not depend on $\var$.
From Assumption~\ref{asslyap}(i)-(ii),
\[ \|f(t)-f_\infty\| \leq C [\E(f(t)) - E(f_\infty)]^{1/3}\]
(here $1/3$ could be $1/2-\var$), so $\|f(t)-f_\infty\|=O(t^{-\infty})$
also. Finally, since $f(t)$ is bounded in all spaces $X^s$ by
Assumption~\ref{asssol}, it follows by interpolation that 
$\|f(t)-f_\infty\|_s = O(t^{-\infty})$ for any $s>0$.

So it all amounts to proving~\eqref{Tleq}. Let $K, K',k,\var_0$
be provided by Theorem~\ref{thmprecise}. (There is no loss of
generality in taking the constants $K$ appearing in (i) and (ii)
to be equal.) Let then $\var_1, K_1, \ell$ be provided
by Lemma~\ref{lemchoice}. For any $\var\leq\min(\var_0,\var_1)$ and
any $t\in [t_0,t_0+T]$ we have
\[ \frac{E}4 \leq \L(f(t)) \leq \frac{5E}4;\qquad
\frac{d}{dt}\, \bigl[ {\cal L}(f(t))\bigr] \leq -K' a_{J-1}E^{1+\var} \leq
-K' K_1 E^{1+(\ell+1)\var}.\]
So
\[ E \geq \L(f(t_0)) - \L(f(t_0+T)) \geq
T\, K'K_1\, E^{1+(\ell+1)\var},\]
and then~\eqref{Tleq} follows with $\lambda=\ell+1$
(which eventually depends only on $J$).
\end{proof}

\section{Compressible Navier--Stokes system}\label{secCNS}

In this section I~start to show how to apply Theorem~\ref{thmmain}
on ``concrete'' examples.

The compressible Navier--Stokes equations take the general form
\begeq\label{CNS}
\begin{cases}
\pa_t\rho + \nabla\cdot (\rho u) = 0\\
\pa_t (\rho u) + \nabla\cdot (\rho u \otimes u) + \nabla p =
\nabla\cdot \tau\\
\pa_t(\rho e) + \nabla\cdot (\rho e u + p u) = \nabla\cdot (\tau u)
- \nabla\cdot q
\end{cases}
\endeq
where $\rho$ is the density, $u$ (vector-valued) is the velocity,
$e$ is the energy, $q$ (vector-valued) is the heat flux, and
$\tau$ (matrix-valued) is the viscous stress. In the case
of perfect gases in dimension $N$, it is natural to use the following
constitutive laws:
\begeq\label{constlaws}
\begin{cases}
\dps p = \rho T \\ \\
\dps e = \frac{|u|^2}2 + \frac{N}2 T \\ \\
\dps \tau = 2\mu\, \{\nabla u\} \\ \\
\dps q = - \frac{N}2 \kappa\, \nabla T,
\end{cases}\endeq
where $T$ is the temperature, $\mu$ is the viscosity,
$\kappa$ is the heat conductivity, and $\{\nabla u\}$ (matrix-valued)
is the traceless symmetric strain:
\[ \{\nabla u\}_{ij} = \frac12 \left(
\derpar{u_i}{x_j} + \derpar{u_j}{x_i}\right) -
\left(\frac{\nabla\cdot u}{N}\right)\, \delta_{ij},\]
and $\delta_{ij}=1_{i=j}$.
Then~\eqref{CNS} takes the form
\begeq\label{CNS'}
\begin{cases}
\dps \pa_t\rho + \nabla\cdot (\rho u) = 0;\\ \\
\dps \pa_t (\rho u) + \nabla\cdot (\rho u \otimes u) + \nabla (\rho T) =
2\mu\, \nabla\cdot \{\nabla u\};\\ \\
\dps \pa_t\left(\rho \frac{|u|^2}2 + \frac{N}2\, \rho T \right) + 
\nabla\cdot \Bigl(\rho \frac{|u|^2}2\, u\,  +\,
\left(\frac{N+2}{2}\right) \rho u T \Bigr) \\
\qquad\qquad\qquad\qquad\qquad\qquad\qquad = 2\mu\, \nabla\cdot
(u\{\nabla u\}) + \frac{N}2\, \kappa \Delta T.
\end{cases}
\endeq
(Note that
\[ \nabla\cdot \{\nabla u\} = \mu \Delta u + 
\mu \left(1-\frac2N\right) \nabla\nabla\cdot u,\]
so in the case considered here, the second Lam\'e coefficient 
is negative and equal to $-(2/N)\mu$, which is the borderline case.)

To avoid discussing boundary conditions I~shall only consider the
case when $x$ varies in the torus $\T^N$.
(Later on, for the Boltzmann equation we'll come to grips 
with boundary conditions a bit more.)

There are $N+2$ conservation laws for~\eqref{CNS}: total mass,
total momentum ($N$ scalar quantities) and total kinetic energy. 
Without loss of generality I~shall assume
\begeq\label{CLCNS} 
\int \rho =1;\qquad 
\int \rho u =0 ;\qquad
\int \rho \frac{|u|^2}{2} + \frac{N}2 \int \rho T = \frac{N}2.
\endeq

There is an obvious stationary state: $(\rho,u,T)\equiv (1,0,1)$.
The goal of this section is the following conditional nonlinear 
stability result. The notation $C^k$ stands for the usual space of
functions whose derivatives up to order $k$ are bounded.

\begin{Thm}[Conditional convergence for compressible Navier--Stokes]
\label{thmstabCNS}
Let $t\to f(t) =(\rho(t),u(t),T(t))$ be a $C^\infty$ solution
of~\eqref{CNS'}, satisfying the uniform bounds
\begeq\begin{cases}
\forall k\in\N \qquad \sup_{t\geq 0}
\Bigl( \|\rho(t)\|_{C^k} + \|u(t)\|_{C^k} + \|T(t)\|_{C^k}\Bigr) < +\infty; \\
\forall t\geq 0,\quad \rho(t) \geq \rho_m>0; \quad
T(t)\geq T_m >0.
\end{cases}
\endeq
Then $\|f(t) - (1,0,1)\|_{C^k}= O(t^{-\infty})$ for all $k$.
\end{Thm}

\begin{proof}[Proof of Theorem~\ref{thmstabCNS}]
Let us check that all assumptions of Theorem~\ref{thmmain}
are satisfied. Assumption~\ref{assfunct} is satisfied with,
say, $X^s = H^s(\T^N; \R\times\R^N\times\R)$, where $H^s$ stands for
the usual $L^2$-Sobolev space of functions with $s$ derivatives in
$L^2$. To fulfill Assumption~\ref{assworkspace}, define
$C_s:= \sup \{ \|f(t)\|_{H^s};\ t\geq 0\}$ and let
\[ X = Y := \Bigl\{ f\in C^\infty(\T^N;\R\times\R^N\times\R);\quad
\forall s,\quad \|f\|_s \leq C_s; \quad
\rho\geq \rho_m; \quad T\geq T_m \Bigr\}.\]
(Note that necessarily $\rho_m\leq 1$, $T_m\leq 1$.)

To check Assumption~\ref{asssol}, rewrite~\eqref{CNS'} in the nonconservative
form
\begeq\label{nonconserv}
\begin{cases}
\dps (\pa_t + u\cdot\nabla) \rho + \rho (\nabla\cdot u) = 0\\ \\
\dps (\pa_t + u\cdot\nabla) u + \nabla T + 
T\left(\frac{\nabla\rho}{\rho}\right)
= \frac{2\mu}{\rho}\, \nabla\cdot\{\nabla u\} \\ \\
\dps (\pa_t + u\cdot\nabla) T + \frac2N\, T(\nabla\cdot u) = 
\frac4N\, \frac{\mu}{\rho}\, |\{\nabla u\}|^2 + 
\frac{\kappa}{\rho}\, \Delta T,
\end{cases}\endeq
and define
\begeq\label{defBCNS}
B f = \Bigl( u\cdot\nabla \rho + \rho(\nabla\cdot u),\
u\cdot\nabla u + \nabla T + T \nabla (\log\rho), \
u\cdot\nabla T + \frac2N\, T (\nabla\cdot u) \Bigr); 
\endeq
\begeq\label{defCCNS}
\C f = \Bigl( 0, \ \frac{2\mu}{\rho}\, \nabla\cdot\{\nabla u\}, \
\frac4N\, \frac{\mu}{\rho}\, |\{\nabla u\}|^2 + \frac{\kappa}{\rho}
\Delta T\Bigr).
\endeq
Then~\eqref{asssol} obviously holds true. 

Assumption~\ref{asseq} is satisfied with $f_\infty=(1,0,1)$. 

As usual in the theory of
viscous compressible flows, an important difficulty to overcome is the
fact that diffusion does not act on the $\rho$ variable. So let
$\Pi_1=\Pi$ be defined by
\[ \Pi (\rho,u,T) = (\rho,0,1).\]
Assumption~\ref{assproj} is obviously satisfied with this choice of
nonlinear projection.

Next, let $\E$ be the negative of the usual entropy for perfect fluids:
\[ \E(\rho,u,T) = \int \rho\log\rho \ - \
\frac{N}2 \int \rho \log T.\]
Taking into account~\eqref{CLCNS},
\[ \E(\rho,u,T) - \E(1,0,1) = 
\int \rho \log \rho + \int \rho \frac{|u|^2}{2}
+ \frac{N}2 \int \rho (T-\log T-1);\]
\[ \E(\rho,u,T) - \E (\Pi(\rho,u,T)) = 
\int \rho \frac{|u|^2}2 + \frac{N}2 \int \rho(T-\log T-1).
\]
Thanks to the uniform bounds from above and below on
$\rho$ and $T$, $\E(f)-\E(f_\infty)$ controls
$\|f-f_\infty\|^2$ from above, and $\E(f)-\E(\Pi f)$ controls
$\|f-\Pi f\|^2$ from above and below; so
Assumption~\ref{asslyap} is satisfied.

It only remains to check Assumption~\ref{asskey}.
By a classical computation, for any $f\in Y$,
\begin{align*}
\E'(f)\cdot (\C f) & = - \left(
2\mu \int \frac{|\{\nabla u\}|^2}{T} \, + \, \kappa
\int \frac{|\nabla T|^2}{T^2} \right)\\
& \leq - K\, \left( \int |\{\nabla u\}|^2 + \int \rho |\nabla T|^2
\right),
\end{align*}
where the last inequality follows again from the lower bound on $T$
and the upper bound on $\rho$.

By Poincar\'e inequality, $\int \rho |\nabla T|^2$ controls
$\int \rho (T- \<T\>_\rho)^2$, where $\<T\>_\rho = \int \rho T$ is
the average of $T$ with respect to $\rho$. In turn, this controls
$\|T-1\|^2 - 2(\<T\>_\rho -1)^2$. Since $\<T\>_\rho - 1 =
(-1/N) \int \rho |u|^2$, we conclude that there are positive
constants $K$ and $C$ such that
\[ \int \rho |\nabla T|^2 \geq K \|T-1\|^2 - C \|u\|^2\]
for all $f\in Y$. On the other hand, by~\cite[Proposition~11]{DV:boltz:05},
\[ \int |\{\nabla u\}|^2 \geq K' \|u\|^2.\]
All in all, there is a positive constant $K$ such that
\[ \int |\{\nabla u\}|^2 + \int \rho |\nabla T|^2 \geq
K \bigl( \|T-1\|^2 + \|u\|^2\bigr) = K \|f-\Pi f\|^2,\]
so Assumption~\ref{asskey}(i) holds true.

By another classical computation, $\E'(f)\cdot (Bf)=0$,
so Assumption~\ref{asskey}(ii) also holds true.

On the range of $\Pi$, the functional derivative $\Pi'$ vanishes
(because $\pa_t\rho=0$ when $u=0$), and $B(\rho,u,T) = 
(0,-\nabla \log\rho, 0)$. Then
\[ (\Id - \Pi)'_{\Pi f} \cdot (B\Pi f) = B \Pi f = 
(0,-\nabla\log \rho, 0).\]
Thus
\[ \Bigl\| (\Id - \Pi)'_{\Pi f} \cdot (B\Pi f)\Bigr\|^2 =
\int |\nabla (\log \rho)|^2,\]
which under our assumptions controls $\int |\nabla \rho|^2$,
and then by Poincar\'e inequality also $\|\rho-1\|^2= \|\Pi f - f_\infty\|^2$.
This establishes Assumption~\ref{asskey}(iii), and then
the conclusion of the theorem follows from Theorem~\ref{thmmain}.
\end{proof}

\section{Weakly self-consistent Vlasov--Fokker--Planck equation}
\label{secVFP}

One of the final goals of the theory which I~have been trying to
start in this memoir is the convergence to equilibrium for the
nonlinear Vlasov--Poisson--Fokker--Planck equation with an external
confinement. This kinetic model, of great importance in plasma physics,
describes the evolution of a cloud of charged particles undergoing
deterministic and random (white noise) forcing, friction,
and influencing each other by means of Coulomb interaction.

Besides the fact that the regularity theory of the 
Vlasov--Poisson--Fokker--Planck equation is still at an early stage
(to say the least), one meets serious difficulties when trying to
apply Theorem~\ref{thmmain} to this model, in particular because the
problem is set in the whole space. So for the moment I~shall be content
to treat a simpler baby problem where (a) the confining potential
is replaced by a periodic boundary condition;
(b) the Coulomb interaction potential is replaced by a {\em small}
and {\em smooth} potential. The smallness assumption is not only a
technical simplification: It will prevent phase transition and guarantee
the uniqueness of equilibrium state.

Even with these simplifications, the problem of convergence to equilibrium
is nontrivial because the model is nonlinear and the diffusion only acts
on the velocity variable. This will be a perfect example of application
of Theorem~\ref{thmmain}.

Here the unknown $f=f(t,x,v)$ is a time-dependent probability
density in phase space ($x\in \T^N$ stands for position and $v\in\R^N$ for
velocity). The equation reads
\begeq\label{wVFP}
\begin{cases}
\dps \pa_t f + v\cdot\nabla_x f + F[f](t,x)\cdot\nabla_v f =
\Delta_v f + \nabla_v\cdot (fv) \\ \\
\dps F[f](t,x) = -\int \nabla W(x-y)\, f(t,y,w)\,dw\,dy.
\end{cases}
\endeq
Here $W\in C^\infty(\T^N)$ is even ($W(-z)=W(z)$), and without loss
of generality $\int W=0$. As we shall see later, if $W$ is small enough
in a suitable sense then the unique equilibrium for~\eqref{wVFP} is
the Maxwellian with constant density:
\[ f_\infty(x,v) = M(v) = \frac{e^{-\frac{|v|^2}{2}}}{(2\pi)^{N/2}}.\]

Since the total mass $\int f(t,x,v)\,dv\,dx$ is preserved with time,
there is an a priori estimate on the force, like
$\|F\|_{C^k} \leq C_k \|W\|_{C^k}$; so there is no real difficulty
in adapting the proofs of regularity for the {\em linear} kinetic
Fokker--Planck equation (see Appendix~\ref{appreg}). In this way one
can establish the existence and uniqueness of a solution as soon as,
say, the initial datum has finite moments of arbitrary order;
and this solution will be smooth for positive times.

The goal of this section is to establish the following
convergence result:

\begin{Thm}[Large-time behavior of the weakly self-consistent
Vlasov--Fokker--Planck equation]
\label{thmVFP}
Let $W\in C^\infty(\T^N)$ satisfy $\int W=0$.
Let $f_0=f_0(x,v)$ be a probability density on $\T^N\times\R^N$,
such that $\int f_0(x,v)|v|^k\,dv\,dx <+\infty$ for all $k\in\N$,
and let $f=f(t,x,v)$ be the unique smooth solution of~\eqref{wVFP}.
Let $\delta$ be so small that
\[ \delta + \frac{\delta^2\, e^\delta}2 < \frac12.\]
If $\max|W|<\delta$ then
\[ \|f(t,\cdot) - M \|_{L^1} = O(t^{-\infty}).\]
\end{Thm}

\begin{Rk} It is not hard to show that the conclusion of Theorem~\ref{thmVFP}
does not hold true without any size condition on $W$, since 
in general~\eqref{wVFP} can admit several stationary states.
In the proof of Theorem~\ref{thmVFP} I~shall show that there is
only one stationary state as soon as $\max|W|<1$; I~don't know
how good this bound is. The assumptions of the theorem are 
satisfied with $\delta=0.38$, which does leaves some margin of
improvement.
\end{Rk}

\begin{proof}[Proof of Theorem~\ref{thmVFP}]
The first step consists in establishing uniform regularity estimates;
I~shall only sketch them very briefly. 

First, one establishes
differential inequalities on the ``regularized'' moments
$M_k(t) = \int f(t,x,v) (1+|v|^2)^{k/2}\,dv\,dx$:
\[ \frac{dM_k}{dt} \leq C M_{k-1} - K M_k,\]
where $C$ and $K$ are positive constants.
(Here the fact that the position space is $\T^N$ induces a 
considerable simplifcation.) Then one deduces easily that
each moment $\int f|v|^k$ remains bounded uniformly in time.

Next, by adapting the arguments in Appendix~\ref{appreg}, one
can prove uniform Sobolev estimates of the form
\[ \forall k\in\N, \quad \forall t_0>0,\qquad
\sup_{t\geq t_0} \|f(t,\cdot)\|_{H^k} < +\infty.\]
These bounds, combined with the moment estimates, imply the boundedness
of the solution $f(t,\cdot)$ in all spaces $X^s$, where $X^k$ is
defined for $k\in\N$ by
\begeq\label{XkVFP} 
\|f\|_{X^k}^2 = \sum_{|\ell|+|m|\leq k}
\int \bigl| \nabla_x^k \nabla_v^m f(x,v)\bigr|^2\,
(1+|v|^2)^k\,dv\,dx.
\endeq
and $X^s$ is defined by interpolation for noninteger $s$.
It is easy to check that these spaces satisfy Assumption~\ref{assfunct}.

Finally, classical methods based on the maximum principle
(as in~\cite[Section~10]{DV:landau:1}) suffice to show that
\[ f(t,x,v) \geq K e^{-a|v|^2}\]
uniformly in $t\geq t_0$, provided that $K$ is small enough and $a$
is large enough. (Here again, the assumption that the position space is $\T^N$
simplifies things quite a bit by allowing $x$ to be treated as a parameter.)
It follows that $\rho(t,x)=\int f(t,x,v)\,dv$ is bounded below 
by a uniform positive constant for $t\geq t_0>0$.

Up to changing the origin of time, we can now assume that $f$ is uniformly
bounded in all spaces $X^s$ and that $\rho$ satisfies a uniform lower bound.
This determines a workspace
\[ X = Y := \Bigl\{ f;\quad \forall s\ 
\|f\|_{X^s}\leq C_s;\quad \rho\geq \rho_m>0\Bigr\},\]
as in Assumption~\ref{assworkspace}.

Then we define
\[ B f = v\cdot\nabla_x f + F[f]\cdot\nabla_v f;\qquad
\C f = \Delta_v f + \nabla_v\cdot (fv); \]
\[ f_\infty = M(v);\qquad  \Pi_1(f) = \Pi(f) = \rho M,\qquad
\rho = \int f\,dv.\]
Assumptions~\ref{asssol}, \ref{asseq}, \ref{assstat} and
\ref{assproj} are readily checked. 

Next let the free energy functional $\E$ be defined by
\[ \E(f) = \int f\log f \,dv\,dx + \int f \frac{|v|^2}{2}\,dv\,dx 
+ \frac12 \int \rho(x)\,\rho(y)\, W(x-y)\,dx\,dy.\]
By standard manipulations,
\begeq\label{EEPiVFP}
\E(f)-\E(\Pi f) = \int_{\T^N\times\R^N} f \log \frac{f}{\rho M};
\endeq
\begeq\label{EEPi'VFP}
\E(\Pi f) - \E(f_\infty) = \int_{\T^N}\rho\log\rho + 
\frac12 \int_{\T^N} \rho(x)\,\rho(y)\,W(x-y)\,dx\,dy.
\endeq

By the Csisz\'ar--Kullback--Pinsker inequality,
$\E(f)-\E(\Pi_1f) \geq (1/2)\|f-\rho M\|_{L^1}^2$; then by interpolation
of $L^2$ between $L^1$ and $H^k$ (as in~\cite[Lemma~10]{DV:boltz:05}),
one deduces
\[ \E(f)-\E(\Pi f) \geq \frac12 \|f-\rho M\|_{L^1}^2 \geq
K\, \|f-\rho M\|_{H^k}^{-\theta} \|f-\rho M\|_{L^2}^{2+\theta},\]
where $\theta$ is arbitrarily small if $k$ is chosen large enough.
This shows that Assumption~\ref{asslyap}(i) is satisfied.

On the other hand, since $\int W=0$,
\begin{align*} \E(\Pi f) - \E(f_\infty) &
= \int \rho \log\rho + \frac12 \int [\rho(x)-1]\,
[\rho(y)-1]\,W(x-y)\,dx\,dy \\
& \geq \frac12 \|\rho-1\|_{L^1}^2 - \frac12 (\max |W|)\,
\|\rho-1\|_{L^1}^2.
\end{align*}
By assumption, $\max |W|<1$; so there is a constant $K>0$ such that
\[ \E(\Pi f) - \E(f_\infty) \geq K \|\rho -1\|_{L^1}^2.\]
By interpolation again, this can be controlled from below by
$\|\rho-1\|_{L^2}^{2+\var}$ for arbitrarily small $\var$, 
and the left inequality in Assumption~\ref{asslyap}(ii) is satisfied.
(This is the first time that we use the smallness assumption on $W$.)
The right inequality in Assumption~\ref{asslyap}(ii) is easy.

By classical computations (see e.g.~\cite[Section~2]{DV:FP:01}),
\[ -\E'(f)\cdot (\C f) = \int f \left| \nabla_v \log\frac{f}{\rho M}\right|^2
\,dv\,dx \geq 2 \int f \log \frac{f}{\rho M} \geq
\frac12 \|f-\rho M\|_{L^1}^2,\]
so there is no difficulty to establish Assumption~\ref{asskey}(i).
Assumption~\ref{asskey}(ii) follows immediately since
$\E'(f)\cdot (Bf)=0$. So it only remains to establish 
Assumption~\ref{asskey}(iii).

As in the example of the compressible Navier--Stokes system,
the functional derivative $\Pi'$ vanishes on the range of $\Pi$, so
\begin{align*} (\Id-\Pi)'_{\Pi f}\cdot (B\Pi f) = 
B\Pi f & = v\cdot\nabla_x (\rho M ) 
- (\nabla W\ast \rho)\cdot\nabla_v(\rho M)\\
& = v\cdot\nabla_x \bigl( \rho + \rho (W\ast \rho)\bigr)M.
\end{align*}
Then
\begin{align*}
\Bigl\| (\Id-\Pi)'_{\Pi f}\cdot (B\Pi f)\Bigr\|^2
& = \int \Bigl| v\cdot\nabla_x \bigl(
\rho + \rho (W\ast\rho) \bigr) \Bigr|^2\, M(v)\,dv\,dx \\
& = \int \Bigl| \nabla_x \bigl( \rho + \rho (\rho\ast W)\bigr) 
\Bigr|^2\,dx \\
& \geq K \int \rho \left| \frac{\nabla\rho}{\rho} + \rho\ast W\right|^2
\,dx,
\end{align*}
where the lower bound on $\rho$ was used in the last inequality.
Let
\[ \mu(x) := \frac{e^{-(\rho\ast W)(x)}}{\int e^{-(\rho\ast W)}}.\]
Since $\mu$ is uniformly bounded from above and below, we can
use a logarithmic Sobolev inequality with reference measure
$\mu(x)\,dx$; so there is a positive constant $K$ such that
\begin{align}
\int \rho \left| \frac{\nabla \rho}{\rho} + \rho\ast W\right|^2\,dx 
& = \int \rho \left| \nabla \log\frac{\rho}{\mu}\right|^2\,dx \nonumber\\
& \geq K \int \rho \log\frac{\rho}{\mu}\,dx \nonumber\\
& = K \left(\int \rho\log\rho + \int \rho (\rho\ast W) -
\log \int e^{-W\ast\rho}\right). \label{pluginW}
\end{align}

By assumption, $\max|W|\leq\delta$; so $|W\ast\rho|\leq\delta$, and
\begin{multline*}
\Bigl| e^{-W\ast\rho} - \bigl(1- W\ast\rho\bigr)\Bigr| \leq
e^{\delta} \, \frac{(W\ast\rho)^2}{2} =
e^\delta\, \frac{[W\ast(\rho-1)]^2}{2} \\
\leq \frac{e^{\delta} (\max |W|)^2}{2}\, \|\rho-1\|_{L^1}^2 
\leq \frac{\delta^2 e^\delta}{2}\, \|\rho-1\|_{L^1}^2.
\end{multline*}
Since $\int (W\ast\rho)=0$, it follows by integration of this bound that
\[ \left| \int e^{-W\ast\rho}\, -1\right| \leq
\frac{\delta^2 e^\delta}2 \, \|\rho-1\|_{L^1}^2.\]
As a consequence,
\[ \log \left(\int e^{-W\ast\rho}\right) \leq
\frac{\delta^2 e^\delta}2 \, \|\rho-1\|_{L^1}^2.\]
From this bound and the inequality $|\int \rho (\rho\ast W)|\leq
\delta \|\rho-1\|_{L^1}^2$ again, we obtain
\begin{align*}
\int \rho\log\rho + \int \rho(\rho\ast W) - \log 
\int e^{-W\ast\rho} & \geq
\frac{\|\rho-1\|_{L^1}^2}{2} \, - \delta \|\rho-1\|_{L^1}^2
- \frac{\delta^2 e^\delta}{2} \|\rho-1\|_{L^1}^2 \\
& \geq \left( \frac12 - \delta - \frac{\delta^2 e^\delta}2\right)
\|\rho-1\|_{L^1}^2.
\end{align*}
By assumption the coefficient in front of $\|\rho-1\|_{L^1}^2$ is
positive, and then we can use interpolation again to get
\[\int \rho\log\rho \ + \int \rho(\rho\ast W) - \log 
\int e^{-W\ast\rho} \geq K_\var \|\rho-1\|_{L^2}^{2+\var}
= K_\var \|\rho M - M \|_{L^2}^{2+\var}.\]
So Assumption~\ref{asskey}(iii) holds. (Here again the
smallness condition was crucially used.) Then all the assumptions
of Theorem~\ref{thmmain} are satisfied, and the conclusion follows
at once.
\end{proof}

\section{Boltzmann equation} \label{secboltz}

This last section is devoted to the Boltzmann equation;
see~\cite{vill:handbook:02} and the references therein for background
and references on this model. I~have personally devoted a considerable
amount of research time on the problem of convergence to equilibrium
for the Boltzmann equation, alone or in collaborations with
Toscani and Desvillettes; a detailed account of this topic can be
found in my lecture notes~\cite{vill:ihp}.

As in Section~\ref{secVFP}
the unknown is a time-dependent probability density
$f=f(t,x,v)$ on the phase space. The variable $x$ will be assumed
to vary in a bounded $N$-dimensional domain $\Om_x$, that will be either
the torus $\T^N$, or a smooth bounded connected open subset of $\R^N$.
The equation reads
\begeq\label{BE}
\begin{cases}
\dps \derpar{f}{t} + v\cdot\nabla_x f = Q(f,f) \\ \\
\dps Q(f,f) = \int_{\R^N} \int_{S^{N-1}}
\Bigl[ f(x,v') f(x,v'_*) - f(x,v) f(x,v_*) \Bigr]\, B(v-v_*,\sigma)\,d\sigma
\, dv_*\\ \\
\dps v' = \frac{v+v_*}{2} + \frac{|v-v_*|}2\,\sigma;\qquad
v'_* = \frac{v+v_*}{2} - \frac{|v-v_*|}2\,\sigma.
\end{cases}
\endeq

Here $B$ is the collision kernel; for simplicity I~shall restrict to the
case $B=|v-v_*|$ (hard spheres interaction), but the analysis works
as soon as Assumptions (5) and (19) in~\cite{DV:boltz:05} are
satisfied, which covers all physically relevant cases that I~know of.

Three kinds of estimates play an important role in the modern theory
of the Boltzmann equation: Sobolev estimates (in $x$ and $v$ variables),
moment estimates and positivity estimates of the form
$f\geq K_0 e^{-A_0|v|^{q_0}}$. At least in some cases, 
the positivity estimates follow from regularity estimates~\cite{mouhot:lb:05},
but I~shall not address this issue here.

To continue the discussion it is necessary to take {\bf boundary conditions}
into account. I~shall consider five cases: (i) periodic boundary conditions;
(ii) bounce-back boundary conditions; (iii) specular reflection
in a nonaxisymmetric domain; (iv) specular reflection in a spherically
symmetric domain; (v) Maxwellian accommodation with constant
wall temperature. Cases (i) to (iii) were already considered
in~\cite{DV:boltz:05}, while cases (iv) and (v) are new and will be
the occasion of interesting developments. 
Specular reflection in a general axisymmetric domain 
(not spherically symmetric) is intermediate between cases (iii) and (iv)
and can probably be treated as a variant, but I~have not tried to do so.
Other conditions could be treated as a variant of (v), 
such as more general accommodation kernels, but they do not seem to cause
any substantial additional difficulty. On the other hand,
the techniques presented here are helpless to treat
accommodation with {\em variable} wall temperature, for which the
collision operator does not vanish; I~shall add a few words about this
issue in the end of the section.

\subsection{Periodic boundary conditions}

In this subsection I~shall consider the Boltzmann 
equation~\eqref{BE} in the position space $\Om_x=\T^N$ 
(the $N$-dimensional torus). Then there are $N+2$ conservation laws:
total mass, total momentum ($N$ components) and total kinetic energy.
Without loss of generality, I~shall assume
\begeq\label{conslawBE}
\int f\,dv\,dx =1;\qquad
\int f v\,dv\,dx = 0;\qquad
\int f |v|^2\,dv\,dx = N.
\endeq
Then the equilibrium state takes the form
\[ f_\infty(x,v) = M(v) = \frac{e^{-\frac{|v|^2}{2}}}{(2\pi)^{N/2}}.\]

Our goal is the next convergence theorem:

\begin{Thm}[Convergence for the Boltzmann equation with 
periodic boundary conditions]
\label{thmcvBEper}
Let $f$ be a solution of~\eqref{BE} in the spatial domain
$\Om_x=\T^N$, satisfying the conservation laws~\eqref{conslawBE},
and the uniform regularity estimates
\begeq\label{uniformregBE}
\begin{cases}
\dps \forall s\geq 0 \qquad
\sup_{t\geq 0}\ \|f(t,\cdot)\|_{H^s(\Om_x\times\R^N_v)} < +\infty; \\ \\
\dps \forall k\geq 0\qquad
\sup_{t\geq 0} \int f(t,x,v)\, |v|^k\,dv\,dx < +\infty;\\ \\
\dps \forall (t,x,v)\in\R_+\times\Om_x\times\R^N_v,\qquad
f(t,x,v) \geq K_0\, e^{-A_0|v|^{q_0}}.
\end{cases}
\endeq
Then
\[ \forall s\geq 0,\qquad 
\bigl\|f(t,\cdot) - M \bigr\|_{H^s} = O(t^{-\infty}).\]
\end{Thm}

\begin{Rk} This theorem is nonempty, in the sense that, when $f_0$ is
very smooth and close to equilibrium in a suitable sense, 
then~\eqref{uniformregBE} holds true. See the discussion in~\cite{DV:boltz:05}
for more information.
\end{Rk}

\begin{proof}[Proof of Theorem~\ref{thmcvBEper}]
Let $f$ satisfy the assumptions of Theorem~\ref{thmcvBEper}.
Let $(X^s)_{s\geq 0}$ be the scale of weighted Sobolev spaces already defined
in the treatment of the Vlasov--Fokker--Planck equation
(recall equation~\eqref{XkVFP}). It follows from the assumptions
that $C_s:=\sup_{t\geq 0} \|f\|_{s}$ is finite for all $s$.
We shall work in the spaces
\[ X := \Bigl\{f;\quad 
\|f\|_s \leq C_s; \ f(x,v) \geq K_0\, e^{-A_0 |v|^{q_0}}\Bigr\};\qquad\]
\[ Y := \Bigl\{f;\quad 
\|f\|_s \leq C'_s; \ f(x,v) \geq K'_0\, e^{-A_0 |v|^{q_0}}\Bigr\},\]
where $C'_s$, $K'_0$ will be determined later on. Then
Assumptions~\ref{assfunct}, \ref{assworkspace} and~\ref{asssol}
are obviously satisfied.

Define 
\[ Bf=v\cdot\nabla_x f;\qquad \C f = Q(f,f).\]
Then Assumption~\ref{asseq}(i) is obviously true,
Assumption~\ref{asseq}(ii) is satisfied since $B$ is linear
continuous $X^{s+1}\to X^s$, and
Assumption~\ref{asseq}(iii) is a consequence of~\cite[eq.~(78)]{DV:boltz:05}.

Assumption~\ref{assstat} holds true with $f_\infty=M$; notice
that both the transport and the collision part vanish on $f_\infty$.

Next, if $f=f(x,v)$ is given, define
\[ \rho = \int f\,dv; \qquad 
u = \frac1{\rho} \int fv\,dv;\qquad
T = \frac1{N\rho}\int f |v-u|^2\,dv,\]
and
\[ M_{\rho\, u\, T} = \frac{\rho(x)\, e^{-\frac{|v-u(x)|^2}{2T(x)}}}
{\bigl[2\pi T(x)\bigr]^{N/2}}.\]
(Note that $M$ depends on $f$ via $\rho,u,T$.) It is easy to derive
uniform estimates of smoothness on $\rho$, $u$ and $T$
in terms of the estimates on $f$; and to derive similarly
strict positivity estimates on $\rho$, $T$:
See~\cite[Proposition~7]{DV:boltz:05}.

Now we can introduce the nonlinear projection operators:
\[ \Pi_1 f = M_{\rho\, u \, T};\qquad
\Pi_2 f = M_{\rho\, 0\, 1};\qquad 
\Pi_3 f = M.\]
By adjusting the constants $C'_s$, $K'_0$, we can ensure that
$\Pi_j(X)\subset Y$. Then the rest of Assumption~\ref{assproj}
follows easily.

The natural Lyapunov functional in the present case is
of course Boltzmann's $H$ functional:
\[ \E(f) = H(f) = \int_{\T^N_x\times\R^N_v} f\log f\,dv\,dx.\]
By standard computations, taking into account~\ref{conslawBE}, we have
\begeq\label{EEPi1BE}
\E(f) - \E(\Pi_1f) = 
\int f \log \frac{f}{M_{\rho u T}};
\endeq
\begeq\label{EEPi'BE}
\E(\Pi_1f) - \E(f_\infty) =
\int \rho\log\rho + \int \rho \frac{|u|^2}{2} + 
\int \rho(T-\log T-1).
\endeq
To find a lower bound on~\eqref{EEPi1BE}, it suffices to
use the Csisz\'ar--Kullback--Pinsker inequality and interpolation,
as we did previously for the Vlasov--Fokker--Planck equation
(recall~\eqref{EEPiVFP}; or~\cite[eq.~(47)${}_1$]{DV:boltz:05}). 
Upper and lower bounds for~\eqref{EEPi'BE} can be obtained as we did 
before for the compressible Navier--Stokes equations. So
Assumption~\ref{asslyap} is satisfied.

Now the crucial step consists in checking Assumption~\ref{asskey}.
By a classical computation,
\[ - H'(f)\cdot (\C f) = \int D\bigl(f(x,\cdot)\bigr)\,dx,\]
where $D(f)$ is Boltzmann's dissipation of information:
\begin{multline*} D(f) = 
\frac14 \int \Bigl( f(v') f(v'_*) - f(v) f(v_*)\Bigr)\,
\Bigl(\log f(v')f(v'_*) - \log f(v)f(v_*)\Bigr)\\ \,B(v-v_*,\sigma)\,
d\sigma\,dv\,dv_*.
\end{multline*}
Known entropy production estimates from~\cite{vill:cer:03}
make it possible to estimate $D(f)$ from below by
$K_\var \bigl[ H(f) - H(M_{\rho\, u\, T})\bigr]^{1+\var}$.
(Such estimates go back to~\cite{TV:entropy:99}; see
also~\cite{vill:ihp} for a detailed account on this problem.)
Then Assumption~\ref{asskey}(i) follows easily,
as in~\cite[Corollary~5]{DV:boltz:05}.

Assumption~\ref{asskey}(ii) is an immediate consequence
of Assumption~\ref{asskey}(i),
since $H'(f)\cdot(Bf) =0$.

It remains to establish Assumption~\ref{asskey}(iii).
For this we use Remark~\ref{rkpractic}. 
According to~\cite[eq.~(69)]{DV:boltz:05},
if $f=M^f_{\rho\, u\, T}$ at time $t=0$ and $f$ evolves according
to $\pa_t f + v\cdot\nabla_x f=0$, then
\[ \left.\frac{d^2}{dt^2} \right|_{t=0} \|f-M^f_{\rho u T}\|
\geq K \left(\int_{\T^N} |\nabla T|^2 + 
\int_{\T^N} |\{\nabla u\}|^2\right),\]
where, as in Section~\ref{secCNS},
\[ \{\nabla u\}_{ij} = \frac12 \left( \frac{\pa u_i}{\pa x_j} 
+ \frac{\pa u_j}{\pa x_i}\right) 
- \frac{1}{N} (\nabla\cdot u)\,\delta_{ij},\]
and $\nabla \cdot u$ is the divergence of $u$.
According to~\cite[Section~IV.2]{DV:boltz:05}, there are
constants $K_1$, $K_2$, $K_3$ only depending on $N$, such that
\[ \begin{cases}
\dps \int_{\T^N} |\nabla T|^2 \geq K_1 \|T-1\|^2 - C \|u\|^2;\\ \\
\dps \int_{\T^N} |\{\nabla u\}|^2 \geq K_2 \|\nabla u\|^2
\geq K_3 \|u\|^2.
\end{cases} \]
This implies
\begin{align*}
\bigl\| (\Id-\Pi_1)'_{\Pi_1f}\cdot (B\Pi_1f)\bigr\|_{L^2}^2 &
\geq K (\|T-1\|^2 + \|u\|^2) \\
& \geq K' \, \|\Pi_1f-\Pi_2f\|^2.
\end{align*}

Next, according to~\cite[eq.~(71)]{DV:boltz:05},
if $f=M_{\rho\, 0\, 1}$ at time $t=0$ and $f$ evolves according
to $\pa_t f + v\cdot\nabla_x f=0$, then
\[ \left.\frac{d^2}{dt^2} \right|_{t=0} \|f-M^f_{\rho 0 1}\|
\geq K \int_{\T^N} |\nabla \rho|^2 .\]
Combining this with a Poincar\'e inequality
(see again~\cite[Section~IV.2]{DV:boltz:05}), we deduce that
\[\bigl\| (\Id-\Pi_2)'_{\Pi_2f}\cdot (B\Pi_2f)\bigr\|_{L^2}^2
\geq K \|\rho-1\|^2 \geq K' \|\Pi_2 f - f_\infty\|_{L^2}^2.\]
(This is in fact as in Section~\ref{secVFP}, if we set $W=0$.)

This concludes the verification of Assumption~\ref{asskey},
and the result follows by an application of Theorem~\ref{thmmain}.
\end{proof}

\begin{Rk} A comparison with the proof of the same result
in~\cite{DV:boltz:05} shows that the crucial functional inequalities
are all the same; but there are essential simplifications
in that (a) it suffices to do the computations for Maxwellian
states (``local equilibrium'' in the language of~\cite{DV:boltz:05});
and especially (b) there is no longer need for the tricky analysis of
the system of differential inequalities. More explicitly,
Sections III.3, V and VI of~\cite{DV:boltz:05} are shortcut by
the use of Theorem~\ref{thmmain}.
\end{Rk}

\subsection{Bounce-back condition} \label{subbb}

Now let $\Om_x$ be a bounded smooth open subset of $\R^N$; up to rescaling
units we may assume that $|\Om_x|=1$ (the Lebesgue measure of the domain
is normalized). In this subsection the boundary condition is of
{\em bounce-back} type:
\begeq\label{bb}
x\in \pa\Om_x \Longrightarrow\qquad
f(x,v) = f(x,-v).
\endeq
A consequence of~\eqref{bb} is that $u=0$ on $\pa\Om_x$ (the mean
velocity vanishes on the boundary).

Now there are only 2 conservation laws: mass and energy.
So, without loss of generality, I~shall assume
\begeq\label{conslawbb} \int f\,dv\,dx =1;\qquad
\int f |v|^2\,dv\,dx = N.
\endeq
The equilibrium is again the steady Maxwellian,
\[ f_\infty(x,v) = M(v) = \frac{e^{-\frac{|v|^2}{2}}}{(2\pi)^{N/2}}.\]

Here is the analogue of Theorem~\ref{thmcvBEper}:

\begin{Thm}[Convergence for the Boltzmann equation with 
bounce-back boundary conditions]
\label{thmcvBEbb}
Let $f$ be a solution of~\eqref{BE} in a smooth bounded connected
spatial domain $\Om_x$, satisfying bounce-back boundary conditions,
the conservation laws~\eqref{conslawbb},
and the uniform regularity estimates~\eqref{uniformregBE}.
Then
\[ \forall s\geq 0,\qquad 
\bigl\|f(t,\cdot) - M \bigr\|_{H^s} = O(t^{-\infty}).\]
\end{Thm}

\begin{proof}[Proof of Theorem~\ref{thmcvBEbb}]
The proof is quite similar to the proof of Theorem~\ref{thmcvBEper},
however the sequence of projection operators is different:
\[ \Pi_1 f = M_{\rho\, u\, T};\qquad
\Pi_2 f = M_{\rho\, u\, \<T\>};\qquad
\Pi_3 f = M_{\rho\, 0\, 1};\qquad
\Pi_4 f = M,\]
where $\<T\>=\int \rho T$ is the average temperature.
According to~\cite[eq.~(70)-(71)]{DV:boltz:05} and a reasoning
similar to the one in the proof of Theorem~\ref{thmcvBEper},
\begeq\begin{cases}
\Bigl\| \bigl(\Id-\Pi_1\bigr)'_{\Pi_1f}\cdot (B\Pi_1 f)\Bigr\|^2\geq
K\|\nabla T\|^2;\\ \\
\Bigl\| \bigl(\Id-\Pi_2\bigr)'_{\Pi_2f}\cdot (B\Pi_2 f)\Bigr\|^2\geq
K\|\nabla^{\rm sym}u\|^2;\\ \\
\Bigl\| \bigl(\Id-\Pi_1\bigr)'_{\Pi_1f}\cdot (B\Pi_1 f)\Bigr\|^2\geq
K\|\nabla \rho\|^2,
\end{cases}\endeq
where $\nabla^{\rm sym} u$ is the symmetrized gradient of $u$,
that is
\[ \bigl(\nabla^{\rm sym}u\bigr)_{ij} = 
\frac12 \left( \derpar{u_i}{x_j} + \derpar{u_j}{x_i}\right).\]
By Poincar\'e inequalities,
\[ \|\nabla T\|^2 \geq K \|T-\<T\>\|^2;\qquad
\|\nabla \rho\|^2 \geq K \|\rho-1\|^2.\]
By the classical Korn inequality, and the Poincar\'e inequality
again (component-wise),
\[ \|\nabla^{\sym} u\|^2 \geq K \|\nabla u\|^2 \geq K' \|u\|^2.\]
These estimates imply $\|(\Id-\Pi_j)'_{\Pi_jf}\cdot (B\Pi_jf)\|^2
\geq K\|\Pi_jf -\Pi_{j+1}f\|^2$ for all $j\in\{1,2,3\}$, so
Assumption~\ref{asskey}(iii) is satisfied in the end.
Then Theorem~\ref{thmmain} applies.
\end{proof}

\subsection{Specular reflection in a nonaxisymmetric domain}

In this subsection the bounce-back boundary condition is replaced
by the {\em specular reflection} condition:
\[ x\in \pa\Om_x \Longrightarrow\qquad
f(x,v) = f(x,R_xv),\qquad R_xv=v-2\<v,n\> n.\]
This condition is more degenerate and the shape of the domain
will influence the form of the equilibrium. For the moment I~shall
assume that the domain is {\em nonaxisymmetric} in dimension $N=3$.
The notation is the same as in Subsection~\ref{subbb}.

\begin{Thm}[Convergence for the Boltzmann equation with 
nonaxisymmetric specular conditions]
\label{thmcvBEspec}
Let $f$ be a solution of~\eqref{BE} in a smooth bounded connected
nonaxisymmetric spatial domain $\Om_x\subset\R^3$, 
satisfying specular boundary condition, the conservation 
laws~\eqref{conslawbb}, and the uniform regularity 
estimates~\eqref{uniformregBE}. Then
\[ \forall s\geq 0,\qquad 
\bigl\|f(t,\cdot) - M \bigr\|_{H^s} = O(t^{-\infty}).\]
\end{Thm}

\begin{proof}[Proof of Theorem~\ref{thmcvBEspec}]
The proof is entirely similar to the proof of Theorem~\ref{thmcvBEbb},
except that the condition $u=0$ on the boundary is replaced by
the weaker condition $u\cdot n=0$, where $n$ is the inner unit normal
to $\Om_x$. Then the classical Korn inequality should be replaced
by the Korn inequality established by Desvillettes and 
myself in~\cite{DV:korn:02}.
\end{proof}

\subsection{Specular reflection in a spherically symmetric domain}

In this subsection $\Om_x$ is a bounded smooth connected spherically
symmetric domain in $\R^3$; so, up to translation, 
$\Om_x$ is either a ball $(|x|<R)$ or a shell $(0<r<|x|<R)$.
Again I~shall assume that $|\Om_x|=1$. I~shall write $N=3$ to
keep track of the role of the dimension in various formulas
(certainly the analysis can be extended to more general domains,
but one has to be careful about the meaning of the conservation of
angular momentum).

Now there are $N+2$ conservation laws: mass, kinetic energy and
{\em angular momentum} ($N$ scalar quantities). Without loss of generality,
I~shall assume
\begin{multline}\label{conslawBEsph} \int f\,dv\,dx=1;\qquad
\int f(x,v)|v|^2\,dv\,dx = N;\qquad \\
\int f(x,v) (v\wedge x)\,dv\,dx = {\bf M} \in \R^N.
\end{multline}
The existence of an equilibrium is not trivial if ${\bf M}\neq 0$,
and the equilibrium does not seem to be explicit. It is a local
Maxwellian with uniform temperature $\theta$, but nonzero
velocity $u_\infty$ and nonhomogeneous density $\rho_\infty$.
The equations determining this equilibrium were studied,
at the beginning of the nineties, by Desvillettes~\cite{desv:eq:90}.
Here I~shall suggest a variational approach to this problem,
by means of the following lemma from elementary calculus of
variations (the proof of which will be only sketched):

\begin{Lem}[stationary solutions in a spherically symmetric domain]
\label{lemstatsph}
Let $\Om_x$ be a spherically symmetric domain in $\R^N$, $N=3$, $|\Om_x|=1$.
Whenever $\rho$ is a nonnegative integrable density on $\Om_x$,
and $m\in L^1(\Om_x;\R^N)$, define
\[ F(\rho,m) = \int \rho\log\rho - \frac{N}2
\log \left(1 - \frac1N \int \frac{|m|^2}{\rho}\right).\]
Then there is a unique $(\rho_\infty,m_\infty)\in C^\infty(\Om_x;
\R_+\times\R^N)$ which minimizes the functional $F$ under the constraints
\begeq\label{constrF}
\int \rho=1;\qquad \int m(x)\wedge x \,dx = {\bf M}.
\endeq
Moreover, $\rho$ is strictly positive; and 
there are an antisymmetric matrix $\Sigma_\infty$
and positive constants $\theta_\infty$ and $Z$ such that for all $x\in\Om_x$,
\[ \frac{m(x)}{\rho(x)} = \Sigma_\infty x;\qquad
\rho_\infty(x) = \frac{e^{\frac{|\Sigma_\infty x|^2}{2\theta_\infty}}}{Z}.\]
\end{Lem}

\begin{proof}[Sketch of proof of Lemma~\ref{lemstatsph}]
Write
\[ \theta = 1 - \frac1N \int \frac{|m|^2}{\rho}, \]
then
\begin{align*}
F(\rho,m) & = \int \rho\log\rho + \frac{N}2 (1-\theta)
+ \frac{N}2 \bigl(\theta-\log\theta-1\bigr)\\
& = \int \rho\log\rho + \int \frac{|m|^2}{2\rho} + \Psi(\theta),
\end{align*}
where $\Psi(\theta) = (N/2)(\theta-\log\theta-1)$.

By a classical computation, $(\rho,m)\longmapsto |m|^2/\rho$ is
convex, so $\theta$ is a concave function of $(\rho,m)$. Moreover,
$\theta$ remains in $(0,1)$, and on that interval $\Psi$ is a convex
{\em decreasing} function of $\theta$. It follows that $\Psi(\theta)$
is a strictly convex function of $(\rho,m)$. So 
\[ F: (\rho,m)\longmapsto \int \rho\log\rho + \int \frac{|m|^2}{2\rho} 
+ \Psi(\theta)\]
is a strictly convex function of $(\rho,m)$. This conclusion does not
change if $F$ is restricted on the domain defined by the {\em linear}
constraints~\eqref{constrF}; so $F$ has at most one minimizer.

The Euler--Lagrange equations for the minimization of $F$ read
\begeq\label{ELsph}
\begin{cases}
\dps\log\rho - \frac{|m|^2}{2\theta\rho^2} = \lambda_0;\\ \\
\dps \frac{N}{\theta} \left(\frac{m_i}{\rho}\right) = 
\var_{ijk}\,\lambda_j x_k,
\end{cases}\endeq
where $(\lambda_j)_{0\leq j\leq N}$ are constants,
$(m_i)_{1\leq i\leq N}$ are the components of $m$,
and $\var_{ijk}$ is defined by the equations
$(a\wedge b)_i = \sum\, \var_{ijk} a_j b_k$.
These equations imply that $m/\rho$ is an antisymmetric linear 
function of $x$. In particular, the minimizer a priori lives in
a finite-dimensional space. The rest of the lemma follows by classical
arguments.
\end{proof}

The goal of the present subsection is the following result:

\begin{Thm}[Convergence for the Boltzmann equation with 
spherically symmetric specular conditions]
\label{thmcvBEsph}
Let $f$ be a solution of~\eqref{BE} in a smooth bounded connected
spherically symmetric spatial domain $\Om_x\subset\R^3$, 
satisfying specular boundary condition, the conservation
laws~\eqref{conslawBEsph} and the uniform regularity 
estimates~\eqref{uniformregBE}. Then
\[ \forall s\geq 0,\qquad 
\bigl\|f(t,\cdot) - f_\infty \bigr\|_{H^s} = O(t^{-\infty}),\]
where 
\[ f_\infty(x,v) = \frac{e^{\frac{|\Sigma_\infty x|^2}{2\theta_\infty}}}{Z}\: 
\frac{e^{-\frac{|v-\Sigma_\infty x|^2}{2\theta_\infty}}}
{(2\pi\theta_\infty)^{3/2}},\]
and the antisymmetric matrix $\Sigma_\infty$, the positive constants
$Z$ and $\theta_\infty$ are provided by Lemma~\ref{lemstatsph}.
\end{Thm}

\begin{Rk} The variable $\theta_\infty$ is the (uniform) equilibrium 
temperature; the velocity field in the stationary state is still rotating,
and the density is lower near the interior of the box.
\end{Rk}

\begin{proof}[Proof of Theorem~\ref{thmcvBEsph}]
The only differences with the previously treated cases lie
in the definition of the projection operators, and the verification
of Assumptions~\ref{asslyap}(ii) and~\ref{asskey}(iii). 

In the present case, let
\[ \Sigma:= \< \nabla^\a u\>, \qquad \theta:=\<T\>_\rho;\]
more explicitly, $\Sigma$ is the average value of the antisymmetric part of 
the matrix-valued field $\nabla u$ (the averaging measure is the 
normalized Lebesgue measure), while $\theta$ is the average value
of the temperature (but now the averaging measure has density $\rho$).
I~shall identify the matrix $\Sigma$ with the velocity field
$x\longmapsto \Sigma x$, and $\theta$ with the constant function $x\to\theta$.
Then the sequence of projection operators is as follows:
\begin{multline*} \Pi_1 f = M_{\rho\, u \, T};\qquad
\Pi_2 f = M_{\rho\, u\, \theta};\qquad
\Pi_3 f = M_{\rho\, \Sigma\, \theta};\\
\Pi_4 f = f_\infty = M_{\rho_{\infty}\, \Sigma_{\infty}\, \theta_\infty}.
\end{multline*}

Once again the Lyapunov functional is
\[ H(f) = \int f\log f.\]
After taking into account the conservation laws~\eqref{conslawBEsph},
one observes that
\begin{multline}\label{HHsph} H(\Pi_1f) - H(f_\infty) =
\left( \int \rho\log\rho + \frac{N}2 
\int \rho (T-\log T-1) + \int \rho \frac{|u|^2}{2} \right) \\
- \left(\int \rho_\infty\log\rho_\infty+
\frac{N}2 \int \rho_\infty(\theta -\log\theta-1)
+ \int \rho_\infty \frac{|u_\infty|^2}{2}\right).
\end{multline}
Let again $\Phi(\theta) = \theta-\log\theta-1$: then by
Jensen's inequality (in quantitative form),
\[ \int \rho \Phi(T) \geq \Phi (\<T\>_\rho) + K \|T-\<T\>_\rho\|^2,\]
where $K$ depends on the bounds on $\rho$ and $T$.
Plugging this in~\eqref{HHsph} and using the same notation as
in Lemma~\ref{lemstatsph}, one obtains the lower bound
\begin{align} H(\Pi_1f) - H(f_\infty) &
= \frac{N}2 \left( \int \rho \Phi(T) - \Phi (\<T\>_\rho)\right)
+ \bigl[ F(\rho,m) - F(\rho_\infty,m_\infty)\bigr] \nonumber\\ 
& \geq K \|T-\<T\>_\rho\|^2 + \bigl[ F(\rho,m) - F(\rho_\infty,m_\infty)\bigr].
\end{align}

The upper bound
\[ H(\Pi_1f) - H(f_\infty)  \leq 
C \|T-\<T\>_\rho\|^2 + \bigl[ F(\rho,m) - F(\rho_\infty,m_\infty)\bigr]\]
is obtained in a similar way. 

So to prove Assumption~\ref{asslyap}(ii), it suffices to check that
\[
K \|\Pi_1 f - f_\infty\|^2 \leq F(\rho,m) - F(\rho_\infty,m_\infty)\leq
C \|\Pi_1 f - f_\infty\|^2;
\]
or, which amounts to the same,
\begin{multline}\label{unifbdsF}
K \Bigl\|(\rho,m)-(\rho_\infty,m_\infty)\Bigr\|^2 
\leq F(\rho,m) - F(\rho_\infty,m_\infty)\\ \leq
C \Bigl\|(\rho,m)-(\rho_\infty,m_\infty)\Bigr\|^2.
\end{multline}
The upper bound is obvious from the definition, the bounds
on $(\rho,m)$ (which follow from the bounds on $f$) and the bounds
on $(\rho_\infty,m_\infty)$. To prove the lower bound, it suffices to
establish the uniform convexity of $F$. Let
\[ f(\rho,m) = \rho\log\rho + \frac{|m|^2}{2\rho}.\]
The Hessian of $f$ has matrix
\[ \left(\begin{array}{cc}\dps \frac{I_N}{\rho} & \dps
-\frac{m}{\rho^2} \\ \\
\dps -\frac{m}{\rho^2} & \dps \frac{|m|^2}{\rho^3} + \frac1{\rho} 
\end{array}\right),\]
where $I_N$ stands for the $N\times N$ identity matrix;
under our assumptions on $\rho$, this Hessian matrix is uniformly
positive, so $f$ is uniformly convex, and the same is true of the functional
$F:(\rho,m)\longmapsto \int f(\rho,m)\,dx + \Psi(\theta)$.
This conclusion does not change when one imposes the linear constraints
constraints~\eqref{constrF}, and the lower bound in~\eqref{unifbdsF} follows.

The last crucial step in the proof consists in the verification of
Assumption~\ref{asskey}(iii). As in the previous subsection,
\begin{align}\label{gainPi1}
\Bigl\| (\Id-\Pi_1)'_{\Pi_1f}\cdot (B\Pi_1f)\Bigr\|^2 
& \geq K \int |\nabla T|^2 + \int |\{\nabla u\}|^2\\
& \geq K' \|T-\<T\>_\rho\|^2\nonumber,
\end{align}
which controls $\|\Pi_1 f- \Pi_2f\|^2$.

Next,
\begin{equation}
\Bigl\| (\Id-\Pi_2)'_{\Pi_2f}\cdot (B\Pi_2f)\Bigr\|^2 
\geq K \int |\nabla^{\rm sym}u|^2
\geq K' \|\nabla u - \Sigma\|^2\nonumber,
\end{equation}
where the second inequality follows from a version of
Korn's inequality~\cite[eq.~(1)]{DV:korn:02}.
Note that $\Sigma=\nabla(\Sigma x)$ (to avoid confusions I~shall now
write $\Sigma x$ for the map $x\to \Sigma x$), 
so one can apply again a Poincar\'e inequality
to obtain in the end
\begeq\label{gainPi2} 
\Bigl\| (\Id-\Pi_2)'_{\Pi_2f}\cdot (B\Pi_2f)\Bigr\|^2
\geq K \|u -\Sigma x\|^2,
\endeq
which controls $\|\Pi_2 f - \Pi_3 f\|^2$.

The gain from $\Pi_3$ is the main novelty.
As a consequence of~\cite[eq.~(65)]{DV:boltz:05},
\begin{align} \label{newPigain}
\frac{1}{M_{\rho\,\Sigma\,\theta}} \Bigl(
\pa_t M_{\rho\,\Sigma\,\theta} + v\cdot\nabla_x 
M_{\rho\,\Sigma\,\theta}\Bigr) = &
\left( \frac{\pa_t\rho + \Sigma\cdot\nabla\rho}{\rho} -
\frac{N}{2}\frac{\pa_t \theta}{\theta}\right) \\
& + \left(\frac{v-u}{\sqrt{\theta}}\right)\cdot
\left( \sqrt{\theta}\, \frac{\nabla\rho}{\rho} + 
\frac{\pa_t\Sigma + (\Sigma\cdot\nabla)\Sigma}{\sqrt{\theta}}
\right) \nonumber\\
& + \sum_i \left(\frac{v_i-u_i}{\sqrt{\theta}}\right)^2\, 
\frac{\pa_t\theta}{\theta}.\nonumber
\end{align}
The first and third lines do not bring any new estimate, so we focus on the
second line. First note that
\[ \bigl((\Sigma\cdot\nabla)\Sigma\bigr)_i(x) = 
\sum_{jk\ell}
(\Sigma_{jk}x_k)\, \pa_j (\Sigma_{i\ell}\,x_\ell) = (\Sigma^2 x)_i.\]
(Do not mistake the symbol of summation with the matrix $\Sigma$.)
Next, the equation for the mean velocity field $u$ is
$\pa_t u + u\cdot \nabla u + \nabla T + T \nabla (\log\rho)
+ (\nabla \cdot D)/\rho$, where $D$ vanishes on the range of $\Pi_1$.
Taking the antisymmetric part of this equation results in
$\pa_t \nabla^\a u =0$, hence $\pa_t\Sigma=0$. The conclusion is that
$\pa_t\Sigma$ vanishes on the range of $\Pi_3$. From all this information,
we deduce that the second line of~\eqref{newPigain} can be
simplified into
\[ (v-u)\cdot \left(\frac{\nabla\rho}{\rho} + \frac{\Sigma^2 x}{\theta}
\right),\]
where again $\Sigma^2 x$ is a shorthand for
the map $x\to \Sigma^2 x$. It follows that
\begeq\label{new3gainPi}
\Bigl\| (\Id-\Pi_3)'_{\Pi_3f}\cdot (B\Pi_3f)\Bigr\|^2
\geq K \Bigl\| \frac{\nabla\rho}{\rho} + 
\frac{\Sigma^2 x}{\theta}\Bigr\|^2.
\endeq

I~shall now show that
\begin{multline}\label{nowshowthat}
\Bigl\|u - \Sigma x\Bigr\|^2 + 
\Bigl\| \frac{\nabla\rho}{\rho} + 
\frac{\Sigma^2 x}{\theta}\Bigr\|^2 \geq
K\: \bigl( \|\rho-\rho_\infty\|^{2}
+ |\theta-\theta_\infty|^{2} \\ + 
\|\Sigma x-\Sigma_\infty x\|^{2}\bigr).
\end{multline}
Since the right-hand side controls $\|\Pi_3 f - \Pi_4f\|^{2+\var}$,
in view of~\eqref{new3gainPi} and~\eqref{gainPi2} this will imply
\begeq\label{IdPi2Pi3}
\Bigl\| (\Id-\Pi_2)'_{\Pi_2f}\cdot (B\Pi_2f)\Bigr\|^2
+ \Bigl\| (\Id-\Pi_3)'_{\Pi_3f}\cdot (B\Pi_3f)\Bigr\|^2\geq 
K \|\Pi_3 f - \Pi_4f\|^{2},
\endeq
completing the verification of Assumption~\ref{asskey}(iii).

Since $|\theta-\theta_\infty|\leq C(\|\rho-\rho_\infty\|
+ \|u-u_\infty\|)$ and $\|\Sigma x - \Sigma_\infty x\|
\leq C \|\Sigma-\Sigma_\infty\| \leq C' \|u-u_\infty\|$,
to establish~\eqref{nowshowthat} it is sufficient to prove
\begeq\label{nowshowthat'}
\Bigl\|u - \Sigma x\Bigr\|^2 + 
\Bigl\| \frac{\nabla\rho}{\rho} + 
\frac{\Sigma^2 x}{\theta}\Bigr\|^2 \geq
K \bigl( \|\rho - \rho_\infty\|^2 + \|u-u_\infty\|^2\bigr).
\endeq
In view of the bounds on $\rho$ and $u$, and the uniform convexity
of $F$ (used above to check Assumption~\ref{asslyap}(ii)),
inequality~\eqref{nowshowthat'} will be a consequence of
\begeq\label{nowshowthat''}
\Bigl\|\nabla u - \Sigma\Bigr\|^2 + 
\Bigl\| \frac{\nabla\rho}{\rho} + 
\frac{\Sigma^2 x}{\theta}\Bigr\|^2 \geq
K \bigl[ F(\rho,m) - F(\rho_\infty,m_\infty)\bigr].
\endeq

The following lemma will be useful:

\begin{Lem} \label{lemPhiK}
Let $\Phi$ be a $K$-uniformly convex function, defined and differentiable
on a convex open subset of a Hilbert space $\H$, and let
$\Lambda:\H\to\R^d$ be a linear map. If $X_\infty$ minimizes
$\Phi$ under the constraints $\Lambda(X)=c$, then 
\[ \Phi(X) - \Phi(X_\infty) \leq (2K)^{-1}
\inf_{\lambda\in (\ker \Lambda)^\bot}\, 
\Bigl\| \grad \Phi(X) + \lambda \Bigr\|^2.\]
\end{Lem}

Postponing the proof of Lemma~\ref{lemPhiK} for the moment,
let us apply it to the uniformly convex functional
\[ \Phi(\rho,m) = \int \rho\log \rho - \frac{N}2
\log \left( 1- \frac{1}{N} \int \frac{|m|^2}{\rho}\right)\]
and the linear map
\[ \Lambda(\rho,m) = \Bigl( \int \rho, \int m\wedge x\Bigr)\in
\R\times\R^N.\]
Then
\[ \grad\Phi = \Bigl( \log\rho - \frac{|m|^2}{2\theta\rho^2},\ 
\frac{m}{\theta\rho}\Bigr),\]
and $(\ker \Lambda)^\bot$ is made of vectors
$\lambda=(\lambda_0,A)$, where $\lambda_0\in\R$ and
$A$ is a (constant!) antisymmetric matrix.
So Lemma~\ref{lemPhiK} implies
\begeq\label{1F}
F(\rho,m) - F(\rho_\infty,m_\infty)
\leq C \inf_{\lambda_0\in\R;\ A^*=-A}\ 
\Bigl( \Bigl\| \log\rho - \frac{|m|^2}{2\theta\rho^2}\Bigr\|^2 
+ \Bigl\| \frac{m}{\theta\rho} - A x\Bigr\|^2 \Bigr),
\endeq

By Poincar\'e inequality, the bounds on $\theta$, and
Korn inequality,
\begeq\label{byPoincKorn} 
\inf_{A^*=-A} \Bigl\| \frac{m}{\theta\rho} - A x\Bigr\|^2 \leq 
C \inf_{A^*=-A} \Bigl\| \nabla (u - Ax)\Bigr\|^2 =
C \Bigl\| \nabla u - \<\nabla^\a u\>\Bigr\|^2.
\endeq

On the other hand,
\begin{align}
\inf_{\lambda_0\in\R}\Bigl\| \log\rho - \frac{|m|^2}{2\theta\rho^2}\Bigr\|^2 
& \leq 2\inf_{\lambda_0} \Bigl\| \log\rho - \frac{|\Sigma x|^2}{2\theta}
- \lambda_0 \Bigr\|^2 
+ 2 \Bigl\| \frac{|\Sigma x|^2}{2\theta}
- \frac{|m|^2}{2\theta\rho^2}\Bigr\|^2 \nonumber\\
& \leq 2\inf_{\lambda_0} \Bigl\| \log\rho - \frac{|\Sigma x|^2}{2\theta}
- \lambda_0 \Bigr\|^2
+ C \Bigl\| u - \Sigma x\|^2.
\label{logrhoa}
\end{align}
By Poincar\'e inequality, the second term in the right-hand side 
of~\eqref{logrhoa} can be bounded by a constant multiple of
$\|\nabla u - \Sigma\|^2$. As for the first term, it can also be bounded
by means of a Poincar\'e inequality:
\begin{align*}
\inf_{\lambda_0} \Bigl\| \log\rho - \frac{|\Sigma x|^2}{2\theta}
- \lambda_0 \Bigr\|^2
& = \Bigl\| \log\rho - \frac{|\Sigma x|^2}{2\theta}
- \left\<\log\rho - \frac{|\Sigma x|^2}{2\theta} \right\>\Bigr\|^2\\
& \leq C \, \Bigl\| \nabla \left(\log\rho - \frac{|\Sigma x|^2}{2\theta}
\right) \Bigr\|^2 \\
& = C \Bigl\| \frac{\nabla\rho}{\rho} + \frac{\Sigma^2 x}{\theta} \Bigr\|^2.
\end{align*}
All in all,
\[ \inf_{\lambda_0\in\R}\Bigl\| \log\rho - \frac{|m|^2}
{2\theta\rho^2}\Bigr\|^2 \leq
C \Bigl\|u - \Sigma x\Bigr\|^2 + 
\Bigl\| \frac{\nabla\rho}{\rho} + 
\frac{\Sigma^2 x}{\theta}\Bigr\|^2.\]
This combined with~\eqref{1F} and~\eqref{byPoincKorn}
concludes the verification of~\eqref{nowshowthat''}.
Then we can apply Theorem~\ref{thmmain} and get
the conclusion of Theorem~\ref{thmcvBEsph}.
\end{proof}

\begin{proof}[Proof of Lemma~\ref{lemPhiK}]
Let $\tilde{\Phi}:X\to \Phi(X_\infty+X)$.
By assumption, 0 is a minimizer of $\tilde{\Phi}$ on $\ker\Lambda$.
Since $\Phi$ is $K$-convex and differentiable, the same is true
of $\tilde{\Phi}$, so that
\[ \tilde{\Phi}(0) \geq \tilde{\Phi}(X) - 
\Bigl\< {\grad}' \tilde{\Phi}(X), \, X\Bigr\> 
+ \frac{K}{2} \|X\|^2,\]
where $\grad'$ stands for the gradient in the space $\ker\Lambda$.
It follows by Young's inequality that
\[ \tilde{\Phi}(X) - \tilde{\Phi}(0) \leq
\frac{\Bigl\| {\grad}' \tilde{\Phi}(X)\Bigr\|^2}{2K}.\]

But $\grad'\tilde{\Phi}$ is nothing but the orthogonal projection of
$\grad\Phi(X)$ (in $\H$) onto $\ker\Lambda$; so
\[ \Bigl\| {\grad}' \tilde{\Phi}(X)\Bigr\|
= \inf_{\lambda\in (\ker\Lambda)^\bot} 
\Bigl\|\grad \tilde{\Phi}(X) +\lambda \Bigr\|.\]
The conclusion of Lemma~\ref{lemPhiK} follows easily.
\end{proof}

\begin{Rk} I~don't know if the term in $\Pi_2$ can be dispended with
in~\eqref{IdPi2Pi3}; in any case this is an example where it is
convenient to have the general formulation of Assumption~\ref{asskey}(iii),
rather than the simplified inequality~\eqref{key'}. In the next
subsection, another example will be presented where this possibility
is crucially used (see Remark~\ref{rkgoodexample}).
\end{Rk}

\subsection{Maxwellian accommodation} \label{secMaxw}

In this subsection $\Om_x$ will again be a bounded smooth connected
open subset of $\R^N$ with unit Lebesgue measure, 
but now the boundary condition will be the
Maxwellian accommodation with a {\em fixed} temperature $T_w$.
Explicitly,
\begeq\label{accomm} 
x\in\pa\Om_x \Longrightarrow\qquad
f^+(x,v) = \left( \int f^-(x,v')\, |v'\cdot n|\,dv'\right)\,
M_w(v),
\endeq
where $f^+$ (resp. $f^-$) stands for the restriction of $f$ to
$\{v\cdot n>0\}$ (resp. $\{v\cdot n<0\}$), $n$ is the {\em inner}
unit normal vector, and $M_w$ is a fixed ``wall'' Maxwellian:
\[ M_w(v) = \frac{e^{-\frac{|v|^2}{2T_w}}}{(2\pi)^{\frac{N-1}2}\,
T_w^{\frac{N+1}2}}.\]
(The analysis would go through if one would impose a more general
condition involving a reflection kernel $C(v'\to v)$, as
in~\cite[Chapter~1]{cer:book:00}.)
An important identity which follows from~\eqref{accomm} is
\begeq\label{fvn}\forall x\in \pa\Om_x,\qquad
\int_{\R^N} f(x,v)\,(v\cdot n) =0;
\endeq
equivalently, the mean velocity satisfies
\begeq\label{uxn}
\forall x\in\pa\Om_x,\qquad u\cdot n=0.
\endeq

In this case there is only one conservation, namely the total mass.
Without loss of of generality, I~shall assume that the solution
is normalized so that
\begeq\label{conslawTw}
\int f\,dv\,dx =1.
\endeq

Then the unique equilibrium is the Maxwellian distribution with
constant temperature equal to the wall temperature:
\begeq\label{f8Tw} 
f_\infty(x,v) = \frac{e^{-\frac{|v|^2}{2T_w}}}{(2\pi)^{N/2}}.
\endeq

\begin{Thm}[Convergence for the Boltzmann equation with 
Maxwellian accommodation]
\label{thmcvBETw}
Let $f$ be a solution of~\eqref{BE} in a smooth bounded connected
spatial domain $\Om_x\subset\R^N$ with $|\Om_x|=1$.
Assume that $f$ satisfies the boundary
condition~\eqref{accomm}, the conservation
laws~\eqref{conslawTw} and the uniform regularity 
estimates~\eqref{uniformregBE}. Then
\[ \forall s\geq 0,\qquad 
\bigl\|f(t,\cdot) - f_\infty \bigr\|_{H^s} = O(t^{-\infty}),\]
where $f_\infty$ is defined by~\eqref{f8Tw}.
\end{Thm}

\begin{proof}[Proof of Theorem~\ref{thmcvBETw}]
The proof follows again the same pattern as in all the previous
theorems in this section. However, the Lyapunov functional is not
Boltzmann's $H$ functional, but a modified version of it:
\[ \E(f) = \int f\log f + \frac{1}{2T_w} \int f|v|^2\,dv\,dx.\]
Moreover, the sequence of nonlinear projection operators will be
\[ \Pi_1 f = M_{\rho\, u\, T};\qquad
\Pi_2 f = M_{\rho\, u\, T_w};\qquad
\Pi_3 f = M_{\rho\, 0\, T_w};\qquad
\Pi_4 f = M_{1\, 0\, T_w}.\]

In particular,
\begin{align*} \E(\Pi_1f) - \E(f_\infty) 
& = \int \rho\log\rho -\frac{N}2 \int \rho\log T
+ \frac1{T_w} \left( \int \rho \frac{|u|^2}{2} +
\frac{N}2 \int \rho T\right) \\
& \qquad\qquad\qquad\qquad\qquad\qquad\qquad -
\frac{N}2 \log T_w + \frac{N}2 \\
& = \int \rho\log\rho + \frac1{T_w} \int \rho\frac{|u|^2}{2}
+ \frac{N}2 \int \rho \left(\frac{T}{T_w}
- \log\frac{T}{T_w} - 1 \right).
\end{align*}
From this it is easy to check Assumption~\ref{asslyap}(ii).

The interesting features of this case reveal themselves when
we try to check Assumption~\ref{asskey}. First, by a classical computation,
\begin{align*}
{\cal D}(f) & = -\E'(f)\cdot (Bf) =
\int_{\Om_x} D(f(x,\cdot))\,dx 
+ \int_{\Om_x\times\R^N} (\log f+1) (-v\cdot\nabla_x f)\\
&\qquad\qquad + \frac{1}{2T_w} \int_{\Om_x\times\R^N} 
(v\cdot\nabla_x f) |v|^2\,dv\,dx\\
& = \int_{\Om_x} D(f(x,\cdot))\,dx - 
\int_{\pa\Om_x\times\R^N} 
\bigl (f\log f (x,v) + f(x,v)\frac{|v|^2}{2T_w}\bigr)
|v\cdot n|\,dv\,dx \\
& = \int_{\Om_x} D(f(x,\cdot))\,dx +
\int_{\pa\Om_x\times\R^N}
f \log \frac{f}{e^{-\frac{|v|^2}{2T_w}}}\,(v\cdot n)\,dv\,dx.
\end{align*}
As before, he first term in the right-hand side is controlled
below by $K[\E(f)-\E(\Pi_1f)]^{1+\var}$. The second term needs
some rewriting.
In view of~\eqref{fvn} and~\eqref{accomm}, we have, with the notation
$\rho-(x) = \int f^-(x,v) |v\cdot n|\,dv$,
\begin{align*} \int_{\pa\Om_x\times\R^N}
- f \log \frac{f}{e^{-\frac{|v|^2}{2T_w}}}\,(v\cdot n)\,dv\,dx
& = - \int f\log \left(\frac{f}{\rho^-\, M_w}\right)\,(v\cdot n)\,dv\,dx\\
& = \int_{v\cdot n<0}  
f\log \left(\frac{f}{\rho^-\, M_w}\right)\,|v\cdot n|\,dv\,dx\\
& = \int_{v\cdot n<0}
f |v\cdot n|\, \log \left(\frac{f\, |v\cdot n|}{\rho_-M_w\, |v\cdot n|}\right)
\,dv\,dx.
\end{align*}
This quantity takes the form of a nonnegative information functional,
as a particular case of the Darroz\`es--Guiraud--Cercignani 
inequality~\cite[Chapter~1]{cer:book:00}.
The Csisz\'ar--Kullback--Pinsker inequality will give an explicit lower
bound: For each $x\in\pa\Om_x$,
\begin{align*}\int_{v\cdot n<0}
f |v\cdot n|\, & 
\log \left(\frac{f\, |v\cdot n|}{\rho_-M_w\, |v\cdot n|}\right)
\,dv \\
& \geq \frac1{2\rho_-(x)} \Bigl\| f|v\cdot n| - \rho_- M_w|v\cdot n|
\Bigr\|_{L^1(\{v\cdot n<0\}; |v\cdot n|\,dv)}^2\\
& = \frac1{2\rho_-(x)} \Bigl\| f|v\cdot n| - \rho_- M_w|v\cdot n|
\Bigr\|_{L^1(|v\cdot n|\,dv)}^2.
%alignement horrible ci-dessus
\end{align*}

After interpolation and use of smoothness bounds, we conclude that
\begeq\label{DTw}
{\cal D}(f) \geq K \bigl[ \E(f) - \E(\Pi_1f)\bigr]^{1+\var}
+ K \bigl\| f - \rho_-M_w\bigr\|^{2+\var}_{L^q(\pa\Om_x\times\R^N;
\, |v\cdot n|\, |v|^q\,dx\,dv)},
\endeq
where $q$ is arbitrarily large and $\var$ is arbitrarily small.
A useful consequence of~\eqref{DTw} is
\begeq\label{DTw'}
{\cal D}(f) \geq K \|T - T_w\|^{2+\var}_{L^q(\pa\Om_x)}
+ K \|u\|^{2+\var}_{L^q(\pa\Om_x)},
\endeq
where again $q$ is arbitrarily large.

The other estimates are similar to the ones in the previous
subsections:
\begeq\label{otherPi}
\begin{cases}
\dps \Bigl\| (\Id-\Pi'_1)_{\Pi_1f}\cdot (B\Pi_1f)\Bigr\|^2
\geq K \bigl(\|\nabla T\|^2 + \|\{\nabla u\}\|^2\bigr);\\ \\
\dps \Bigl\| (\Id-\Pi'_2)_{\Pi_2f}\cdot (B\Pi_2f)\Bigr\|^2
\geq K \|\nabla^{\rm sym} u\|^2 ;\\ \\
\dps \Bigl\| (\Id-\Pi'_3)_{\Pi_3f}\cdot (B\Pi_3f)\Bigr\|^2
\geq K \|\nabla \rho\|^2.
\end{cases}\endeq

Thanks to~\eqref{DTw'} and~\eqref{otherPi}${}_1$,
\begin{align*}
{\cal D}(f) + \Bigl\| (\Id-\Pi'_1)_{\Pi_1f}\cdot (B\Pi_1f)\Bigr\|^2
& \geq K \Bigl( \|\nabla T\|_{L^2(\Om_x)}^2 + 
\|T-T_w\|_{L^q(\pa\Om_x)}^{2+\var}\Bigr) \\
& \geq K' \| T-T_w\|_{L^p(\Om_x)}^{2+\var},
\end{align*}
where the latter inequality comes from, say, the trace Sobolev inequality
if, say, $p=(2N)/(N-2)$ and $q=2(N-1)/(N-2)$.
(If $N=2$ a slightly different argument based on a variant of
the Moser--Trudinger inequality can be used to give the same result.) 
After interpolation one concludes that
\begeq\label{Dgreaterthan}
{\cal D}(f) + \Bigl\| (\Id-\Pi'_1)_{\Pi_1f}\cdot (B\Pi_1f)\Bigr\|^2
\geq K \| T-T_w\|^{2+\var} \geq K' \|\Pi_1 f - \Pi_2 f\|^{2+\var'}.
\endeq

Next, if $\Om_x$ is not axisymmetric, then the boundary 
condition~\eqref{uxn}, the Korn inequality from~\cite{DV:korn:02}
and the Poincar\'e inequality imply
\begin{align}\label{ifOmnotaxi} 
\Bigl\| (\Id-\Pi'_2)_{\Pi_2f}\cdot (B\Pi_2f)\Bigr\|^2
& \geq K \|\nabla^{\rm sym} u\|^2 \\
& \geq K'\|u\|^2 \geq K'' \|\Pi_2 f - \Pi_3 f\|^{2+\var}.\nonumber
\end{align}

If $\Om_x$ is axisymmetric, the previous argument breaks down,
but we can use~\eqref{DTw'} and replace~\eqref{ifOmnotaxi} by 
\begin{align}\label{ifOmaxi}
{\cal D}(f) + \Bigl\| (\Id-\Pi'_2)_{\Pi_2f}\cdot (B\Pi_2f)\Bigr\|^2
& \geq K \|\nabla^{\rm sym} u\|^2 + \|u\|_{L^2}^{2+\var}\\
& \geq K'\|u\|^{2+\var} \geq K'' \|\Pi_2 f - \Pi_3 f\|^{2+\var'},\nonumber
\end{align}
where the but-to-last inequality follows from a
{\em trace Korn inequality} (Proposition~\ref{thmkorn} in
Appendix~\ref{toolbox}).

Finally, 
\begin{align*} 
\Bigl\| (\Id-\Pi'_3)_{\Pi_3f}\cdot (B\Pi_3f)\Bigr\|^2
& \geq K \|\nabla\rho\|^2 \\
& \geq K \|\rho - 1\|^2 \geq K' \|\Pi_3 f - \Pi_4 f\|^{2+\var'}.
\end{align*}

Then Assumption~\ref{asskey}(iii) is satisfied, and one can
use Theorem~\ref{thmmain} to prove Theorem~\ref{thmcvBETw}.
\end{proof}

\begin{Rk}\label{rkgoodexample}
 This is an example where the range of $\Pi_1$ is much
larger than the set where the dissipation ${\cal D}$ vanishes.
Trying to devise a projection operator onto the space where
${\cal D}$ vanishes gives rise to a horrendous nonlocal variational 
problem whose solution is totally unclear. On the other hand,
inequalities~\eqref{Dgreaterthan} and~\eqref{ifOmaxi} 
would be false without the contribution of ${\cal D}(f)$. 
In this example we see that the possibility
to use the generalized condition appearing in Assumption~\ref{asskey}(iii), 
rather than the simplified condition~\eqref{key'}, leads to a great
flexibility.
\end{Rk}

\subsection{Further comments}

In many important situations (variable wall temperature,
evaporation problems, etc.), one is led to study {\em non-Maxwellian}
stationary solutions of the Boltzmann equation; then there is usually
no variational principle for these solutions, and the mere existence
of stationary solutions is a highly nontrivial problem, 
see e.g.~\cite{ark:stat:94,arknouri:stat:02,arknouri:statcyl:05}.

From the technical point of view, the non-Maxwellian nature of
the stationary state means that
if one defines $B=v\cdot\nabla_x$ (transport operator)
and $\C=Q$ (collision operator), then the equations $Bf_\infty=0$ 
and $\C f_\infty=0$ cease to hold. No need to say, Theorem~\ref{thmmain}
collapses, and it is quite hard to figure out how to save it.

There is a thin analogy with the (linear) problem of the oscillator
chain considered in Subsection~\ref{applosc} in the case when
the two temperatures are not equal; in that case a change of reference 
measure, based on Proposition~\ref{rulesdiffeq}(ii), was at least
able to reduce the problem to one of the type $A^*A+B$, $B^*=-B$.
By analogy, one could imagine that a first step to come to grips with
the quantitative analysis of stability for non-Maxwellian stationary
solutions consists in re-defining the ``antisymmetric'' and the
``diffusive'' parts of the Boltzmann equation by performing some
change of reference measure. Even this first step is nontrivial.

\part*{Appendices}

\setcounter{Thm}{0}
\def\theThm{\hspace*{0mm}A.\arabic{Thm}}
\def\thesection{\hspace*{0mm}A.\arabic{section}}

This last part is devoted to some technical appendices used throughout
the memoir, some of them with their own interest.

In Appendix~\ref{appproduct} I have gathered some sufficient conditions 
for a probability measure to admit a Poincar\'e inequality. After
recalling some well-known criteria for Poincar\'e inequality
in $\R^n$, I~shall prove some useful results about tensor products;
they might belong to folklore in certain mathematical circles, but I~am
not aware of any precise reference.

Appendices~\ref{appuniqueFP} and~\ref{appreg} are devoted to some
properties of the linear (kinetic) Fokker--Planck equation.
First in Appendix~\ref{appuniqueFP} I~shall prove a uniqueness
theorem; the method is quite standard, although computations are a bit tricky.
Appendix~\ref{appreg} is much more original and could be considered
as a research paper on its own right: There I~shall present a new 
strategy to get hypoelliptic regularization estimates. The method
has the advantage to be very elementary, to avoid fractional derivatives 
as well as localization, and to yield optimal exponents of
decay in short time. As Nash's theory of elliptic regularity, it is based on
differential equations satisfied by certain functionals of the solutions.
The results are nonstandard in several respects: They are global,
directly yield pointwise in time estimates, and apply for initial
data that do not lie in an $L^2$-type space.
I developed the method during a stay in Reading University,
from January to March 2003; thanks are due to Mike Cullen for his
hospitality. 

A closely related, but somewhat simpler strategy was found independently
and almost simultaneously by Fr\'ed\'eric H\'erau. I~shall explain his method
in Subsection~\ref{subvariants}, and develop it into
an abstract theorem of global regularization applying to the same kind 
of operators that have been considered in Part~\ref{partAAB} of this memoir.
This extension grew out from discussions with Denis Serre.

Finally, in Appendix~\ref{toolbox} I~gathered various technical lemmas
and functional inequalities which are used throughout the memoir.
I~draw the attention of the reader to the ``distorted Nash inequality''
appearing in Lemma~\ref{nash}, which might have an interesting role
to play in the future for ``global'' hypoelliptic regularization estimates.

\vfill\pagebreak

\section{Some criteria for Poincar\'e inequalities}
\label{appproduct}

To begin with, I~shall recall a popular and rather general criterion 
for Poincar\'e inequalities in $\R^n$.

\begin{Thm} \label{thmDS}
Let $V\in C^2(\R^n)$, such that $e^{-V}$ is a probability density
on $\R^n$. If
\begeq\label{DS} 
\frac{|\nabla V(x)|^2}2- \Delta V(x) \xrightarrow[|x|\to\infty]{} +\infty,
\endeq
then $\mu$ satisfies a Poincar\'e inequality.
\end{Thm}

\begin{proof}
The key estimate can be found in Deuschel
and Stroock~\cite[Proof of Theorem~6.2.21]{DS:LD:89}:
If $w = |\nabla V|^2/2 - \Delta V$, then for any $h\in C^1_c(\R^n)$,
\begeq\label{intwh2}\int w h^2\,d\mu \leq 4 \int |\nabla h|^2\,d\mu.
\endeq

Let $R_0>0$ be large enough that $w(|x|)\geq 0$ for $|x|\geq R_0$.
For $R>R_0$, define $\var(R):= [\inf \{w(|x|); |x|\geq R\}]^{-1}$;
then $\var(R)\to 0$ as $R\to\infty$.
So it follows from~\eqref{intwh2} that
\begeq\label{intxgR} 
\int_{|x|\geq R} h^2\,d\mu \leq \var(R) \left[
4 \int |\nabla h|^2\,d\mu - (\inf w) \int h^2\,d\mu\right].
\endeq

Now let $h\in C^1(\R^n,\R)$ with $\int h\,d\mu=0$. 
For any $R>0$, let $B_R$ be the ball of radius $R$ in $\R^n$, 
and let $\mu_R$ be the restriction of $\mu$ to $B_R$ (normalized
to be a probability measure). Since $B_R$ is bounded, $\mu_R$ satisfies 
a Poincar\'e inequality with a constant $P(R)$ depending on $R$, so
\[ \int h^2\,d\mu_R \leq P(R) \int |\nabla h|^2\,d\mu_R + \left( \int h\,d\mu_R \right)^2.\]
Of course $\mu_R$ has density $(\mu[B_R])^{-1}e^{-V(x)}1_{|x|\leq R}$. If $R$
is large enough, then $\mu[B_R]\geq 1/2$, so
\begeq\label{intermcritPoinc} 
\int_{|x|\leq R} h^2\,e^{-V} \leq P(R) \int_{|x|\leq R} 
|\nabla h|^2\,e^{-V} + 2\left( \int_{|x|\leq R} h\,e^{-V} \right)^2.
\endeq

Since $\int h e^{-V}=0$ and $e^{-V}$ is a probability density,
\begeq\label{estimreste}
\left( \int_{|x|\leq R} h\,e^{-V} \right)^2  
= \left( \int_{|x|>R} h\,e^{-V}\right)^2 \leq \int_{|x|>R} h^2 e^{-V}.
\endeq
Plugging this into~\eqref{intermcritPoinc}, one deduces that
\begeq\label{intermcritPoinc2} 
\int h^2\,e^{-V} \leq P(R) \int |\nabla h|^2\,e^{-V} + 
3 \int_{|x|> R} h^2 \,e^{-V}.
\endeq
Combining this with~\eqref{intxgR}, we recover
\[ \int h^2\,e^{-V} \leq [P(R) + 12\var(R)] \int |\nabla h|^2\,e^{-V} 
-3(\inf w)\var(R) \int h^2 \,e^{-V}.\]
So, if $R$ is large enough that $3(\inf w)\var(R)>-1$, one has
\[ \int h^2\, e^{-V} \leq \left(\frac{P(R)+12\var(R)}
{1+3 (\inf w)\var(R)}\right) \int |\nabla h|^2\, e^{-V}. \]
This concludes the proof of Theorem~\ref{thmDS}.
\end{proof}

The sequel of this Appendix is devoted to Poincar\'e inequalities
in product spaces.
It is well-known that ``spectral gap inequalities tensorize'', in the following sense:
If each $L_\ell$ ($\ell=1,2$) is a nonnegative operator on a Hilbert space $\H_\ell$, admitting a
spectral gap $\kappa_\ell$, then $L_1\otimes I + I\otimes L_2$ admits a spectral gap
$\kappa=\min (\kappa_1,\kappa_2)$. Now the goal is to extend this
result in a form which allows multipliers. 
I shall start with an abstract theorem and then particularize it.

\begin{Thm}\label{thmcoertens} 
For $\ell=1,2$, let $L_\ell$ be a nonnegative unbounded operator on a Hilbert space $\H_\ell$,
admitting a finite-dimensional kernel, and a spectral gap $\kappa_\ell>0$. 
Let $\ov{M}$ be a nonnegative unbounded operator acting on $\H_2$, 
whose restriction to the kernel $\K_2$ of $L_2$ is bounded and coercive. 
Then the unbounded operator 
\[ L = L_1\otimes \ov{M} + I\otimes L_2\]
admits a spectral gap $\kappa>0$. More precisely, 
for any nonnegative operator $M\leq \ov{M}$,
such that the restriction $M|_{\K_2}$ of $M$ to $\K_2$ satisfies
$\lambda I \leq M \leq \Lambda I$, one has
\[ \kappa \geq \min \left ( \frac{\kappa_2}2, \ \frac{\kappa_2}{16} \frac{\lambda}{\Lambda^2},\
\frac{\kappa_1}2 \lambda \right). \]
\end{Thm}

\begin{Thm} \label{thmSGproduct}

(i) For $\ell=1,2$, let $(X_\ell,\mu_\ell)$ be a probability space, and let
$L_\ell$ be a nonnegative operator on $\H_\ell=L^2(\mu_\ell)$, whose kernel
is made of constant functions, admitting a spectral gap $\kappa_\ell>0$. 
Let $\ov{m}$ be a nonnegative measurable function on $X_2$, 
which does not vanish $\mu_2$-almost everywhere, 
and $\ov{M}$ be the multiplication operator by $\ov{m}$.
Then the unbounded operator 
\[ L = L_1\otimes \ov{M} + I\otimes L_2\]
admits a spectral gap $\kappa>0$. More precisely, for any nonnegative function
$m\leq \ov{m}$, lying in $L^2(\mu_2)$,
\[ \kappa \geq \min \left (\frac{\kappa_2}2,\ 
\frac{\kappa_2}{16} \frac{\|m\|_{L^1}^2}{\|m\|_{L^2}^2},\
\frac{\kappa_1}2 \|m\|_{L^1} \right). \]
\sm

(ii) More generally, for each $\ell\in \{1,\ldots,N\}$, 
let $(X_\ell, \mu_\ell)$ be a probability space,
and let $L_\ell$ ($1\leq \ell\leq N$) be a nonnegative symmetric operator 
on $\H_\ell=L^2(\mu_\ell)$, whose kernel is made of constant 
functions, admitting a spectral gap $\kappa_\ell$;
let $\ov{m}_\ell$ be a nonnegative measurable
function on $X_{\ell+1}\times\ldots\times X_N$, which does not vanish 
$\mu_{\ell+1}\otimes\ldots\otimes \mu_N$-almost everywhere,
and let $\ov{M}_\ell$ be the associated multiplication operator. Then
the linear operator
\[ L = \sum_{\ell=1}^N I^{\otimes \ell-1}\otimes L_\ell\otimes \ov{M}_\ell \]
admits a spectral gap.
\end{Thm}

\begin{Ex} Let $\mu$ and $\nu$ be two probability measures on $\R$, each satisfying
a Poincar\'e inequality. Equip $\R^2$ with the tensor measure 
$\mu\otimes\nu(dx\,dy)=\mu(dx)\,\nu(dy)$. Then
\[ L = -(\pa_x^*\pa_x + x^2 \pa_y^*\pa_y) \]
is \emph{coercive} on $L^2(\mu\otimes\nu)/\R$.
\end{Ex}

\begin{proof}[Proof of Theorem~\ref{thmcoertens}]
Let $P_\ell$ be the orthogonal projection on $(\ker L_\ell)^\bot$ in
$\H_\ell$. The spectral gap assumption means $L_\ell\geq \kappa_\ell\, P_\ell$. Let $M$
be the multiplication operator by $m$, then $\ov{M}\geq M$.

When applied to nonnegative operators, tensorization preserves the order:
when $A\geq A'\geq 0$ and $B\geq B'\geq 0$, one has $A\otimes B\geq
A\otimes B'\geq A'\otimes B'$. Thus,
\[ L_1\otimes \ov{M} + I\otimes L_2 \geq \kappa_1 P_1\otimes M + 
\kappa_2 I\otimes P_2. \]
So it is sufficient to prove the theorem when $L_\ell = P_\ell$ and
$\ov{m}=m\in L^2(\mu_2)$.

Let $(e_i^1)_{i\geq 0}$ be an orthonormal basis for $\H_1$, such that
$(e_i^1)_{i\leq k_1}$ is an orthonormal basis of $\K_1:=\ker L_1$; and let
$(e_j^2)_{j\geq 0}$ be an orthonormal basis for $\H_2$, such that
$(e_j^2)_{j\leq k_2}$ is an orthonormal basis of $\K_2:=\ker L_2$. Then
$(e_i^1\otimes e_j^2)_{i,j\geq 0}$ is an orthonormal basis for $\H$.
Moreover, the kernel of $L$ is the vector space generated by
$(e_i^1\otimes e_j^2)_{i,j\in \K}$, where 
$\K:=\{(i,j);\ i\leq k_1, j\leq k_2\}$. So the goal is to prove
\begin{multline*} 
f =\sum_{(i,j)} c_{ij} e_i^1 \otimes e_j^2 \Longrightarrow \\ \qquad
\kappa_1\ip{(P_1\otimes M) f}{f} + \kappa_2 \ip{(I\otimes P_2) f}{f} \geq
\kappa \sum_{(i,j)\notin \K} c_{ij}^2. 
\end{multline*}

First of all,
\begeq\label{IP2} \ip{(I\otimes P_2) f}{f} = \sum_{(i,j)} c_{ij}\, e_i^1\otimes P_2 e_j^2
= \sum_{i\geq 0;\,j\geq k_2+1} c_{ij}\, e_i^1 \otimes e_j^2. 
\endeq
Next,
\begin{align*} 
\ip{(P_1\otimes M) f}{f} & = 
\sum_{(i,i',j,j')} c_{ij} c_{i'j'} \ip{P_1e_i^1}{e_{i'}^1} 
\ip{Me_j^2}{e_{j'}^2} \\
& = \sum_{i\geq k_1+1;\: j,j'\geq 0} c_{ij} c_{ij'} \ip{Me_j^2}{e_{j'}^2}\\
& = \sum_{i\geq k_1+1;\: j,j'\geq k_2+1} c_{ij} c_{i j'} \ip{Me_j^2}{e_{j'}^2}
+ 2 \sum_{i\geq k_1+1;\: j\geq k_2+1;\: j'\leq k_2} c_{ij} c_{i j'} 
\ip{Me_j^2}{e_{j'}^2} \\
  & \hspace*{70mm}+ \sum_{i\geq k_1+1;\: j,j'\leq k_2} c_{ij} c_{i j'} 
\ip{Me_j^2}{e_{j'}^2}.
\end{align*}
We shall estimate these three sums one after the other:

- The first sum might be rewritten as
\[ \sum_{i\geq k_1+1} \Bigl\<M \Bigl( \sum_{j\geq k_2+1} c_{ij} e_j^2\Bigr),\
\Bigl(\sum_{j\geq k_2+1} c_{ij} e_j^2\Bigr)\Bigr\>,\]
and is therefore nonnegative. 

- Similarly, the third sum might be rewritten as
\[ \sum_{i\geq k_1+1} \Bigl\<M \Bigl( \sum_{j\leq k_2} c_{ij} e_j^2\Bigr),\
\Bigl(\sum_{j\leq k_2} c_{ij} e_j^2\Bigr)\Bigr\>,\]
which can be bounded below by
\[ \lambda \sum_{i\geq k_1+1;\, j\leq k_2} c_{ij}^2.\]

- Finally, by applying the inequality $\|Me_j\|\leq \Lambda$ ($j\leq k_2$)
and the Cauchy--Schwarz inequality twice, one can bound the second sum from below by
\[ -2 \sum_{j\leq k_2} \sum_{i\geq k_1+1} c_{ij} \|Me_j^2\|\, 
\left\|\sum_{j'\geq k_2+1} c_{ij'} e_{j'}^2 \right\| \] 
\[ \geq
-2 \Lambda\ \sum_{j\leq k_2} \sqrt{\sum_{i\geq k_1+1} c_{ij}^2}\:
\sqrt{\sum_{i\geq k_1+1} \Bigl\| \sum_{j'\geq k_2+1} c_{ij'} e_{j'}^2 \Bigr\|^2}\]
\[ \geq -2\Lambda \sqrt{\sum_{i\geq k_1+1;\, j\leq k_2} c_{ij}^2}\:
\sqrt{\sum_{i\geq k_1+1;\: j\geq k_2+1} c_{ij}^2}.\]

All in all,
\begin{align*} 
\ip{(P_1\otimes M) f}{f} & \geq
\lambda \sum_{i\geq k_1+1;\: j\leq k_2} c_{ij}^2
- 2\Lambda\: \sqrt{\sum_{i\geq k_1+1;\: j\leq k_2} c_{ij}^2}\:
\sqrt{\sum_{i\geq k_1+1;\: j\geq k_2+1} c_{ij}^2} \\
& \geq \frac{\lambda}2 \sum_{i\geq k_1+1;\: j\leq k_2} c_{ij}^2
- \frac{4\Lambda^2}{\lambda} \sum_{i\geq k_1+1;\: j\geq k_2+1} c_{ij}^2.
\end{align*}

Combining this with~\eqref{IP2}, we see that for all $\theta\in [0,1]$,
\begin{align*} \ip{Lf}{f} & \geq \kappa_2 \ip{(I\otimes P_2)f}{f} +
\kappa_1 \theta \ip{(P_1\otimes M)f}{f} \\
& \geq \kappa_2 \sum_{i\geq 0;\, j\geq k_2+1} c_{ij}^2
+ \frac{\kappa_1\theta \lambda}2 \sum_{i\geq k_1+1;\: j\leq k_2} c_{ij}^2
- \frac{4 \kappa_1 \theta \Lambda^2}{\lambda} 
\sum_{i\geq k_1+1;\: j\geq k_2+1} c_{ij}^2\\
& \geq \kappa \left(\sum_{i\geq 0;\: j\geq k_2+1} c_{ij}^2 + 
\sum_{i\geq k_1+1;\: j\leq k_2} c_{ij}^2 \right),
\end{align*}
with
\[ \kappa:= \min \left( \kappa_2 - \frac{4 \kappa_1 \theta \Lambda^2}{\lambda},
\ \frac{\kappa_1\theta \lambda}2\right).\]
To conclude the proof of Theorem~\ref{thmcoertens}, it suffices to choose
\[ \theta:= \min \left( 1, \,
\frac{\kappa_2 \lambda}{8\kappa_1 \Lambda^2}\right).\]
\end{proof}

\begin{proof}[Proof of Theorem~\ref{thmSGproduct}]
Let $M$ be the multiplication operator by $m$. The restriction of $M$ to
constant functions is obviously coercive with constant $\lambda:= \int m\,d\mu_2$,
and $M$ is bounded by $\|m\|_{L^1}I$. Then (i) follows by a direct application
of Theorem~\ref{thmcoertens}. After that, statement (ii) follows from (i) 
by a simple induction on $N$.
\end{proof}

\section{Well-posedness for the Fokker--Planck equation} \label{appuniqueFP}

The goal of this Appendix is the following uniqueness theorem:

\begin{Thm} \label{thmuniqueFPf'}
With the notation of Theorem~\ref{thmuniqueFP}, for any
$f_0\in L^2((1+E)\,dv\,dx)$, the Fokker--Planck equation~\eqref{kFP} admits
at most one distributional solution $f=f(t,x,v) \in C(\R_+;{\cal D}'(\R^n_x\times\R^n_v))
\cap L^\infty_\loc(\R_+;L^2((1+E)\,dv\,dx)) \cap L^2_\loc (\R_+; H^1_v(\R^n_x\times\R^n_v))$,
such that $f(0,\cdot)=f_0$.
\end{Thm}

\begin{Rk} The a priori estimates
\[ \frac{d}{dt} \int f^2\,dv\,dx = - 2\int |\nabla_v f|^2\,dv\,dx + 2n \int f^2\,dv\,dx \]
\[ \frac{d}{dt} \int f^2 E\,dv\,dx = - 2 \int |\nabla_v f|^2\,dv\,dx
+ n \int f^2 (1+E)\,dv\,dx\]
allow to prove {\em existence} of a solution, too, for an initial
datum $f_0\in L^2((1+E)\,dv\,dx)$; but this is not what we are after
here. (Actually, an existence theorem can be established under much more
general assumptions.)
\end{Rk}

Before going on with the argument, I should explain why the uniqueness statement
in Theorem~\ref{thmuniqueFPf} implies the one in Theorem~\ref{thmuniqueFP}.
In that case, Proposition~\ref{rulesdiffeq}(iii) can be applied even if
$\nabla V$ is only continuous: indeed, the differential operator
$\nabla V(x)\cdot\nabla_v$ always makes distributional sense. So, if
$h$ is any solution of~\eqref{kFPh}, satisfying the assumptions of
Theorem~\ref{thmuniqueFP}, then $f:=h\rho_\infty$ defines a
solution of~\eqref{kFP}, and it also satisfies the assumptions in
Theorem~\ref{thmuniqueFPf}, in view of the inequalities
\[ \int f^2 (1+E)\,dv\,dx \leq \int f^2 e^E\,dv\,dx = \int h^2 e^{-E}\,dv\,dx,\]
\begin{align*} \int |\nabla_v f|^2\,dv\,dx & = 
\int |\nabla_v (\rho_\infty h)|^2\,dv\,dx \\
& \leq 2 \int |\nabla_v h|^2\,\rho_\infty^2\,dv\,dx 
+ 2 \int h^2 |\nabla_v \rho_\infty|^2\,dv\,dx\\
& \leq C \sup_{x,v} \Bigl[ (1+|v|^2)e^{-V(x)} e^{-\frac{|v|^2}2} \Bigr]
\left( \int |\nabla_v h|^2\,d\mu + \int h^2\,d\mu\right).
\end{align*}

\begin{proof}[Proof of Theorem~\ref{thmuniqueFPf'}]
By linearity, it is enough to prove
\[ \|f(T,\cdot)\|_{L^2} \leq e^{CT} \|f(0,\cdot)\|_{L^2},\]
which will also yield short-time stability.
So let $f$ solve the Fokker--Planck equation in distribution sense,
and let $T>0$ be an arbitrary time.

For any $\varphi\in C^\infty(t,x,v)$, compactly supported in $(0,T)\times\R^n_x\times\R^n_v$,
one has
\[ \iint f \Bigl( \partial_t \varphi + v\cdot\nabla_x\varphi 
-\nabla V(x)\cdot\nabla_v \varphi + \Delta_v\varphi - v\cdot\nabla_v\varphi\Bigr)
\,dt\,dv\,dx =0.\]
Since $f\in L^\infty([0,T]; L^2(\R^n_x\times\R^n_v))$, a standard approximation procedure
shows that
\begin{multline} \label{fphifphi}
\int f(T,\cdot)\varphi(T,\cdot)\,dv\,dx - \int f(0,\cdot)\varphi(0,\cdot)\,dv\,dx =\\
\iint f \Bigl( \partial_t \varphi + v\cdot\nabla_x\varphi 
-\nabla V(x)\cdot\nabla_v \varphi + \Delta_v\varphi - v\cdot\nabla_v\varphi\Bigr)
\,dt\,dv\,dx 
\end{multline}
for all $\varphi\in C^1((0,T); C^2_c(\R^n_x\times\R^n_v))\cap C([0,T];C^2_c(\R^n_x\times\R^n_v))$,
where $\iint$ stands for the integral over $[0,T]\times\R^n_x\times\R^n_v$.

Let $\chi,\eta$ be $C^\infty$ functions on $\R^n$ with $0\leq \chi\leq 1$,
$\chi(x)\equiv 1$ for $|x|\leq 1$, $\chi(x)\equiv 0$ for $|x|\geq 2$,
$\eta\geq 0$, $\int \eta=1$, $\eta$ radially symmetric, $\eta(x)\equiv 0$
for $|x|\geq 1$. With the notation $\var=(\var_1,\var_2)$,
$\delta=(\delta_1,\delta_2)$. Define
\[ \chi_\var(x,v) = \chi(\var_1 x)\,\chi(\var_2 v), \qquad
\eta_\delta(x,v)=\eta(\delta_1 x)\,\eta(\delta_2 v).\]
In words: $\chi_\var$ is a family of smooth cut-off functions, and $\eta_\delta$
is a family of mollifiers. (The introduction of $\eta_\delta$ is the main
modification with respect to the argument in~\cite[Proposition~5.5]{helfnier:witten:05}.)

Define now
\[ f_{\var,\delta}:= (f\chi_\var)\ast\eta_\delta,\qquad
\varphi_{\var,\delta} := \chi_\var ((f\chi_\var)\ast\eta_\delta\ast\eta_\delta).\]
The goal is of course to let $\delta\to 0$, $\var \to 0$ in a suitable way.

Since $\eta$ is radially symmetric, the identity $\int g(f\ast\eta) = \int (g\ast\eta)f$
holds true. So, for any $t\in [0,T]$,
\begeq\label{fta}\int f(t,\cdot)\varphi_{\var,\delta}(t,\cdot)\,dv\,dx = 
\int f_{\var,\delta}(t,\cdot)^2\,dv\,dx.
\endeq
Similarly,
\begeq\label{ftb} \int f\partial_t \varphi_{\var,\delta}\,dv\,dx =
\int f_{\var,\delta}\, \partial_t f_{\var,\delta}\,dv\,dx 
= \frac12\:\frac{d}{dt} \int f_{\var,\delta}^2\,dv\,dx.
\endeq
By combining~\eqref{fta} and~\eqref{ftb}, we get
\[ \iint f\partial_t \varphi_{\var,\delta}\,dv\,dx\,dt =
\frac12 \left( \int f(T,\cdot)\varphi_{\var,\delta}(T,\cdot)\,dv\,dx 
- \int f(0,\cdot)\varphi_{\var,\delta}(0,\cdot) \,dv\,dx\right). \]
So, by plugging $\varphi=\varphi_{\var,\delta}$ into~\eqref{fphifphi},
one obtains
\begeq \label{fphi'}
\frac12 \left( \int f(T,\cdot)\varphi(T,\cdot)\,dv\,dx \right.\left.  
- \int f(0,\cdot)\varphi(0,\cdot) \,dv\,dx\right)\endeq
\begeq= \iint f \chi_\var 
\bigl(v\cdot\nabla_x - \nabla V(x)\cdot\nabla_v 
   + \Delta_v -v\cdot\nabla_v\bigr) 
( f_{\var,\delta}\ast \eta_\delta)\, dv\,dx\,dt \label{fphi1}
\endeq
\begeq +\iint f \Bigl[ \bigl(v\cdot\nabla_x - \nabla V(x)\cdot\nabla_v 
+ \Delta_v -v\cdot\nabla_v\bigr) \chi_\var\Bigr] 
(f_{\var,\delta}\ast\eta_\delta)\,dv\,dx\,dt \label{fphi2}
\endeq
\begeq
+ 2 \iint f\chi_\var \nabla_v\chi_\var\cdot\nabla_v 
(f_{\var,\delta}\ast\eta_\delta)
\,dv\,dx\,dt \label{fphi3}. 
\endeq

For any given $\var>0$, all the functions involved are restricted to a compact set $K_\var$
in the variable $X=(x,v)$, uniformly in $\delta\leq 1$.
Now use the identities $\nabla(g\ast\eta)=(\nabla g)\ast\eta$,
$\Delta(g\ast\eta) = (\Delta g)\ast\eta$, $\int g(h\ast\eta) =\int h(g\ast\eta)$
to rewrite~\eqref{fphi1} as
\begin{multline} \label{fphi'1}
\iint f_{\var,\delta} \bigl(v\cdot\nabla_x 
- \nabla V(x)\cdot\nabla_v + \Delta_v - v\cdot\nabla_v\bigr)
f_{\var,\delta}\,dv\,dx\,dt \\
+ \iint f\chi_\var \Bigl[ \xi\cdot \nabla (f_{\var,\delta}\ast\eta_\delta)
- (\xi\cdot\nabla f_{\var,\delta})\ast\eta_\delta\Bigr],
\end{multline}
where $\xi$ is a temporary notation for the vector field
$(v,-\nabla V(x)-v)$. By integration by parts, the first integral
in~\eqref{fphi'1} can be rewritten as
\begeq\label{fphi''1} 
- \iint |\nabla_v f_{\var,\delta}|^2 - \frac12 \iint (v\cdot\nabla_v)f_{\var,\delta}^2
= - \iint |\nabla_v f_{\var,\delta}|^2 + \frac{n}2 \iint f_{\var,\delta}^2.
\endeq

Now we should estimate
\begeq\label{Biglxi}
\Bigl\| \xi\cdot \nabla (f_{\var,\delta}\ast\eta_\delta) - 
(\xi\cdot\nabla f_{\var,\delta})\ast\eta_\delta \Bigr\|_{L^2} 
= \left\| \int [\xi(X)-\xi(Y)]\cdot \nabla f_{\var,\delta}(Y)\,
\eta_\delta(Y-X) \right\|_{L^2(dX)}.
\endeq

We shall estimate the contributions of $v\cdot\nabla_x$, $v\cdot\nabla_v$
and $\nabla_x V\cdot\nabla_v$ separately. First, with obvious notation,
\begin{align*}
\bigl\| [v\cdot\nabla_x, \eta_\delta\ast] \, f_{\var,\delta}\bigr\|_{L^2} 
& = \Bigl\| \int (v-w)\cdot\nabla_x f_{\var,\delta}(y,w)\,
\eta_{\delta_1}(x-y)\,\eta_{\delta_2}(v-w)\,dw\,dy\Bigr\|_{L^2}\\
& = \Bigl\| \int f_{\var,\delta}(y,w)\,(v-w)\cdot \nabla_x\eta_{\delta_1}(x-y)\,
\eta_{\delta_2}(v-w)\,dw\,dy \Bigr\|_{L^2}.
\end{align*}
Inside the integral, one has $|v-w|\leq \delta_2$, $|x-y|\leq\delta_1$,
and also $|\nabla_x\eta_{\delta_1}|= O(\delta_1^{-(n+1)})$,
$\eta_{\delta_2}=O(\delta_2^{-n})$; so, all in all,
\begin{align*}
\bigl\| [v\cdot\nabla_x, \eta_\delta\ast] \, f_{\var,\delta}\bigr\|_{L^2} 
& \leq C \frac{\delta_2}{\delta_1} \ \Bigl\|
\frac{1}{\delta_1^n\,\delta_2^n} \int f_{\var,\delta}(y,w)\,
1_{|x-y|\leq \delta_1} 1_{|v-w|\leq\delta_2}\,dw\,dy \Bigr\|_{L^2}\\
& \leq C \frac{\delta_2}{\delta_1}\|f_{\var,\delta}\|_{L^2}\\ 
& \leq C \frac{\delta_2}{\delta_1}\|f\|_{L^2},
\end{align*}
where the last two inequalities follow from Young's convolution inequality.

Next,
\[ \bigl\| [v\cdot\nabla_v, \, \eta_\delta\ast] f_{\var,\delta} \bigr\|_{L^2}
= \Bigl\| \int (v-w)\cdot \nabla_v f_{\var,\delta} (y,w)\,\eta_{\delta_1}(x-y)\,
\eta_{\delta_2}(v-w)\,dw\,dy \Bigr\|.\]
Using the fact that $|v-w|\leq \delta_2$ inside the integral and applying
Young's convolution inequality as before, we find
\[ \bigl\| [v\cdot\nabla_v, \, \eta_\delta\ast] f_{\var,\delta}\bigr\|_{L^2}
\leq C\, \delta_2 \|\nabla_v f_{\var,\delta}\|_{L^2}
\leq C \,\delta_2\|\nabla_v f\|_{L^2}.\]

Finally,
\begin{align*}
\bigl\| [\nabla_x V\cdot\nabla_v,\, & 
\eta_\delta\ast] f_{\var,\delta} \bigr\|_{L^2} \\
%horreur alignement ci-dessus 
& = \Bigl\| \int [\nabla V(x) -\nabla V(y)]\cdot\nabla_v f_{\var,\delta}(y,w)\,
\eta_{\delta_1}(x-y)\,\eta_{\delta_2}(v-w)\,dw\,dy\Bigr\|_{L^2}\\
& \leq C \ \sup\Bigl\{ |\nabla V(x) - \nabla V(y)|;\
|x-y|\leq \delta_1;\ x,y\in K_{\var_1}\Bigr\} \ \|\nabla_v f_{\var,\delta}\|_{L^2}\\
& \leq C\, \theta_{\var_1}(\delta_1)\, \|\nabla_v f_{\var,\delta}\|_{L^2},
\end{align*}
where $\theta_\var$ stands for the modulus of continuity of $\xi$ on the compact
set $K_\var$. In all these estimates, the $L^2$ norm was taken with respect to 
all variables $t,x,v$. The conclusion is that the $L^2$ norm in~\eqref{Biglxi}
is bounded like
\begeq\label{fphi''12} 
O \left( \frac{\delta_2}{\delta_1} \|f\|_{L^2}^2 + 
\delta_2 \|\nabla_v f\|_{L^2} + \theta_{\var_1}(\delta_1)\|\nabla_v f\|_{L^2}\right).
\endeq

Next, since $\|\nabla_v\chi_\var\|_{L^\infty} \leq C \var_2$, it is possible
to bound~\eqref{fphi3} by
\begin{multline}\label{fphi''3} 
C \var_2 \|f\chi_\var\|_{L^2} 
\|\nabla_v (f_{\var,\delta}\ast\eta_\delta)\|_{L^2}
\leq C\var_2 \|f\|_{L^2} (\|\nabla_v f\|_{L^2} + \|f\|_{L^2}).
\end{multline}

Finally, the terms in the integrand of~\eqref{fphi2}
can be bounded with the help of the inequalities
\begin{multline*} |v\cdot \nabla_x \chi_\var(x,v)| \leq C |v|\var_1,\qquad
|\nabla V(x)\cdot\nabla_v\chi_\var(x,v)| \leq C \var_2 M(\var_1^{-1}), \\
|\Delta_v\chi_\var(x,v)| \leq C\, \var_2^2, \qquad
|v\cdot\nabla_v \chi_\var(x,v)| \leq C |v|\var_2, 
\end{multline*}
where $M(R):= \sup \{ |\nabla V(x)|;\: |x|\leq 2R\}$.
Then, by Cauchy--Schwarz again, \eqref{fphi2} can be bounded by
\[C \Bigl [ \var_1 + \var_2 ( 1+M(\var_1^{-1}))\Bigr]
\sqrt{\iint f^2 |v|^2\,dv\,dx\,dt}\ \sqrt{\iint (f_{\var,\delta}\ast\eta_\delta)^2\,dv\,dx\,dt}.\]
Since $|v|^2 \leq 2E - 2 (\inf V)$, in the end~\eqref{fphi2} is controlled by
\begeq\label{fphi''2} C \Bigl [ \var_1 + \var_2 ( 1+M(\var_1^{-1}))\Bigr]\
\sqrt{\iint f^2 E\,dv\,dx\,dt} \sqrt{\iint f^2\,dv\,dx\,dt}. 
\endeq

By plugging the bounds~\eqref{fphi''1}, \eqref{fphi''2} and~\eqref{fphi''3} 
into~\eqref{fphi'}, we conclude that
\begin{multline} \label{mult12}
\frac12 \left( \int f_{\var,\delta}^2 (T,x,v)\,dv\,dx -
\int f_{\var,\delta}^2(0,x,v)\,dv\,dx \right) \\
\leq - \iint |\nabla_v f_{\var,\delta}|^2 (t,x,v)\,dv\,dx\,dt 
+ \frac{n}2 \iint f_{\var,\delta}^2 (t,x,v)\,dv\,dx\,dt \\
+ C \Bigl( \frac{\delta_2}{\delta_1}\|f\|_{L^2}^2
+ \delta_2 \|f\|_{L^2}\,\|\nabla_v f_{\var,\delta}\|_{L^2}
+ \theta_{\var_1}(\delta_1)
\|\nabla_v f_{\var,\delta}\|_{L^2}\|f_{\var,\delta}\|_{L^2} \\
+ C (\var_1 + \var_2 M(\var_1^{-1})) \|f(1+E)\|_{L^2}^2,
\end{multline}
where all the $L^2$ norms in the right-hand side are with respect
to $dv\,dx\,dt$. Now let $\delta_2\to 0$, then $\delta_1\to 0$,
then $\var_2\to 0$, then $\var_1\to 0$, then $\delta\to 0$:
all the error terms in the right-hand side of~\eqref{mult12} vanish
in the limit, and $f_{\var,\delta}$ converges to $f$ almost everywhere 
and in $L^2(dv\,dx\,dt)$. So
\begin{align*}
\int f^2(T,x,v)\,dv\,dx & \leq
\liminf \int f_{\var,\delta}^2 (T,x,v)\,dv\,dx \\
& \leq \liminf \left[ \int f_{\var,\delta}^2 (0,x,v)\,dv\,dx +
\frac{n}2 \iint f_{\var,\delta}^2 (t,x,v)\,dv\,dx\,dt \right]\\
& = \int f^2(0,x,v)\,dv\,dx + \frac{n}2 \iint f^2(t,x,v)\,dv\,dx\,dt.
\end{align*}
By Gronwall's lemma,
\[ \|f(t,\cdot)\|_{L^2(\R^n\times\R^n)} \leq e^{\frac{nt}{4}} 
\|f(0,\cdot)\|_{L^2(\R^n\times\R^n)},\]
which concludes the argument.
\end{proof}

\begin{Rk} Just as in~\cite[Proposition~5.5]{helfnier:witten:05}, 
the particular structure of the Fokker--Planck equation was used
in the estimate $|\nabla V(x)\cdot\nabla_v\chi_\var|\leq \var_2 M(\var_1^{-1})$.
It would be interesting to understand to what extent this computation
can be generalized to larger classes of linear equations, and whether
this has anything to do with the hypoelliptic structure.
\end{Rk}

\section{Some methods for global hypoellipticity} \label{appreg}

This Appendix is devoted to various regularization estimates for the
Fokker--Planck equation. 
I~shall consider only two particular cases (those which were used in
the present paper): First, the $L^2(\mu)\to H^1(\mu)$ regularization for
the Fokker--Planck equation in the form~\eqref{kFPh'}; secondly,
the $M\to H^k$ regularization for the Fokker--Planck equation
in the form~\eqref{kFPf'} (Here $M$ is the space of bounded measures,
and $H^k$ is the {\em non-weighted} Sobolev space of order $k$.)

\subsection{From weighted $L^2$ to weighted $H^1$}\label{weigL2H1}

In the sequel $V$, is a $C^2$ potential on $\R^n$, bounded below, 
$\gamma(v) = (2\pi)^{-n/2} e^{-|v|^2/2}$ is the standard Gaussian,
and $\mu(dx\,dv)=\gamma(v) e^{-V(x)}\,dv\,dx$ 
stands for the equilibrium measure associated 
with the Fokker--Planck equation~\eqref{kFPh'} 
(it might have finite or infinite mass).
Apart from that, the only regularity assumption is the 
existence of a constant $C$ such that
\begeq \label{hypV}
|\nabla^2 V| \leq C ( 1+ |\nabla V|).
\endeq
As we shall see, this is sufficient to get estimate~\eqref{estregFP},
independently of the fact that $e^{-V}$ satisfies the Poincar\'e 
inequality~\eqref{Poincx} or not.

\begin{Thm} \label{thmhypFP} Let $V$ be a $C^2$ function on $\R^n$,
bounded below and satisfying~\eqref{hypV}.
Then, solutions of the Fokker--Planck equation~\eqref{kFPh'} with initial datum $h_0$ satisfy
\[ 0\leq t\leq 1 \Longrightarrow\qquad
\|\nabla_x h(t,\cdot)\|_{L^2(\mu)} + \sum_{k=1}^3 \|\nabla_v^k h(t,\cdot)\|_{L^2(\mu)}
\leq \frac{C}{t^{3/2}} \|h_0\|_{L^2(\mu)}\]
for some constant $C$, only depending on $n$ and the constant
$C$ appearing in~\eqref{hypV}.
\end{Thm}

\begin{Rk}
These estimates seem to be new. The proof can be adapted to cover the case 
of $L^1$ initial data, at the price of a deterioration of the exponents. 
I~shall explain this later on.
\end{Rk}

\begin{Rk}\label{H1H3interp}
Theorem~\ref{thmhypFP} shows that (with obvious notation)
$e^{-tL}$ maps $L^2$ into $H^1_x\cap H^3_v$ with norm $O(t^{-3/2})$. It also maps
$L^2$ into $L^2$ with norm $O(1)$; so, by interpolation, it maps $L^2$ into
$H^{\alpha}_x\cap H^{3\alpha}_v$ with norm $O(t^{-3\alpha/2})$, for all $\alpha\in [0,1]$.
Since $\int_0^1 e^{-t(1+L)}\,dt$ is a parametrix for $(I+L)^{-1}$,
and $t^{-\beta}$ is integrable at $t=0$ for $\beta<1$, one can deduce
a ``stationary'' hypoelliptic regularity estimate {\em \`a la Kohn}:
\begeq\label{Kohn}
\|h\|_{H_x^\alpha(\mu)}  + \|h\|_{H_v^{3\alpha}(\mu)} \leq 
C \bigl(\|h\|_{L^2(\mu)} + \|Lh\|_{L^2(\mu)}\bigr), \qquad \forall \alpha<2/3.
\endeq
With a much more refined analysis, it is actually possible to catch the
optimal exponent $\alpha=2/3$ in the above estimate. (This realization
came after discussions with Christ.) I~shall not develop this tricky
issue here.
\end{Rk}

\begin{proof}[Proof of Theorem~\ref{thmhypFP}]
As a consequence of Theorem~\ref{thmuniqueFP} and an 
approximation argument which is omitted here, 
it is sufficient to prove this theorem for smooth, rapidly decaying solutions.
So I~shall not worry about technical justification of the manipulations below.
Also, $C$ will stand for various constants which only depend on
$n$ and the constant in~\eqref{hypV}.

The following estimates will be used several times. 
As a consequence of Lemma~\ref{lemD2Vbdd} in Appendix~\ref{toolbox},
for each $v$,
\[ \int_{\R^n} |\nabla V(x)|^2 g^2(x,v) e^{-V(x)}\,dx \leq C 
\left ( \int g^2(x,v) e^{-V(x)}\,dx + \int |\nabla_x g(x,v)|^2 e^{-V(x)}\,dx
\right);\]
by integrating this with respect to $\gamma(v)\,dv$ one obtains
\begeq\label{D2Vbdd1}
\int_{\R^n\times\R^n} |\nabla V|^2 g^2\,d\mu \leq C 
\left ( \int g^2\,d\mu + \int |\nabla_x g|^2\,d\mu \right).
\endeq
Similarly,
\begeq\label{D2Vbdd2}
\int_{\R^n\times\R^n} |v|^2 g^2\,d\mu \leq C 
\left ( \int g^2\,d\mu + \int |\nabla_v g|^2\,d\mu \right).
\endeq

Now we turn to the main part of the argument, which can be decomposed into
four steps.
\med

{\bf Step~1:} {\em ``Energy'' estimate in $H^1_x$ and $H^3_v$ norms
combined}.
\sm

To avoid heavy notation, I shall use symbolic matrix
notation which should be rather self-explanatory, and write
\[ L = v\cdot\nabla_x - \nabla V(x)\cdot\nabla_v -\Delta_v - v\cdot\nabla_v.\]

By differentiating the equation once with respect to $x$, and three times with respect to $v$,
one finds
\begeq \label{1x}
\left(\derpar{}{t} +L\right) \nabla_x h = \nabla^2_x V(x)\cdot \nabla_v h,
\endeq
\begeq \label{3v}
\left(\derpar{}{t} + L\right) \nabla_v^3 h = -3 \nabla^2_v\nabla_x h -3 \nabla^3_v h.
\endeq
After taking the scalar product of~\eqref{1x} by $\nabla_x h$ and integrating against
$\mu$, we get
\begeq\label{1xi}
\frac12\: \frac{d}{dt}\int |\nabla_x h|^2\,d\mu + \int |\nabla_v\nabla_x h|^2\,d\mu
= \int(\nabla^2_x V)\nabla_v h \cdot \nabla_x h \,d\mu.
\endeq
Similarly, from~\eqref{3v} it follows that
\begeq\label{3vi}
\frac12\:\frac{d}{dt} \int |\nabla_v^3 h|^2 \,d\mu + \int |\nabla_v^4 h|^2\,d\mu
= -3 \int |\nabla_v^3 h|^2\,d\mu - 3 \int \nabla_v^3 h\cdot \nabla_v^2\nabla_x h\,d\mu.
\endeq

Let us bound the right-hand side of~\eqref{1xi}. Since $(\nabla_v)^*=-\nabla_v +v$
(where the $*$ is for the adjoint in $L^2(\mu)$), one has
\[ \int(\nabla^2_x V)\nabla_v h \cdot \nabla_x h \,d\mu
= - \int (\nabla^2_x V) h\cdot\nabla_x\nabla_v h\,d\mu 
- \int \bigl\<(\nabla^2_x V) h v, \nabla_xh\bigr\>\,d\mu. \] 
By Cauchy--Schwarz and Young's inequality,
\[ - \int (\nabla^2_x V) h\cdot\nabla_x\nabla_v h\,d\mu 
\leq \int|\nabla^2_x V|^2 h^2\,d\mu + \frac14 \int |\nabla_x\nabla_v h|^2\,d\mu.\]
Thanks to~\eqref{D2Vbdd1}, this can be bounded by
\[ C\left(\int |\nabla_xh|^2\,d\mu +
\int h^2\,d\mu\right) + \frac14 \int |\nabla_x\nabla_v h|^2\,d\mu.\]

By Cauchy--Schwarz inequality again,
\begeq\label{interm11} - \int \bigl\<(\nabla^2_x V) h v,\, \nabla_xh\bigr\>\,d\mu
\leq \sqrt{\int |\nabla^2_x V|^2 h^2\,d\mu}\sqrt{\int |v|^2 |\nabla_x h|^2\,d\mu}.
\endeq
In view of~\eqref{D2Vbdd2},
\begeq\label{interm22} \frac12 \int |v|^2 |\nabla_x h|^2\,d\mu
\leq C \left(\int |\nabla_xh|^2\,d\mu + \int |\nabla_v\nabla_xh|^2\,d\mu\right).
\endeq
By~\eqref{interm11}, \eqref{interm22} and Young's inequality,
there is a constant $C$ such that
\begin{multline*} - \int \bigl\<(\nabla^2_x V) h v,\, \nabla_xh\bigr\>\,d\mu
\leq C \left( \int |\nabla^2_x V|^2 h^2\,d\mu 
+ \int |\nabla_xh|^2\,d\mu \right) \\
+ \frac14 \int |\nabla_v\nabla_xh|^2\,d\mu.
\end{multline*}

All in all,
\begeq\label{step1.1}
\frac12\: \frac{d}{dt}\int |\nabla_x h|^2\,d\mu + \frac12\int |\nabla_v\nabla_x h|^2\,d\mu
\leq C \left(\int |\nabla^2_x V|^2\,h^2\,d\mu + \int |\nabla_x h|^2\,d\mu\right).
\endeq

The right-hand side in~\eqref{3vi} is estimated in a similar way:
\[ -3\int \nabla_v^3 h\cdot \nabla_v^2\nabla_x h\,d\mu
= 3\int \nabla_v^4 h\cdot \nabla_v\nabla_x h\,d\mu
- 3 \int \nabla_v^3 h \cdot v \nabla_v\nabla_x h\,d\mu;\]
then on one hand
\[ 3\int \nabla_v^4 h\cdot \nabla_v\nabla_x h\,d\mu
\leq \frac14 \int |\nabla_v^4 h|^2\,d\mu + 9 \int |\nabla_v\nabla_x h|^2\,d\mu;\]
on the other hand, again by~\eqref{D2Vbdd2},
\begin{align*} -3 \int \nabla_v^3 h \cdot v \nabla_v \nabla_x h\,d\mu 
& \leq 3 \sqrt{\int |v|^2 |\nabla_v^3h|^2\,d\mu}\
\sqrt{\int |\nabla_v\nabla_x h|^2\,d\mu} \\
& \leq C\sqrt{\int |\nabla_v^4h|^2\,d\mu + \int |\nabla_v^3h|^2\,d\mu}\
\sqrt{\int |\nabla_v\nabla_x h|^2\,d\mu}\\
& \leq \frac14 \int |\nabla_v^4h|^2\,d\mu + \frac14 \int |\nabla_v^3h|^2\,d\mu
+ C \int |\nabla_v\nabla_x h|^2\,d\mu.
\end{align*}
So there is a constant $C$ such that
\begin{multline}\label{step1.2}
\frac12\:\frac{d}{dt} \int |\nabla_v^3 h|^2 \,d\mu + \frac12 \int |\nabla_v^4 h|^2\,d\mu
\leq C\left( \int |\nabla_v\nabla_x h|^2\,d\mu + \int |\nabla_v^3 h|^2\,d\mu \right.\\
\left.+ \int |\nabla_v^4 h|^2\,d\mu\right).
\end{multline}

As a consequence of~\eqref{step1.1} and~\eqref{step1.2}
it is possible to find numerical constants $a,K,C>0$ (only depending on
$n$ and $C$ in~\eqref{hypV}) such that
\begin{multline}\label{step1} 
\frac{d}{dt} \left( \int |\nabla_xh|^2\,d\mu + a \int |\nabla^3_vh|^2\,d\mu\right)
\leq -K \left( \int |\nabla_v^4h|^2\,d\mu + \int |\nabla_v\nabla_x h|^2\,d\mu\right) \\
+ C \left(\int h^2\,d\mu + \int |\nabla_xh|^2\,d\mu + \int |\nabla_v^3h|^2\,d\mu\right).
\end{multline}
This concludes the first step.
\med

{\bf Step~2:} {\em Time-behavior of the mixed derivative.}

In this step I~shall focus on the mixed derivative integral
$\int \nabla_x h\cdot\nabla_v h\,d\mu$. By differentiating the equation with respect
to $x$ and multiply by $\nabla_vh$, differentiating the equation with respect to $v$
and multiply by $\nabla_xh$, then using the chain rule and the identity
$F\Delta_v G + G\Delta_v F = \Delta_v (FG) - 2 \nabla_vF\cdot\nabla_v G$,
one easily obtains
\begin{multline*} \left( \derpar{}{t} + L\right) (\nabla_x h\cdot\nabla_v h) = 
\ip{\nabla^2_x V\cdot \nabla_vh}{\nabla_vh} \\- 2 \nabla_v\nabla_x h\cdot \nabla^2_vh\,d\mu
- |\nabla_xh|^2 - \nabla_xh\cdot\nabla_vh.
\end{multline*}
After integration against $\mu$, this yields
\begin{multline}\label{afteripmu} \frac12\:\frac{d}{dt} \int \nabla_xh\cdot\nabla_vh\,d\mu
= \int \ip{\nabla^2_x V\cdot \nabla_vh}{\nabla_vh}\,d\mu
- 2 \int\nabla_v\nabla_x h\cdot \nabla^2_vh\,d\mu \\
- \int |\nabla_xh|^2\,d\mu - \int \nabla_x h\cdot\nabla_vh\,d\mu.
\end{multline}
The first term in the right-hand side need some rewriting:
Since $(\nabla_v)^*=-\nabla_v +v$,
\[ \int \ip{\nabla^2_x V\cdot \nabla_vh}{\nabla_vh}\,d\mu
= - \int \nabla^2_x V h \nabla^2_v h\,d\mu - \int h \ip{\nabla^2_x V v}{\nabla_v h}\,d\mu\]
\[ \leq
\sqrt{\int |\nabla^2_x V|^2 h^2\,d\mu}\ \sqrt{\int |\nabla^2_vh|^2\,d\mu}
+ \sqrt{\int |\nabla^2_x V|^2h^2\,d\mu}\ \sqrt{\int |v|^2 |\nabla_vh|^2\,d\mu}. \]
With the help of Young's inequality and~\eqref{D2Vbdd2} again, this can be bounded by
\[ \var \int |\nabla^2_x V|^2 h^2\,d\mu + C_\var \left(
+ \int |\nabla_v^2h|^2\,d\mu + \int |\nabla_v h|^2\,d\mu\right).\]
By Lemma~\ref{lemD2Vbdd}, if $\var$ is small enough then this is bounded by
\[ \frac14 \left( \int |\nabla_x h|^2\,d\mu + \int h^2\,d\mu\right)
+ C \left(\int |\nabla_v h|^2\,d\mu + \int |\nabla_v^2h|^2\,d\mu \right).\]

Now for the second term in the right-hand side of~\eqref{afteripmu},
we just write
\[ - 2 \int\nabla_v\nabla_x h \nabla^2_vh\,d\mu
\leq \int |\nabla_v \nabla_x h|^2\,d\mu + \int |\nabla_v^2 h|^2\,d\mu.\]

Summarizing all the above computations: There is a numerical constant $C$,
only depending on $n$ and $C$ in~\eqref{hypV}, such that
\begin{multline}
\label{step2} \frac{d}{dt} \int \nabla_x h\cdot\nabla_v h \,d\mu
\leq -\frac12 \int |\nabla_x h|^2\,d\mu \\ + C 
\left( \int h^2\,d\mu + \int |\nabla_v h|^2\,d\mu
+ \int |\nabla_v^2 h|^2\,d\mu + \int |\nabla_x \nabla_v h|^2\,d\mu\right).
\end{multline}
This concludes the second step of the proof.

\begin{Rk} We could also have conducted the computations in the following way:
\[ - 2 \int\nabla_v\nabla_x h \nabla^2_vh\,d\mu
= 2 \int \nabla_xh \cdot\nabla^3_vh\,d\mu
- 2\int \nabla_x h \cdot v \nabla^2_v h\,d\mu.\]
Then on one hand,
\[ 2 \int \nabla_xh \cdot\nabla^3_vh\,d\mu
\leq \frac14 \int |\nabla_x h|^2\,d\mu + 4 \int |\nabla^3_vh|^2\,d\mu;\]
on the other hand, just as before,
\begin{align*} - 2\int \nabla_x h \cdot v \nabla^2_v h\,d\mu
& \leq \frac14 \int |\nabla_xh|^2\,d\mu + 4 \int |v|^2|\nabla^2_vh|^2\,d\mu \\
& \leq \frac14 \int |\nabla_xh|^2\,d\mu
+ C \left( \int |\nabla^2_v h|^2\,d\mu + \int |\nabla^3_vh|^2\,d\mu\right).
\end{align*}
By doing so, we would have obtained the same result as~\eqref{step2},
except that the integral $\int |\nabla_x\nabla_v h|^2\,d\mu$
would be replaced by $\int |\nabla_v^3 h|^2\,d\mu$. Then the
rest of the proof would have worked through.
\end{Rk}
\med

{\bf Step~3:} {\em Interpolation inequalities}
\sm

If $h$ is a function of $v$, lying in $L^2(\gamma)$, one can write
$h=\sum_k c_k H_k$, where $H_k$ are normalized Hermite polynomials and $k$
are multi-indices in $\N^n$; then
\begin{multline*} 
\int h^2\,d\gamma = \sum_k c_k^2, \qquad \int |\nabla_v h|^2\,d\mu
= \sum |k|^2 c_k^2, \qquad \\
\int |\nabla_v^2 h|^2\,d\mu = \sum |k|^4 c_k^2, \qquad
\text{etc.} 
\end{multline*}
(here $|k|^2 = \sum k_\ell^2$, $1\leq \ell\leq n$ and $|k|^4 = (|k|^2)^2$). 
Then, by H\"older's inequality (in the $k$ variable), 
one can prove interpolation inequalities such as
\[ \int |\nabla_v h|^2\,d\gamma \leq C \left(\int h^2\,d\gamma\right)^{2/3}
\left(\int |\nabla_v^3 h|^2\,d\gamma\right)^{1/3}. \]
Now if $h=h(x,v)$ is a function of both variables $x$ and $v$, one can apply
the previous inequality to $h(x,\cdot)$ for each $x$, then integrate
with respect to $e^{-V}(x)\,dx$, and apply H\"older's inequality in the $x$ variable,
to find
\[ \int |\nabla_v h|^2\,d\mu \leq C \left(\int h^2\,d\mu\right)^{2/3}
\left(\int |\nabla_v^3 h|^2\,d\mu\right)^{1/3}. \]
Similarly,
\[ \int |\nabla_v^j h|^2\,d\mu \leq C \left(\int h^2\,d\mu\right)^{1-(j/4)}
\left(\int |\nabla_v^4 h|^2\,d\mu\right)^{j/4}, \qquad 1\leq j\leq 3.\]
\med

{\bf Step~4:} {\em Conclusion}
\sm

Now we can turn to the proof of estimate~\eqref{estregFP}. Without loss of
generality, assume $\int h^2\,d\mu=1$ at time~0. Then, since this quantity
is nonincreasing with time, $\int h^2(t,\cdot)\,d\mu \leq 1$ for all $t\geq 0$. 
By combining the results of Steps~1, 2 and~3, we discover that the quantities
\begin{multline*} 
X:= \int |\nabla_x h|^2\,d\mu, \qquad Y_j:= \int |\nabla_v^j h|^2\,d\mu
\qquad (0\leq j\leq 4), \\ 
{\cal M}:= \int \nabla_x h \cdot\nabla_v h, 
\qquad W=\int |\nabla_x \nabla_v h|^2
\end{multline*}
viewed as functions of $t$, solve the system of differential inequalities
\begeq\label{system} \begin{cases} \dps\frac{d}{dt} (X+aY_3) \leq - K (Y_4+W) + C (1+X+Y_3);\\ \\
\dps \frac{d}{dt} {\cal M} \leq -K X + C (1+Y_1+Y_2+W) \\ \\
|{\cal M}| \leq \sqrt{X\,Y_1}; \qquad Y_1\leq C Y_2^{1/2}\leq C' Y_3^{1/3}\leq C'' Y_4^{1/4}
\end{cases} 
\endeq

It is a consequence of Lemma~\ref{lemsystem} in Appendix~\ref{toolbox}
that solutions of~\eqref{system} satisfy
\[ 0\leq t\leq 1 \Longrightarrow\qquad X(t) + Y(t)\leq \frac{A}{t^3} \]
for some computable constant $A$. As a consequence, for $0\leq t\leq 1$,
\[ \int |\nabla_x h|^2\,d\mu = O(t^{-3}), \qquad
\int |\nabla_v^3 h|^2\,d\mu = O(t^{-3}). \]
Then, by interpolation $\int |\nabla_v h|^2 = O(t^{-1})$, 
$\int |\nabla_v^2 h|^2 = O(t^{-2})$.
This concludes the proof of~\eqref{estregFP}.
\end{proof}

\subsection{Variants} \label{subvariants}

Here I~studied the regularization effect by means of a
{\em system} of differential inequalities.
It is natural to ask whether one can do the same with
just one differential inequality. The answer is affirmative:
It is possible to use a trick similar to the one in
the proof of Theorem~\ref{thmsimple},
that is, add a carefully chosen lower-order term which is derived
from the mixed derivative $\int \nabla_x h\cdot\nabla_v h$.

A first possibility is to consider the Lyapunov functional
\begin{multline*} {\cal E}(h) = 
\int h^2\,d\mu + a \int |\nabla_x h|^2\,d\mu
+ 2b \int \nabla_x (D_x^{1/3}h)\cdot \nabla_v (D_x^{1/3}h)\,d\mu \\
+ c \int |\nabla_v^3 h|^2\,d\mu,
\end{multline*}
where $D_x=(-\Delta_x)^{1/2}$. Then by using computations similar
to the ones in Subsection~\ref{weigL2H1}, plus estimates on the commutator
$[D_x^{1/3},\nabla V]$, one can establish the following
a priori estimate along the Fokker--Planck equation:
For well-chosen positive constants $a,b,c$,
\[ \frac{d}{dt}\, {\cal E}(h) \leq -K {\cal E}(h)^{4/3},\qquad
h=e^{-tL}h_0.\]
The desired result follows immediately.

One drawback of this method is the introduction of fractional derivatives.
There is a nice variant due to H\'erau~\cite{herau:FP} in which one
avoids this by using powers of $t$:
\[ {\cal F}(t,h) = \int h^2\,d\mu + a t \int |\nabla_v h|^2\,d\mu
+ 2 b t^2 \int \nabla_vh\cdot\nabla_x h + c t^3 \int |\nabla_x h|^2\,d\mu.\]
Then one can estimate the time-derivative of ${\cal F}(t,e^{-tL}h_0)$
by means of computations similar to those in Subsection~\ref{weigL2H1},
and the inequalities
\[ t\, \left | \int (\nabla_v h\cdot\nabla_x h)\,d\mu \right|\leq
C \int |\nabla_v h|^2\,d\mu + \var t^2 \int |\nabla_xh|^2\,d\mu;\]
\[ t^2\, \left| \int \nabla_x \nabla_v h\cdot\nabla_v^2 h\,d\mu \right|
\leq C t\, \int |\nabla_v^2 h|^2\,d\mu + \var t^3\,
\int |\nabla_x \nabla_vh|^2\,d\mu.\]
In the end, if $a,b,c$ are well-chosen, one obtains,
with the shorthand $h=e^{-tL}h_0$,
\begin{multline*} \frac{d}{dt} {\cal F}(t,h) \leq
-K \left( \int |\nabla_v h|^2\,d\mu + 
t \int |\nabla_v^2 h|^2\,d\mu\right. + t^2 \int |\nabla_x h|^2\,d\mu
\\\left. + t^3 \int |\nabla_x \nabla_v h|^2\,d\mu\right).
\end{multline*}
It follows that ${\cal F}(t,h)$ is nonincreasing, and
therefore
\[ \int |\nabla_v h|^2\,d\mu = O(t^{-1}),\qquad
\int |\nabla_x h|^2\,d\mu = O(t^{-3}).\]
The conclusion is not so strong as the one we had before, since
we only have estimates on the first-order derivative in $v$.
But the exponents are again optimal, and it is possible to
adapt the method and recover estimates on higher-order derivatives.
Furthermore, estimates on $\int |\nabla_x h|^2$ and $\int |\nabla_v h|^2$
are exactly what is needed for Theorem~\ref{thmFPL2} to apply.

H\'erau's method lends itself very well to an abstract treatment.
For instance, let us consider an abstract operator $L=A^*A+B$,
satisfying Assumptions~(i)--(iii) in Theorem~\ref{thmsimple},
then the following decay rates (in general optimal) can be proven, 
at least formally:
\begeq\label{AetLh} 
\|A e^{-tL}h\| = O(t^{-1/2});\qquad \|C e^{-tL}h\| = O(t^{-3/2}).
\endeq
To show this, introduce
\[ {\cal F}(t,h):= \|e^{-tL}h\|^2 + at \|Ae^{-tL}h\|^2 + 2bt^2 \<Ae^{-tL}h,
Ce^{-tL}h\> + ct^3 \|Ce^{-tL}h\|^2.\]
Then, we can perform computations similar to the ones in 
Subsection~\ref{proofsimple}, except that now there are extra terms
coming from the time-dependence of the coefficients $a,b,c$.
Writing $h$ for $e^{-tL}h$, we have, if $a,b/a,c/b,c^2/b,b^2/ac$ are
small enough:
\begin{multline} \label{newlineregul}
\frac{d{\cal F}(t,h)}{dt} \leq
- \kappa \Bigl( \|Ah\|^2 + at \|A^2h\|^2 + bt^2 \|Ch\|^2 + ct^3 \|CAh\|^2\Bigr)
\\ + a \|Ah\|^2 + 4 bt \<Ah, Ch\> + 3 ct^2 \|Ch\|^2,
\end{multline}
were $\kappa$ is a positive number.
When $0\leq t\leq 1$, the positive terms in the right-hand side 
of~\eqref{newlineregul} can all be controlled by the
negative terms if $a$, $b$ and $c/b$ are small enough.
Then
\[ \frac{d{\cal F}(t,h)}{dt} \leq - K (\|Ah\|^2 + \|A^2h\|^2 + \|Ch\|^2).\]
In particular, ${\cal F}$ is a nonincreasing function of $t$, and
then the desired bounds $\|Ah\|^2= O(t^{-1})$, $\|Ch\|^2=O(t^{-3})$ follow
(as well as the bound $\|A^2h\|^2= O(t^{-2})$).

The very same scheme of proof allows to establish
a regularization theorem similar to Theorem~\ref{hypocomult}:

\begin{Thm}\label{hypoellmult}
Let $\H$ be a Hilbert space, $A:\H\to\H^n$ and $B:\H\to \H$ be unbounded
operators, $B^*=-B$, let $L:=A^*\! A+B$. Assume the existence of $\Nc\in\N$ and
(possibly unbounded) operators $C_0, C_1, \ldots, 
C_{\Nc+1}$, $R_1,\ldots,R_{\Nc+1}$
and $Z_1,\ldots,Z_{\Nc+1}$ such that
\[ C_0=A, \qquad [C_j,B] = Z_{j+1}C_{j+1} + R_{j+1}\quad (0\leq j\leq \Nc), \qquad C_{\Nc+1}=0,\]
and, for all $k\in \{0,\ldots, \Nc\}$,
\sm

(i) $[A,C_k]$ is bounded relatively to $\{C_j\}_{0\leq j\leq k}$ 
and $\{C_jA\}_{0\leq j\leq k-1}$;
\sm

(ii) $[C_k,A^*]$ is bounded relatively to $I$ and $\{C_j\}_{0\leq j\leq k}$;
\sm

(iii) $R_k$ is bounded relatively to $\{C_j\}_{0\leq j\leq k-1}$ and $\{C_jA\}_{0\leq j\leq k-1}$.
\sm

(iv) There are positive constants $\lambda_j$, $\Lambda_j$ such that
$\lambda_j I \leq Z_j \leq \Lambda_j I$.
\sm

Then the following bound holds true along the semigroup $e^{-tL}$:
\[ \forall k\in \{0,\ldots,\Nc\}, \quad \|e^{-tL}h\|
\leq C\, \frac{\|h\|}{t^{k+\frac12}},\]
where $C$ is a constant only depending on the constants appearing implicitly
in Assumptions (i)--(iv).
\end{Thm}

\begin{Rk} A reasoning similar to Remark~\ref{H1H3interp} shows that
the exponents $1/(k+1/2)$ cannot be improved. Indeed, in H\"ormander's
theory, the weight attributed to the commutator $C_k$ would be $2k+1$,
and the regularity estimates established by 
Rothschild and Stein~\cite{rothschildstein:nilp:76},
which are optimal in general, provide regularization by an order
$2/(2k+1)=1/(k+1/2)$.
\end{Rk}

\begin{Rk} I~shall show below how H\'erau's method can be adapted 
to yield regularization from $L\log L$ initial datum. On the other hand, 
it is not clear that it can be used to establish regularization from measure 
initial data.
\end{Rk}

\subsection{Higher regularity from measure initial data} \label{measureSob}

Now I~shall explain how to extend the previous results by
(a) establishing Sobolev regularity of higher order, (b) removing
the assumption of $L^2$ integrability for the initial datum.

I~shall only consider the case when $\nabla V$ is Lipschitz and
has all its derivatives uniformly bounded. There are three motivations
for these restrictions: (i) even if they are far from optimal, they
will simplify the presentation quite a bit; (ii) they ensure the
uniqueness of the solution of the Fokker--Planck equation starting
from a measure initial datum; (iii) the theorems of convergence
to equilibrium studied in the present paper use the Lipschitz
regularity of $\nabla V$ anyway.

As before, the equation under study is
\begeq\label{FPagain}
\partial_t f + v\cdot\nabla_x f - \nabla V(x)\cdot\nabla_v f
= \Delta_v f + v\cdot\nabla_v f + nf.
\endeq
This equation admits a unique solution as soon as $f_0$ is a
finite nonnegative measure (say a probability measure) with
finite energy, and it is easy to prove the propagation of regularity
and of moment bounds.

So to establish regularization in higher-order Sobolev space
$H^k_x(H^\ell_v)$ it is enough to prove, for smooth and
rapidly decaying solutions, an a priori estimate like
\[ \|f_t\|_{H^k_xH^\ell_v(\R^n_x\times\R^n_v)} \leq 
\frac{C}{t^{-{1/\kappa}}},\]
with constants $C$ and $\kappa$ that do not depend on the regularity
of $f_0$.

In the sequel, $C$ and $K$ will stand for various constants depending
only on $n$ and $V$. As in Subsection~\ref{weigL2H1} the a priori
estimate is divided into four steps. The conservation of mass
(that is, the preservation of $\int f\,dx\,dv$) along equation~\eqref{FPagain}
will be used several times.
\med

{\bf Step~1:} {\em ``Energy'' estimate in higher order Sobolev spaces.}

Let $k$ and $\ell$ be given integers ($k$ will be the regularity
in $x$ and $\ell$ the regularity in $v$).
Computations similar to those in Subsection~\ref{weigL2H1}
(differentiating the equation and integrating) yield
\begin{multline*}
\frac{d}{dt}
\int |\nabla_x^k \nabla_v^\ell f|^2\,dx\,dv \ \leq
-K \int |\nabla_x^k \nabla_v^{\ell+1} f|^2\,dx\,dv
+ C \int |\nabla_x^k \nabla_v^\ell f|^2\,dx\,dv \\
+ C \int |\nabla_x^{k+1}\nabla_v^{\ell-2} f|^2\,dx\,dv 
+ C \sum_{1\leq i\leq k} \int |\nabla_x^{k-i}\nabla_v^\ell f|^2
|\nabla_x^{i+1}V|^2\,dx\,dv.
\end{multline*}
By assumption $|\nabla_x^iV|$ is bounded for any $i$, so
the above equation reduces to
\begin{multline*}
\frac{d}{dt} \int |\nabla_x^k \nabla_v^\ell f|^2
\leq -K \int |\nabla_x^k \nabla_v^{\ell+1} f|^2 + 
C \sum_{j\leq k} \int |\nabla_x^j \nabla_v^\ell f|^2 \\
+ C \int |\nabla_x^{k+1} \nabla_v^{\ell-2} f|^2.
\end{multline*}

Then one can repeat the computation with $(k,\ell)$ replaced
by $(k+1,\ell-2)$ and then $(k+2,\ell-4)$, etc. By an easy induction,
for a given integer $m$, we can find positive constants 
$K,C, a_0=1,a_1,\ldots,a_m$ such that
\begin{multline*} \frac{d}{dt} \sum_{k=0}^m a_k
\int |\nabla_x^k \nabla_v^{3(m-k)} f|^2\,dx\,dv\leq
- K \sum_{k=0}^m \int |\nabla_x^k \nabla_v^{3(m-k)+1} f|^2\,dx\,dv \\
+ C \sum_{k=0}^m \sum_{j\leq k} \int |\nabla_x^j \nabla_v^{3(m-k)}f|^2\,dx\,dv.
\end{multline*}

Repeating the same operation for lower order terms (that is,
decreasing $m$), for each couple of nonnegative integers 
$(k,\ell)$ with $3k+\ell\leq 3m$ we can find a positive constant 
$a_{k,\ell}$ such that
\begin{multline*} \frac{d}{dt}
\sum_{3k+\ell\leq 3m} a_{k,\ell} \int |\nabla_x^k \nabla_v^\ell f|^2\,dx\,dv
\leq -K \int |\nabla_v^{3m+1} f|^2\,dx\,dv  \\
+ C \sum_{3k+\ell\leq 3m} \int |\nabla_x^k\nabla_v^\ell f|^2\,dx\,dv.
\end{multline*}
Then we can define an ``energy functional'' of order $m$, which
controls the $L^2$-regularity of $f$ up to order $m$ in $x$ and
$3m$ in $v$:
\begeq\label{Em}
{\cal E}_m(f) = \sum_{3k+\ell\leq 3m} a_{k,\ell}
\int |\nabla_x^k \nabla_v^\ell f|^2\,dx\,dv.
\endeq
(Recall, to avoid any confusion, that
$\int |\nabla_x^k \nabla_v^\ell f|^2$ is the sum of all terms
$\int (\pa_{x_1}^{k_1}\ldots\pa_{x_n}^{k_n}\pa_{v_1}^{\ell_1}\ldots
\pa_{v_n}^{\ell_n} f)^2$ with $k_1+\ldots+k_n =k$, 
$\ell_1+\ldots+\ell_n=\ell$.)

Then the a priori estimate on the Fokker--Planck equation~\eqref{FPagain} 
can be recast as
\begeq\label{dEmdt}
\frac{d}{dt} {\cal E}_m(f) \leq -K \int |\nabla_v^{3m+1}f|^2\,dx\,dv
+ C {\cal E}_m(f).
\endeq

The important terms in ${\cal E}_m$ are the extreme ones, that is
for $(k,\ell)=(m,0)$, $(0,3m)$ or $(0,0)$. All the other ones can
be controlled by these three extremal terms; to see this, it suffices
to apply H\"older's inequality in Fourier space: Denoting by
$\xi$ the conjugate variable to $x$ and by $\eta$ the conjugate variable
to $v$, one has
\begin{align*}
& \int |\nabla_x^k\nabla_v^\ell f|^2\,dx\,dv 
= C \int |\xi|^{2k} |\eta|^{2\ell}|\hat{f}|\,d\xi\,d\eta\\
& \qquad 
\leq C \left(\int |\xi|^{2m}|\hat{f}|^2\,d\xi\,d\eta\right)^{\frac{k}{m}}
\left(\int |\eta|^{6m}|\hat{f}|^2\,d\xi\,d\eta\right)^{\frac{\ell}{3m}}
\left(\int |\hat{f}|^2\,d\xi\,d\eta\right)^{1-\left(\frac{k}{m} + 
\frac{\ell}{3m}\right)}\\
& \qquad = C \left(\int |\nabla_x^m f|^2\,dx\,dv\right)^{\frac{k}{m}}
\left(\int |\nabla_v^{3m} f|^2\,dx\,dv\right)^{\frac{\ell}{3m}}
\left(\int f^2\,dx\,dv\right)^{{1-\left(\frac{k}{m} + 
\frac{\ell}{3m}\right)}}.
\end{align*}

It follows easily that there are positive constants $K,C$ such that
\begin{multline}\label{KEC} 
K \left( \int |\nabla_x^m f|^2 + \int |\nabla_v^{3m} f|^2 + \int f^2\right)
\leq {\cal E}_m(f) \\ \leq
C \left( \int |\nabla_x^m f|^2 + \int |\nabla_v^{3m} f|^2 + \int f^2\right).
\end{multline}
\med

{\bf Step~2:} {\em Mixed derivatives}

Now define the higher order mixed derivative functional
\begeq\label{Mm}
{\cal M}_m(f) = \int \nabla_x^m f\cdot\nabla_x^{m-1}\nabla_v f
= \sum_{1\leq i_1,\ldots, i_m \leq n} \int \frac{\pa^mf}{\pa x_{i_1} \ldots
\pa_{x_{i_m}}}\, \frac{\pa^mf}{\pa x_{i_1} \ldots
\pa_{x_{i_{m-1}}} \pa_{v_{i_m}}}.
\endeq

By computations in the same style as in Step~2 of Subsection~\ref{weigL2H1},
one can establish
\[ \frac{d}{dt} {\cal M}_m(f) \leq
-K \int |\nabla_x^m f|^2\,dx\,dv 
+ C \sum_{k<m,\ 3k+\ell\leq 3m} \int |\nabla_x^k\nabla_v^\ell f|^2.\]
Each of the terms appearing in the latter sum can then be estimated
by elementary interpolation inequalities as in Step~1: If $k<m$ then
\[ \int |\nabla_x^k \nabla_v^\ell f|^2 \leq
\var \int |\nabla_x^m f|^2 + C \left(
\int |\nabla_v^{3m} f|^2 + \int f^2\right),\]
where $\var$ is an arbitrarily small positive number. The conclusion is that
\begeq\label{ddtMm}
\frac{d}{dt} {\cal M}_m(f) \leq -K \int |\nabla_x^m f|^2\,dx\,dv\
+ C \left( \int |\nabla_v^{3m} f|^2\,dx\,dv
+ \int f^2\,dx\,dv\right).
\endeq
\med

{\bf Step~3:} {\em Interpolation inequalities.}
\sm

There are two things to check: (i) that ${\cal M}_m$ is ``much smaller''
than ${\cal E}_m$, and (ii) that $\int |\nabla_v^{3m}f|^2$
is ``much smaller'' than $\int |\nabla_v^{3m+1}f|^2$. The difficulty
is that we cannot just use interpolation in $L^2$-type spaces.
In replacement, we shall use the anisotropic {\em Nash}-type interpolation
inequality exposed in Appendix~\ref{toolbox}.

First, by Cauchy--Schwarz,
\[ |{\cal M}_m(f)| \leq 
\left( \int |\nabla_x^m f|^2\right)^{\frac12}
\left( \int |\nabla_x^{m-1}\nabla_v f|^2\right)^{\frac12}.\]
Then the second term is estimated thanks to Lemma~\ref{nash}
with $\lambda=m-1$, $\mu=1$, $\lambda'=m$, $\mu'=3m$:
\[ \int |\nabla_x^{m-1}\nabla_v f|^2
\leq \left( \int |\nabla_x^m f|^2 + \int |\nabla_v^{3m}f|^2\right)^{1-\theta}
\left(\int f\right)^{2\theta},\]
where $\theta=2/(3m+6n)$ is a positive number. 
Since the mass $\int f$ is preserved under the time-evolution by
the Fokker--Planck equation, we arrive at the estimate
\begeq\label{Mdelta}
|{\cal M}_\delta(f)| \leq C {\cal E}_m(f)^{1-\delta},
\endeq
where $\delta=\theta/2$ is a positive constant.

Next, apply Lemma~\ref{nash} again with $\lambda=0$,
$\lambda'=m$, $\mu=3m$, $\mu'=3m+1$. Noting that
$(\lambda/\lambda') + (\mu/\mu')=3m/(3m+1)<1$, we see that there
exists $\theta\in (0,1)$ such that
\[
\int |\nabla_v^{3m}f|^2\,dx\,dv
\leq C \left( \int |\nabla_x^m f|^2 + \int |\nabla_v^{3m+1} f|^2\,dx\,dv
\right)^{1-\theta}\left(\int f\,dx\,dv \right)^{2\theta}.
\]
The same estimate holds true for $\int f^2$ (this can be treated
by the usual Nash inequality), and then one can use the fact that
$\int f$ is preserved by the Fokker--Planck equation, to obtain the
a priori estimate
\begeq\label{3m+1}
\int |\nabla_v^{3m}f|^2\,dx\,dv + \int f^2\,dx\,dv
\leq C \left( \int |\nabla_x^m f|^2 + \int |\nabla_v^{3m+1} f|^2\,dx\,dv
\right)^{1-\theta}
\endeq
\med

{\bf Step~4:} {\em Conclusion}

Equations~\eqref{KEC}, \eqref{Mdelta}, \eqref{dEmdt}, \eqref{3m+1}
and \eqref{ddtMm} together show that we can apply Lemma~\ref{lemsystem}
with ${\cal E}={\cal E}_m$, ${\cal M}={\cal M}_m$, $X=\int |\nabla_x^m f|^2$,
$Y=\int |\nabla_v^{3m} f|^2 + \int f^2$, $Z=\int |\nabla^{3m+1}f|^2$.
Thus there are constants $C$ and $\kappa$ such that
${\cal E}_m(f_t) \leq C/t^{1/\kappa}$. This concludes the proof
of the a priori estimate.

\subsection{Regularization in an $L\log L$ context}

If the initial datum is assumed to have finite entropy, then
H\'erau's method can be adapted to yield the regularization in
Fisher information sense, with exponents that are likely to be optimal.
Here is a rather general result in this direction, under the
same assumptions as Theorem~\ref{thmLlogL}:

\begin{Thm} \label{hypoellLog}
Let $E\in C^2(\R^N)$, such that $e^{-E}$ is rapidly decreasing,
and $\mu(dX) = e^{-E(X)}\,dX$ is a probability measure on $\R^N$. 
Let $(A_j)_{1\leq j\leq m}$ and $B$ be
first-order derivation operators with smooth coefficients.
Denote by $A_j^*$ and $B^*$ their respective adjoints in $L^2(\mu)$, 
and assume that $B^*=-B$. Denote by $A$ the collection $(A_1,\ldots,A_m)$, 
viewed as an unbounded operators whose range is made of functions valued in $\R^m$.
Define 
\[ L = A^*A + B = \sum_{j=1}^m A_j^*A_j + B,\]
and assume that $e^{-tL}$ defines a well-behaved semigroup on a
suitable space of positive functions
(for instance, $e^{-tL}h$ and $\log (e^{-tL}h)$ are $C^\infty$ and all
their derivatives grow at most polynomially if $h$ is itself $C^\infty$ 
with all derivatives bounded, and $h$ is bounded below by a positive
constant).

Next assume the existence of $\Nc\geq 1$,
derivation operators $C_0,\ldots,C_{\Nc+1}$ and $R_1,\ldots,R_{\Nc+1}$,
and vector-valued functions $Z_1,\ldots,Z_{\Nc+1}$ (all of them with
$C^\infty$ coefficients, growing at most polynomially, as their
partial derivatives) such that
\[ C_0 =A, \qquad [C_j,B] = Z_{j+1}\,C_{j+1} + R_{j+1}\quad
(0\leq j\leq \Nc),\qquad C_{\Nc+1}=0,\]
and

(i) $[A,C_k]$ is pointwise bounded relatively to $A$;

(ii) $[C_k,A^*]$ is pointwise bounded relatively to 
$I, \{C_j\}_{0\leq j\leq k}$;

(iii) $R_k$ is pointwise bounded with respect to $\{C_j\}_{0\leq j\leq k-1}$;

(iv) there are positive constants $\lambda_j,\Lambda_j$ such that
$\lambda_j \leq Z_j\leq \Lambda_j$;

(v) $[A,C_k]^*$ is pointwise bounded relatively to $I,A$.
\sm

Then the following bound holds true: With the notation $h(t) = e^{-tL}h_0$,
\[ \forall k\in \{0,\ldots,\Nc\}, \quad 
\int h(t) \bigl|C_k\log h(t)\bigr|^2\,d\mu \leq
C\: \frac{\dps \int h_0\log h_0\,d\mu}{t^{2k+1}},\]
where $C$ is a constant only depending on the constants appearing implicitly
in Assumptions (i)--(v).
\end{Thm}

\begin{proof}
The proof is patterned after the proofs of Theorems~\ref{thmLlogL}
and~\ref{hypoellmult}: Write $u=\log h$, $f=e^{-E}h$, and introduce the 
Lyapunov functional
\[ {\cal F}(t,h) = \int f u \ +\ \sum_{k=0}^{\Nc}
\Bigl( a_k t^{2k+1} \int f|C_ku|_m^2 \: + \:
2 b_k t^{2k+2} \int f\<C_ku, C_{k+1}u\>_m\Bigr).\]
The computations for $d{\cal F}/dt$ are the same as in the proof
of Theorem~\ref{thmLlogL}, except that now there are additional terms
caused by the explicit dependence on $t$. So
\begin{multline} \frac{d{\cal F}(t,h(t)}{dt} \leq \\ -K 
\left( \int f|Au|^2 + \sum_k a_k t^{2k+1}\int f|C_kAu|^2 
+ \sum_k b_k t^{2k+2} \int f|C_{k+1}u|^2\right)\\
+ \sum_k (2k+1)a_k t^{2k} \int f|C_ku|^2 + 
\sum_k (2k+2) b_k t^{2k+1} \int f\<C_ku, C_{k+1}u\>.
\end{multline}
Obviously, the additional terms can be controlled by the ones in
the first line of the right-hand side, provided that
$a_{k}/b_{k-1}$ and $b_kt^{2k+1}/\sqrt{(b_{k-1}t^{2k})(b_kt^{2k+1})}$
are small enough; the second condition reduces to $b_k/b_{k-1}$ small
enough. These conditions have been enforced in the proof of
Theorem~\ref{thmLlogL}. So all in all,
${\cal F}(t,h(t))$ is a nonincreasing function of $t$, and the conclusion
follows immediately.
\end{proof}

\section{Toolbox} \label{toolbox}

The following elementary lemma is used in the proofs 
of Theorems~\ref{thmsimple}, \ref{thmmultsimple} and~\ref{thmCjA}.

\begin{Lem}\label{ll} 
Let $\delta>0$ and $u_0>0$ be given. Then it is always possible
to choose positive numbers $u_1, u_2, \ldots, u_N$ in such a way that
\[ \begin{cases} \forall k\in \{0,\ldots, N-1\}, \qquad
u_{k+1} \leq \delta\, u_k; \\
\forall k\in \{1,\ldots, N-1\}, \qquad u_k^2 \leq \delta\, u_{k-1}\, u_{k+1}.
\end{cases}\]
\end{Lem}

\begin{proof} Without loss of generality, assume $u_0=1$.
Set $m_0=0, m_1=1$; by induction, it is possible to pick up positive
numbers $m_k$ such that
\[ m_{k+1} \in (m_k, 2 m_k - m_{k-1}).\]
The resulting sequence will be increasing and satisfy $m_k > (m_{k-1}+ m_{k+1})/2$. 
Next set $u_k = \var^{m_k}$; for $\var$ small enough, 
the desired inequalities are satisfied.
\end{proof}

The next lemma is used in the proof of Theorem~\ref{thmprecise};
it is a kind of nonlinear analogue of Lemma~\ref{ll}.

\begin{Lem} \label{lemchoice}
Let $K, \ov{E}, k>0, J$ be given. Then there exists constants
$\var_1=\var_1(J)>0$, $\ell=\ell(J,k)>0$ and $K_1=K_1(K,\ov{E},k,J)>0$
with the following property: For any $\var\in (0,\var_1)$
and $E\in (0,\ov{E})$, there exist coefficients $a_1,\ldots,a_{J-1}>0$
satisfying
\begeq\label{condcoeff}
\begin{cases} 
1= a_0 \geq a_1 \geq a_2\geq\ldots \geq a_{J-1};\\ \\
a_1 \leq K E^\var; \\
\forall j\in\{1,J-1\}, \quad \dps \frac{a_j^2}{a_{j-1}} \leq K\, 
a_{J-1}^{1+\var}\, E^{k\var};\\ \\
a_{J-1} \geq K_1\, E^{\ell \var}.
\end{cases}
\endeq
\end{Lem}

\begin{proof}[Proof of Lemma~\ref{lemchoice}]
Without loss of generality we may assume that $K$ is 
bounded above by $m:=\min(1,(\ov{E})^{-k})$;
otherwise, just replace $K$ by $m$.

We shall choose the coefficients $a_j$ in such a way that the
inequality in the third line of~\eqref{condcoeff} holds as an equality.
For $j=J-1$ this gives
\[ a_{J-1}^2 = K a_{J-1}^{1+\var} a_{J-2} E^{k\var},\]
hence
\[ a_{J-2} = a_{J-1} \, (K E^{k\var} a_{J-1}^{\var})^{-1}.\]
Then the equality
\[ \frac{a_j}{a_{J-1}} = \left(\frac{a_{j+1}}{a_{J-1}}\right)^2\,
\bigl(K E^{k\var} a_{J-1}^\var\bigr)^{-1}\]
yields, by decreasing induction,
\[ a_j = a_{J-1}\, (K E^{k\var} a_{J-1}^\var)^{-\alpha_j},\]
where $\alpha_j$ is defined by the (decreasing) recursion relation
\[ \alpha_{J-2} =1,\qquad \alpha_{j-1} = 2 \alpha_j +1.\]

The sequence $(a_j)_{1\leq j\leq J-1}$ so defined is
nonincreasing if $KE^{k\var} a_{J-1}\leq 1$. From the bound
$K\leq \min(1,(\ov{E})^{-k})$, we know that
$KE^{k\var}\leq 1$ as soon as $\var\leq \var_1\leq 1$
($\var_1$ to be chosen later), and $a_{J-1}\leq 1$.

Then $\alpha_1=2^{J-2}-1$ is a positive integer depending only on $J$, and
\[ a_1 = a_{J-1}\, (KE^{k\var} a_{J-1}^\var)^{-\alpha_1}=
a_{J-1}^{1-\alpha_1 \var} (K E^{k\var})^{-\alpha_1}.\]
If $\var\leq \var_1:=1/(2\alpha_1)$, then $a_{J-1}$ appears
in the right-hand side with a positive exponent $1-\alpha_1\var\,\in(1/2,1)$.
Also $\var_1\leq 1$, as assumed before.

To make sure that the first condition in~\eqref{condcoeff} is
fulfilled, we impose
\[ a_{J-1}^{1-\alpha_1 \var} (K E^{k\var})^{-\alpha_1}
= K E^{\var},\]
that is
\[ a_{J-1} = \bigl[ K^{1+\alpha_1} E^{(1+k\alpha_1)_\var}\bigr]
^{\frac1{1-\alpha_1}}\var.\]
Up to decreasing $K$ again, we may assume that the quantity
inside square brackets is bounded by~1; this also implies that
$a_{J-1}\leq 1$, as assumed before. Then, since
$1/(1-\alpha_1\var)\leq 1/(1-\alpha_1\var_1)=2$, one has
\[ a_{J-1} \geq \bigl[ K^{1+\alpha_1} E^{(1+k\alpha_1) \var}\bigr]^2,\]
and the lemma follows upon choosing
$K_1 = K^{2(1+\alpha_1)}$, $\ell=2(1+k\alpha_1)$.
\end{proof}

The next lemma, used to check~\eqref{condthmsimpleFP} in Section~\ref{secFP},
states that $|\nabla^2 V|$ defines a bounded operator $H^1(e^{-V})\to L^2(e^{-V})$
as soon as $|\nabla^2 V|$ is dominated by $|\nabla V|$.

\begin{Lem} \label{lemD2Vbdd}
Let $V$ be a $C^2$ function on $\R^n$, satisfying~\eqref{condD2V}.
Then, for all $g\in H^1(e^{-V})$,
\sm

(i) $\dps \int_{\R^n} |\nabla V|^2\, g^2\, e^{-V} \leq 8(1+\sqrt{n}C)^2
\left ( \int_{\R^n} g^2\,e^{-V} + \int_{\R^n}|\nabla g|^2\,e^{-V} \right);$
\sm

(ii) $\dps \int_{\R^n} |\nabla^2 V|^2\,g^2\,e^{-V} \leq 16\, C^2(1+\sqrt{2n}C)^2
\left ( \int_{\R^n} g^2\,e^{-V} + \int_{\R^n}|\nabla g|^2\,e^{-V} \right).$
\end{Lem}

%\begin{Rk} In the particular case when $V$ is quadratic, i.e. $e^{-V}$ is the 
%normalized centered Gaussian $\gamma$, then computations simplify into
%\begeq\label{lemD2Vbddg}
%\int_{\R^n} |x|^2 g^2(x)\,\gamma(dx) \leq 2 \max(n,2)
%\left( \int_{\R^n} g^2\,d\gamma + \int_{\R^n} |\nabla g|^2\,d\gamma \right). 
%\endeq
%\end{Rk}

\begin{proof}[Proof of Lemma~\ref{lemD2Vbdd}]
By a density argument, we may assume that $g$ is smooth and decays fast enough
at infinity. Then, by the identity $\nabla (e^{-V})=-(\nabla V)e^{-V}$ and
an integration by parts,
\begin{multline*} \int |\nabla V|^2 g^2\, e^{-V}= - \int g^2 \nabla V\cdot\nabla (e^{-V})
= - \int \nabla \cdot (g^2 \nabla V)\,e^{-V} \\ 
= - \int g^2 (\Delta V) e^{-V} - 2 \int g (\nabla g\cdot \nabla V)\, e^{-V}. 
\end{multline*}
By Cauchy--Schwarz inequality,
\begin{multline}\label{afterCS}
\int |\nabla V|^2 g^2\, e^{-V} \leq \sqrt{\int g^2 (\Delta V)^2\, e^{-V}}
\sqrt{\int g^2\,e^{-V}} \\ + 2 \sqrt{\int |\nabla V|^2 g^2\,e^{-V}}
\sqrt{\int |\nabla g|^2\,e^{-V}}.
\end{multline}
Since, by~\eqref{Poincx},
\[ (\Delta V)^2 \leq n |\nabla^2 V|^2 \leq n C^2(1+|\nabla V|)^2\leq 2nC^2 (1+|\nabla V|^2),\]
it follows from~\eqref{afterCS} that
\begin{multline*}
\int |\nabla V|^2 g^2\, e^{-V} \leq \sqrt{2n}C\ \sqrt{\int g^2\, e^{-V} 
+ \int |\nabla V|^2 g^2\, e^{-V}} \sqrt{\int g^2\,e^{-V}} \\
+ 2 \sqrt{\int |\nabla V|^2 g^2\,e^{-V}} \sqrt{\int |\nabla g|^2\,e^{-V}}
\end{multline*}
\begin{multline*} \leq \sqrt{2n}C \int g^2\, e^{-V} + \sqrt{2n}C 
\sqrt{\int |\nabla V|^2 g^2\, e^{-V}}
\sqrt{\int g^2\,e^{-V}} \\ + 2 \sqrt{\int |\nabla V|^2 g^2\,e^{-V}}
\sqrt{\int |\nabla g|^2\,e^{-V}}.
\end{multline*}
Thanks to Young's inequality, this can be bounded by
\begin{multline*} 
\sqrt{2n}C\int g^2\,e^{-V} + \left( \frac14 \int |\nabla V|^2 g^2\,e^{-V} 
+ 2nC^2 \int g^2\,e^{-V}\right) \\
+ \left( \frac14 \int |\nabla V|^2 g^2\,e^{-V} + 4 \int |\nabla g|^2\, e^{-V} \right).
\end{multline*}
All in all,
\begin{multline*}
\int |\nabla V|^2 g^2\, e^{-V} \leq \frac12 \int |\nabla V|^2 g^2\, e^{-V}
+ (\sqrt{2n}C + 2nC^2) \int g^2\, e^{-V} \\
+ 4 \int |\nabla g|^2 \,e^{-V},
\end{multline*}
so
\begeq\label{intermnablaV}
\int |\nabla V|^2 g^2\, e^{-V} \leq 
2 (\sqrt{2n}C + 2nC^2) \int g^2\, e^{-V} + 8 \int |\nabla g|^2 \,e^{-V},
\endeq
This easily leads to statement (i) after crude upper bounds.

To prove statement (ii), start again from~\eqref{intermnablaV} and
apply~\eqref{condD2V} again, in the form $|\nabla^2V|^2 \leq 2 C(1+|\nabla V|^2)$:
the desired conclusion follows at once.
\end{proof}

Next, we shall study an interpolation inequality ``in Nash style''.
First recall the classical Nash inequality~\cite{nash:58} in $\R^n_x$:
If $f$ is a nonnegative function of $x\in\R^n$, then
\[ \int_{\R^n} f^2\,dx \leq C(n)\, \left( \int_{\R^n} |\nabla_x f|^2\,dx\right)^{1-\theta}
\left(\int_{\R^n} f\,dx\right)^{2\theta},\]
where
\[ \theta = \frac{2}{n+2}.\]
It is easy to generalize this inequality for higher order,
or fractional derivatives: If $D=(-\Delta)^{1/2}$, and $0\leq \lambda<\lambda'$, then
\[ \int_{\R^n} |D_x^\lambda f|^2\,dx \leq C(n,\lambda,\lambda')\,
\left( \int_{\R^n} |D_x^{\lambda'} f|^2\,dx\right)^{1-\theta}
\left(\int_{\R^n} f\,dx\right)^{2\theta},\]
where now
\[ \theta = \frac{2(\lambda'-\lambda)}{n+2\lambda'}.\]

The next lemma generalizes this to functions which depend on two variables,
$x$ and $v$, and allows different orders of derivations in these variables.
The symbol $D$ will again stand for $(-\Delta)^{1/2}$.

\begin{Lem} \label{nash}
Let $f=f(x,v)$ be a nonnegative (smooth, rapidly decaying) function on $\R^n_x\times\R^n_v$.
Let $\lambda,\lambda',\mu,\mu'$ be four nonnegative numbers with
$\lambda',\mu'>0$. If
\[ \frac{\lambda}{\lambda'} + \frac{\mu}{\mu'} <1,\]
then there is a constant $C=C(n,\lambda,\mu,\lambda',\mu')$ such that
\begeq\label{lmlm'}
\int |D_x^\lambda D_v^\mu f|^2\,dx\,dv \leq
C \left( \int |D_x^{\lambda'} f|^2\,dx\,dv + 
\int |D_v^{\mu'}f|^2\,dx\,dv \right)^{1-\theta} \left(\int f\right)^{2\theta},
\endeq
where
\[ \theta = \frac{1-\left(\frac{\lambda}{\lambda'} + \frac{\mu}{\mu'}\right)}
{1 + \frac{n}{2} \left( \frac1{\lambda'} + \frac1{\mu'}\right)}.\]
\end{Lem}

\begin{proof}[Proof of Lemma~\ref{nash}]
The strategy here will be the same as in the classical proof 
(actually due to Stein) of Nash's inequality: Go to Fourier space and
separate according to high and low frequencies, then optimize.
I~shall denote by $\hat{f}$ the Fourier transform of $f$,
by $\xi$ the Fourier variable that is dual to $x$,
and by $\eta$ the variable that is dual to $v$.
So the inequality to prove is
\begeq\label{toprove} \int |\xi|^{2\lambda} |\eta|^{2\mu} |\hat{f}|^2\,d\xi\,d\eta \leq
C \left( \int |\xi|^{2\lambda'} |\hat{f}|^2\,d\xi\,d\eta +
\int |\eta|^{2\mu'}|\hat{f}|^2\,d\xi\,d\eta \right)^{1-\theta} 
\|\hat{f}\|_{L^\infty}^{2\theta}.
\endeq

First start with the case $\lambda=0$, and separate the integral in the
left-hand side of~\eqref{toprove} in three parts:
\begin{multline*} \int (\ldots)\,d\xi\,d\eta \ = \
\int_{|\xi|\leq R, \ |\eta|\leq S} (\ldots)\,d\xi\,d\eta +
\int_{|\xi|>R,\ |\eta|\leq S} (\ldots)\,d\xi\,d\eta \\
+ \int_{|\xi|>R,\ |\eta|>S} (\ldots)\,d\xi\,d\eta,
\end{multline*}
where $R$ and $S$ are positive numbers that will be chosen later on.

Then,
\begin{align} \nonumber
\int_{|\xi|\leq R,\ |\eta|\leq S}
|\eta|^{2\mu} |\hat{f}(\xi,\eta)|^2\,d\eta\,d\xi 
& \leq S^{2\mu} \, \vol (|\xi|\leq R)\, \vol (|\eta|\leq S)\,
\|\hat{f}\|_{L^\infty}^2\\
& \leq C_n R^n S^{2\mu+n}\|\hat{f}\|_{L^\infty}^2,
\label{firstintxv}
\end{align}
where $C_n$ only depends on $n$, and $\vol$ is a notation for the Lebesgue
volume in $\R^n$.

Next
\begeq\label{secondintxv}
\int_{|\xi|> R,\ |\eta|\leq S}
|\eta|^{2\mu} |\hat{f}(\xi,\eta)|^2\,d\xi\,d\eta
\leq \frac{S^{2\mu}}{R^{2\lambda'}}
\int |\xi|^{2\lambda'} |\hat{f}(\xi,\eta)|^2\,d\xi\,d\eta.
\endeq

Finally,
\begeq\label{thirdintxv}
\int_{|\xi|> R,\ |\eta|> S}
|\eta|^{2\mu} |\hat{f}(\xi,\eta)|^2\,d\xi\,d\eta
\leq \frac1{S^{2(\mu'-\mu)}}
\int |\eta|^{2\mu'} |\hat{f}(\xi,\eta)|^2\,d\xi\,d\eta.
\endeq

Choose $R$ and $S$ such that $S^{2\mu}/R^{2\lambda'}= 1/S^{2(\mu'-\mu)}$, 
i.e. $R=S^{\mu'/\lambda'}$. This yields a bound like
\[ C_n S^{2\mu+n\left(1 + \frac{\mu'}{\lambda'}\right)}
\|\hat{f}\|_{L^\infty}^2 + 
S^{2(\mu-\mu')} \left(\int |\xi|^{2\lambda'}|\hat{f}|^2 +
\int |\eta|^{2\mu'} |\hat{f}|^2\right).\]
Then the result follows by optimization in $S$.
\med

By symmetry, the same argument works for the case when $\mu=0$.
Now for the general case, we first choose $p$ and $q$ such that
$p^{-1}+q^{-1}=1$ and $p^{-1}\geq \lambda/\lambda'$, $q^{-1}\geq \mu/\mu'$,
and apply H\"older's inequality with conjugate exponents $p$ and $q$:
\begeq\label{beforeHolder}
\int |\xi|^{2\lambda} |\eta|^{2\mu}
|\hat{f}|^2\,d\xi\,d\eta  \leq
\left(\int |\xi|^{2\lambda p} |\hat{f}|^2\,d\xi\,d\eta\right)^{\frac1{p}}
\left( \int|\eta|^{2\mu q} |\hat{f}|^2\,d\xi\,d\eta\right)^{\frac1{q}}.
\endeq
Then we apply to the integrals in the right-hand side of~\eqref{beforeHolder}
the results obtained before for $\lambda=0$ and $\mu=0$:
\[ \int |\xi|^{2\lambda p} |\hat{f}|^2
\leq C \left( \int |\xi|^{2\lambda'} |\hat{f}|^2
+ \int |\eta|^{2\mu'} |\hat{f}|^2 \right)^{1-\theta_1}
\|\hat{f}\|_{L^\infty}^{2\theta_1},\]
and
\[ \int |\eta|^{2\mu q} |\hat{f}|^2
\leq C \left( \int |\xi|^{2\lambda'} |\hat{f}|^2
+ \int |\eta|^{2\mu'} |\hat{f}|^2 \right)^{1-\theta_2}
\|\hat{f}\|_{L^\infty}^{2\theta_2},\]
where
\[ \theta_1 = \frac{\lambda' - \lambda p}{\lambda + \frac{n}2
\left( 1 + \frac{\lambda'}{\mu'}\right)},\qquad
\theta_2 = \frac{\mu' - \mu q}{\mu + \frac{n}2
\left(1 + \frac{\mu'}{\lambda'}\right)}.\]
After some calculation, one finds
\[ \frac{\theta_1}{p} + \frac{\theta_2}{q} = 
\frac{1- \left(\frac{\lambda}{\lambda'} + \frac{\mu}{\mu'}\right)}
{1 + \frac{n}2 \left( \frac1{\lambda'} + \frac1{\mu'}\right)},\]
and the result follows.
\end{proof}

The next technical lemma in this Appendix is an estimate about a system
of differential inequalities. The system may look very particular,
but I~believe that it arises naturally in many problems of hypoelliptic
regularization. In any case, this system is used
in both subsections~\ref{weigL2H1} and~\ref{measureSob} of
Appendix~\ref{appreg}.

\begin{Lem} \label{lemsystem}
Let ${\cal E}$, $X$, $Y$, $Z$ and ${\cal M}$ be continuous functions of $t\in [0,1]$,
with ${\cal E}, X,Y,Z\geq 0$, such that
\begeq\label{syst1}
K (X+Y) \leq {\cal E} \leq C (X+Y),
\endeq
\begeq\label{syst2}
|{\cal M}| \leq C {\cal E}^{1-\delta},
\endeq
\begeq\label{syst3}
\frac{d{\cal E}}{dt} \leq - K Z + C {\cal E},
\endeq
\begeq\label{syst4}
Y \leq C (X+Z)^{1-\theta},
\endeq
\begeq\label{syst5}
\frac{d{\cal M}}{dt} \leq -KX + C (Y+Z),
\endeq
where $C,K$ are positive constants, and $\delta,\theta$ are
real numbers lying in $(0,1)$. Then
\[ {\cal E}(t) \leq \frac{\ov{C}}{t^{1/\kappa}},\qquad
\kappa= \min \left(\delta, \, \frac{\theta}{1-\theta}\right),\]
where $\ov{C}$ is an explicit constant which only depends on
$C,K,\theta,\delta$.
\end{Lem}

\begin{proof}[Proof of Lemma~\ref{lemsystem}]
Let $\tilde{\cal E}(t)=e^{-Ct}{\cal E}(t)$; then
$\tilde{\cal E}$ satisfies estimates similar to ${\cal E}$,
except that equation~\eqref{syst3} becomes $d\tilde{\cal E}/dt\leq -K Z$.
In the sequel I~shall keep the notation ${\cal E}$ for $\tilde{\cal E}$,
so this just amounts to replacing~\eqref{syst3} by
\begeq\label{syst3'}
\frac{d{\cal E}}{dt} \leq - K Z.
\endeq
In particular, ${\cal E}$ is nonincreasing.

Now let $E>0$, and let $I\subset [0,1]$ be the time-interval where 
$(E/2) \leq {\cal E}(t)\leq E$. The goal is to show that
the length $|I|$ of $I$ is bounded like $O(E^{-\kappa})$
for some $\kappa>0$. If that is the case, then the conclusion
follows. Indeed, let $E_0>0$ be given, and let $T$ be the first
time $t$ such that ${\cal E}(t)\leq E_0$, then
\[ T \leq C' \sum_{n\geq 1} E_0^{-n\kappa} \leq C'' E_0^{-\kappa};\]
so $E_0\leq T^{-1/\kappa}$. (Here as in the sequel,
$C$, $C'$, $C''$ stand for various constants that only depend on
the constants $C$ and $K$ appearing in the statement of the lemma.)

If $E\leq 1$ then the conclusion obviously holds true.
So we might assume that $E\geq 1$.

It follows by integration of~\eqref{syst3'} over $I$ that
\begeq\label{intIZE}
\int_I Z(t)\,dt \leq E -\frac{E}2 = \frac{E}2.
\endeq

By integrating~\eqref{syst4}, we find
\begin{align*}
\int_I Y(t)\,dt & \leq C \int_I \bigl[ X(t) + Z(t)\bigr]^{1-\theta}\,dt\\
& \leq C' \left( \int_I X(t)^{1-\theta}\,dt + 
   \int_I Z(t)^{1-\theta}\,dt\right)\\
& \leq C' \left( |I|\: \Bigl[\sup_I X(t)^{1-\theta}\Bigr]
   + \left(\int_I Z(t)\,dt \right)^{1-\theta} |I|^\theta\right).
\end{align*}
To estimate the first term inside the parentheses, note that 
$X\leq C {\cal E}\leq CE$; to bound the second term, use~\eqref{intIZE}.
The result is
\begeq\label{intIY} \int_I Y(t)\,dt \leq C \Bigl( |I| E^{1-\theta} + 
E^{1-\theta} |I|^\theta\Bigr) \leq C' |I|^\theta E^{1-\theta},
\endeq
where the last inequality follows from $|I|\leq |I|^\theta$.
(Note indeed that $|I|\leq 1$ and $\theta<1$.)

Next, integrate inequality~\eqref{syst5} over $I=[t_1,t_2]$, to get
\begin{align}
K \int_I X(t)\,dt & \leq |{\cal M}(t_1)| + |{\cal M}(t_2)|
+ C \int_I [Y(t)+Z(t)]\,dt \noindent \\
& \leq 2 \sup_{t\in I} |{\cal M}(t)| + C\left(\int_I Y(t)\,dt +
\int_I Z(t)\,dt\right).
\label{intIX} 
\end{align}

Also, since ${\cal E}\geq E/2$ on $I$, we have
\begeq\label{TE2}
\frac{|I|\,E}{2} \leq \int_I {\cal E}(t)\,dt
\leq C \left( \int_I X(t)\,dt + \int_I Y(t)\,dt\right),
\endeq
where the last inequality follows from~\eqref{syst1}.

The combination of~\eqref{intIX} and~\eqref{TE2} implies
\[ \frac{|I|\,E}2 \leq C \Bigl ( \sup_{t\in I} |{\cal M}(t)| + 
\int_I Y(t)\,dt + \int_I Z(t)\,dt\Bigr).\]
To estimate the first term inside the brackets, use~\eqref{syst4};
to estimate the second one, use~\eqref{intIY}; to estimate
the third one, use~\eqref{intIZE}. The result is
\begeq\label{TE22} 
|I|\,E \leq C (E^{1-\delta} + |I|^\theta E^{1-\theta} + E).
\endeq

Now we can conclude, separating three cases according to which of
the three terms in the right-hand side of~\eqref{TE22} is largest:

- If it is $E^{1-\delta}$, then
$|I|\,E \leq 3 C E^{1-\delta}$, so $|I|\leq 3 C E^{-\delta}$;

- If it is $|I|^\theta E^{1-\theta}$, then
$|I|\,E \leq 3 C |I|^\theta E^{1-\theta}$, so
$|I| \leq (3C)^{\frac1{1-\theta}} E^{-\frac{\theta}{1-\theta}}$;

- If it is $E$, then $|I|\leq 3C$.
\sm

In any case, there is an estimate like $|I|\leq \ov{C} E^{-\kappa}$,
where $\kappa$ is as in the statement of the lemma.
So the proof is complete.
\end{proof}

The final result in this appendix is a variation of
the usual Korn inequality, used in Subsection~\ref{secMaxw}.

\begin{Prop}[trace Korn inequality] \label{thmkorn}
Let $\Om$ be a smooth bounded connected open subset of $\R^N$.
Then there is a constant $C=C(\Om)$ such that for any vector field
$u\in H^1(\Om;\R^N)$, tangent to the boundary $\pa\Om$,
\begeq\label{ineqkorntrace} 
\|\nabla u\|^2_{L^2(\Om)} \leq C \bigl( \|\nabla^{\rm sym} u\|_{L^2(\Om)}
+ \|u\|_{L^2(\Om)}\bigr),
\endeq
where $\nabla^{\rm sym}u$ stands for the symmetric part of the
matrix-valued field $\nabla u$.
\end{Prop}

\begin{proof} By density, we may assume that $u$ is smooth.
According to~\cite[eq.~(39)-(42)]{DV:korn:02}, if $u$ is tangent to the
boundary, then
\[ \int_\Om |\nabla^{\rm sym} u|^2 = 
\int_\Om |\nabla^\a u|^2 + \int_\Om (\nabla\cdot u)^2 - 
\int_{\pa\Om} \II_\Om (u,u),\]
where $\nabla^\a u$ stands for the antisymmetric part of $\nabla u$,
and $\II_\Om$ for the second fundamental form of the domain $\Om$.
It follows that
\[ \int_\Om |\nabla^\a u|^2 \leq \int_\Om |\nabla^{\rm sym} u|^2
+ C \int_{\pa\Om} |u|^2,\]
where $C=\max_{\pa\Om} \|\II_\Om\|$. Inequality~\eqref{ineqkorntrace}
follows immediately.
\end{proof}

%\backmatter

%\bibliographystyle{acm} 	

%\bibliography{biblio,mybiblio}	

\def\cprime{$'$} \def\cprime{$'$} \def\cprime{$'$} \def\cprime{$'$}
  \def\cprime{$'$}

\end{document}